\documentclass[aos,preprint]{imsart}

\RequirePackage{amsthm,amsmath,amsfonts,amssymb}
\RequirePackage[authoryear]{natbib}

\RequirePackage[OT1]{fontenc}

\usepackage[linesnumbered, ruled,vlined]{algorithm2e}
\usepackage{url}

\usepackage{soul}
\usepackage{mathabx}

\RequirePackage[authoryear]{natbib}
  
\usepackage{graphicx}

\usepackage{tikz}
\usetikzlibrary{calc,fadings,decorations.pathreplacing, angles,quotes, shadings}
\usepackage{ragged2e}
\usepackage{lmodern}
\usetikzlibrary{arrows}
\tikzstyle{block}=[draw opacity=1,line width=1.4cm]
\usepackage{pstricks,pst-node,pst-text,pst-3d}
\usepackage{amsmath}
\usepackage{latexsym}
\usepackage{amsmath}
\usepackage{graphicx}
\usepackage{epic}
\usepackage{eepic}
\usepackage{paralist}
\usepackage{caption}
\usepackage{subcaption}

\startlocaldefs
\numberwithin{equation}{section}
\theoremstyle{plain}

\endlocaldefs

\newif\ifsupplement
\supplementfalse

\usepackage{xr}
\definecolor{AliceBlue}{rgb}{0.94,0.97,1.00}
\definecolor{AntiqueWhite1}{rgb}{1.00,0.94,0.86}
\definecolor{AntiqueWhite2}{rgb}{0.93,0.87,0.80}
\definecolor{AntiqueWhite3}{rgb}{0.80,0.75,0.69}
\definecolor{AntiqueWhite4}{rgb}{0.55,0.51,0.47}
\definecolor{AntiqueWhite}{rgb}{0.98,0.92,0.84}
\definecolor{BlanchedAlmond}{rgb}{1.00,0.92,0.80}
\definecolor{BlueViolet}{rgb}{0.54,0.17,0.89}
\definecolor{CadetBlue1}{rgb}{0.60,0.96,1.00}
\definecolor{CadetBlue2}{rgb}{0.56,0.90,0.93}
\definecolor{CadetBlue3}{rgb}{0.48,0.77,0.80}
\definecolor{CadetBlue4}{rgb}{0.33,0.53,0.55}
\definecolor{CadetBlue}{rgb}{0.37,0.62,0.63}
\definecolor{CornflowerBlue}{rgb}{0.39,0.58,0.93}
\definecolor{DarkBlue}{rgb}{0.00,0.00,0.55}
\definecolor{DarkCyan}{rgb}{0.00,0.55,0.55}
\definecolor{DarkGoldenrod1}{rgb}{1.00,0.73,0.06}
\definecolor{DarkGoldenrod2}{rgb}{0.93,0.68,0.05}
\definecolor{DarkGoldenrod3}{rgb}{0.80,0.58,0.05}
\definecolor{DarkGoldenrod4}{rgb}{0.55,0.40,0.03}
\definecolor{DarkGoldenrod}{rgb}{0.72,0.53,0.04}
\definecolor{DarkGray}{rgb}{0.66,0.66,0.66}
\definecolor{DarkGreen}{rgb}{0.00,0.39,0.00}
\definecolor{DarkGrey}{rgb}{0.66,0.66,0.66}
\definecolor{DarkKhaki}{rgb}{0.74,0.72,0.42}
\definecolor{DarkMagenta}{rgb}{0.55,0.00,0.55}
\definecolor{DarkOliveGreen1}{rgb}{0.79,1.00,0.44}
\definecolor{DarkOliveGreen2}{rgb}{0.74,0.93,0.41}
\definecolor{DarkOliveGreen3}{rgb}{0.64,0.80,0.35}
\definecolor{DarkOliveGreen4}{rgb}{0.43,0.55,0.24}
\definecolor{DarkOliveGreen}{rgb}{0.33,0.42,0.18}
\definecolor{DarkOrange1}{rgb}{1.00,0.50,0.00}
\definecolor{DarkOrange2}{rgb}{0.93,0.46,0.00}
\definecolor{DarkOrange3}{rgb}{0.80,0.40,0.00}
\definecolor{DarkOrange4}{rgb}{0.55,0.27,0.00}
\definecolor{DarkOrange}{rgb}{1.00,0.55,0.00}
\definecolor{DarkOrchid1}{rgb}{0.75,0.24,1.00}
\definecolor{DarkOrchid2}{rgb}{0.70,0.23,0.93}
\definecolor{DarkOrchid3}{rgb}{0.60,0.20,0.80}
\definecolor{DarkOrchid4}{rgb}{0.41,0.13,0.55}
\definecolor{DarkOrchid}{rgb}{0.60,0.20,0.80}
\definecolor{DarkRed}{rgb}{0.55,0.00,0.00}
\definecolor{DarkSalmon}{rgb}{0.91,0.59,0.48}
\definecolor{DarkSeaGreen1}{rgb}{0.76,1.00,0.76}
\definecolor{DarkSeaGreen2}{rgb}{0.71,0.93,0.71}
\definecolor{DarkSeaGreen3}{rgb}{0.61,0.80,0.61}
\definecolor{DarkSeaGreen4}{rgb}{0.41,0.55,0.41}
\definecolor{DarkSeaGreen}{rgb}{0.56,0.74,0.56}
\definecolor{DarkSlateBlue}{rgb}{0.28,0.24,0.55}
\definecolor{DarkSlateGray1}{rgb}{0.59,1.00,1.00}
\definecolor{DarkSlateGray2}{rgb}{0.55,0.93,0.93}
\definecolor{DarkSlateGray3}{rgb}{0.47,0.80,0.80}
\definecolor{DarkSlateGray4}{rgb}{0.32,0.55,0.55}
\definecolor{DarkSlateGray}{rgb}{0.18,0.31,0.31}
\definecolor{DarkSlateGrey}{rgb}{0.18,0.31,0.31}
\definecolor{DarkTurquoise}{rgb}{0.00,0.81,0.82}
\definecolor{DarkViolet}{rgb}{0.58,0.00,0.83}
\definecolor{DeepPink1}{rgb}{1.00,0.08,0.58}
\definecolor{DeepPink2}{rgb}{0.93,0.07,0.54}
\definecolor{DeepPink3}{rgb}{0.80,0.06,0.46}
\definecolor{DeepPink4}{rgb}{0.55,0.04,0.31}
\definecolor{DeepPink}{rgb}{1.00,0.08,0.58}
\definecolor{DeepSkyBlue1}{rgb}{0.00,0.75,1.00}
\definecolor{DeepSkyBlue2}{rgb}{0.00,0.70,0.93}
\definecolor{DeepSkyBlue3}{rgb}{0.00,0.60,0.80}
\definecolor{DeepSkyBlue4}{rgb}{0.00,0.41,0.55}
\definecolor{DeepSkyBlue}{rgb}{0.00,0.75,1.00}
\definecolor{DimGray}{rgb}{0.41,0.41,0.41}
\definecolor{DimGrey}{rgb}{0.41,0.41,0.41}
\definecolor{DodgerBlue1}{rgb}{0.12,0.56,1.00}
\definecolor{DodgerBlue2}{rgb}{0.11,0.53,0.93}
\definecolor{DodgerBlue3}{rgb}{0.09,0.45,0.80}
\definecolor{DodgerBlue4}{rgb}{0.06,0.31,0.55}
\definecolor{DodgerBlue}{rgb}{0.12,0.56,1.00}
\definecolor{FloralWhite}{rgb}{1.00,0.98,0.94}
\definecolor{ForestGreen}{rgb}{0.13,0.55,0.13}
\definecolor{GhostWhite}{rgb}{0.97,0.97,1.00}
\definecolor{GreenYellow}{rgb}{0.68,1.00,0.18}
\definecolor{HotPink1}{rgb}{1.00,0.43,0.71}
\definecolor{HotPink2}{rgb}{0.93,0.42,0.65}
\definecolor{HotPink3}{rgb}{0.80,0.38,0.56}
\definecolor{HotPink4}{rgb}{0.55,0.23,0.38}
\definecolor{HotPink}{rgb}{1.00,0.41,0.71}
\definecolor{IndianRed1}{rgb}{1.00,0.42,0.42}
\definecolor{IndianRed2}{rgb}{0.93,0.39,0.39}
\definecolor{IndianRed3}{rgb}{0.80,0.33,0.33}
\definecolor{IndianRed4}{rgb}{0.55,0.23,0.23}
\definecolor{IndianRed}{rgb}{0.80,0.36,0.36}
\definecolor{LavenderBlush1}{rgb}{1.00,0.94,0.96}
\definecolor{LavenderBlush2}{rgb}{0.93,0.88,0.90}
\definecolor{LavenderBlush3}{rgb}{0.80,0.76,0.77}
\definecolor{LavenderBlush4}{rgb}{0.55,0.51,0.53}
\definecolor{LavenderBlush}{rgb}{1.00,0.94,0.96}
\definecolor{LawnGreen}{rgb}{0.49,0.99,0.00}
\definecolor{LemonChiffon1}{rgb}{1.00,0.98,0.80}
\definecolor{LemonChiffon2}{rgb}{0.93,0.91,0.75}
\definecolor{LemonChiffon3}{rgb}{0.80,0.79,0.65}
\definecolor{LemonChiffon4}{rgb}{0.55,0.54,0.44}
\definecolor{LemonChiffon}{rgb}{1.00,0.98,0.80}
\definecolor{LightBlue1}{rgb}{0.75,0.94,1.00}
\definecolor{LightBlue2}{rgb}{0.70,0.87,0.93}
\definecolor{LightBlue3}{rgb}{0.60,0.75,0.80}
\definecolor{LightBlue4}{rgb}{0.41,0.51,0.55}
\definecolor{LightBlue}{rgb}{0.68,0.85,0.90}
\definecolor{LightCoral}{rgb}{0.94,0.50,0.50}
\definecolor{LightCyan1}{rgb}{0.88,1.00,1.00}
\definecolor{LightCyan2}{rgb}{0.82,0.93,0.93}
\definecolor{LightCyan3}{rgb}{0.71,0.80,0.80}
\definecolor{LightCyan4}{rgb}{0.48,0.55,0.55}
\definecolor{LightCyan}{rgb}{0.88,1.00,1.00}
\definecolor{LightGoldenrod1}{rgb}{1.00,0.93,0.55}
\definecolor{LightGoldenrod2}{rgb}{0.93,0.86,0.51}
\definecolor{LightGoldenrod3}{rgb}{0.80,0.75,0.44}
\definecolor{LightGoldenrod4}{rgb}{0.55,0.51,0.30}
\definecolor{LightGoldenrodYellow}{rgb}{0.98,0.98,0.82}
\definecolor{LightGoldenrod}{rgb}{0.93,0.87,0.51}
\definecolor{LightGray}{rgb}{0.83,0.83,0.83}
\definecolor{LightGreen}{rgb}{0.56,0.93,0.56}
\definecolor{LightGrey}{rgb}{0.83,0.83,0.83}
\definecolor{LightPink1}{rgb}{1.00,0.68,0.73}
\definecolor{LightPink2}{rgb}{0.93,0.64,0.68}
\definecolor{LightPink3}{rgb}{0.80,0.55,0.58}
\definecolor{LightPink4}{rgb}{0.55,0.37,0.40}
\definecolor{LightPink}{rgb}{1.00,0.71,0.76}
\definecolor{LightSalmon1}{rgb}{1.00,0.63,0.48}
\definecolor{LightSalmon2}{rgb}{0.93,0.58,0.45}
\definecolor{LightSalmon3}{rgb}{0.80,0.51,0.38}
\definecolor{LightSalmon4}{rgb}{0.55,0.34,0.26}
\definecolor{LightSalmon}{rgb}{1.00,0.63,0.48}
\definecolor{LightSeaGreen}{rgb}{0.13,0.70,0.67}
\definecolor{LightSkyBlue1}{rgb}{0.69,0.89,1.00}
\definecolor{LightSkyBlue2}{rgb}{0.64,0.83,0.93}
\definecolor{LightSkyBlue3}{rgb}{0.55,0.71,0.80}
\definecolor{LightSkyBlue4}{rgb}{0.38,0.48,0.55}
\definecolor{LightSkyBlue}{rgb}{0.53,0.81,0.98}
\definecolor{LightSlateBlue}{rgb}{0.52,0.44,1.00}
\definecolor{LightSlateGray}{rgb}{0.47,0.53,0.60}
\definecolor{LightSlateGrey}{rgb}{0.47,0.53,0.60}
\definecolor{LightSteelBlue1}{rgb}{0.79,0.88,1.00}
\definecolor{LightSteelBlue2}{rgb}{0.74,0.82,0.93}
\definecolor{LightSteelBlue3}{rgb}{0.64,0.71,0.80}
\definecolor{LightSteelBlue4}{rgb}{0.43,0.48,0.55}
\definecolor{LightSteelBlue}{rgb}{0.69,0.77,0.87}
\definecolor{LightYellow1}{rgb}{1.00,1.00,0.88}
\definecolor{LightYellow2}{rgb}{0.93,0.93,0.82}
\definecolor{LightYellow3}{rgb}{0.80,0.80,0.71}
\definecolor{LightYellow4}{rgb}{0.55,0.55,0.48}
\definecolor{LightYellow}{rgb}{1.00,1.00,0.88}
\definecolor{LimeGreen}{rgb}{0.20,0.80,0.20}
\definecolor{MediumAquamarine}{rgb}{0.40,0.80,0.67}
\definecolor{MediumBlue}{rgb}{0.00,0.00,0.80}
\definecolor{MediumOrchid1}{rgb}{0.88,0.40,1.00}
\definecolor{MediumOrchid2}{rgb}{0.82,0.37,0.93}
\definecolor{MediumOrchid3}{rgb}{0.71,0.32,0.80}
\definecolor{MediumOrchid4}{rgb}{0.48,0.22,0.55}
\definecolor{MediumOrchid}{rgb}{0.73,0.33,0.83}
\definecolor{MediumPurple1}{rgb}{0.67,0.51,1.00}
\definecolor{MediumPurple2}{rgb}{0.62,0.47,0.93}
\definecolor{MediumPurple3}{rgb}{0.54,0.41,0.80}
\definecolor{MediumPurple4}{rgb}{0.36,0.28,0.55}
\definecolor{MediumPurple}{rgb}{0.58,0.44,0.86}
\definecolor{MediumSeaGreen}{rgb}{0.24,0.70,0.44}
\definecolor{MediumSlateBlue}{rgb}{0.48,0.41,0.93}
\definecolor{MediumSpringGreen}{rgb}{0.00,0.98,0.60}
\definecolor{MediumTurquoise}{rgb}{0.28,0.82,0.80}
\definecolor{MediumVioletRed}{rgb}{0.78,0.08,0.52}
\definecolor{MidnightBlue}{rgb}{0.10,0.10,0.44}
\definecolor{MintCream}{rgb}{0.96,1.00,0.98}
\definecolor{MistyRose1}{rgb}{1.00,0.89,0.88}
\definecolor{MistyRose2}{rgb}{0.93,0.84,0.82}
\definecolor{MistyRose3}{rgb}{0.80,0.72,0.71}
\definecolor{MistyRose4}{rgb}{0.55,0.49,0.48}
\definecolor{MistyRose}{rgb}{1.00,0.89,0.88}
\definecolor{NavajoWhite1}{rgb}{1.00,0.87,0.68}
\definecolor{NavajoWhite2}{rgb}{0.93,0.81,0.63}
\definecolor{NavajoWhite3}{rgb}{0.80,0.70,0.55}
\definecolor{NavajoWhite4}{rgb}{0.55,0.47,0.37}
\definecolor{NavajoWhite}{rgb}{1.00,0.87,0.68}
\definecolor{NavyBlue}{rgb}{0.00,0.00,0.50}
\definecolor{OldLace}{rgb}{0.99,0.96,0.90}
\definecolor{OliveDrab1}{rgb}{0.75,1.00,0.24}
\definecolor{OliveDrab2}{rgb}{0.70,0.93,0.23}
\definecolor{OliveDrab3}{rgb}{0.60,0.80,0.20}
\definecolor{OliveDrab4}{rgb}{0.41,0.55,0.13}
\definecolor{OliveDrab}{rgb}{0.42,0.56,0.14}
\definecolor{OrangeRed1}{rgb}{1.00,0.27,0.00}
\definecolor{OrangeRed2}{rgb}{0.93,0.25,0.00}
\definecolor{OrangeRed3}{rgb}{0.80,0.22,0.00}
\definecolor{OrangeRed4}{rgb}{0.55,0.15,0.00}
\definecolor{OrangeRed}{rgb}{1.00,0.27,0.00}
\definecolor{PaleGoldenrod}{rgb}{0.93,0.91,0.67}
\definecolor{PaleGreen1}{rgb}{0.60,1.00,0.60}
\definecolor{PaleGreen2}{rgb}{0.56,0.93,0.56}
\definecolor{PaleGreen3}{rgb}{0.49,0.80,0.49}
\definecolor{PaleGreen4}{rgb}{0.33,0.55,0.33}
\definecolor{PaleGreen}{rgb}{0.60,0.98,0.60}
\definecolor{PaleTurquoise1}{rgb}{0.73,1.00,1.00}
\definecolor{PaleTurquoise2}{rgb}{0.68,0.93,0.93}
\definecolor{PaleTurquoise3}{rgb}{0.59,0.80,0.80}
\definecolor{PaleTurquoise4}{rgb}{0.40,0.55,0.55}
\definecolor{PaleTurquoise}{rgb}{0.69,0.93,0.93}
\definecolor{PaleVioletRed1}{rgb}{1.00,0.51,0.67}
\definecolor{PaleVioletRed2}{rgb}{0.93,0.47,0.62}
\definecolor{PaleVioletRed3}{rgb}{0.80,0.41,0.54}
\definecolor{PaleVioletRed4}{rgb}{0.55,0.28,0.36}
\definecolor{PaleVioletRed}{rgb}{0.86,0.44,0.58}
\definecolor{PapayaWhip}{rgb}{1.00,0.94,0.84}
\definecolor{PeachPuff1}{rgb}{1.00,0.85,0.73}
\definecolor{PeachPuff2}{rgb}{0.93,0.80,0.68}
\definecolor{PeachPuff3}{rgb}{0.80,0.69,0.58}
\definecolor{PeachPuff4}{rgb}{0.55,0.47,0.40}
\definecolor{PeachPuff}{rgb}{1.00,0.85,0.73}
\definecolor{PowderBlue}{rgb}{0.69,0.88,0.90}
\definecolor{RosyBrown1}{rgb}{1.00,0.76,0.76}
\definecolor{RosyBrown2}{rgb}{0.93,0.71,0.71}
\definecolor{RosyBrown3}{rgb}{0.80,0.61,0.61}
\definecolor{RosyBrown4}{rgb}{0.55,0.41,0.41}
\definecolor{RosyBrown}{rgb}{0.74,0.56,0.56}
\definecolor{RoyalBlue1}{rgb}{0.28,0.46,1.00}
\definecolor{RoyalBlue2}{rgb}{0.26,0.43,0.93}
\definecolor{RoyalBlue3}{rgb}{0.23,0.37,0.80}
\definecolor{RoyalBlue4}{rgb}{0.15,0.25,0.55}
\definecolor{RoyalBlue}{rgb}{0.25,0.41,0.88}
\definecolor{SaddleBrown}{rgb}{0.55,0.27,0.07}
\definecolor{SandyBrown}{rgb}{0.96,0.64,0.38}
\definecolor{SeaGreen1}{rgb}{0.33,1.00,0.62}
\definecolor{SeaGreen2}{rgb}{0.31,0.93,0.58}
\definecolor{SeaGreen3}{rgb}{0.26,0.80,0.50}
\definecolor{SeaGreen4}{rgb}{0.18,0.55,0.34}
\definecolor{SeaGreen}{rgb}{0.18,0.55,0.34}
\definecolor{SkyBlue1}{rgb}{0.53,0.81,1.00}
\definecolor{SkyBlue2}{rgb}{0.49,0.75,0.93}
\definecolor{SkyBlue3}{rgb}{0.42,0.65,0.80}
\definecolor{SkyBlue4}{rgb}{0.29,0.44,0.55}
\definecolor{SkyBlue}{rgb}{0.53,0.81,0.92}
\definecolor{SlateBlue1}{rgb}{0.51,0.44,1.00}
\definecolor{SlateBlue2}{rgb}{0.48,0.40,0.93}
\definecolor{SlateBlue3}{rgb}{0.41,0.35,0.80}
\definecolor{SlateBlue4}{rgb}{0.28,0.24,0.55}
\definecolor{SlateBlue}{rgb}{0.42,0.35,0.80}
\definecolor{SlateGray1}{rgb}{0.78,0.89,1.00}
\definecolor{SlateGray2}{rgb}{0.73,0.83,0.93}
\definecolor{SlateGray3}{rgb}{0.62,0.71,0.80}
\definecolor{SlateGray4}{rgb}{0.42,0.48,0.55}
\definecolor{SlateGray}{rgb}{0.44,0.50,0.56}
\definecolor{SlateGrey}{rgb}{0.44,0.50,0.56}
\definecolor{SpringGreen1}{rgb}{0.00,1.00,0.50}
\definecolor{SpringGreen2}{rgb}{0.00,0.93,0.46}
\definecolor{SpringGreen3}{rgb}{0.00,0.80,0.40}
\definecolor{SpringGreen4}{rgb}{0.00,0.55,0.27}
\definecolor{SpringGreen}{rgb}{0.00,1.00,0.50}
\definecolor{SteelBlue1}{rgb}{0.39,0.72,1.00}
\definecolor{SteelBlue2}{rgb}{0.36,0.67,0.93}
\definecolor{SteelBlue3}{rgb}{0.31,0.58,0.80}
\definecolor{SteelBlue4}{rgb}{0.21,0.39,0.55}
\definecolor{SteelBlue}{rgb}{0.27,0.51,0.71}
\definecolor{VioletRed1}{rgb}{1.00,0.24,0.59}
\definecolor{VioletRed2}{rgb}{0.93,0.23,0.55}
\definecolor{VioletRed3}{rgb}{0.80,0.20,0.47}
\definecolor{VioletRed4}{rgb}{0.55,0.13,0.32}
\definecolor{VioletRed}{rgb}{0.82,0.13,0.56}
\definecolor{WhiteSmoke}{rgb}{0.96,0.96,0.96}
\definecolor{YellowGreen}{rgb}{0.60,0.80,0.20}
\definecolor{aliceblue}{rgb}{0.94,0.97,1.00}
\definecolor{antiquewhite}{rgb}{0.98,0.92,0.84}
\definecolor{aquamarine1}{rgb}{0.50,1.00,0.83}
\definecolor{aquamarine2}{rgb}{0.46,0.93,0.78}
\definecolor{aquamarine3}{rgb}{0.40,0.80,0.67}
\definecolor{aquamarine4}{rgb}{0.27,0.55,0.45}
\definecolor{aquamarine}{rgb}{0.50,1.00,0.83}
\definecolor{azure1}{rgb}{0.94,1.00,1.00}
\definecolor{azure2}{rgb}{0.88,0.93,0.93}
\definecolor{azure3}{rgb}{0.76,0.80,0.80}
\definecolor{azure4}{rgb}{0.51,0.55,0.55}
\definecolor{azure}{rgb}{0.94,1.00,1.00}
\definecolor{beige}{rgb}{0.96,0.96,0.86}
\definecolor{bisque1}{rgb}{1.00,0.89,0.77}
\definecolor{bisque2}{rgb}{0.93,0.84,0.72}
\definecolor{bisque3}{rgb}{0.80,0.72,0.62}
\definecolor{bisque4}{rgb}{0.55,0.49,0.42}
\definecolor{bisque}{rgb}{1.00,0.89,0.77}
\definecolor{black}{rgb}{0.00,0.00,0.00}
\definecolor{blanchedalmond}{rgb}{1.00,0.92,0.80}
\definecolor{blue1}{rgb}{0.00,0.00,1.00}
\definecolor{blue2}{rgb}{0.00,0.00,0.93}
\definecolor{blue3}{rgb}{0.00,0.00,0.80}
\definecolor{blue4}{rgb}{0.00,0.00,0.55}
\definecolor{blueviolet}{rgb}{0.54,0.17,0.89}
\definecolor{blue}{rgb}{0.00,0.00,1.00}
\definecolor{brown1}{rgb}{1.00,0.25,0.25}
\definecolor{brown2}{rgb}{0.93,0.23,0.23}
\definecolor{brown3}{rgb}{0.80,0.20,0.20}
\definecolor{brown4}{rgb}{0.55,0.14,0.14}
\definecolor{brown}{rgb}{0.65,0.16,0.16}
\definecolor{burlywood1}{rgb}{1.00,0.83,0.61}
\definecolor{burlywood2}{rgb}{0.93,0.77,0.57}
\definecolor{burlywood3}{rgb}{0.80,0.67,0.49}
\definecolor{burlywood4}{rgb}{0.55,0.45,0.33}
\definecolor{burlywood}{rgb}{0.87,0.72,0.53}
\definecolor{cadetblue}{rgb}{0.37,0.62,0.63}
\definecolor{chartreuse1}{rgb}{0.50,1.00,0.00}
\definecolor{chartreuse2}{rgb}{0.46,0.93,0.00}
\definecolor{chartreuse3}{rgb}{0.40,0.80,0.00}
\definecolor{chartreuse4}{rgb}{0.27,0.55,0.00}
\definecolor{chartreuse}{rgb}{0.50,1.00,0.00}
\definecolor{chocolate1}{rgb}{1.00,0.50,0.14}
\definecolor{chocolate2}{rgb}{0.93,0.46,0.13}
\definecolor{chocolate3}{rgb}{0.80,0.40,0.11}
\definecolor{chocolate4}{rgb}{0.55,0.27,0.07}
\definecolor{chocolate}{rgb}{0.82,0.41,0.12}
\definecolor{coral1}{rgb}{1.00,0.45,0.34}
\definecolor{coral2}{rgb}{0.93,0.42,0.31}
\definecolor{coral3}{rgb}{0.80,0.36,0.27}
\definecolor{coral4}{rgb}{0.55,0.24,0.18}
\definecolor{coral}{rgb}{1.00,0.50,0.31}
\definecolor{cornflowerblue}{rgb}{0.39,0.58,0.93}
\definecolor{cornsilk1}{rgb}{1.00,0.97,0.86}
\definecolor{cornsilk2}{rgb}{0.93,0.91,0.80}
\definecolor{cornsilk3}{rgb}{0.80,0.78,0.69}
\definecolor{cornsilk4}{rgb}{0.55,0.53,0.47}
\definecolor{cornsilk}{rgb}{1.00,0.97,0.86}
\definecolor{cyan1}{rgb}{0.00,1.00,1.00}
\definecolor{cyan2}{rgb}{0.00,0.93,0.93}
\definecolor{cyan3}{rgb}{0.00,0.80,0.80}
\definecolor{cyan4}{rgb}{0.00,0.55,0.55}
\definecolor{cyan}{rgb}{0.00,1.00,1.00}
\definecolor{darkblue}{rgb}{0.00,0.00,0.55}
\definecolor{darkcyan}{rgb}{0.00,0.55,0.55}
\definecolor{darkgoldenrod}{rgb}{0.72,0.53,0.04}
\definecolor{darkgray}{rgb}{0.66,0.66,0.66}
\definecolor{darkgreen}{rgb}{0.00,0.39,0.00}
\definecolor{darkgrey}{rgb}{0.66,0.66,0.66}
\definecolor{darkkhaki}{rgb}{0.74,0.72,0.42}
\definecolor{darkmagenta}{rgb}{0.55,0.00,0.55}
\definecolor{darkolive}{rgb}{0.33,0.42,0.18}
\definecolor{darkorange}{rgb}{1.00,0.55,0.00}
\definecolor{darkorchid}{rgb}{0.60,0.20,0.80}
\definecolor{darkred}{rgb}{0.55,0.00,0.00}
\definecolor{darksalmon}{rgb}{0.91,0.59,0.48}
\definecolor{darksea}{rgb}{0.56,0.74,0.56}
\definecolor{darkslate}{rgb}{0.18,0.31,0.31}
\definecolor{darkslate}{rgb}{0.18,0.31,0.31}
\definecolor{darkslate}{rgb}{0.28,0.24,0.55}
\definecolor{darkturquoise}{rgb}{0.00,0.81,0.82}
\definecolor{darkviolet}{rgb}{0.58,0.00,0.83}
\definecolor{deeppink}{rgb}{1.00,0.08,0.58}
\definecolor{deepsky}{rgb}{0.00,0.75,1.00}
\definecolor{dimgray}{rgb}{0.41,0.41,0.41}
\definecolor{dimgrey}{rgb}{0.41,0.41,0.41}
\definecolor{dodgerblue}{rgb}{0.12,0.56,1.00}
\definecolor{firebrick1}{rgb}{1.00,0.19,0.19}
\definecolor{firebrick2}{rgb}{0.93,0.17,0.17}
\definecolor{firebrick3}{rgb}{0.80,0.15,0.15}
\definecolor{firebrick4}{rgb}{0.55,0.10,0.10}
\definecolor{firebrick}{rgb}{0.70,0.13,0.13}
\definecolor{floralwhite}{rgb}{1.00,0.98,0.94}
\definecolor{forestgreen}{rgb}{0.13,0.55,0.13}
\definecolor{gainsboro}{rgb}{0.86,0.86,0.86}
\definecolor{ghostwhite}{rgb}{0.97,0.97,1.00}
\definecolor{gold1}{rgb}{1.00,0.84,0.00}
\definecolor{gold2}{rgb}{0.93,0.79,0.00}
\definecolor{gold3}{rgb}{0.80,0.68,0.00}
\definecolor{gold4}{rgb}{0.55,0.46,0.00}
\definecolor{goldenrod1}{rgb}{1.00,0.76,0.15}
\definecolor{goldenrod2}{rgb}{0.93,0.71,0.13}
\definecolor{goldenrod3}{rgb}{0.80,0.61,0.11}
\definecolor{goldenrod4}{rgb}{0.55,0.41,0.08}
\definecolor{goldenrod}{rgb}{0.85,0.65,0.13}
\definecolor{gold}{rgb}{1.00,0.84,0.00}
\definecolor{gray0}{rgb}{0.00,0.00,0.00}
\definecolor{gray100}{rgb}{1.00,1.00,1.00}
\definecolor{gray10}{rgb}{0.10,0.10,0.10}
\definecolor{gray11}{rgb}{0.11,0.11,0.11}
\definecolor{gray12}{rgb}{0.12,0.12,0.12}
\definecolor{gray13}{rgb}{0.13,0.13,0.13}
\definecolor{gray14}{rgb}{0.14,0.14,0.14}
\definecolor{gray15}{rgb}{0.15,0.15,0.15}
\definecolor{gray16}{rgb}{0.16,0.16,0.16}
\definecolor{gray17}{rgb}{0.17,0.17,0.17}
\definecolor{gray18}{rgb}{0.18,0.18,0.18}
\definecolor{gray19}{rgb}{0.19,0.19,0.19}
\definecolor{gray1}{rgb}{0.01,0.01,0.01}
\definecolor{gray20}{rgb}{0.20,0.20,0.20}
\definecolor{gray21}{rgb}{0.21,0.21,0.21}
\definecolor{gray22}{rgb}{0.22,0.22,0.22}
\definecolor{gray23}{rgb}{0.23,0.23,0.23}
\definecolor{gray24}{rgb}{0.24,0.24,0.24}
\definecolor{gray25}{rgb}{0.25,0.25,0.25}
\definecolor{gray26}{rgb}{0.26,0.26,0.26}
\definecolor{gray27}{rgb}{0.27,0.27,0.27}
\definecolor{gray28}{rgb}{0.28,0.28,0.28}
\definecolor{gray29}{rgb}{0.29,0.29,0.29}
\definecolor{gray2}{rgb}{0.02,0.02,0.02}
\definecolor{gray30}{rgb}{0.30,0.30,0.30}
\definecolor{gray31}{rgb}{0.31,0.31,0.31}
\definecolor{gray32}{rgb}{0.32,0.32,0.32}
\definecolor{gray33}{rgb}{0.33,0.33,0.33}
\definecolor{gray34}{rgb}{0.34,0.34,0.34}
\definecolor{gray35}{rgb}{0.35,0.35,0.35}
\definecolor{gray36}{rgb}{0.36,0.36,0.36}
\definecolor{gray37}{rgb}{0.37,0.37,0.37}
\definecolor{gray38}{rgb}{0.38,0.38,0.38}
\definecolor{gray39}{rgb}{0.39,0.39,0.39}
\definecolor{gray3}{rgb}{0.03,0.03,0.03}
\definecolor{gray40}{rgb}{0.40,0.40,0.40}
\definecolor{gray41}{rgb}{0.41,0.41,0.41}
\definecolor{gray42}{rgb}{0.42,0.42,0.42}
\definecolor{gray43}{rgb}{0.43,0.43,0.43}
\definecolor{gray44}{rgb}{0.44,0.44,0.44}
\definecolor{gray45}{rgb}{0.45,0.45,0.45}
\definecolor{gray46}{rgb}{0.46,0.46,0.46}
\definecolor{gray47}{rgb}{0.47,0.47,0.47}
\definecolor{gray48}{rgb}{0.48,0.48,0.48}
\definecolor{gray49}{rgb}{0.49,0.49,0.49}
\definecolor{gray4}{rgb}{0.04,0.04,0.04}
\definecolor{gray50}{rgb}{0.50,0.50,0.50}
\definecolor{gray51}{rgb}{0.51,0.51,0.51}
\definecolor{gray52}{rgb}{0.52,0.52,0.52}
\definecolor{gray53}{rgb}{0.53,0.53,0.53}
\definecolor{gray54}{rgb}{0.54,0.54,0.54}
\definecolor{gray55}{rgb}{0.55,0.55,0.55}
\definecolor{gray56}{rgb}{0.56,0.56,0.56}
\definecolor{gray57}{rgb}{0.57,0.57,0.57}
\definecolor{gray58}{rgb}{0.58,0.58,0.58}
\definecolor{gray59}{rgb}{0.59,0.59,0.59}
\definecolor{gray5}{rgb}{0.05,0.05,0.05}
\definecolor{gray60}{rgb}{0.60,0.60,0.60}
\definecolor{gray61}{rgb}{0.61,0.61,0.61}
\definecolor{gray62}{rgb}{0.62,0.62,0.62}
\definecolor{gray63}{rgb}{0.63,0.63,0.63}
\definecolor{gray64}{rgb}{0.64,0.64,0.64}
\definecolor{gray65}{rgb}{0.65,0.65,0.65}
\definecolor{gray66}{rgb}{0.66,0.66,0.66}
\definecolor{gray67}{rgb}{0.67,0.67,0.67}
\definecolor{gray68}{rgb}{0.68,0.68,0.68}
\definecolor{gray69}{rgb}{0.69,0.69,0.69}
\definecolor{gray6}{rgb}{0.06,0.06,0.06}
\definecolor{gray70}{rgb}{0.70,0.70,0.70}
\definecolor{gray71}{rgb}{0.71,0.71,0.71}
\definecolor{gray72}{rgb}{0.72,0.72,0.72}
\definecolor{gray73}{rgb}{0.73,0.73,0.73}
\definecolor{gray74}{rgb}{0.74,0.74,0.74}
\definecolor{gray75}{rgb}{0.75,0.75,0.75}
\definecolor{gray76}{rgb}{0.76,0.76,0.76}
\definecolor{gray77}{rgb}{0.77,0.77,0.77}
\definecolor{gray78}{rgb}{0.78,0.78,0.78}
\definecolor{gray79}{rgb}{0.79,0.79,0.79}
\definecolor{gray7}{rgb}{0.07,0.07,0.07}
\definecolor{gray80}{rgb}{0.80,0.80,0.80}
\definecolor{gray81}{rgb}{0.81,0.81,0.81}
\definecolor{gray82}{rgb}{0.82,0.82,0.82}
\definecolor{gray83}{rgb}{0.83,0.83,0.83}
\definecolor{gray84}{rgb}{0.84,0.84,0.84}
\definecolor{gray85}{rgb}{0.85,0.85,0.85}
\definecolor{gray86}{rgb}{0.86,0.86,0.86}
\definecolor{gray87}{rgb}{0.87,0.87,0.87}
\definecolor{gray88}{rgb}{0.88,0.88,0.88}
\definecolor{gray89}{rgb}{0.89,0.89,0.89}
\definecolor{gray8}{rgb}{0.08,0.08,0.08}
\definecolor{gray90}{rgb}{0.90,0.90,0.90}
\definecolor{gray91}{rgb}{0.91,0.91,0.91}
\definecolor{gray92}{rgb}{0.92,0.92,0.92}
\definecolor{gray93}{rgb}{0.93,0.93,0.93}
\definecolor{gray94}{rgb}{0.94,0.94,0.94}
\definecolor{gray95}{rgb}{0.95,0.95,0.95}
\definecolor{gray96}{rgb}{0.96,0.96,0.96}
\definecolor{gray97}{rgb}{0.97,0.97,0.97}
\definecolor{gray98}{rgb}{0.98,0.98,0.98}
\definecolor{gray99}{rgb}{0.99,0.99,0.99}
\definecolor{gray9}{rgb}{0.09,0.09,0.09}
\definecolor{gray}{rgb}{0.75,0.75,0.75}
\definecolor{green1}{rgb}{0.00,1.00,0.00}
\definecolor{green2}{rgb}{0.00,0.93,0.00}
\definecolor{green3}{rgb}{0.00,0.80,0.00}
\definecolor{green4}{rgb}{0.00,0.55,0.00}
\definecolor{greenyellow}{rgb}{0.68,1.00,0.18}
\definecolor{green}{rgb}{0.00,1.00,0.00}
\definecolor{grey0}{rgb}{0.00,0.00,0.00}
\definecolor{grey100}{rgb}{1.00,1.00,1.00}
\definecolor{grey10}{rgb}{0.10,0.10,0.10}
\definecolor{grey11}{rgb}{0.11,0.11,0.11}
\definecolor{grey12}{rgb}{0.12,0.12,0.12}
\definecolor{grey13}{rgb}{0.13,0.13,0.13}
\definecolor{grey14}{rgb}{0.14,0.14,0.14}
\definecolor{grey15}{rgb}{0.15,0.15,0.15}
\definecolor{grey16}{rgb}{0.16,0.16,0.16}
\definecolor{grey17}{rgb}{0.17,0.17,0.17}
\definecolor{grey18}{rgb}{0.18,0.18,0.18}
\definecolor{grey19}{rgb}{0.19,0.19,0.19}
\definecolor{grey1}{rgb}{0.01,0.01,0.01}
\definecolor{grey20}{rgb}{0.20,0.20,0.20}
\definecolor{grey21}{rgb}{0.21,0.21,0.21}
\definecolor{grey22}{rgb}{0.22,0.22,0.22}
\definecolor{grey23}{rgb}{0.23,0.23,0.23}
\definecolor{grey24}{rgb}{0.24,0.24,0.24}
\definecolor{grey25}{rgb}{0.25,0.25,0.25}
\definecolor{grey26}{rgb}{0.26,0.26,0.26}
\definecolor{grey27}{rgb}{0.27,0.27,0.27}
\definecolor{grey28}{rgb}{0.28,0.28,0.28}
\definecolor{grey29}{rgb}{0.29,0.29,0.29}
\definecolor{grey2}{rgb}{0.02,0.02,0.02}
\definecolor{grey30}{rgb}{0.30,0.30,0.30}
\definecolor{grey31}{rgb}{0.31,0.31,0.31}
\definecolor{grey32}{rgb}{0.32,0.32,0.32}
\definecolor{grey33}{rgb}{0.33,0.33,0.33}
\definecolor{grey34}{rgb}{0.34,0.34,0.34}
\definecolor{grey35}{rgb}{0.35,0.35,0.35}
\definecolor{grey36}{rgb}{0.36,0.36,0.36}
\definecolor{grey37}{rgb}{0.37,0.37,0.37}
\definecolor{grey38}{rgb}{0.38,0.38,0.38}
\definecolor{grey39}{rgb}{0.39,0.39,0.39}
\definecolor{grey3}{rgb}{0.03,0.03,0.03}
\definecolor{grey40}{rgb}{0.40,0.40,0.40}
\definecolor{grey41}{rgb}{0.41,0.41,0.41}
\definecolor{grey42}{rgb}{0.42,0.42,0.42}
\definecolor{grey43}{rgb}{0.43,0.43,0.43}
\definecolor{grey44}{rgb}{0.44,0.44,0.44}
\definecolor{grey45}{rgb}{0.45,0.45,0.45}
\definecolor{grey46}{rgb}{0.46,0.46,0.46}
\definecolor{grey47}{rgb}{0.47,0.47,0.47}
\definecolor{grey48}{rgb}{0.48,0.48,0.48}
\definecolor{grey49}{rgb}{0.49,0.49,0.49}
\definecolor{grey4}{rgb}{0.04,0.04,0.04}
\definecolor{grey50}{rgb}{0.50,0.50,0.50}
\definecolor{grey51}{rgb}{0.51,0.51,0.51}
\definecolor{grey52}{rgb}{0.52,0.52,0.52}
\definecolor{grey53}{rgb}{0.53,0.53,0.53}
\definecolor{grey54}{rgb}{0.54,0.54,0.54}
\definecolor{grey55}{rgb}{0.55,0.55,0.55}
\definecolor{grey56}{rgb}{0.56,0.56,0.56}
\definecolor{grey57}{rgb}{0.57,0.57,0.57}
\definecolor{grey58}{rgb}{0.58,0.58,0.58}
\definecolor{grey59}{rgb}{0.59,0.59,0.59}
\definecolor{grey5}{rgb}{0.05,0.05,0.05}
\definecolor{grey60}{rgb}{0.60,0.60,0.60}
\definecolor{grey61}{rgb}{0.61,0.61,0.61}
\definecolor{grey62}{rgb}{0.62,0.62,0.62}
\definecolor{grey63}{rgb}{0.63,0.63,0.63}
\definecolor{grey64}{rgb}{0.64,0.64,0.64}
\definecolor{grey65}{rgb}{0.65,0.65,0.65}
\definecolor{grey66}{rgb}{0.66,0.66,0.66}
\definecolor{grey67}{rgb}{0.67,0.67,0.67}
\definecolor{grey68}{rgb}{0.68,0.68,0.68}
\definecolor{grey69}{rgb}{0.69,0.69,0.69}
\definecolor{grey6}{rgb}{0.06,0.06,0.06}
\definecolor{grey70}{rgb}{0.70,0.70,0.70}
\definecolor{grey71}{rgb}{0.71,0.71,0.71}
\definecolor{grey72}{rgb}{0.72,0.72,0.72}
\definecolor{grey73}{rgb}{0.73,0.73,0.73}
\definecolor{grey74}{rgb}{0.74,0.74,0.74}
\definecolor{grey75}{rgb}{0.75,0.75,0.75}
\definecolor{grey76}{rgb}{0.76,0.76,0.76}
\definecolor{grey77}{rgb}{0.77,0.77,0.77}
\definecolor{grey78}{rgb}{0.78,0.78,0.78}
\definecolor{grey79}{rgb}{0.79,0.79,0.79}
\definecolor{grey7}{rgb}{0.07,0.07,0.07}
\definecolor{grey80}{rgb}{0.80,0.80,0.80}
\definecolor{grey81}{rgb}{0.81,0.81,0.81}
\definecolor{grey82}{rgb}{0.82,0.82,0.82}
\definecolor{grey83}{rgb}{0.83,0.83,0.83}
\definecolor{grey84}{rgb}{0.84,0.84,0.84}
\definecolor{grey85}{rgb}{0.85,0.85,0.85}
\definecolor{grey86}{rgb}{0.86,0.86,0.86}
\definecolor{grey87}{rgb}{0.87,0.87,0.87}
\definecolor{grey88}{rgb}{0.88,0.88,0.88}
\definecolor{grey89}{rgb}{0.89,0.89,0.89}
\definecolor{grey8}{rgb}{0.08,0.08,0.08}
\definecolor{grey90}{rgb}{0.90,0.90,0.90}
\definecolor{grey91}{rgb}{0.91,0.91,0.91}
\definecolor{grey92}{rgb}{0.92,0.92,0.92}
\definecolor{grey93}{rgb}{0.93,0.93,0.93}
\definecolor{grey94}{rgb}{0.94,0.94,0.94}
\definecolor{grey95}{rgb}{0.95,0.95,0.95}
\definecolor{grey96}{rgb}{0.96,0.96,0.96}
\definecolor{grey97}{rgb}{0.97,0.97,0.97}
\definecolor{grey98}{rgb}{0.98,0.98,0.98}
\definecolor{grey99}{rgb}{0.99,0.99,0.99}
\definecolor{grey9}{rgb}{0.09,0.09,0.09}
\definecolor{grey}{rgb}{0.75,0.75,0.75}
\definecolor{honeydew1}{rgb}{0.94,1.00,0.94}
\definecolor{honeydew2}{rgb}{0.88,0.93,0.88}
\definecolor{honeydew3}{rgb}{0.76,0.80,0.76}
\definecolor{honeydew4}{rgb}{0.51,0.55,0.51}
\definecolor{honeydew}{rgb}{0.94,1.00,0.94}
\definecolor{hotpink}{rgb}{1.00,0.41,0.71}
\definecolor{indianred}{rgb}{0.80,0.36,0.36}
\definecolor{ivory1}{rgb}{1.00,1.00,0.94}
\definecolor{ivory2}{rgb}{0.93,0.93,0.88}
\definecolor{ivory3}{rgb}{0.80,0.80,0.76}
\definecolor{ivory4}{rgb}{0.55,0.55,0.51}
\definecolor{ivory}{rgb}{1.00,1.00,0.94}
\definecolor{khaki1}{rgb}{1.00,0.96,0.56}
\definecolor{khaki2}{rgb}{0.93,0.90,0.52}
\definecolor{khaki3}{rgb}{0.80,0.78,0.45}
\definecolor{khaki4}{rgb}{0.55,0.53,0.31}
\definecolor{khaki}{rgb}{0.94,0.90,0.55}
\definecolor{lavenderblush}{rgb}{1.00,0.94,0.96}
\definecolor{lavender}{rgb}{0.90,0.90,0.98}
\definecolor{lawngreen}{rgb}{0.49,0.99,0.00}
\definecolor{lemonchiffon}{rgb}{1.00,0.98,0.80}
\definecolor{lightblue}{rgb}{0.68,0.85,0.90}
\definecolor{lightcoral}{rgb}{0.94,0.50,0.50}
\definecolor{lightcyan}{rgb}{0.88,1.00,1.00}
\definecolor{lightgoldenrod}{rgb}{0.93,0.87,0.51}
\definecolor{lightgoldenrod}{rgb}{0.98,0.98,0.82}
\definecolor{lightgray}{rgb}{0.83,0.83,0.83}
\definecolor{lightgreen}{rgb}{0.56,0.93,0.56}
\definecolor{lightgrey}{rgb}{0.83,0.83,0.83}
\definecolor{lightpink}{rgb}{1.00,0.71,0.76}
\definecolor{lightsalmon}{rgb}{1.00,0.63,0.48}
\definecolor{lightsea}{rgb}{0.13,0.70,0.67}
\definecolor{lightsky}{rgb}{0.53,0.81,0.98}
\definecolor{lightslate}{rgb}{0.47,0.53,0.60}
\definecolor{lightslate}{rgb}{0.47,0.53,0.60}
\definecolor{lightslate}{rgb}{0.52,0.44,1.00}
\definecolor{lightsteel}{rgb}{0.69,0.77,0.87}
\definecolor{lightyellow}{rgb}{1.00,1.00,0.88}
\definecolor{limegreen}{rgb}{0.20,0.80,0.20}
\definecolor{linen}{rgb}{0.98,0.94,0.90}
\definecolor{magenta1}{rgb}{1.00,0.00,1.00}
\definecolor{magenta2}{rgb}{0.93,0.00,0.93}
\definecolor{magenta3}{rgb}{0.80,0.00,0.80}
\definecolor{magenta4}{rgb}{0.55,0.00,0.55}
\definecolor{magenta}{rgb}{1.00,0.00,1.00}
\definecolor{maroon1}{rgb}{1.00,0.20,0.70}
\definecolor{maroon2}{rgb}{0.93,0.19,0.65}
\definecolor{maroon3}{rgb}{0.80,0.16,0.56}
\definecolor{maroon4}{rgb}{0.55,0.11,0.38}
\definecolor{maroon}{rgb}{0.69,0.19,0.38}
\definecolor{mediumaquamarine}{rgb}{0.40,0.80,0.67}
\definecolor{mediumblue}{rgb}{0.00,0.00,0.80}
\definecolor{mediumorchid}{rgb}{0.73,0.33,0.83}
\definecolor{mediumpurple}{rgb}{0.58,0.44,0.86}
\definecolor{mediumsea}{rgb}{0.24,0.70,0.44}
\definecolor{mediumslate}{rgb}{0.48,0.41,0.93}
\definecolor{mediumspring}{rgb}{0.00,0.98,0.60}
\definecolor{mediumturquoise}{rgb}{0.28,0.82,0.80}
\definecolor{mediumviolet}{rgb}{0.78,0.08,0.52}
\definecolor{midnightblue}{rgb}{0.10,0.10,0.44}
\definecolor{mintcream}{rgb}{0.96,1.00,0.98}
\definecolor{mistyrose}{rgb}{1.00,0.89,0.88}
\definecolor{moccasin}{rgb}{1.00,0.89,0.71}
\definecolor{navajowhite}{rgb}{1.00,0.87,0.68}
\definecolor{navyblue}{rgb}{0.00,0.00,0.50}
\definecolor{navy}{rgb}{0.00,0.00,0.50}
\definecolor{oldlace}{rgb}{0.99,0.96,0.90}
\definecolor{olivedrab}{rgb}{0.42,0.56,0.14}
\definecolor{orange1}{rgb}{1.00,0.65,0.00}
\definecolor{orange2}{rgb}{0.93,0.60,0.00}
\definecolor{orange3}{rgb}{0.80,0.52,0.00}
\definecolor{orange4}{rgb}{0.55,0.35,0.00}
\definecolor{orangered}{rgb}{1.00,0.27,0.00}
\definecolor{orange}{rgb}{1.00,0.65,0.00}
\definecolor{orchid1}{rgb}{1.00,0.51,0.98}
\definecolor{orchid2}{rgb}{0.93,0.48,0.91}
\definecolor{orchid3}{rgb}{0.80,0.41,0.79}
\definecolor{orchid4}{rgb}{0.55,0.28,0.54}
\definecolor{orchid}{rgb}{0.85,0.44,0.84}
\definecolor{palegoldenrod}{rgb}{0.93,0.91,0.67}
\definecolor{palegreen}{rgb}{0.60,0.98,0.60}
\definecolor{paleturquoise}{rgb}{0.69,0.93,0.93}
\definecolor{paleviolet}{rgb}{0.86,0.44,0.58}
\definecolor{papayawhip}{rgb}{1.00,0.94,0.84}
\definecolor{peachpuff}{rgb}{1.00,0.85,0.73}
\definecolor{peru}{rgb}{0.80,0.52,0.25}
\definecolor{pink1}{rgb}{1.00,0.71,0.77}
\definecolor{pink2}{rgb}{0.93,0.66,0.72}
\definecolor{pink3}{rgb}{0.80,0.57,0.62}
\definecolor{pink4}{rgb}{0.55,0.39,0.42}
\definecolor{pink}{rgb}{1.00,0.75,0.80}
\definecolor{plum1}{rgb}{1.00,0.73,1.00}
\definecolor{plum2}{rgb}{0.93,0.68,0.93}
\definecolor{plum3}{rgb}{0.80,0.59,0.80}
\definecolor{plum4}{rgb}{0.55,0.40,0.55}
\definecolor{plum}{rgb}{0.87,0.63,0.87}
\definecolor{powderblue}{rgb}{0.69,0.88,0.90}
\definecolor{purple1}{rgb}{0.61,0.19,1.00}
\definecolor{purple2}{rgb}{0.57,0.17,0.93}
\definecolor{purple3}{rgb}{0.49,0.15,0.80}
\definecolor{purple4}{rgb}{0.33,0.10,0.55}
\definecolor{purple}{rgb}{0.63,0.13,0.94}
\definecolor{red1}{rgb}{1.00,0.00,0.00}
\definecolor{red2}{rgb}{0.93,0.00,0.00}
\definecolor{red3}{rgb}{0.80,0.00,0.00}
\definecolor{red4}{rgb}{0.55,0.00,0.00}
\definecolor{red}{rgb}{1.00,0.00,0.00}
\definecolor{rosybrown}{rgb}{0.74,0.56,0.56}
\definecolor{royalblue}{rgb}{0.25,0.41,0.88}
\definecolor{saddlebrown}{rgb}{0.55,0.27,0.07}
\definecolor{salmon1}{rgb}{1.00,0.55,0.41}
\definecolor{salmon2}{rgb}{0.93,0.51,0.38}
\definecolor{salmon3}{rgb}{0.80,0.44,0.33}
\definecolor{salmon4}{rgb}{0.55,0.30,0.22}
\definecolor{salmon}{rgb}{0.98,0.50,0.45}
\definecolor{sandybrown}{rgb}{0.96,0.64,0.38}
\definecolor{seagreen}{rgb}{0.18,0.55,0.34}
\definecolor{seashell1}{rgb}{1.00,0.96,0.93}
\definecolor{seashell2}{rgb}{0.93,0.90,0.87}
\definecolor{seashell3}{rgb}{0.80,0.77,0.75}
\definecolor{seashell4}{rgb}{0.55,0.53,0.51}
\definecolor{seashell}{rgb}{1.00,0.96,0.93}
\definecolor{sienna1}{rgb}{1.00,0.51,0.28}
\definecolor{sienna2}{rgb}{0.93,0.47,0.26}
\definecolor{sienna3}{rgb}{0.80,0.41,0.22}
\definecolor{sienna4}{rgb}{0.55,0.28,0.15}
\definecolor{sienna}{rgb}{0.63,0.32,0.18}
\definecolor{skyblue}{rgb}{0.53,0.81,0.92}
\definecolor{slateblue}{rgb}{0.42,0.35,0.80}
\definecolor{slategray}{rgb}{0.44,0.50,0.56}
\definecolor{slategrey}{rgb}{0.44,0.50,0.56}
\definecolor{snow1}{rgb}{1.00,0.98,0.98}
\definecolor{snow2}{rgb}{0.93,0.91,0.91}
\definecolor{snow3}{rgb}{0.80,0.79,0.79}
\definecolor{snow4}{rgb}{0.55,0.54,0.54}
\definecolor{snow}{rgb}{1.00,0.98,0.98}
\definecolor{springgreen}{rgb}{0.00,1.00,0.50}
\definecolor{steelblue}{rgb}{0.27,0.51,0.71}
\definecolor{tan1}{rgb}{1.00,0.65,0.31}
\definecolor{tan2}{rgb}{0.93,0.60,0.29}
\definecolor{tan3}{rgb}{0.80,0.52,0.25}
\definecolor{tan4}{rgb}{0.55,0.35,0.17}
\definecolor{tan}{rgb}{0.82,0.71,0.55}
\definecolor{thistle1}{rgb}{1.00,0.88,1.00}
\definecolor{thistle2}{rgb}{0.93,0.82,0.93}
\definecolor{thistle3}{rgb}{0.80,0.71,0.80}
\definecolor{thistle4}{rgb}{0.55,0.48,0.55}
\definecolor{thistle}{rgb}{0.85,0.75,0.85}
\definecolor{tomato1}{rgb}{1.00,0.39,0.28}
\definecolor{tomato2}{rgb}{0.93,0.36,0.26}
\definecolor{tomato3}{rgb}{0.80,0.31,0.22}
\definecolor{tomato4}{rgb}{0.55,0.21,0.15}
\definecolor{tomato}{rgb}{1.00,0.39,0.28}
\definecolor{turquoise1}{rgb}{0.00,0.96,1.00}
\definecolor{turquoise2}{rgb}{0.00,0.90,0.93}
\definecolor{turquoise3}{rgb}{0.00,0.77,0.80}
\definecolor{turquoise4}{rgb}{0.00,0.53,0.55}
\definecolor{turquoise}{rgb}{0.25,0.88,0.82}
\definecolor{violetred}{rgb}{0.82,0.13,0.56}
\definecolor{violet}{rgb}{0.93,0.51,0.93}
\definecolor{wheat1}{rgb}{1.00,0.91,0.73}
\definecolor{wheat2}{rgb}{0.93,0.85,0.68}
\definecolor{wheat3}{rgb}{0.80,0.73,0.59}
\definecolor{wheat4}{rgb}{0.55,0.49,0.40}
\definecolor{wheat}{rgb}{0.96,0.87,0.70}
\definecolor{whitesmoke}{rgb}{0.96,0.96,0.96}
\definecolor{white}{rgb}{1.00,1.00,1.00}
\definecolor{yellow1}{rgb}{1.00,1.00,0.00}
\definecolor{yellow2}{rgb}{0.93,0.93,0.00}
\definecolor{yellow3}{rgb}{0.80,0.80,0.00}
\definecolor{yellow4}{rgb}{0.55,0.55,0.00}
\definecolor{yellowgreen}{rgb}{0.60,0.80,0.20}
\definecolor{yellow}{rgb}{1.00,1.00,0.00}

\usepackage{lscape}

\usepackage{pgfplots}
\usepackage{dsfont}    

\usepackage{mathtools}

\newtheorem{Lem}{Lemma}[section] 
\newtheorem{Rem}{Remark}[section]
\newtheorem{proposition}{Proposition}[section]

\newtheorem{Example}{Example}[section]
\newtheorem{definition}{Definition}[section]
\newtheorem{assumption}{Assumption}

\DeclareMathOperator*{\argmin}{arg\,min}

\newcommand{\mockalph}[1]{}

\newcommand{\mbf}[1]{\mbox{\boldmath$#1$\unboldmath}}

\newcommand{\cqfd}{\hfill $\square$}

\newcommand{\R}{\mathbb R}

\newcommand{\n}{^{(n)}}

\newcommand{\Xb}{\mathbf{X}}
\newcommand{\Sb}{\mathbf{S}}

\newcommand{\Ib}{\mathbf{I}}

\newcommand{\Qb}{\mathbf{Q}}
\newcommand{\Rb}{\mathbf{R}}

\newcommand{\Fb}{\mathbf{F}}

\newcommand{\ub}{\ensuremath{\mathbf{u}}}

\newcommand{\vb}{\ensuremath{\mathbf{v}}}
\newcommand{\xb}{\ensuremath{\mathbf{x}}}
\newcommand{\yb}{\ensuremath{\mathbf{y}}}
\newcommand{\zb}{\ensuremath{\mathbf{z}}}

\newcommand{\Ub}{\ensuremath{\mathbf{U}}}

\newcommand{\Yb}{\ensuremath{\mathbf{Y}}}

\newcommand{\wb}{\ensuremath{\mathbf{w}}}

\newcommand{\thetab}{{\pmb \theta}}

\newcommand{\mub}{{\pmb \mu}}

\newcommand{\Gammab}{{\pmb \Gamma}}
\newcommand{\kappab}{{\pmb \kappa}}

\newcommand{\nub}{{\pmb \nu}}

\newcommand{\xib}{{\pmb \xi}}

\newcommand{\pr}{^{\prime}}

\newcommand{\ny}{n\rightarrow\infty}

\begin{document}
  
  \setlength{\abovedisplayskip}{5pt}
\setlength{\belowdisplayskip}{5pt}
\setlength{\abovedisplayshortskip}{5pt}
\setlength{\belowdisplayshortskip}{5pt}

\begin{frontmatter}
\title{Quantiles 
 and Quantile Regression \\ on Riemannian Manifolds\\  a  measure-transportation-based approach}
\runtitle{Quantiles on Riemannian Manifolds}

%

\begin{aug}
\author[A,B]{\fnms{Marc} \snm{Hallin}\ead[label=e1]{mhallin@ulb.ac.be}}
\and
\author[C]{\fnms{Hang} \snm{Liu}\ead[label=e2]{hliu01@ustc.edu.cn}}

\address[A]{ECARES {\it and Department of Mathematics, Universit\'e Libre de Bruxelles, Belgium}\vspace{-2mm}}

\address[B]{{\it Institute of Information Theory and Automation, Czech Academy of Sciences, Prague, Czech Republic}\vspace{1mm} \\ {\rm\printead{e1}}}

\address[C]{{\it Faculty of Business in SciTech, School of Management, University of Science and Technology of China} {\rm\printead{e2}}}
\end{aug}


\runauthor{M. Hallin and H. Liu}

\begin{abstract}
Increased attention has been given recently to the statistical analysis of variables with values on nonlinear  manifolds. A natural but nontrivial problem in that context is the definition of quantile concepts.  We are proposing a solution   for compact Riemannian manifolds without boundaries; typical examples are polyspheres, hyperspheres, and toro\"{\i}dal manifolds equipped with their Riemannian metrics.  Our concept of quantile function comes along with a concept of distribution function and, in the empirical case,  ranks and signs. The absence of a canonical ordering is offset by resorting to the data-driven ordering induced by optimal transports. Theoretical properties, such as the uniform convergence of the empirical distribution and conditional (and unconditional) quantile functions and distribution-freeness of ranks and signs, are established.  Statistical inference applications, from  goodness-of-fit to distribution-free rank-based testing,  are without number. Of particular importance is  the case of quantile regression with directional or  toro\"{\i}dal multiple output, which is given special attention in this paper.  Extensive simulations are carried out to illustrate these novel concepts.   
\end{abstract}



\end{frontmatter}

\section{Measure transportation on Riemannian manifolds}\hspace{-3mm}

Let $({\cal M}, g)$ denote a $p$-dimensional manifold  equipped\footnote{Appendix A, where we refer to for definitions, 
 is providing a summary of all Riemannian geometry concepts  used throughout the paper, such as {\it Riemannian metric, manifold smoothness, geodesic distance,} etc.  For  
 details, see   \cite{petersen2006riemannian}, \cite{jost2008riemannian}, and \cite{lee2012smooth}}    with  {\it Riemannian metric}~$g$. We throughout assume that this manifold is  connected, compact, $C^4$-smooth, and has no boundary. Typical examples are the hypersphere or $p$-sphere~$\mathcal{S}^{p}$  in~${\mathbb R}^{p+1}$,  the~$p$-torus~$\mathcal{T}^{p}$ 
(the product $\mathcal{S}^{1}\!\times \ldots\times\! \mathcal{S}^{1}$ of  $p$ circles), and the  polysphere $\boldsymbol{\mathcal{S}}^{p_1\ldots p_k}$ (the pro\-duct~$\mathcal{S}^{p_1} \times \cdots \times \mathcal{S}^{p_k}$ of  $k$ hyperspheres). While the $1$- and~$2$-spheres~$\mathcal{S}^{1}$ and~$\mathcal{S}^{2}$ have countless applications in directional statistics (see, for instance, \cite{Ley2017a}),   the tori $\mathcal{T}^{2}$ and $\mathcal{T}^{3}$ and the more complex polyspheres   are  given increasing attention due to their applications in protein bioinformatics \citep{Mardia1999a, PatrangenaruEllington15, Ley2017a,Ley2018, Eltzner2018, Ameijeiras-Alonso2020a, Pewsey2021,Gonzalezetal21}  and the skeleton approach to shape analysis \citep{Eltzner2015, Eltzner2017, Pizer2017, Pizer2020, polyspheres2022}. Parallel to this application-oriented litera\-ture, a more abstract   mathematical approach also has been developed, e.g., in \cite{Bhattacharya2003, Bhattacharya2005} and \cite{BB08, bhattacharya2012nonparametric}. Measure transportation methods in the area are more recent and mainly focused the hypersphere---see \cite{HLV2022} where ranks and quantiles are introduced for $\mathcal{S}^{p}$-valued observations. In the present paper, we are extending that approach  to more general Riemannian manifolds, with emphasis on the definitions of quantile contours and regions, along with the corresponding notions of ranks and signs, and develop  a method of quantile regression for~$\mathcal{S}^{p}$- or~$\mathcal{T}^{p}$-valued responses and covariates in  general Polish metric spaces. 


Denote by $d(\yb_1, \yb_2)$ the {\it geodesic distance} between $\yb_1$ and  $\yb_2$ in $\mathcal{M}$, by ${\frak F}_{\mathcal{M}}$  the  corresponding Borel $\sigma$-field on~$\mathcal{M}$, and  by ${{\rm P}_1 }$ and ${{\rm P}_2 }$ two probability measures on $(\mathcal{M}, {\frak F}_{\mathcal{M}})$. Consider the set ${\mathcal T}({{\rm P}_1 ,Q})$ of all measurable maps~$\bf T$ from~$(\mathcal{M}, {\frak F}_{\mathcal{M}})$ to $(\mathcal{M}, {\frak F}_{\mathcal{M}})$ (call them {\it transport maps}) that are
 {\it pushing~${{\rm P}_1 }$ forward to ${{\rm P}_2 }$}.\footnote{That is, such that that~$({\bf T}\#{{\rm P}_1 })(V) \coloneqq {{\rm P}_1 }({\bf T}^{-1}(V)) ={{\rm P}_2 }(V)$
for any~$V \in {\frak F}_{\mathcal{M}}$.}   Measure transportation, in its original form  initiated 
  (for~$\mathcal{M} =~\!{\mathbb R}^p$) by \cite{monge1781}, aims at minimizing, over~${\bf T}\in{\mathcal T}({{\rm P}_1 ,Q})$  
  the expected {\it transportation cost}
\begin{equation}\label{cost}
{\rm C}_{\rm M}({\bf T}) \coloneqq \int_{\mathcal{M}} c(\yb, {\bf {\bf T}}(\yb)) {\rm d}{{\rm P}_1 }(\yb)    
\end{equation}
where~$c(\yb, \zb)$ is the cost of transporting $\yb$ to $\zb$. Throughout, 
 we consider the cost function~$c(\yb, \zb) = d^2(\yb, \zb)/2$.


Because  $\mathcal{T}({{\rm P}_1 }, {{\rm P}_2 })$ is not convex, Monge's problem is an uneasy one. 
A more tractable formulation was proposed by \cite{kantorovich} who, rather than   transport maps,   considers {\it transport plans}.  A transport plan is an element $\gamma$ of the set $ \Gamma({{\rm P}_1 }, {{\rm P}_2 })$ of all proba\-bility distributions over $(\mathcal{M}, {\frak F}_{\mathcal{M}})\!\times\! (\mathcal{M}, {\frak F}_{\mathcal{M}})$ with marginals ${\rm P}_1 $ and~${\rm P}_2 $ (a convex set). Kantorovich's problem consists in minimizing, over~$\gamma\in \Gamma({{\rm P}_1 }, {{\rm P}_2 })$,  the total transportation cost
\begin{equation}\label{costK}
    {\rm C}_{\rm K}(\gamma) \coloneqq \int_{\mathcal{M}\times \mathcal{M}} c(\yb, \zb){\rm d} \gamma(\yb, \zb). 
\end{equation}
The  dual of this minimization problem is   the maximization  over $(\psi, \phi)\! \in\! C(\mathcal{M})\! \times~\!\! C(\mathcal{M}) $~of
$$ \int_{\mathcal{M}} \phi(\zb)  {\rm d}{{\rm P}_2 } - \int_{\mathcal{M}} \psi(\yb)  {\rm d}{{\rm P}_1 }
$$
 where $C(\mathcal{M})$ denotes the set of continuous  functions from~$\mathcal{M}$ to $\mathbb R$  satisfying 
$$\phi(\zb)-\psi(\yb) \leq c(\yb, \zb)\quad\text{ for all }\quad (\yb, \zb) \in \mathcal{M}\! \times\! \mathcal{M}.$$

It turns out that the solution to Kantorovich's problem is a transport plan of the form~$\gamma = ({\bf I}_{\rm d} \times {\bf F}) \# {{\rm P}_1 }$ where ${\bf F}$ is a solution to Monge's problem (see Proposition~\ref{keylemma} below). Moreover, ${\bf F}$ takes the form ${\bf F}(\yb) = \exp_{\yb}(\nabla \psi (\yb))$ (see Appendix A for a definition) where $(\psi,\phi)$ is the solution to Kantorovich's dual problem  and satisfies 
\begin{equation}\label{cconv}
\psi(\yb)=\!\! \sup_{\zb \in   \mathcal{M}} \{ \phi(\zb) -c(\yb, \zb) \} \ \ \text{ and }\ \ 
 \phi(\zb)=\!\! \inf_{\yb \in   \mathcal{M}} \{ \psi(\yb) + c(\yb, \zb) \}. 
\end{equation}
Functions $\psi$ and $\phi$  satisfying \eqref{cconv} are called {\it $c$-convex} and {\it $c$-concave}, respectively;  $\phi$ (resp.\  $\psi$) is called the \emph{c-transform} of $\psi$ (resp.\  $\phi$). 

Another important concept in this  
 context is {\it $c$-cyclical monotonicity}. 
\begin{definition}{\rm 
A subset $\Omega \subseteq {\cal M} \times {\cal M}$ is called {{\it $c$-cyclically monotone}}  if, denoting by $\Sigma(k)$  the set of permutations of~$\{1,\ldots,k\}$, for all $k \in {\mathbb N}$, all~$\sigma \in \Sigma(k)$, and all~$(\yb_1, \zb_1), \ldots, (\yb_k, \zb_k)\in \Omega$,  }
 $\sum_{i=1}^k c(\yb_i, \zb_i) \leq \sum_{i=1}^k c(\yb_{\sigma(i)}, \zb_i).$
\end{definition}

Let $\mathfrak{P}$ denote the set of probability measures that are  absolutely continuous with respect to the Riemannian volume measure. The following proposition, which provides 
 key ingredients for our further theoretical developments,  summarizes some of the main  results 
   on measure transportation on manifolds.
\begin{proposition}\label{keylemma} 
Let ${\rm P}_1  \in \mathfrak{P}$ and ${\rm P}_2 $ be two probability measures on~$\mathcal{M}$. Then,  letting  $c(\yb, \zb) = d^2(\yb, \zb)/2$,
\begin{compactenum}
  \item[(i)] a transport plan $\gamma\in\Gamma({{\rm P}_1 },{{\rm P}_2 })$ is the solution (is minimizing the transportation cost~$C_{\rm K}$  in \eqref{costK} over $\Gamma ({{\rm P}_1 }, {{\rm P}_2 })$) of the Kantorovich problem if and only if it is supported on a   $c$-cyclically monotone subset of~$ {\cal M}\!  \times\! {\cal M}$; 
 \item[(ii)] this solution $\gamma$ of the Kantorovich problem exists, is ${\rm P}_1 $-a.s.  unique, and is of  the   form~$({\bf I}_{\rm d}~\!\!\times~\!\!{\bf F}) \# {{\rm P}_1 }$ where ${\bf F}\in{\mathcal T}({{\rm P}_1 },{{\rm P}_2 })$ is the~${\rm P}_1 $-a.s. unique solution of the corresponding Monge problem (minimizing  $C_{\rm M}$   in \eqref{cost}  over ${\mathcal T} ({{\rm P}_1 }, {{\rm P}_2 })$); 
 \item[(iii)] there exist $c$-convex differentiable mappings
   $\psi$  from~$\mathcal{M}$ to $\R$ such\linebreak  that~${\bf F}(\yb)  = \exp_{\yb}(\nabla \psi (\yb))$.\smallskip
   \end{compactenum}
   
\noindent If, moreover, ${\rm P}_2  \in \mathfrak{P}$, then\smallskip
 
 \begin{compactenum}
 \item[(iv)]  ${\bf Q}(\xb)\coloneqq  \exp_{\xb}(\nabla \phi (\xb))$, with $\phi$ the $c$-transform of $\psi$, belongs to ${\mathcal T}({{\rm P}_2 }, {{\rm P}_1 })$, is the~${\rm P}_2 $-a.s. unique minimizer of $C_{\rm M}$ in~\eqref{cost} over ${\mathcal T} ({{\rm P}_2 }, {{\rm P}_1 })$  (i.e., is the~${\rm P}_2 $-a.s.\  unique solution of the  corresponding  Monge problem),  and  satisfies 
 $${\bf Q}({\bf F}(\xb)) = \xb\quad {\rm P}_1 \text{-a.s.\quad and \quad} {\bf F}({\bf Q}(\xb)) = \xb\quad {\rm P}_2 \text{-a.s.} $$
   \end{compactenum}
\end{proposition}
\noindent Under part (iv) of this proposition, $\bf Q$ thus is a ${\rm P}_1 $-a.s. inverse of $\bf F$. 

Necessity and sufficiency  in Part~(i) were established by \cite{Rusch1996} and \citet[Theorem~1]{SCHACHERMAYERTEICHMANN2008}, respectively. Parts (ii), (iii), and~(iv) are proved in \cite{McCann2001}.

\section{Quantile functions, quantile regions, and quantile contours}\label{sec.DistQuan}

In this section, building on the results of Proposition~\ref{keylemma}, we propose   concepts of   distribution and quantile functions for a random variable (r.v.) $\Yb$ with values in $( {\cal M}, {\frak F}_{\mathcal{M}})$. To this end, we need the continuity of the transport maps involved, which relies on the regularity of the cost function~$c$.


\subsection{Regularity conditions}

\color{black}
Letting $\text{\rm Cut}(\yb)$ denote the {\it cut locus}\footnote{See Appendix~\ref{App.Rieman} for  definitions.} of $\yb \in {\cal M}$, define~$\text{\rm Cut}(\mathcal{M}) \coloneqq  \bigcup\left\{\{\yb\}\times \text{\rm Cut}(\yb) \vert\, \yb\in{\mathcal M}\right\}$. Write $\nabla_{\yb} c(\yb, \zb)$ for the gradient$^3$ computed at $\Yb = \yb$ 
 of $c(\cdot, \zb): {\cal M} \rightarrow \mathbb{R}, \Yb \mapsto c(\Yb, \zb)$.

We assume that the following regularity conditions (similar to \cite{loeper2011}) hold for the restriction to~${\cal M}_*^2 \coloneqq  \{ (\yb, \zb), \yb \in \mathcal{M} , \zb \in \mathcal{M} \setminus {\rm Cut}(\yb) \}  $ of the cost function $c$; \linebreak for~$c=d^2/2$, this assumption holds if ${\cal M}$ is $C^4$-smooth. 

\begin{assumption}\label{ass.cost}
The cost function $c$  is such that  (i) $(\yb, \zb)\mapsto c(\yb, \zb)$ is $C^4$-smooth on ${\cal M}_*^2$; (ii) for all $\yb \in {\cal M}$, 
the mapping $\zb \mapsto -\nabla_{\yb} c(\yb, \zb)$ is injective on ${\cal M}\setminus\text{\rm Cut}(\yb)$; 
(iii) $\det((c_{i,j})_{1\leq i, j\leq p}) \neq 0$
for all $(\yb, \zb) \in {\cal M}_*^2$, where $c_{i,j}$ is defined below.
\end{assumption}



Assumption~\ref{ass.cost} allows us to define (Definition~\ref{costseccurv} below) the notion of {\it cost-sectional curvature} on ${\cal M}_*^2$, which was first proposed by \cite{ma2005regularity} for the Euclidean space and later on extended by, e.g., \cite{loeper2009} and \cite{KimMcCann2010} to Riemannian manifolds. Non-negativity of the cost-sectional curvature is a necessary condition 
 for the continuity of the optimal transport. Combined with the regularity conditions on $c$, ${\rm P}_1$, and ${{\rm P}_2 }$, it is also a sufficient condition for  the $C^1$ or $C^{\infty}$ continuity of the optimal transport.

Write (under Assumption~\ref{ass.cost}) 
 $c_{i_1 \ldots i_a, \cdot} (\yb, \zb)$, $c_{_\cdot ,\, j_1 \ldots j_b} (\yb, \zb)$, and~$c_{i_1 \ldots i_a,\, j_1 \ldots j_b} (\yb, \zb)$,  respectively,  for the derivatives of $(\yb ,\zb )\mapsto c(\yb ,\zb )=d^2(\yb, \zb)/2$ with respect to~$y_{i_1}, \ldots , y_{i_a}$, to~$z_{j_1},\ldots , z_{j_b}$, and to~$y_{i_1}, \ldots , y_{i_a}, z_{j_1},\ldots , z_{j_b}$ evaluated at~$(\yb ,\zb)\in {\cal M}_*^2$.   
%
%
%
%
Still under Assumption~\ref{ass.cost}, denote by~$c^{r,s}$  the~$(r,s)$ entry of the inverse of the matrix~$\nabla^2_{\yb \zb} c(\yb, \zb) = (c_{i,j})$. Finally, let~$T_{\yb}{\cal M}$ stand for the {\it tangent space} of ${\cal M}$ at $\yb$.  The {\it cost-sectional curvature} {at $(\yb ,\zb) \in{\cal M}_*^2$} is defined as follows  \citep{LoeperVillani2010}.

\begin{definition}\label{costseccurv}
Under Assumption~\ref{ass.cost},   the {\it cost-sectional curvature 
   {at~$(\yb, \zb) \in {\cal M}_*^2$}} is the  mapping (from  $T_{\yb}{\cal M} \times T_{\zb}{\cal M}$   to $\mathbb R$)
\begin{equation}\label{costcurv}
{\mathfrak C}(\yb, \zb) :\, ({\pmb{\xi}}, {\pmb{\eta}})\mapsto {\mathfrak C}(\yb, \zb) ({\pmb{\xi}}, {\pmb{\eta}}) \coloneqq \frac{3}{2} \sum_{ijk\ell r s  {=1}}^{{p}}(c_{ij,r} c^{r,s} c_{s,k\ell} - c_{ij,k\ell}) \xi^{i} \xi^j \eta^{k} \eta^{\ell}
\end{equation}
 where $\xi^i$ and $\eta^i$ denote the $i$-th components of $\pmb{\xi}$ and~${\pmb{\eta}}$,~respectively. 
\end{definition}


On this cost-sectional curvature, we make the following nonnegativity assumption.
\begin{assumption}\label{ass.costCurve}
For any $(\yb, \zb) \in {\cal M}_*^2$ and any $(\pmb{\xi}, {\pmb{\eta}}) \in T_{\yb}{\cal M} \times T_{\zb}{\cal M}$ such that~$\pmb{\xi} \bot {\pmb{\eta}}$, the cost-sectional curvature satisfies 
${\mathfrak C}(\yb, \zb)
({\pmb{\xi}}, {\pmb{\eta}}) \geq 0$.
\end{assumption}

\begin{Rem}
{\rm 
This  uniform nonnegativity requirement  on ${\mathfrak C}(\yb, \zb)$ 
is necessary  if $\psi$ in Proposition~\ref{keylemma}  is to be $C^1$-smooth. 
For $\zb = \yb$, the expression~\eqref{costcurv} of ${\mathfrak C}(\yb, \zb)
({\pmb{\xi}}, {\pmb{\eta}})$ 
  coincides (up to a positive scalar factor) with the Riemannian sectional curvature for~$c=d^2/2$ (see Theorem~3.8 in \citet{loeper2009}). Hence, if the sectional curvature of a Riemannian manifold is not everywhere nonnegative, 
Assumption~\ref{ass.costCurve} cannot hold, and the optimal transport maps between ${\rm P}_1 $ and ${\rm P}_2 $ might not be continuous.
}
\end{Rem}

Denote by $(\exp_{\yb})^{-1}: {\cal M}\setminus \text{\rm Cut}(\yb) \rightarrow T_{\yb}{\cal M}$ the inverse of the exponential mapping~$\exp_{\yb}$. 
For $c=d^2/2$ and~${\cal M}_*^2$ satisfying Assumption~\ref{ass.cost}, we impose  on ${\cal M}\setminus \text{\rm Cut}(\yb)$ the following assumption (which is similar to the convex support assumption in $\mathbb{R}^d$ made in \citet{delbarrio2020regularity}).

\begin{assumption}\label{ass.convex}
For all $\yb \in \mathcal{M}$, the set $(\exp_{\yb})^{-1}({\cal M}\setminus \text{\rm Cut}(\yb)) \subset T_{\yb}{\cal M}$ is convex.

\end{assumption}

In order 
 to establish the continuity 
  of the optimal transport map~$\Fb$, we need to assume that its graph $\{(\yb, \Fb(\yb)), \yb \in {\cal M}\}$   lies in a subdomain of~${\cal M}_*^2$, hence remains {\it uniformly}  away from the cut locus.  

\begin{assumption}\label{Ass.FyCut}
Let ${{\rm P}_1 } \in \mathfrak{P}$  have bounded density and  ${{\rm P}_2 }\in \mathfrak{P}$ have nonvanishing density on ${\cal M}$. There exists  $\delta>0$ such that the graph $\{(\yb, \Fb(\yb)), \yb \in \mathcal{M}\} $ lies in $ {\cal D}_{\delta}$, where~${\cal D}_{\delta} \coloneqq \{(\yb, \zb)\in \mathcal{M}\times \mathcal{M}: d(\text{\rm Cut}(\yb), \zb) \geq \delta\}$.
\end{assumption}


Verification of Assumptions~\ref{ass.costCurve} and \ref{Ass.FyCut} is difficult even for regular manifold such as hyperspheres. Substantial efforts to provide sufficient conditions have been made in the recent literature on measure transportation, including \cite{LoeperVillani2010}, \cite{loeper2011}, \cite{FKM2013}, \cite{Ge2021}, and several others.  {Three examples  of  practical interest  and for which these assumptions are satisfied are as follows (obviously, Examples~\ref{ex.sphere} and \ref{ex.torus} are special cases of Example~\ref{ex.product}).}

\begin{Example}\label{ex.sphere}
The unit hypersphere $\mathcal{S}^{p} \coloneqq  \{\yb \in {\mathbb R}^{p+1}: 
 \yb^\top \yb =1\}$ equipped with the geodesic distance $d_{\rm s}(\yb, \zb)\coloneqq  \cos^{-1}(\yb^\top \zb)$.
\end{Example}

\begin{Example}\label{ex.torus}
The torus~$\mathcal{T}^{p}$, which is the product $\mathcal{S}^{1}\!\times \ldots\times\! \mathcal{S}^{1}$ of  $p$ circles, equipped with the geodesic distance $${d_{\rm t}(\yb, \zb)\coloneqq \left(\sum_{i=1}^{p} d_{\rm s}^2(\yb_i, \zb_i)\right)^{1/2}}\coloneqq \left(\sum_{i=1}^{p} (\cos^{-1}(\yb_i^\top \zb_i))^2\right)^{1/2},$$ where $\yb = (\yb_1^\top\!, \ldots, \yb_{p}^\top)^\top\!\!$ (resp.,~$\zb = (\zb_1^\top\!, \ldots, \zb_{p}^\top)^\top\!\!$) with $\yb_i\! \in\! \mathcal{S}^{1}\!$  (resp.,~$\zb_i \!\in~\!\!\mathcal{S}^{1}\!$). For~$p=2$, this yields the familiar 2-torus $\mathcal{T}^2 = \mathcal{S}^1\times \mathcal{S}^1$  in ${\mathbb{R}^3}$.
\end{Example}

\begin{Example}\label{ex.product}
The polysphere $\mathcal{S}^{p_1} \times \cdots \times \mathcal{S}^{p_k}$ ($p_1 + \ldots + p_k = p$) equipped with the geodesic distance $${d_{\rm p}(\yb, \zb)\coloneqq \left(\sum_{i=1}^k d_{\rm s}^2(\yb_i, \zb_i)\right)^{1/2}}\coloneqq \left(\sum_{i=1}^k (\cos^{-1}(\yb_i^\top \zb_i))^2\right)^{1/2},$$ where $\yb = (\yb_1^\top\!, \ldots, \yb_k^\top)^\top\!\!$ (resp.,~$\zb = (\zb_1^\top\!, \ldots, \zb_k^\top)^\top\!\!$) with $\yb_i\! \in\! \mathcal{S}^{p_i}\!$  (resp.,~$\zb_i \!\in~\!\!\mathcal{S}^{p_i}\!$).
\end{Example}


Cylindrical product spaces of the form ${\cal M} = \mathcal{S}^{p_1} \times \mathbb{R}^{p_2}$ or ${\cal M} = \mathcal{S}^{p_1} \times {\Omega}$ with compact~$\Omega \subset \mathbb{R}^{p_2}$ also  have been studied in the literature;  according to Remark 3.10 in  \citet{KimMcCann2010}, Assumption~\ref{ass.costCurve} holds for such ${\cal M}$. However, to the best of our knowledge, the validity of Assumption~\ref{Ass.FyCut} has not been established; hence the continuity of the optimal transport is not guaranteed. 
Further examples, such as $M$ being {\it locally nearly spherical}  (in the sense of sectional curvature; see \cite{Ge2021}) or  a complex projective space, have not been considered from a statistical perspective, but their study is  beyond the scope of this paper.


In order to prove the continuity of $\Fb$, we make a very mild assumption on the cut locus, which is clearly satisfied by the manifolds listed in the above examples.

\begin{assumption}\label{ass.Cutyy0}
The Riemannian manifold $({\cal M}, g)$ is such  that for any $r>0$ and~$\yb, \yb_0 \in {\cal M}$ with $d(\yb, \yb_0) < r$, there exist 
$C>0$ such that $d(\xb, {\rm Cut}(\yb_0))< Cr$ for any $\xb \in {\rm Cut}(\yb)$.
\end{assumption}

Under the above regularity assumptions, we establish the continuity of $\Fb$, which is stated as follows; see Appendix~\ref{App.proof} for the proof.
\begin{proposition}\label{Prop.cont}
Let ${{\rm P}_1 }, {{\rm P}_2 } \in \mathfrak{P}$,  be such that ${{\rm P}_1 }$  has bounded density and  ${{\rm P}_2 }$ has nonvanishing density on ${\cal M}$. For~$c=d^2/2$, let Assumptions~\ref{ass.cost}--\ref{ass.Cutyy0} hold. 
 Then, the optimal transport map~$\Fb$ between ${{\rm P}_1 }$ and ${{\rm P}_2 }$ (see Proposition~\ref{keylemma}) is continuous. 
\end{proposition}
%
\subsection{Distribution and quantile functions}\label{sec.FbQb}

Proposition~\ref{keylemma} allows us to define the distribution function of a variable $\Yb \sim {\rm P}$ with values in ${\cal M}$ as the (${\rm P}$-a.s.\ unique) map that pushes ${\rm P}$ forward to some adequate reference probability measure ${{\rm P}_2 }$. Motivated by \cite{Hallin2021} and \cite{HLV2022}, we choose ${{\rm P}_2 }$ as  the uniform  distribution~${\rm U}_{\mathcal M}$ over $\mathcal M$.
  This choice is quite natural in view of the fact that, for $\Yb \sim {\rm P}$, $\Fb(\Yb)$ is uniformly distributed over ${\cal M}$, which is consistent with the fact that the classical distribution function $F$ of a real-valued r.v. $Y\sim{\rm P}$    is a probability integral transformation  pushing $\rm P$  forward to the uniform ${\rm U}_{[0,1]}$ over~$[0, 1]$. 
 \begin{definition} \label{defdistribution}
{\rm Call   {\it (Riemannian) distribution function of} ${\rm P} \in \mathfrak{P}$ (equivalently, of $\Yb\sim{\rm P}\in \mathfrak{P}$)  the ${\rm U}_{\mathcal M}$-a.s.~unique optimal transport map $\Fb$ from $ \mathcal{M}$ to $ \mathcal{M}$  such that~$\Fb\# {\rm P} ={\rm U}_{\mathcal M}$. }
  \end{definition}
We then have the following essential result (see Appendix~\ref{App.proof} for the proof).
%
%
%
%
%
%
%
%
%
%
\begin{proposition}\label{ProHom}  
{Let ${\rm P} \in \mathfrak{P}$ have bounded density on ${\cal M}$. For~$c=d^2/2$, let Assumptions~\ref{ass.cost}--\ref{ass.Cutyy0} hold. Then, the distribution  function $\Fb$ of $\Yb\sim{\rm P}$ is a homeomorphism from $\mathcal{M}^a$ to $\mathcal{M}^b$, where $\mathcal{M}^a$ and $\mathcal{M}^b$ are two compact subsets of $\mathcal{M}$ such that $\mathcal{M}\setminus \mathcal{M}^a$ is of ${\rm P}$-measure zero and $\mathcal{M}\setminus \mathcal{M}^b$ is of ${\rm U}_{\mathcal M}$-measure (hence also Riemannian volume measure) zero.}
\end{proposition}
It  follows  that the inverse $\Qb \coloneqq \Fb^{-1}$ of $\Fb$ is ${\rm U}_{\mathcal M}$-a.s. well defined and a homeomorphism too; as the inverse of a distribution function, $\Qb $  qualifies as a {\it quantile function}.

\begin{definition} \label{defQuantile}
{\rm  {Let ${\rm P}$ and $\mathcal{M}$ satisfy the conditions of Proposition~\ref{ProHom}}.\linebreak   Call~$\Qb \coloneqq~\!\Fb^{-1}$ (equivalently, the~${\rm U}_{\mathcal M}$-a.s.~unique optimal transport map $\Qb$ from $ \mathcal{M}$ to~$ \mathcal{M}$  such that $\Qb\# {\rm U}_{\mathcal M} ={\rm P}$) the {\it (Riemannian) quantile function} of ${\rm P} \in \mathfrak{P}$  (equivalently, of $\Yb\sim{\rm P}$).}
\end{definition}

\subsection{Quantile regions and  contours} Once a concept of quantile function $\bf Q$ is available, quantile regions and contours for $\rm P$ are naturally obtained as the $\bf Q$-images  of   quantile regions and contours associated with  ${\rm U}_{\mathcal M}$. This, {\it mutatis mutandis}, is the approach taken in \cite{chernoetal17}  and \cite{Hallin2021} for~${\mathbb R}^p$, and in  \cite{HLV2022} for the unit sphere ${\mathcal S}^p$. 

\subsubsection{Quantile regions and contours in ${\mathbb R}^p$}\label{subsec.Rp}
 In  \cite{Hallin2021},  the domain of the quantile functions $\Qb$ with range ${\mathbb R}^p$ is the open unit ball~${\mathbb S}^p$, with the spherical uniform~${\rm U}_p$ as the reference distribution ${\rm U}_{\mathcal M}$. This spherical uniform~${\rm U}_p$  is the product of a uniform over   the interval $[0,1)$ and a uniform over the hypersphere~${\mathcal S}^{p-1}$. That factorization   induces a collection 
  of nested closed balls of the form 
 $${\mathbb C}^{{\rm U}}_{ {\mathbb R}^p}(\tau)\coloneqq \tau\overline{\mathbb S}^p =
  \{t \vb \vert\, \vb\in{\mathcal S}^{p-1},\ 
   t\in[0,\tau]\}, \quad \tau\in[0,1)$$
(where $\overline{\mathbb S}^p$ stands for the closed unit ball)  with boundaries the $(p-1)$-spheres  
$${\mathcal C}^{{\rm U}}_{ {\mathbb R}^p}(\tau)\coloneqq \tau{\mathcal S}^{p-1} =
  \{\tau \vb \vert\, \vb\in{\mathcal S}^{p-1}\}, \quad \tau\in[0,1).$$
Note that $\{t\vb \vert t\in{\mathbb R}_+\}$ is the geodesic with direction $\vb$ starting from the origin. 

By construction, each ${\mathbb C}^{{\rm U}}_{ {\mathbb R}^p}(\tau)$  has ${\rm U}$-probability content $\tau$ and  enjoys the same   invariance properties as~${\rm U}$ under orthogonal transformations:  the collection over $\tau$ of these regions  
 therefore constitutes a natural family of center-outward quantile regions of order $\tau\in[0,1)$ for~${\rm U}_p$, with center-outward quantile contours~${\mathcal C}^{{\rm U}}_{ {\mathbb R}^p}(\tau)$. \linebreak 
The~$\Qb$-images~${\mathbb C}^{{\rm P}}_{ {\mathbb R}^p}(\tau)\coloneqq \Qb\big({\mathbb C}^{{\rm U}}_{ {\mathbb R}^p}(\tau)\big)$ and~${\mathcal C}^{{\rm P}}_{ {\mathbb R}^p}(\tau)\coloneqq \Qb\big({\mathcal C}^{{\rm U}}_{ {\mathbb R}^p}(\tau)\big)$  then constitute  the collection, centered at ${\mathbb C}^{{\rm P}}_{ {\mathbb R}^p}(0)\coloneqq \Qb (\{{\boldsymbol 0}\})$, of center-outward   quantile regions and contours of order $\tau\in [0,1)$ of~$\rm P$. We refer to \cite{Hallin2021} for details.

 
\subsubsection{Quantile regions and contours on the hypersphere ${\mathcal S}^p$}\label{232sec}  To define quantile regions on~${\mathcal M}=~\!{\mathcal S}^p$, \cite{HLV2022}  similarly factorize the uniform~${\rm U}_{\mathcal M}$ over the unit hypersphere into the product of a uniform on an interval $[0,\pi]$ of {\it latitudes} and a uniform ${\rm U}_{{\mathcal S}^{p-1}}$ over a $(p-1)$-sphere~${\mathcal S}^{p-1}$ of {\it longitudes} 
by selecting a {\it pole}  (a point in~${\mathcal S}^p$). That pole  characterizes a system of (hyper)meridians (geodesics starting from the pole) and parallels ($(p-1)$-spheres orthogonal to the meridians), among which an equatorial  $(p-1)$-sphere (the set of points with latitude $\pi/2$). \cite{HLV2022}  suggest selecting that pole as the~$\bf F$-image ${\bf F}(\thetab)={\bf F}(\thetab_{\rm Fr})$ of an element $\thetab_{\rm Fr}$ of the Fr\'echet mean set (see~\eqref{def.Frechet} for the definition) of $\rm P$, but other choices  ${\bf F}(\thetab)$,  $\thetab\in{\mathcal S}^p$ are possible. That factorization induces a collection 
 of nested {\it spherical caps} centered at ${\bf F}(\thetab )$, of the form 
\begin{equation}\label{sphericalpole}
{\mathbb C}^{\rm U, {\mathcal S}^p}_{{\bf F}(\thetab )}(\tau)\! \coloneqq\! \big\{\!\exp_{{\bf F}(\thetab )}(t\vb){\big\vert}\, { t\in[0,s(\tau)],}    \ \vb \in T_{{\bf F}(\thetab )}{\cal S}^{p}  \, \text{and}  \,  \Vert \vb \Vert = \pi\big\}
\end{equation}
where $\exp_{{\bf F}(\thetab)}(t\vb)$, $t\in[0,1]$ is a geodesic  starting from ${\bf F}(\thetab )$ and moving towards~$- {\bf F}(\thetab )$ while the mapping $s: [0,1] \rightarrow [0,1]$ is such that ${\mathbb C}^{\rm U, {\mathcal S}^p}_{{\bf F}(\thetab )}(\tau)$ has ${\rm U}_{\cal M}$-probability content~$\tau$, that is,
\begin{equation}\label{stauSphere}
	{\rm U}_{\cal M}\Big[{\mathbb C}^{\rm U, {\mathcal S}^p}_{{\bf F}(\thetab )}(\tau)\Big] =
		{\int_{0}^{s(\tau)} (\sin(\pi t))^{p-1} {\rm d}t }/{\int_{0}^1 (\sin(\pi t))^{p-1} {\rm d}t}  = \tau ;
	\end{equation}
note that $s(0)=0$ and $s(1)=1$. 	
 The boundary ${\mathcal C}^{\rm U, {\mathcal S}^p}_{{\bf F}(\thetab )}(\tau)$ of~${\mathbb C}^{\rm U, {\mathcal S}^p}_{{\bf F}(\thetab )}(\tau) $  is the  $(p-1)$ sphere with   latitude 
$s(\tau)\pi$ and~${\mathbb C}^{\rm U, {\mathcal S}^p}_{{\bf F}(\thetab )}(\tau)$ is the enclosed   ball. Since $\tau = 1/2$ yields $s(1/2) = 1/2$ in~\eqref{stauSphere} (latitude $\pi/2$),
 	 ${\mathcal C}^{\rm U, {\mathcal S}^p}_{{\bf F}(\thetab )}(1/2)$ is the equatorial $(p-1)$-sphere associated with the chosen pole ${\bf F}(\thetab )$. 


An equally valid alternative (which is not formally covered in \cite{HLV2022}), involving the same 
$(p-1)$-spheres ${\mathcal C}^{\rm U, {\mathcal S}^p}_{{\bf F}(\thetab )}(\tau)$, would be a collection of nested {\it equatorial strips} $ {\mathbb C}^{\rm U, {\mathcal S}^p}_{{\bf F}(\thetab )}((1+\tau)/2)\setminus {\mathbb C}^{\rm U, {\mathcal S}^p}_{{\bf F}(\thetab )}((1-\tau)/2)$ centered at the equator~{${\mathcal C}^{\rm U, {\mathcal S}^p}
 _{{\bf F}(\thetab )}(1/2)
$} (a submanifold ${{\cal M}_0}$ of ${\cal M}={\cal S}^p$) rather than  the pole and limited by the~$(p-1)$-spheres~{${\mathcal C}^{\rm U, {\mathcal S}^p}_{{\bf F}(\thetab )}((1+\tau)/2)$ and~${\mathcal C}^{\rm U, {\mathcal S}^p}_{{\bf F}(\thetab )}((1-\tau)/2)$}, respectively. These equa\-torial strips also can be defined as 
\begin{align}\nonumber
{\mathbb C}^{\rm U\updownarrows}_{{\cal M}_0}(\tau)\coloneqq 
\big\{\exp_{\yb}(t\vb){\big\vert}\, t\in [0, 1 - 2s((1-\tau)/2)],&\, \vb \in T_{\yb}{\cal S}^{p},   \vb 
\perp T_{\yb}{\mathcal C}
^{{\rm U}, {\mathcal S}^p}_{{\bf F}(\thetab )}(1/2)  \\ & \, \text{and}  \,  \Vert \vb \Vert = \pi/2
\ \text{for} \ \yb \in {\mathcal C}^{{\rm U}, {\mathcal S}^p}_{{\bf F}(\thetab )}(1/2)
\big\}\label{sphericalstrip}
\end{align}
where  $\exp_{\yb}(t\vb)$ is a geodesic curve originating in  $\yb\in{{\cal M}_0}={\mathcal C}^{\rm U, {\mathcal S}^p}_{{\bf F}(\thetab )}(1/2)$  orthogonally to~${\mathcal C}^{\rm U, {\mathcal S}^p}_{{\bf F}(\thetab )}(1/2)$. 


This choice between pole-centered caps and equator-centered strips as quantile regions rests with the analyst and depends on his/her view on the data-generating process and/or the objectives of the analysis. Once that choice is made, the quantile regions of order $\tau$ for~$\rm P$   are the $\bf Q$-images 
\begin{equation}\label{capsdef}{\mathbb C}^{\rm P,  {\mathcal S}^p}_{\thetab }(\tau)\coloneqq {\bf Q}\left({\mathbb C}^{\rm U,  {\mathcal S}^p}_{{\bf F}(\thetab )}(\tau)\right),\qquad \tau\in[0,1]
\end{equation}
or
\begin{equation}\label{stripsdef} {\mathbb C}_{\thetab }^{\rm P{ \updownarrows}}(\tau)\coloneqq {\bf Q}\left({\mathbb C}^{\rm U\updownarrows}_{{\cal M}_0}(\tau)
\right)
={\bf Q}\left( {\mathbb C}^{\rm U,  {\mathcal S}^p}_{{\bf F}(\thetab )}\left(\frac{1+\tau}{2}\right) \setminus {\mathbb C}^{\rm U,  {\mathcal S}^p}_{{\bf F}(\thetab )}\left(\frac{1-\tau}{2}\right)\right)
,\qquad \tau\in[0,1]
\end{equation} 
of the uniform ones; their boundaries 
$${\mathcal C}^{\rm P,  {\mathcal S}^p}_{\thetab }(\tau)\coloneqq {\bf Q}\big({\mathcal C}^{\rm U,  {\mathcal S}^p}_{{\bf F}(\thetab )}(\tau)\big)\quad\text{ and }\quad {\mathcal C}^{\rm P  \updownarrows}_{\thetab}(\tau)\coloneqq \Qb\Big({\mathcal C}^{\rm U,  {\mathcal S}^p}_{{\bf F}(\thetab )}({(1-\tau})/{2}) \bigcup {\mathcal C}^{\rm U,  {\mathcal S}^p}_{{\bf F}(\thetab )} (({1+\tau})/{2})
\Big)$$ are the corresponding quantile contours. {Note that here, unlike in Section~\ref{subsec.Rp}, quantile regions of order $\tau = 1$ are well defined thanks to the fact that $\rm P$ is compactly supported. 
}

We now proceed with extending this approach to more general manifolds.

\subsubsection{Quantile regions and contours on ${\cal M}$}

 Throughout this section, we consider the general case of a manifold $\cal M$ with distribution and quantile functions $\Fb$ and $\Qb$, respectively,  satisfying the assumptions of Proposition~\ref{ProHom}. In order to define quantile regions and  contours on~${\cal M}$,  we need, as in $\mathbb{R}^d$ and ${\mathcal S}^p$,   a pre-specified central region---a submanifold   ${{\cal M}_0} \subset {\cal M}$ (equipped with the   Riemannian metric inherited from $\cal M$) with~${\rm U}_{\cal M}$-probability zero (hence,~${{\cal M}_0}= \partial {{\cal M}_0}$)  and carrying a {\it normal vector field}. 
 
Recall that  a {\it normal vector set of ${\cal M}_0$} at $\yb\in{\cal M}_0$ is a set  ${V}^{\perp}_{\yb}\coloneqq \{\vb_{\yb}\}$ of  vectors~$\vb_{\yb}\in~\!T_{\yb}{\cal M}$ such that $\vb_{\yb} \perp T_{\yb}  {{\cal M}_0}$. 
  A vector field~$V^{\perp}_{{\cal M}_0} \coloneqq \bigcup_{\yb \in {\cal M}_0} V^{\perp}_{\yb}$ is a {\it normal vector field} of ${{\cal M}_0}$ if  ${V}^{\perp}_{\yb}$ is a vector normal set at each~$\yb \in~\!{{\cal M}_0}$.  In case ${{\cal M}_0}$ is a singleton~$\{\yb\}$, any set of vectors  in $T_{\yb}{\cal M}$ is a normal vector field of ${{\cal M}_0}$.\smallskip

Let us make the following assumption on the existence and properties of~${V}^{\perp}_{{\cal M}_0}$.
 
\begin{assumption}\label{ass.Wincrease}
The submanifold ${\cal M}_0$ admits a normal vector field  ${V}^{\perp}_{{\cal M}_0}$  such that the regions 
 \begin{equation}\label{defgenregion} 
 {\mathbb C}^{\rm U}_{{\cal M}_0}(\tau) \coloneqq \left\{\gamma^{\vb}_\yb (t)\coloneqq \exp_\yb (t\vb)\vert t\in [0, s(\tau)], \yb\in{{\cal M}_0}, \vb\in V_\yb^\perp \subseteq {V}^{\perp}_{{\cal M}_0}
  \right\}, \quad \tau\in [0,1]
  \end{equation}  
  satisfy ${\rm U}_{\cal M}\left( {\mathbb C}^{\rm U}_{{\cal M}_0}(\tau)\right) =\tau$ for some strictly increasing $s:\, [0,1]\rightarrow  [0,1]$.\footnote{This also holds for ${\mathcal M}={\mathbb R}^p$, which is not compact, with, however, $s:\, [0,1)\rightarrow  [0,1)$. }
\end{assumption}

Under Assumption~\ref{ass.Wincrease}, the regions $ {\mathbb C}^{\rm U}_{{\cal M}_0}(\tau)$, $\tau\in [0, 1]$ defined in \eqref{defgenregion} are closed, connected, and nested, with ${\rm U}_{\cal M}$-probability content $\tau$, hence constitute a natural collection of quantile regions for ${\rm U}_{\cal M}$, with quantile contours 
  \begin{equation}\label{defgencontour} 
 {\mathcal C}^{\rm U}_{{\cal M}_0}(\tau) \coloneqq \left\{\gamma^{\vb}_\yb (s(\tau ))=\exp_\yb (s(\tau)\vb) \vert 
  \yb\in{{\cal M}_0}, \vb\in V_\yb^\perp \subseteq {V}^{\perp}_{{\cal M}_0}
  \right\}\quad \tau\in [0,1].
  \end{equation} 
Moreover,  Assumption~\ref{ass.Wincrease} implies $s(0)=0$,   $s(1)=1$, and~${\mathbb C}^{\rm U}_{{\cal M}_0}(\tau =1) = \cal M$.
  
 {\it  Equatorial strips,}  of the form
  \begin{equation*}
{{\mathbb C}^{\rm U\updownarrows}_{{\cal M}_0}(\tau)} \!\coloneqq\! \left\{\gamma^{\vb}_\yb (t) \vert\, s((1-\tau)/2)\leq t\leq s((\tau +1)/2), \yb\in{{\cal M}_0}, \vb\in V_\yb^\perp \subseteq {V}^{\perp}_{{\cal M}_0}
  \right\},\  \tau\in [0,1]
  \end{equation*} 
   also can be constructed; they are centered about ${\mathcal C}^{\rm U}_{{\cal M}_0}(1/2)$, which then plays the role of an equatorial submanifold of ${\cal M}$;    their boundaries are 
\begin{equation}\label{contoursupdown}
{{\mathcal C}^{\rm U\updownarrows}_{{\cal M}_0}(\tau)}\coloneqq {\mathcal C}^{\rm U}_{{\cal M}_0}((1+\tau )/2) \bigcup {\mathcal C}^{\rm U}_{{\cal M}_0}((1-\tau )/2),\quad \tau\in [0,1].
\end{equation}

 Finally, the  quantile regions ${{\mathbb C}^{\rm P}_{{\cal M}_0}(\tau)} $ and~${{\mathbb C}^{\rm P\updownarrows}_{{\cal M}_0}(\tau)} $  and the quantile  contours  ${\mathcal C}^{\rm P}_{{\cal M}_0}(\tau)$ and ${{\mathcal C}^{\rm U\updownarrows}_{{\cal M}_0}(\tau)}$ 
 of~$\rm P$ are obtained as the   $\Qb$-images of their~${\rm U}_{\cal M}$ counterparts.
 
\begin{definition}\label{def.Ctau}
	{\rm 
		Call 
		$\mathbb{C}^{\rm P}_{{\mathcal M}_0}(\tau) \coloneqq  {\Qb}(\mathbb{C}_{{\cal M}_0}^{\rm U}(\tau))$ and $\mathcal{C}^{\rm P}_{{\mathcal M}_0}(\tau)  \coloneqq {\Qb}(\mathcal{C}_{{\cal M}_0}^{\rm U}(\tau))
		$
		   {\it (Riemannian) quantile region} and {\it quantile contour},  respectively, {\it of order}~$\tau \in [0, 1]$ of $\rm P$; similarly call ${{\mathbb C}^{\rm P\updownarrows}_{{\cal M}_0}(\tau)} \coloneqq \Qb( {{\mathbb C}^{\rm U\updownarrows}_{{\cal M}_0}(\tau)} )$ a {\it (Riemannian) equatorial strip} of order $\tau$,  and   ${\Qb}({\cal M}_0)$ a {\it (Riemannian) median} of ${\rm P}$}. 
\end{definition}
 
 \begin{Rem}\label{Rem.suff} 
 		{\rm
 			Let  ${\cal M}_0$ admit a normal vector field  ${V}^{\perp}_{{\cal M}_0}$ such that  
 			\begin{equation}\label{eq.C*}
 			{\mathcal C}^{\rm U}_* \coloneqq \left\{\gamma^{\vb}_\yb (s_*)=\exp_\yb (s_*\vb) \vert 
 			\yb\in{{\cal M}_0}, \vb\in V_\yb^\perp \subseteq {V}^{\perp}_{{\cal M}_0}
 			\right\}
 			\end{equation}
 			is  a submanifold of dimension $(p-1)$ for some $s_*\in[0,1]$.  A sufficient condition for Assumption~\ref{ass.Wincrease} to hold is that a uniform variable  $\Ub\sim~\!{\rm U}_{\cal M}$ can be written under the form~$\Ub = \exp_{\Yb_{\Ub}}(w_{\Ub}\vb_{\Yb_{\Ub}})$, where $w_{\Ub}\sim {\rm U}_{[0,1]}$ and $\Yb_{\Ub}\sim {\rm U}_{{\mathcal C}^{\rm U}_*}$ are mutually independent and $\vb_{\Yb_{\Ub}} \in V_{\Yb_{\Ub}}^\perp \subseteq {V}^{\perp}_{{\mathcal C}^{\rm U}_*}$ for some normal vector field ${V}^{\perp}_{{\mathcal C}^{\rm U}_*}$  determined by~${V}^{\perp}_{{\cal M}_0}$. 
 		} 
 		{\rm This condition is satisfied for~${\cal M}=~\!{\mathbb R}^p$ (Section~\ref{subsec.Rp}) with~${{\cal M}_0}=\{{\boldsymbol 0}\}$, ${ V}^{\perp}_{{\cal M}_0}$ the~$(p-1)$-sphere~${\mathcal S}^{p-1}$ of the radii of the unit ball,~$s(t)=t$, $t\in [0,1)$, and $ {\mathcal C}^{\rm U}_{*} = {s_*\mathcal S}^p$\linebreak  for any~$s_*\in (0,1)$.  It is satisfied also for~${\cal M}={\mathcal S}^p$ (Section~\ref{232sec}) with   
 			\begin{enumerate}
 				\item[(a)] either ${{\cal M}_0} = \{\Fb (\thetab)\}$ for some pole $\thetab \in {\cal S}^{p}$,  ${V}^{\perp}_{{\cal M}_0} =  \big\{\!\vb \in T_{{\bf F}(\thetab )}{\cal S}^{p}  \,\big\vert \, \Vert \vb \Vert = \pi \big\}$,
 			$s(t)$ defined in~\eqref{stauSphere}, and ${\mathcal C}^{\rm U}_{*}$ 
 			the equatorial $(p-1)$-sphere associated with $\Fb (\thetab)$, yielding the pole-centered spheri\-cal caps~\eqref{sphericalpole},
 			\item[(b)] or  ${\cal M}_0$ 
 			the equatorial $(p-1)$-sphere associated with some   $\Fb (\thetab)$, 		   $${V}^{\perp}_{{\cal M}_0} =  \big\{\!\vb \in T_{\yb}{\cal S}^{p} \big\vert\, \vb \perp T_{\yb} {\cal M}_0 ,\,  \Vert \vb \Vert = \pi/2, \ \text{and} \ \yb \in {\cal M}_0\big\},$$
 			with	    $s(t) = 1 - 2\tilde{s}((1-t)/2)$ where $\tilde{s}$ satisfies \eqref{stauSphere}, ${\mathcal C}^{\rm U}_{*}={\mathcal M}_0$, and $s_* = 0$,  
 			yielding the equatorial strips  \eqref{sphericalstrip}.
 		\end{enumerate}
 }\end{Rem}

\begin{Rem}
	{\rm 		  
	  One may choose ${\mathcal M}_0$ and $V^\perp_{{\mathcal M}_0}$ in Assumption~\ref{ass.Wincrease} such that $\mathbb{C}_{{{\cal M}_0}}^{\rm U}(\tau)$ enjoys, as much as possible, the symmetry properties  of ${\rm U}_{\mathcal M}$.  For instance, the normal vector fields  in Sections~\ref{subsec.Rp} and \ref{232sec} are chosen to be rotationally symmetric over~${{\cal M}_0}$, hence maintaining the rotational symmetries of the uniform distributions over~$\mathbb{S}^p $ and~$\mathcal{S}^p$, respectively; see  the case of the torus  in  Section~\ref{233sec} for a further example.}
	  \end{Rem}


\begin{Rem}
	{\rm 
The choice of ${\cal M}_0$ depends on the analyst and determines the topo\-logical properties of the quantile regions (see Lemma~\ref{lem.homotopy} below). One possible choice (cf. Remark~\ref{Rem.suff}(a)) for~${\cal M}={\mathcal S}^p$   is~${{\cal M}_0} = \{ \Fb(\thetab_{\rm Fr})\}$, where $\thetab_{\rm Fr}$  is a {\it  Fr\'echet mean}, that is, belongs to the   {\it  Fr\'echet mean set}  
\begin{equation}\label{def.Frechet}
	{\cal A}_{\rm Fr}\coloneqq  {\rm argmin}_{\yb \in  \mathcal{M}} {\rm E}_{{\rm P}}[c(\Yb, \yb)].
\end{equation}
of~$\Yb\!\sim~\!{\rm P}$ (${\cal A}_{\rm Fr}$  is non-empty since ${\cal M}$ is compact). The ${\rm U}_{{\mathcal M}}$-quantile regions~$\mathbb{C}_{{{\cal M}_0}}^{\rm U}(\tau)$ 
 can be defined  from the set of geodesics connec\-ting~$\Fb(\thetab_{\rm Fr})$ and $\text{\rm Cut}(\Fb(\thetab_{\rm Fr}))$\footnote{see Appendix~A for the definition.} (playing the role of extremes of ${\rm U}_{\cal M}$), namely,\vspace{2mm}
 
 $\mathbb{C}_{{{\cal M}_0}}^{\rm U}(\tau) \coloneqq \Big\{\gamma^{\vb}_{\Fb(\thetab_{\rm Fr})}(t) \! \coloneqq\!  \exp_{\Fb(\thetab_{\rm Fr})}(t\vb){\big\vert}\, t \in [0, \tau], \ \vb \in T_{\Fb(\thetab_{\rm Fr})}{\cal M}$\vspace{-4.5mm}

\begin{equation}\label{def.Vs}
\hspace{55mm}\text{such that} \, \exp_{\Fb(\thetab_{\rm Fr})}(\vb) \in \text{\rm Cut}(\Fb(\thetab_{\rm Fr})) \Big\}. 
\end{equation}
Assumption~\ref{ass.Wincrease} then follows from the Hopf-Rinow Theorem (see, e.g., Chapter~5 in \cite{petersen2006riemannian}). {The resulting quantile regions $\mathbb{C}_{{{\cal M}_0}}^{\rm P}(\tau)={\Qb}(\mathbb{C}_{{{\cal M}_0}}^{\rm U}(\tau))$  of $\rm P$  are centered about the pole~$\thetab_{\rm Fr}$.}
}
\end{Rem}


Lemma~\ref{lem.homotopy} states that ${{\cal M}_0}$ and  $\mathbb{C}_{{{\cal M}_0}}^{\rm U}(\tau )$ are topologically similar. \vspace{-2mm} 
 
\begin{Lem}\label{lem.homotopy}
	Under Assumption~\ref{ass.Wincrease},
	\begin{compactenum}
		\item[(i)] for any $\tau \in [0, 1]$, ${\cal M}_0 = \partial{{\cal M}_0}$ and $\mathbb{C}_{{{\cal M}_0}}^{\rm U}(\tau )$ 
		  are homotopy equivalent;
		\item[(ii)] 
		for any $\tau_1, \tau_2 \in (0, 1)$, $\mathbb{C}_{{\cal M}_0}^{\rm U}(\tau_1)$ is  diffeomorphic to $\mathbb{C}_{{\cal M}_0}^{\rm U}(\tau_2)$.
	\end{compactenum}
\end{Lem}
See Appendix~\ref{App.proof} for the proof.\smallskip

Proposition~\ref{prop.FQC} summarizes the main properties of $\Fb$, $\Qb$, $\mathbb{C}_{{\mathcal M}_0}^{\rm P}(\tau)
$, 
 and~$\mathcal{C}^{\rm P}_{{\mathcal M}_0}(\tau)
$, which follow  from 
 Proposition~\ref{ProHom} and Lemma~\ref{lem.homotopy}. Details are left to the reader.

\begin{proposition} \label{prop.FQC}
	{Let $\Yb\sim {\rm P}$ where ${\rm P} \in \mathfrak{P}$, with distribution and quantile functions~$\Fb$ and~$\Qb$, respectively, has bounded density on ${\cal M}$. For~$c=d^2/2$, let Assumptions~\ref{ass.cost}--\ref{ass.Cutyy0} hold}. 
	Then, for any $\tau \in [0, 1]$,
	\begin{compactenum}
		\item[(i)] $\Fb$ entirely characterizes ${\rm P}$ and  $\Fb(\Yb)\sim\Fb\#{\rm P} =~\!{\rm U}_{\cal M}$;
		\item[(ii)] $\Qb$ entirely characterizes ${\rm P}$ and $\Qb(\Ub)\sim\Qb\#{\rm U}_{\cal M}=~\!{\rm P}$  for $\Ub\sim {\rm U}_{\cal M}$;
\item[(iii)] under Assumption~\ref{ass.Wincrease}, for any $\tau \in [0, 1]$,
      \begin{compactenum}
	 \item[(a)] 	irrespective of $\rm P \in  \mathfrak{P}$, the $\rm P$-probability  ${\rm P}\!\left({\mathbb C}_{{\mathcal M}_0}^{\rm P}(\tau )
		\right)$ of ${\mathbb C}_{{\mathcal M}_0}^{\rm P}(\tau )$ and the $\rm P$-proba\-bility~${\rm P}\!\big({{\mathbb C}^{\rm P\updownarrows}_{{\cal M}_0}(\tau)} 
		\big)$ of ${\mathbb C}^{\rm P\updownarrows}_{{\cal M}_0}(\tau)$ are equal to $\tau$;
		\item[(b)] the quantile regions~${\mathbb C}_{{\mathcal M}_0}^{\rm P}(\tau)$ are closed, and ${\mathbb C}_{{\mathcal M}_0}^{\rm P}(\tau)\bigcap {\cal M}^a$ is homotopy equivalent to $\Qb({{\cal M}_0})\bigcap {\cal M}^a$ with $\mathcal{M}^a\subset\mathcal{M}$ of ${\rm P}$-measure one   defined in Proposition~\ref{ProHom};
		\item[(c)] for any $\tau_1, \tau_2 \in (0, 1)$, $\mathbb{C}_{{\mathcal M}_0}^{\rm P}(\tau_1)\cap {\cal M}^a$ is  homeomorphic to $\mathbb{C}_{{\mathcal M}_0}^{\rm P}(\tau_2)\cap {\cal M}^a$.
	 \end{compactenum}
	\end{compactenum}
\end{proposition}

Next we illustrate how    quantile regions 
 ${\mathbb C}_{{\mathcal M}_0}^{\rm P}(\tau)$, ${\mathbb C}_{{\mathcal M}_0}^{\rm P\updownarrows}(\tau)$, and the contours ${\mathcal C}_{{\mathcal M}_0}^{\rm P}(\tau)$ look like on the torus and the polysphere. 

\subsubsection{Quantile regions and contours on the 2-torus $\mathcal{T}^{2}\!$}\label{233sec}  The quantile regions 
in~${\mathbb R}^d$ and the spherical-cap-based ones \eqref{capsdef} 
 in ${\mathcal S}^p$ are homotopic to one point, while the equator-centered strips \eqref{stripsdef} 
 are homotopic to the equator $\Qb ({\mathcal M}_0) = {\mathcal C}_\thetab ^{{\rm P},{\mathcal S}^p}(1/2)$.  The homotopy structure of  the torus is more complex.  To start with, consider the 2-torus~$\mathcal{M}=\mathcal{T}^{2}= \mathcal{S}^{1}_1\!\times\mathcal{S}^{1}_2$. Some connected regions there are homotopic to one point, some others to $\mathcal{S}^{1}_1$, and  some others to~$\mathcal{S}^{1}_2$.  Letting  $\Fb(\thetab) = (t_1, t_2)^\top$,  $t_i\in {\mathcal S}^1_i$, $i=1,2$, meaningful quantile regions can be constructed for each of these three cases.\smallskip

 \begin{compactenum} 
 \item[(a)] Choosing ${{\cal M}_0} = \{\Fb(\thetab
 )\}$, ${V}^{\perp}_{{\cal M}_0} = \big(\{\pm\pi\}\times [-\pi, \pi)\big) \bigcup \big([-\pi, \pi)\times \{\pm\pi\}\big)$ 
(i.e. the boundary of a square centered at $\Fb (\thetab)$ in $T_{\Fb(\thetab )}\mathcal{T}^{2}$)  and $s({\tau}) = \tau^{1/2}$ yields a collection of nested squares 
 \begin{align*}
 \mathbb{C}_{{{\cal M}_0}}^{\rm U}({\tau}) &= \left\{\exp_{\Fb(\thetab)} (t\vb)\vert t\in [0, s(\tau)], \vb\in {V}^{\perp}_{{\cal M}_0}\right\} \\
 &= \big[t_1-{\pi\tau^{1/2}} , t_1 +\pi\tau^{1/2} \big]\times \big[ t_2- \pi\tau^{1/2}, t_2 + {\pi\tau^{1/2}}\big]
 \end{align*}
 centered at and homotopic to~$\Fb(\thetab)$,   with ${\rm U}_{\mathcal M}$-probability content $\tau\in[0,1]$. Call  $\thetab$ (with, e.g., $\thetab = \thetab_{\rm Fr}$) a {\it toroidal pole} and the associated regions $\mathbb{C}_{{{\cal M}_0}}^{\rm U}({\tau})$ {\it toroidal caps}. This choice leads to  quantile regions of the form~$\mathbb{C}_{{\mathcal M}_0}^{\rm P}(\tau) \coloneqq \Qb(\mathbb{C}_{{{\cal M}_0}}^{\rm U}({\tau}))$ that, similar to the quantile regions considered in Section~\ref{232sec} for~${\mathcal S}^p$, are  homotopic to a point---the pole~$\thetab$.
 
 \item[(b)] Letting ${{\cal M}_0} = \{t_1\}\times\mathcal{S}^{1}_2$, ${V}^{\perp}_{{\cal M}_0} = \{\vb_{\yb}{\big\vert}\,  \vb_{\yb} \in \{(\pi, 0)^\top, (-\pi, 0)^\top\}, \yb \in {{\cal M}_0}\}$, and $s({\tau}) = \tau$ yields a family of nested regions 
$\mathbb{C}_{{{\cal M}_0}}^{{\rm U}\updownarrows}({\tau}) = [t_1 - {\pi\tau}, t_1+ {\pi\tau}]\times  \mathcal{S}^{1}_2
$ 
centered at and homotopic to~$\{t_1\}\times\mathcal{S}^{1}_2$ (playing for ${\rm U}_{\mathcal M}$ the role of a {\it toroidal equator}), with ${\rm U}_{\mathcal M}$-probability content $\tau\in[0, 1]$; call them {\it toroidal equatorial  strips}.   The corresponding family of quantile regions $\mathbb{C}_{{\mathcal M}_0}^{\rm P \updownarrows}(\tau) \coloneqq  \Qb(\mathbb{C}_{{{\cal M}_0}}^{\rm U \updownarrows}({\tau}))$ are homotopic to 
${\bf Q}(\{t_1\}\times \mathcal{S}^{1}_2)$, and have the same interpretation as the quantile regions  described in Section~\ref{232sec} for~${\mathcal S}^p$.
 \item[(b')]
 Of course, $\mathcal{S}^{1}_2$ and $\mathcal{S}^{1}_1$, in the 2-torus $\mathcal{T}^{2}$, play symmetric roles;  exchanging them   yields (with obvious notation) quantile regions  of the form 
  $${\mathbb C}_{{\mathcal M}_0}^{\rm P \updownarrows}(\tau)\coloneqq \Qb(\mathcal{S}^{1}_1 \times  [t_2 - {\pi\tau}, t_2+ {\pi\tau}]),\quad \tau\in [0,1]$$ 
which are homotopic to~${\bf Q}(\mathcal{S}^{1}_1\times\{t_2\})$.
 \end{compactenum}

\subsubsection{Quantile regions and contours on the p-torus $\mathcal{T}^{p}$}\label{ptorus}  More cases are possible in  the $p$-torus~$\mathcal{T}^{p}=\mathcal{S}^{1}_1\!\times \ldots\times\! \mathcal{S}^{1}_p$ 
 with $p\geq 3$.
\smallskip

 \begin{compactenum} 
 \item[(a)] Letting ${{\cal M}_0} = \Fb(\thetab) = (t_1, \ldots, t_p)^\top$ (e.g., with~$\thetab = \thetab_{\rm Fr}$), 
 $${V}^{\perp}_{{\cal M}_0} = \{(v_1, \ldots, v_p)^\top{\big\vert}\, v_k = \pm\pi, v_j \in [-\pi, \pi] \,\, \text{for} \,\,  j \neq k, k = 1, \ldots, p\}$$ (i.e. the boundary of a hypercube centered at $\Fb(\thetab)$ in $T_{\Fb(\thetab)}\mathcal{T}^{p}$), and $s({\tau}) = \tau^{1/p}$
  yields a collection of nested regions (call them toroidal caps)  
 $$\mathbb{C}_{{{\cal M}_0}}^{\rm U}({\tau}) \coloneqq  \bigtimes_{i=1}^p
 \big[t_i-\pi\tau^{1/p}, t_i + \pi\tau^{1/p}\big]$$
centered at and homotopic to~$\Fb(\thetab)$ (which plays the role of a {\it toroidal pole} for ${\rm U}_{{\mathcal M}}$),   with~${\rm U}_{\mathcal M}$-probability content $\tau\in[0,1]$. These  toroidal caps in turn define a collection of quantile regions $\mathbb{C}_{{\mathcal M}_0}^{\rm P}(\tau) \coloneqq \Qb(\mathbb{C}_{{{\cal M}_0}}^{\rm U}({\tau}))$ with $\rm P$-probability content $\tau$ which enjoy the same interpretation as the pole-centered quantile regions in Section~\ref{233sec} (a).\smallskip

 \item[(b)] Noting that, for each $k\in\{1,\ldots,p\}$, the $p$-torus $\mathcal{T}^{p}=\mathcal{S}^{1}_1\!\times \ldots\times\! \mathcal{S}^{1}_p$  with $p\geq 3$ can be rewritten~as 
 $\mathcal{T}^{p} 
 = \mathcal{S}^{1}_k\!\times\mathcal{S}^{1}_1\!\times \cdots\times\! \mathcal{S}^{1}_{k-1}\!\times\mathcal{S}^{1}_{k+1}\!\times\cdots\times\! \mathcal{S}^{1}_p\eqqcolon \mathcal{S}^{1}_k\!\times\mathcal{T}^{p-1}_{\lceil k}$ and letting~${{\cal M}_0} = \{t_k\}\times\mathcal{T}^{p-1}_{\lceil k}$ with 
  $${V}^{\perp}_{{\cal M}_0} = \{(v_1, \ldots, v_p)^\top{\big\vert}\, v_k = \pm\pi, v_j = 0 \,\, \text{for} \,\,  j \neq k, k = 1, \ldots, p\}$$  
  yields, for  $s({\tau}) = \tau$, a family of nested regions
  $$\mathbb{C}_{{{\cal M}_0}}^{\rm U \updownarrows}({\tau}) = [t_k - {\pi\tau}, t_k+ {\pi\tau}]\times \mathcal{T}^{p-1}_{\lceil k},\quad \tau\in[0,1]$$
  centered at and homotopic to~$\{t_k\}\times\mathcal{T}^{p-1}_{\lceil k}$ (playing the role of a {\it toroidal equator}), with~${\rm U}_{\mathcal M}$-probability content $\tau$: call them {\it$(p-1)$-toroidal equatorial  strips}. The quantile regions $\mathbb{C}_{{\mathcal M}_0}^{\rm P \updownarrows}(\tau) \coloneqq  \Qb(\mathbb{C}_{{{\cal M}_0}}^{\rm U \updownarrows}({\tau}))$ and their boundaries, have the same interpretation 
   and enjoy the same properties as the family of nested {\it equatorial strips} centered at the equator and their~$(p-1)$-hypersphere boundaries described in Section~\ref{232sec}.

 \item[(c)] Note that, for $\mathcal{M} = \mathcal{T}^p$ with $p\geq 3$, ``mixed''  cases involving   poles {\it and }equators, also are possible---choosing, for instance,   
  ${{\cal M}_0} = \{t_{i_1}\}\times \cdots \times \{t_{i_k}\} \times \mathcal{T}^{p-k}_{\lceil i_1, \ldots, i_k}$, with~$\{i_1, \ldots, i_k\} \subseteq \{1, \ldots, p\}$ and $$\mathcal{T}^{p-k}_{\lceil i_1, \ldots, i_k} \coloneqq \mathcal{S}^{1}_{j_{1}} \times \cdots \mathcal{S}^{1}_{j_{p-k}},\quad  \{{j_{1}}, \ldots, j_{p-k}\} = \{1, \ldots, p\}\setminus \{i_1, \ldots, i_k\}.$$
Details are left to the reader.
   \end{compactenum}
  
  
\subsubsection{Quantile regions and contours on the polysphere $\boldsymbol{\mathcal{S}}^{p_1\ldots p_k}$}\label{polysphere}  The case of the polysphere follows along similar lines. For any $j\!\in\!\{1,\ldots, k\}$,  $\boldsymbol{\mathcal{S}}^{p_1\ldots p_k}\!\coloneqq{\mathcal{S}}^{p_1} \times \cdots \times~\!{\mathcal{S}}^{p_k}$ rewrites as 
 $\boldsymbol{\mathcal{S}}^{p_1\ldots p_k}
= \mathcal{S}^{p_j} \times  \mathcal{S}^{p_1}\times \cdots \times  \mathcal{S}^{p_{j-1}}\times  \mathcal{S}^{p_{j+1}}\times \cdots\times \mathcal{S}^{p_k}\eqqcolon  \mathcal{S}^{p_j} \times \boldsymbol{\mathcal{S}}^{p_1\ldots  p_k}_{\lceil p_j}.$ 

\begin{compactenum}
 \item[(a)] Let $\boldsymbol{s}_{j}$, $j=1,\ldots ,k$ denote the $ \mathcal{S}^{p_j}$ components of $\Fb(\thetab )$ (with, e.g., $\thetab=\thetab_{\rm Fr}$). Choosing 
 ${{\cal M}_0}$ as the  sub-polysphere ${{\cal M}_0}\coloneqq \boldsymbol{s}_{j} \times  \boldsymbol{\mathcal{S}}^{p_1\ldots  p_k}_{\lceil p_j}$, 
 $${V}^{\perp}_{{\cal M}_0} = \{(\boldsymbol{0}^\top, \ldots, \boldsymbol{0}^\top, (\vb^{{j}})^\top, \boldsymbol{0}^\top, \ldots, \boldsymbol{0}^\top)^\top\big\vert \vb^{j} \in T_{\boldsymbol{s}_{{j}}}{\cal S}^{{j}} \ \text{and} \ \Vert\vb^{j}\Vert = \pi\},$$ 
and denoting by~$s({\tau})$  the solution to \eqref{stauSphere} yields the nested  {\it polyspherical caps}
 \begin{align}\nonumber 
 {\mathbb C}^{\rm U}_{{\cal M}_0}({\tau}) &= \{\exp_{\boldsymbol{s}_j}(\tau\vb){\big\vert}\, \vb \in {V}^{\perp}_{{\cal M}_0} \} \\ 
 &=\{\exp_{\boldsymbol{s}_{{j}}}(\tau\vb^j) {\big\vert}\,  \vb^{j} \in T_{\boldsymbol{s}_{{j}}}{\cal S}^{{j}}  \ \text{and} \ \Vert\vb^{j}\Vert = \pi\}   
 \times  \boldsymbol{\mathcal{S}}^{p_1\ldots  p_k}_{\lceil p_j}, \label{polycapsdef}
 \end{align}
 centered at $\boldsymbol{s}_{j} \times  \boldsymbol{\mathcal{S}}^{p_1\ldots  p_k}_{\lceil p_j}$. The quantile regions $\mathbb{C}_{{\mathcal M}_0}^{\rm P}(\tau) \coloneqq \Qb(\mathbb{C}_{{{\cal M}_0}}^{\rm U}(s_{\tau}))$, with ${\rm P}$-probability contents $\tau\in [0,1]$, have the same interpretation and enjoy the same properties as the $\Qb$-images of  the spherical caps \eqref{capsdef} defined  on ${\mathcal S}^p$. 
 
 \item[(b)] Similarly to Section~\ref{232sec}, replacing the polyspherical caps ${\mathbb C}^{\rm U}_{{\cal M}_0}(\tau)$ in \eqref{polycapsdef} with {\it polyspherical  equatorial strips} of the form ${\mathbb C}^{\rm U }_{{\cal M}_0}((1+\tau)/2) \setminus {\mathbb C}^{\rm U}_{{\cal M}_0}((1-\tau)/2)$   yields   quantile regions  homotopic to $\boldsymbol{\mathcal{S}}^{p_1\ldots  p_k}_{\lceil p_j}$ ($k$ possible choices for $j$). 
 
  \item[(c)] Mixed types of submanifolds, involving   poles {\it and} sub-polyspheres, also are possible. Details are left to the reader.
\end{compactenum}

\section{Empirical distribution and quantile functions, ranks, and signs}\label{sec.ranksign}

In this section, we introduce the empirical versions $\Fb\n$ and $\Qb\n$ of $\Fb$ and $\Qb$ defined in Section~\ref{sec.FbQb}, along with the corresponding empirical quantile regions and contours, ranks,  and signs. 
Throughout, let~$c=d^2/2$ and denote by $\Yb_1\n, \Yb_2\n, \ldots, \Yb_n\n$   an $n$-tuple of i.i.d.\ copies of $\Yb\sim {\rm P}$, with distribution and quantile functions $\Fb$ and $\Qb$.

\subsection{Empirical distribution and quantile  functions}\label{Sec31}

Let $\mathfrak{G}\n=~\!\{{\scriptstyle{\mathfrak{G}}}\n_1,\ldots,{\scriptstyle{\mathfrak{G}}}\n_n\}$ denote a sequence of   grids of~$n$ distinct points on $\cal M$ such that ${\rm P}^{{\mathfrak{G}}\n} \!\!\!\coloneqq \sum_{i=1}^n \delta_{{\scriptstyle{\mathfrak{G}}}\n_i}$, where~$\delta_{\yb}$ stands for the Dirac measure at $\yb \in~\!{\cal M}$, converges weakly to ${\rm U}_{\mathcal M}$ as~$n\to~\!\infty$:  with a convenient abuse of language, we say that the grids $\mathfrak{G}\n$ are con\-verging weakly to~${\rm U}_{\mathcal M}$.
 Since~$\Fb$ is   minimizing the expected transport cost between~$\rm P$ and~${\rm U}_{\mathcal M}$, we quite natu\-rally  define the {\it empirical (Riemannian) distribution function}~$\Fb\n$ as the optimal transport map pushing
 ~${\rm P}\n\!\! \coloneqq \sum_{i=1}^n \delta_{\Yb_i}$ forward to~${\rm P}^{{\mathfrak{G}}\n}\!\!$.   More precisely,~$\Fb\n$ is mapping the sample $\Yb\n \!\! \coloneqq\{\Yb_1\n, \ldots, \Yb_n\n\}$ to~$\{\Fb\n(\Yb_1\n), \ldots, \Fb\n(\Yb_n\n)\}=~\!\mathfrak{G}\n$ 
in such a way that  (an optimal coupling problem)
\begin{equation}\label{eq.OptCouple}
\sum_{i=1}^n c(\Yb_i\n , \Fb^{(n)}(\Yb_i\n)) = \underset{\pi\in \Pi\n }{\min}  \sum_{i=1}^n  c(\Yb_i\n , {\scriptstyle{\mathfrak{G}}}\n_{\pi (i)}),
\end{equation}
where $\Pi\n$ denotes the set of all permutations $\pi$ of the integers $\{1,\ldots,n\}$.  Similarly define the {\it empirical (Riemannian) quantile function} $\Qb\n$ as the inverse of $\Fb\n$ (which is a bijective mapping). 

The following result establishes the uniform strong convergence of~$\Fb\n$ and $\Qb\n$   to~$\Fb$ and $\Qb$; see Appendix~\ref{App.proof} for the proof.


\begin{proposition}\label{Gliv}
{\rm (Glivenko-Cantelli).} {Let ${\rm P} \in \mathfrak{P}$ have bounded density on ${\cal M}$ and  let Assumptions~\ref{ass.cost}--\ref{ass.Cutyy0} hold}. 
Then, assuming that ${\rm P}^{{\mathfrak{G}}\n}$ converges weakly to ${\rm U}_{\cal M}$ as $n\to\infty$,
$${\max}_{1 \leq i \leq n} d(\Fb^{(n)}(\Yb_i^{(n)}), \Fb(\Yb_i^{(n)})) \rightarrow 0 \  {\text{and} \   {\max}_{1 \leq i \leq n} d(\Qb^{(n)}({\scriptstyle{\mathfrak{G}}}\n_i), \Qb({\scriptstyle{\mathfrak{G}}}\n_i)) \rightarrow 0 \ \ \text{a.s.}}$$
as $\ny$.
\end{proposition}

\subsection{Riemannian empirical quantile regions and contours,  ranks, and signs}\label{subsec.empQuan}

Unlike the real line,   Riemannian manifolds are not canonically ordered, which makes the definition  of ranks and signs problematic. The definition of ranks, actually, is directly related to the definition of empirical quantile contours: the rank $R\n_i$ of the observation $\Yb\n_i$, indeed, is larger than  the rank $R\n_j$ of the observation $\Yb\n_j$ if and only if the order of the empirical quantile contour running through $\Yb\n_i$ is larger than the order of the empirical quantile contour running through $\Yb\n_j$, and two observations have the same rank iff they belong to the same empirical quantile contour. We therefore first proceed with the construction of empirical quantile contours, a condition that requies a structured grid.

In ${\mathbb R}^p$ ($p>1$), \cite{Hallin2021} have shown  how    measure-transportation-based empirical quantile contours inducing a data-driven ordering of the observations, hence  ranks,   can be obtained by imposing a regular structure on the grid~$\mathfrak{G}\n$ used to construct the empirical distribution function. Similar ideas have been implemented by \cite{HLV2022} in ${\mathcal S}^p$ (directional data) and also apply here. 

In  ${\mathbb R}^p$, the target uniform distribution of the theoretical center-outward distribution~$\Fb$ is the uniform over the unit ball, centered at the origin. That centering, in ${\mathbb R}^p$, has no impact on $\Fb$ due to    shift invariance: denoting by $\Fb_{{\boldsymbol \mu}}$ the center-outward distribution of~$\Yb + {\boldsymbol \mu}$, indeed,  $\Fb_{{\boldsymbol \mu}}(\yb + {\boldsymbol \mu}) = \Fb_{\boldsymbol 0} (\yb)$ for any $\yb$ and ${\boldsymbol \mu}$ in ${\mathbb R}^p$. As a result, the same grid centered at the origin  can be used to define consistent empirical quantile regions and contours, irrespective of the distribution. This is no longer the case on manifolds, where the definition of empirical quantile regions requires the choice of a  distribution-specific submanifold~${\mathcal M}_0$.  The latter typically depends on the underlying distribution~$\rm P$ with distribution function $\Fb$ via  the image $\Fb (\thetab)$ of some   $\thetab\in{\mathcal M}$ (such as a Fr\' echet mean of~$\rm P$) determining a pole or an equator: see Sections~\ref{233sec},~\ref{ptorus}, and~\ref{polysphere} for examples. We henceworth emphasize this fact by  writing ${\mathcal M}_{\Fb(\thetab)}$ instead of ${\mathcal M}_0$ and make the following assumption.  

{\begin{assumption}\label{ass.M0}
{\rm The submanifold ${\cal M}_0$ used in the definition of quantile regions and contours is of the form ${\cal M}_0={\cal M}_{\Fb (\thetab)}$ for some 
  $\thetab\in{\cal M}$ where  
(i) the mapping~${\bf f}\mapsto~\!{\cal M}_{\bf f}$, ${\bf f}\in~\!{\cal M}$ is continuous in the sense that, denoting by $d_H$ the Hausdorff distance\footnote{Recall that the {\it Hausdorff distance} between two subsets $A$ and $B$ of  ${\cal M}$ (equipped with the distance~$d$)
 is defined as 
$d_{\rm H}({A}, {B}) \coloneqq \inf
\{\epsilon\geq 0: {A}\subset { B}^{\epsilon}  \text{ and }  {B}\subset {A}^{\epsilon}\}$,  
where $A^{\epsilon}$ and~$B^{\epsilon}$ denote the {\it $\epsilon$-fattenings} (the {\it $\epsilon$-fattening} of   ${C} \subset  \mathcal{M}$ is defined as
	${C}^{\epsilon} \coloneqq \bigcup_{\yb \in {C}} \{\zb \in \mathcal{M}  : d(\zb, \yb)  \leq ~\!\epsilon\})$ of $A$ and $B$,  respectively,  Note that continuity also implies measurability of ${\bf f}\mapsto{\cal M}_{\bf f}$ from $\cal M$ equipped with ${\mathfrak F}_{\cal M}$ to ${\mathfrak F}_{\cal M}$ equipped with the Borel $\sigma$-field induced by the Hausdorff distance.} 
 bet\-ween subsets of $\cal M$ equipped with the Riemannian distance~$d\vspace{-0.75mm}$, $d({\bf f}\n, {\bf f})\rightarrow 0$ as $n\to\infty$  implies~$d_H({\cal M}_{{\bf f}\n}, {\cal M}_{\bf f})\to 0$ as $n\to\infty$;
 (ii) a strongly consistent estimator ${\hat\thetab}\n$of $\thetab$ is available. 
 }
\end{assumption}

This assumption is satisfied, for the manifolds~$\cal M$ and  the various choices of ${\cal M}_0$ described in Section~\ref{232sec} (the $p$-sphere), \vspace{-.5mm}
Section~\ref{ptorus} (the $p$-torus), and Section~\ref{polysphere} (the polysphere), with~${\hat\thetab}\n$an  empirical Fr\' echet mean\footnote{An {\it empirical Fr\' echet mean}  $\hat{\thetab}_{\rm Fr}\n\vspace{-0.8mm}$ is an element of 
	${\cal A}_{\rm Fr}\n \coloneqq  {\rm argmin}_{\yb \in  \mathcal{M}} \sum_{j=1}^n c(\Yb_j\n, \yb)$, the {\it empi\-rical Fr\'echet mean set}.  For any $\epsilon >0$ and $n$ sufficiently large,  ${\cal A}_{\rm Fr}\n$,  with probability one, is a subset  of the $\epsilon$-{\it fattening} ${\cal A}^{\epsilon}_{\rm Fr}$ of ${\cal A}_{\rm Fr}\vspace{-0.8mm}$ (see, e.g., \citet[Theorem~2.2]{bhattacharya2012nonparametric}). Hence, without loss of generality, we can assume that~$\hat{\thetab}_{\rm Fr}\n$ is a consistent estimator of~$\thetab_{\rm Fr}$.
}~$\hat{\thetab}\n_{\rm Fr}$.

Denote by $\mathfrak{G}\n=~\!\{{\scriptstyle{\mathfrak{G}}}\n_1,\ldots,{\scriptstyle{\mathfrak{G}}}\n_n\}$ the grid from which $\Fb\n$ and $\Qb\n$ have been constructed as in Section~\ref{Sec31}. Let ${\scriptstyle\mathfrak{G}}\n_* \in \mathfrak{G}\n$ be the gridpoint such that $$\Qb\n({\scriptstyle\mathfrak{G}}\n_*)={\argmin}_{\Yb\n_i \in \Yb\n} d(\Yb\n_i, \hat\thetab\n);$$ ${\scriptstyle\mathfrak{G}}\n_*$, thus, is the $\Fb\n$-image of the observation which is closest, in the $d$-distance,~to~${\hat\thetab}\n\!$.

The following Lemma, which is a special case of Lemma~\ref{lem.Mxhat} in Section~\ref{sec.QuanReg},
states that~$\widehat{\cal M}\n_0\coloneqq {\cal M}_{{\scriptstyle\mathfrak{G}}\n_*}$ is a consistent estimator of ${\cal M}_{\Fb(\thetab)}$.

\begin{Lem}\label{lem.Mhat}
	Let ${\rm P} \in \mathfrak{P}$ have bounded density on ${\cal M}$ and  let Assumptions~\ref{ass.cost}--\ref{ass.M0} hold. Then $d_{\rm H}(\widehat{\cal M}\n_0, {\cal M}_{\Fb(\thetab)}) \rightarrow 0$ a.s. as $n \rightarrow \infty$.
\end{Lem}



As in 
\cite{HLV2022}, still under Assumption~\ref{ass.M0},  the regular grid then can constructed via the following two steps. \vspace{2mm}

 \noindent{\bf Step 1.} First,  determine a consistent empirical version~$\widehat{\cal M}\n_0$ of ${\cal M}_0$. This requires  a first empirical transport based, e.g., on a grid $\mathfrak{G}\n_0$ of $n$ independent copies of ${\bf U}\sim {\rm U}_{\mathcal M}$, yielding, as in Section~\ref{Sec31},  empirical distribution  and quantile functions $\Fb\n$ and $\Qb\n$. Under Assumption~\ref{ass.M0}, $\widehat{\cal M}_{0}\n$ then converges a.s. to  ${\cal M}_{0}$ in the Hausdorff distance. For each $n$,  the manifold $\widehat{\cal M}\n_0$   determines a vector field~$V^{\perp}_{\widehat{\cal M}\n_0}$ and  nested   contours ${\mathcal C}^{\rm U}_{\widehat{\cal M}\n_0} (\tau)\vspace{-1mm}$ of the form  \eqref{defgencontour} for~$\tau\in[0,1]$, each of which is a (random) submanifold homotopic to~$\widehat{\cal M}\n_0$,  the sequence of which converges a.s., in~$d_H$ distance, to ${\mathcal C}^{\rm U}_{{\cal M}_0} (\tau)$ as $n\to\infty$.\vspace{2mm}

\noindent{\bf Step 2(a).} 
Next, factorize $n$ into $n = n_R n_S+n_0$, with~$n_R\to\infty$ and~$n_S \rightarrow~\!\infty$ as~$n  \rightarrow \infty$  and $0 \leq n_0 < n_S$.\footnote{Imposing, as in $\mathbb{R}^d$,  $0 \leq n_0 <\min (n_R, n_S)$ may not be possible: for instance, if ${\cal M}_0$ is of dimension $(p-1)$, we should have $n_0 = O(n^{(p-1)/p})$ and $n_R = O(n^{1/p})$.}  Select $n_R$ contours, of the form   ${\mathcal C}^{\rm U}_{\widehat{\cal M}\n_0} (r/(n_R + 1))$, $r=1,\ldots,n_R$. 


 \noindent{\bf Step 2(b).} On each contour ${\mathcal C}^{\rm U}_{\widehat{\cal M}\n_0} (r/(n_R + 1))\vspace{-1mm}$ (a submanifold of $\cal M$), construct  a grid~${\mathfrak{G}}^{(n_S)}_{\widehat{\cal M}_0\n}(r)$ of~$n_S$ points i.i.d.\footnote{As random elements, these points are not living on the same probability space as the observations.} with uniform distribution ${\rm U}_{{\mathcal C}^{\rm U}_{\widehat{\cal M}\n_0} (r/(n_R + 1))}$  (for $p=2$, the dimension of~${\mathcal C}^{\rm U}_{\widehat{\cal M}\n_0} (r/(n_R~\!+~\!1))\vspace{-0mm}$ is one and an equispaced grid can be considered~instead). 

\noindent{\bf Step 2(c).} On $\widehat{\cal M}\n_0$, construct  a grid $\mathring{\mathfrak{G}}^{(n_0)}_{\widehat{\cal M}_0\n}$ of~$n_0$ points i.i.d. (in case ${\cal M}_0$ contains a single point, consider $n_0$ copies of that unique point;  in that case,   $n_0 < \min\{n_R, n_S\}$ can be imposed).

  
\noindent{\bf Step 2(d).}  The final regular grid ${\mathfrak{G}}\n_{\widehat{\cal M}_0\n}$ that we obtain is
\begin{equation}\label{eq.Gn}
	{\mathfrak{G}}\n_{\widehat{\cal M}_0\n} \coloneqq \mathring{\mathfrak{G}}^{(n_0)}_{\widehat{\cal M}_0\n} \bigcup {\mathfrak{G}}^{(n_S)}_{\widehat{\cal M}_0\n}(1) \bigcup \cdots \bigcup {\mathfrak{G}}^{(n_S)}_{\widehat{\cal M}_0\n}(n_R).
\end{equation}


Denote by $\Fb\n_*$  the optimal transport---see \eqref{eq.OptCouple}---from the sample to this grid and by~$\Qb\n_*$ its inverse.\footnote{Note that, as soon as $n_0>1$, $\Qb\n_*$ is a set-valued function.} The following Lemma, which follows from Lemma~\ref{lem.Mhat}, shows that the grid ${\mathfrak{G}}\n_{\widehat{\cal M}_0\n}\vspace{-0.6mm}$ satisfies the assumptions of Proposition~\ref{Gliv}  
so that~$\Fb\n_*$ and~$\Qb\n_*$, just as~$\Fb\n$ and~$\Qb\n$, enjoy the Glivenko-Cantelli property. 

\begin{Lem}\label{gridconvlem}   
Let ${\rm P} \in \mathfrak{P}$ have bounded density on ${\cal M}$ and  let Assumptions~\ref{ass.cost}--\ref{ass.M0} hold.
Then the grid ${{\mathfrak{G}}}\n_{\widehat{\cal M}_0\n}$  converges weakly to ${\rm U}_{\mathcal M}$ as $n\rightarrow \infty$.
\end{Lem}

Based on $\Fb\n_*$ and~$\Qb\n_*$,  define  empirical quantile contours, regions, ranks, and signs as follows.

\begin{definition}\label{defRanks}   
{\rm Call  
{\it empirical quantile contour of order $r/{(n_R +1)}$} ($r = 1, \ldots, n_R$ if $n_0 = 0$, and~$r = 0, \ldots, n_R$ if $n_0 \neq 0$) 
 the collection  of observations 
 \[
{\mathcal C}\n(r/{(n_R +1)}) \coloneqq 
\begin{cases}
	\Qb\n_*\left(
	\mathring{\mathfrak{G}}^{(n_0)}_{\widehat{\cal M}_0\n}
	\right)  & \qquad  r = 0\\
	\Qb\n_*\left(
	{\mathfrak{G}}^{(n_S)}_{\widehat{\cal M}_0\n}(r)
	\right) & \qquad  r = 1, \ldots, n_R
\end{cases}
,
\]
 and  {\it empirical quantile region of order $r/{(n_R +1)}$}
 the collection  of observations  
 $${\mathbb C}\n(r/{(n_R +1)})\coloneqq \Qb\n_*\left(\mathring{\mathfrak{G}}^{(n_0)}_{\widehat{\cal M}_0\n} \bigcup {\mathfrak{G}}^{(n_S)}_{\widehat{\cal M}_0\n}(1) \bigcup \cdots \bigcup {\mathfrak{G}}^{(n_S)}_{\widehat{\cal M}_0\n}(r)
\right).$$
Call  {\it rank} 
of $\Yb_i\n$, $i=1,\ldots,n$ the integer 
 \[
R_i\n\!\! = R\n(\Yb_i\n) \coloneqq 
\begin{cases}
	0 & \quad \text{if}\  \  \Yb_i\n \in \Qb\n_*\left(\mathring{\mathfrak{G}}^{(n_0)}_{\widehat{\cal M}_0\n}\right) \\
	k & \quad \text{if}\  \  \Yb_i\n \in {\mathcal C}\n(k/{(n_R +1)}), k = 1, \ldots, n_R
\end{cases}
\]
and   {\it sign} $\Sb_i\n$ of $\Yb_i\n$ 
 \begin{equation}\label{eq.viSi}
 \Sb_i\n\!\! = \Sb\n(\Yb_i\n) \coloneqq 
 \begin{cases}
 	{\bf 0} & \quad \text{if}\  \  \Yb_i\n \in \Qb\n_*\left(\mathring{\mathfrak{G}}^{(n_0)}_{\widehat{\cal M}_0\n}\right) \\
 	{\vb^{(n)}_i}/{\Vert \vb^{(n)}_i \Vert}  & \quad \text{otherwise}
 \end{cases}
 \end{equation}
with $\vb\n_i\in V^{\perp}_{\widehat{\cal M}\n_0}$ such that $\exp_{\yb}(t\vb\n_i) =  \Fb\n_*(\Yb_i\n)$ for some $t\in (0, 1]$ and $\yb \in  \widehat{\cal M}\n_0$.

	}
	\end{definition}

Step 2 of the construction of ${{\mathfrak{G}}}\n_{\widehat{\cal M}_0\n}$ simplifies in case the factorization mentioned in Remark~\ref{Rem.suff} holds. More precisely, Step 2(b) can be replaced with the following Step~\!2(b$^\prime$).\smallskip

\noindent{\bf Step 2(b$^\prime$).} Select, among the contours ${\mathcal C}^{\rm U}_{\widehat{\cal M}\n_0} (\tau)\vspace{-1mm}$, $\tau\in [0,1]$ considered in Step~1,  a contour~$\widehat{\mathcal C}^{\rm U}_* \coloneqq {\mathcal C}^{\rm U}_{\widehat{\cal M}\n_0} (\tau_*)$ of dimension $(p-1)$ with $\tau_* \in [0, 1]$ satisfying $s(\tau_*) = s_*\vspace{-1mm}$ for~$s(\cdot)$ given in Assumption~\ref{ass.Wincrease} and $s_*$ in \eqref{eq.C*}. On this  contour $\widehat{\mathcal C}^{\rm U}_*$, construct  a grid of~$n_S$ i.i.d.\ points with uniform distribution ${\rm U}_{\widehat{\mathcal C}^{\rm U}_*}$ (for $p=2$, an equispaced grid can be considered instead) and the~$n_S$ geodesics $\gamma^{(n_S)}_i$, $i = 1, \ldots, n_S$ running through these $n_S$ points, where $\gamma^{(n_S)}_i: [0, 1] \rightarrow {\cal M}$ satisfies $(\gamma^{(n_S)}_i)\pr(0)  \in V^\perp_{\widehat{\mathcal M}\n_0}$\footnote{See Appendix~\ref{App.Rieman} for the definition of the tangent vector $\gamma\pr$.}. The $n$-points grid ${\mathfrak{G}}\n_{\widehat{\cal M}_0\n}$ to be used is the set of $n_Rn_S$ intersections between these $n_S$ geodesics and the $n_R$ contours obtained in Step~2(a), along with the $n_0$ gridpoints in Step~2(c). By construction, the empirical distribution over ${\mathfrak{G}}\n_{\widehat{\cal M}_0\n}$ converges weakly to ${\rm U}_{\mathcal M}$. The   independence between~$w_{\Ub}$ and $\Yb_{\Ub}$ in  $\Ub = \exp_{\Yb_{\Ub}}(w_{\Ub}\vb_{\Yb_{\Ub}})  \sim {\rm U}_{\mathcal M}$  
 then entails  (see Proposition~\ref{PropDistFree}(vi))  the independence of the signs and ranks of Definition~\ref{defRanks}.


Proposition~\ref{prop.ContAsymptotic} establishes the uniform consistency (in  Hausdorff  distance) of the empirical quantile contours $\mathcal{C}\n$ and regions $\mathbb{C}\n$, which follows as a particular case of more general consistency results for quantile regression (see Proposition~\ref{prop.QwAsymptotic} in Section~\ref{sec.QuanReg}).

\begin{proposition}\label{prop.ContAsymptotic}
Let ${\rm P} \in \mathfrak{P}$ have bounded density on ${\cal M}$. For~$c=d^2/2$, let Assumptions~\ref{ass.cost}--\ref{ass.M0} hold. Let $\{\Yb\n_1, \ldots, \Yb\n_n\}$ be i.i.d.\ with distribution~${\rm P}$ over~$\mathcal{M}$. Then, for any $\epsilon>0$, as $n \rightarrow \infty$,
$${\rm P}\left\lbrace \max_{r \in {\cal R}_n} d_{\rm H}\left(\mathbb{C}\n(r/(n_R+1)),  \mathbb{C}^{\rm P}_{\mathcal{M}_0}(r/(n_R+1))\right) > \epsilon \right\rbrace \rightarrow 0$$
and
$${\rm P}\left\lbrace \max_{r \in {\cal R}_n} d_{\rm H}\left(\mathcal{C}\n(r/(n_R+1)),  \mathcal{C}^{\rm P}_{\mathcal{M}_0}(r/(n_R+1))\right) > \epsilon \right\rbrace \rightarrow 0,$$
where
\begin{equation*} 
	{\cal R}_n \coloneqq 
	\begin{cases}
		\{1, \ldots, n_R\} & \text{if} \ \ n_0 = 0 \\
		\{0, \ldots, n_R\} & \text{if} \ \ n_0 \neq 0
	\end{cases}
	.
\end{equation*}
\end{proposition}

\subsection{Distribution-freeness and ancillarity}

We now establish the main finite-sample properties of the ranks and quantiles just defined. These properties are essential in the construction of distribution-free testing procedures such as rank tests, R-estimators,~etc. 


Denote by ${\rm P}_{\Yb}\n$ the joint distribution of the sample~$\Yb\n = \big( \Yb_1\n, ..., \Yb_n\n\big)$ and by~$\Yb\n_{(\cdot)}=\big( \Yb_{(1)}\n, ..., \Yb_{(n)}\n\big)$  the {\it order statistic} of $\Yb\n$,  obtained by ordering the observations via some arbitrary but fixed  criterion; e.g.,  for~${\cal M}$ in Examples~\ref{ex.sphere}-\ref{ex.product}, by increasing values of first components. 
The order statistic gene\-rates~the sigma-field of permutationally invariant measurable functions of~$\Yb\n$ and  it readily follows from the factorization criterion that it is sufficient for  the  family~$\mathfrak{P}\n \coloneqq  \{{\rm P}_{\Yb}\n \vert \,   {\rm P} \in~\!\mathfrak{P}\}$---actually, it is minimal sufficient: see Proposition~\ref{PropDistFree} (i).  Let~$\Fb\n(\Yb\n)\! \coloneqq\!  (\Fb\n(\Yb_1\n), \ldots, \Fb\n(\Yb_n\n))$ where $\Fb\n$ denotes the empirical distribution function defined in Section~\ref{Sec31}:\linebreak then,~$\Fb\n(\Yb\n)$ is uniformly distributed over the $n!$ permutations of the  grid $\mathfrak{G}\n$. 

Assuming  that  $\widehat{\cal M}\n_0$ is $\Yb\n_{(\cdot)}$-measurable---which can be taken care of via Rao-Blackwellization), ${\Fb_*}\n(\Yb\n) \coloneqq  ({\Fb_*}\n(\Yb_1\n), \ldots, {\Fb_*}\n(\Yb_n\n))$ with ${\Fb_*}\n$ defined in Section~\ref{subsec.empQuan} is, conditionally on $\widehat{\cal M}\n_0$,  uniformly distributed over the  permutations (with repetitions in case ${\cal M}_0$ is a singleton and $n_0>1$) of the  regular grid~$\mathfrak{G}\n_{\widehat{\cal M}_0\n}$.  The resulting vector of ranks $\Rb\n\coloneqq (R\n_1, \dots,R\n_n)$ thus is conditionally, hence also unconditionally, uniform over the permutations with repetitions of~$\big\{0\text{ ($n_0$ copies), }1\text{ ($n_S$ copies), }\ldots, n_R \text{ ($n_S$ copies)}\big\}$. The vector~${\bf S}\n \coloneqq ({\bf S}\n_1,\ldots,{\bf S}\n_{n})$ of signs similarly enjoys distribution-freeness. 
Let $\tilde{n}_0$ denote the number of repetitions in $\mathring{\mathfrak{G}}^{(n_0)}_{\widehat{\cal M}_0\n}$ (for instance, $\tilde{n}_0 = n_0$ if ${\cal M}_0$ is a singleton; $\tilde{n}_0 = 1$ if ${\cal M}_0$ is of dimension $(p-1)$). 
The following proposition  (see Appendix~\ref{App.proof} for a proof) provides a precise summary of these  properties.

\begin{proposition}\label{PropDistFree}
	Let $\{\Yb_1\n, ..., \Yb_n\n\}$ be i.i.d.\  with distribution  ${\rm P} \in \mathfrak{P}$. Then,  for all ${\rm P} \in \mathfrak{P}$,
	\begin{compactenum}
		\item[(i)] the order statistic $\Yb_{(\cdot)}\n$ is sufficient and complete, hence minimal sufficient;
		\item[(ii)]  $\Fb\n(\Yb\n)$ 
		is distribution-free and uniformly distributed over the $n!$ permutations  of the   grid $\mathfrak{G}\n_0$;
		\item[(iii)] assuming that $\hat{\thetab}\n$ is a symmetric function of $\Yb\n$ (that is, is $\Yb\n_{(\cdot)}$-measurable), conditionally on $\widehat{\cal M}_0\n$, ${\Fb}\n_*(\Yb_n\n))$ is uniformly distributed over  the $n!/\tilde{n}_0!$ permutations (the $\tilde{n}_0$ copies in $\mathring{\mathfrak{G}}^{(n_0)}_{\widehat{\cal M}_0\n}$ being counted as $\tilde{n}_0\vspace{-1mm}$ indistinguishable points) of the grid~$\mathfrak{G}\n_{\widehat{\cal M}\n_0}$;
		\item[(iv)] ${\bf R}\n$ is distribution-free and  uniformly distributed over the $n!/n_0!(n_S!)^{n_R}$ permutations with repetitions of 
		$\big\{0\text{ ($n_0$ copies), }1\text{ ($n_S$ copies), }\ldots, n_R \text{ ($n_S$ copies)}\big\}$; 
		the rank~$R\n_i$, $i=1,\ldots,n$   takes value~0 with probability $n_0/n$ and value $j$  with proba\-bility $n_S/n$, $j=1,\ldots,n_R$.		
\item[(v)] letting $\mathcal{I} \coloneqq \left\{i \vert i \in \{1, \ldots,  n\} \ \text{and} 
\  \Yb_i\n \notin \Qb\n_*\left(\mathring{\mathfrak{G}}^{(n_0)}_{\widehat{\cal M}_0\n}\right) \right\}$, the vector of signs~${\bf S}\n $ is, conditionally on $\widehat{\cal M}_0\n$, uniformly distributed over the $n!/n_0!$ permutations with repetitions of 
 $\big\{{\boldsymbol 0}\text{ ($n_0$ copies) }, \vb\n_i/\Vert \vb\n_i \Vert \ \vert \ i \in \mathcal{I} \big\}$
with $\vb\n_i$ as in \eqref{eq.viSi};
\item[(vi)] for $n_0 = 0$, if ${{\mathfrak{G}}}\n_{\widehat{\cal M}_0\n}$ is constructed as in Step 2($b\pr$), the vector $(R_1\n, \ldots, R_n\n)$ of ranks  and  the vector $({\bf S}_1\n, \ldots, {\bf S}_n\n)$ of signs are mutually  independent;
\item[(vii)] the $\sigma$-field generated by $(R\n_1,\ldots,R\n_n)$ and $({\bf S}_1\n, \ldots, {\bf S}_n\n)$ is maximal ancillary for the family of distributions $\mathfrak{P}\n$.		
	\end{compactenum}
\end{proposition}

Note that part (vi) of Proposition~\ref{PropDistFree} only holds for   $n_0=0$. For $n_0>0$, indeed,~$R\n_i =~\!0$ automatically implies ${\Sb\n_i={\boldsymbol 0}}$. However, since
$n_0/n \rightarrow~\!0$ as~$n\to\infty$, this   dependence between ranks and signs induced by the $n_0$ observations mapped to~$\hat{\thetab}_{\bf U}\n$ is   asymptotically negligible.  A tie-breaking procedure similar to that proposed in \cite{Hallin2021} easily allows for removing this restriction; details are omitted.

\subsection{Ranks and empirical quantiles: graphical illustrations}\label{subsec.SimContour}

In this section, we provide graphic  illustrations,  for the 2-torus $\mathcal{T}^2 = {\cal S}^1 \times {\cal S}^1$ and the 2-sphere ${\cal S}^2$,  of the grid~$\mathfrak{G}\n_{\widehat{\cal M}_0\n}$ and, from simulated samples from various types of distributions, 
of  the concepts of ranks, signs, and empirical quantile contours and regions. We illustrate both the cases that ${\cal M}_0$ is a singleton and is of dimension $(p-1)$,  in the toroidal/spherical cap- and equatorial strip-type quantile regions, respectively.

\subsubsection{Cap-type quantile regions}
We first Illustrate how regular grid can be constructed via the two-steps approach in Section~\ref{subsec.empQuan} when ${\cal M}_0$ contains only one point.  To improve readability, and since the objective is to provide  clear pictures of the grids, unreasonably small values of~$n$ and~$n_R$ have been chosen 
 in Fi\-gure~\ref{Fig:Grid}, namely,  $n=121$, $n_0=1$, $n_R = 3$, $n_S=40$.
Plots of ${\mathfrak{G}}\n_{\widehat{\cal M}_0\n}\vspace{-1.5mm}$ for $\mathcal{T}^2$ and  ${\cal S}^2$ are shown in the first row of
Figure~\ref{Fig:Grid}. The purple point represents the single element ${\scriptstyle\mathfrak{G}}\n_*$ of~$\widehat{\cal M}_0\n$, and the associated $n_R$ contours ${\mathcal C}^{\rm U}_{\widehat{\cal M}\n_0} (r/(n_R + 1))$ for $r = 1, 2, 3$ are shown in red, blue and green, respectively. 
Two distinct representations of $\mathcal{T}^2$ are provided: in the left panel, the  flat~$[-\pi, \pi)^2$ square representation  centered at ${\scriptstyle\mathfrak{G}}\n_*\!$,  in which   opposite edges (i.e., ${\scriptstyle\mathfrak{G}}\n_*$'s cut locus)  are identified; in the central panel,   the familiar representation in~${\mathbb R}^3$.

\begin{figure}[htbp!]
	\centering
	\begin{tabular}{ccc}
		\hspace{-10mm}    \includegraphics[scale=0.28]{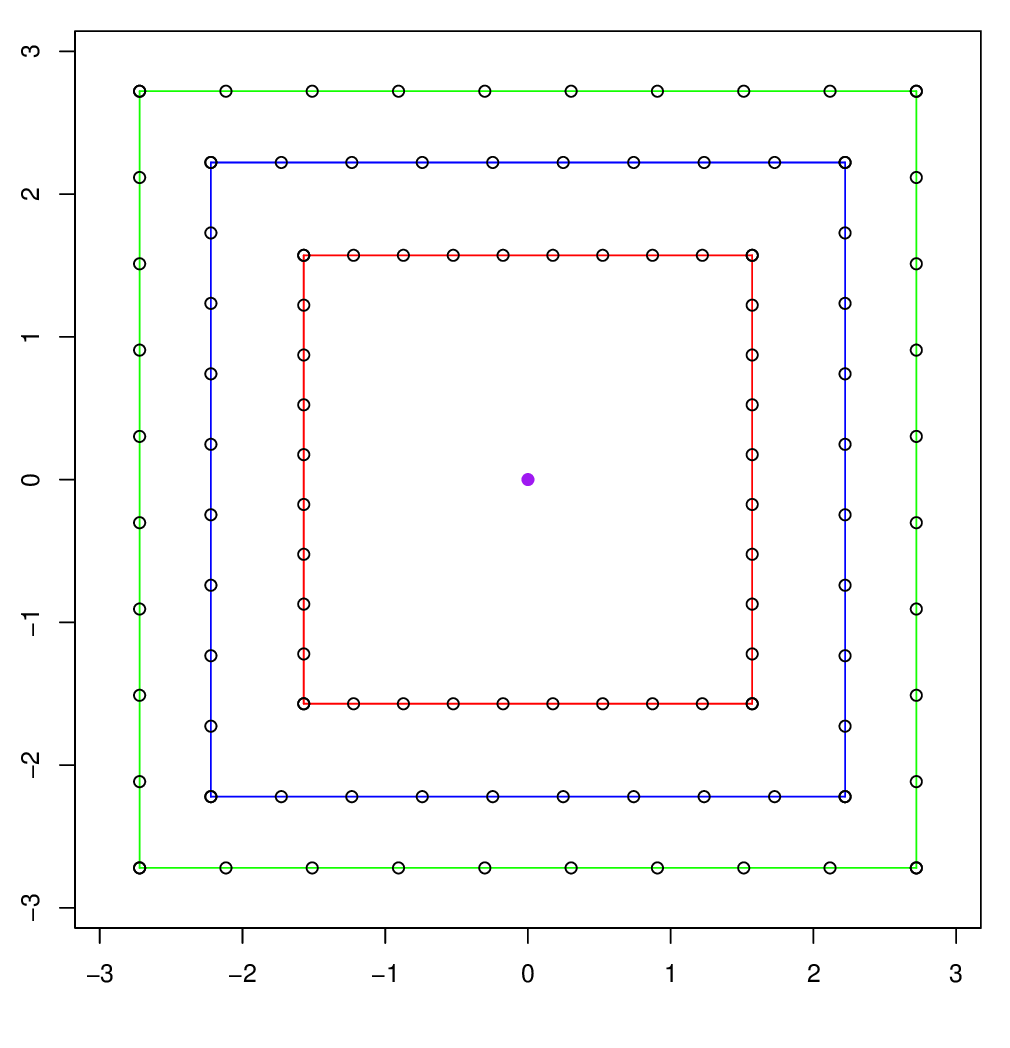} \vspace{-5mm}&
		\raisebox{7mm}{ \hspace{-8mm}\includegraphics[scale=0.15]{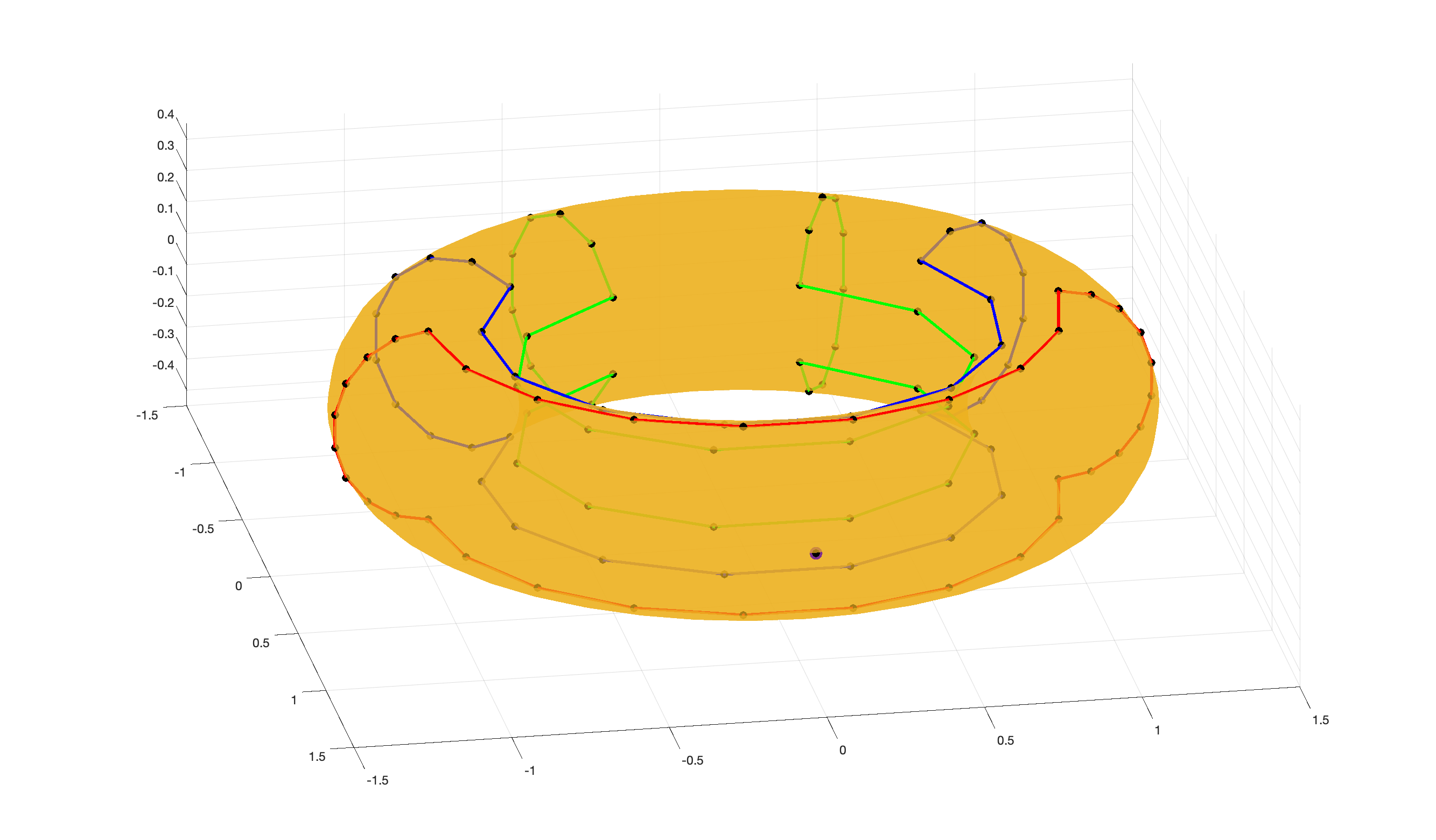}}\vspace{-5mm}&\raisebox{-3mm}{\hspace{-22mm} \includegraphics[trim=50 20 40 20,clip, width=0.45\textwidth]{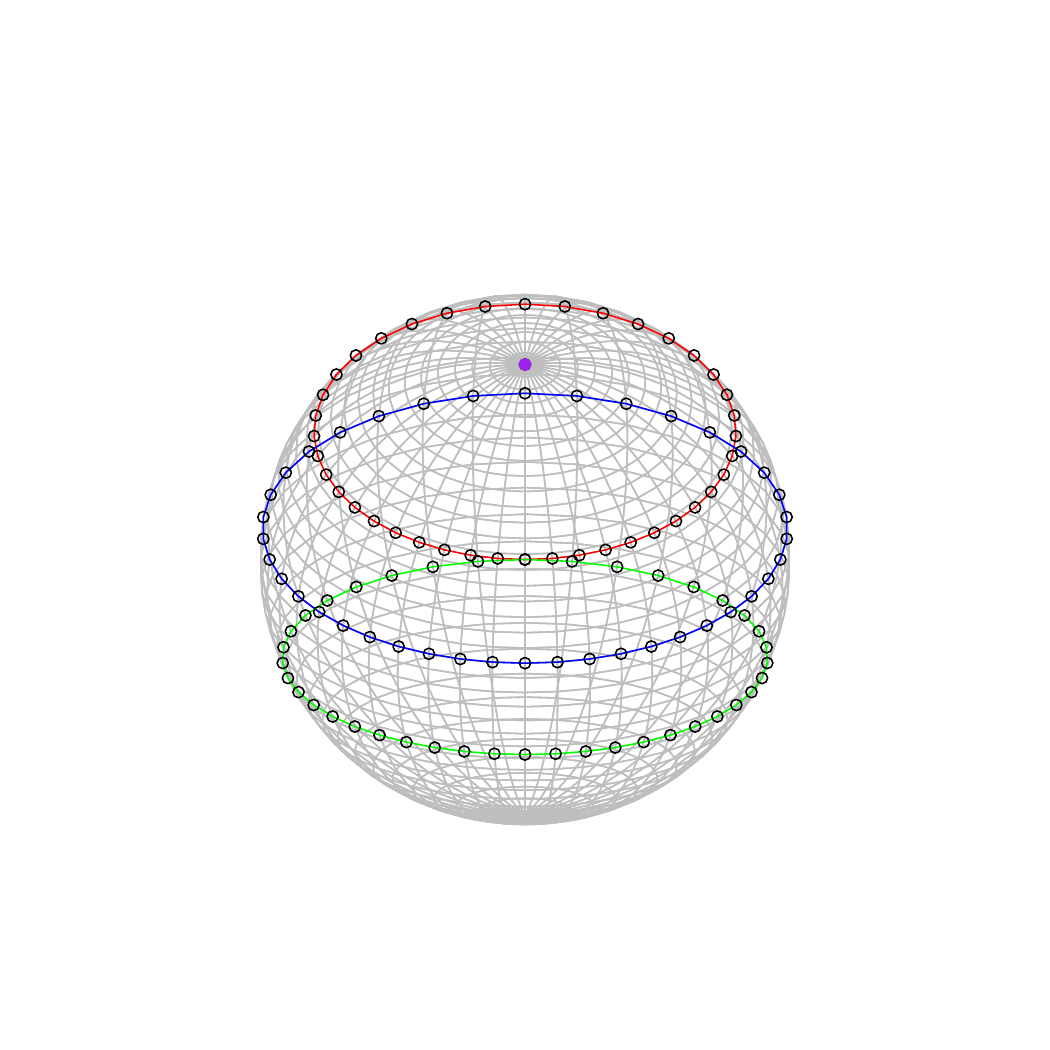}} \vspace{-3mm}
		\\
		
		\hspace{-10mm}        \includegraphics[scale=0.28]{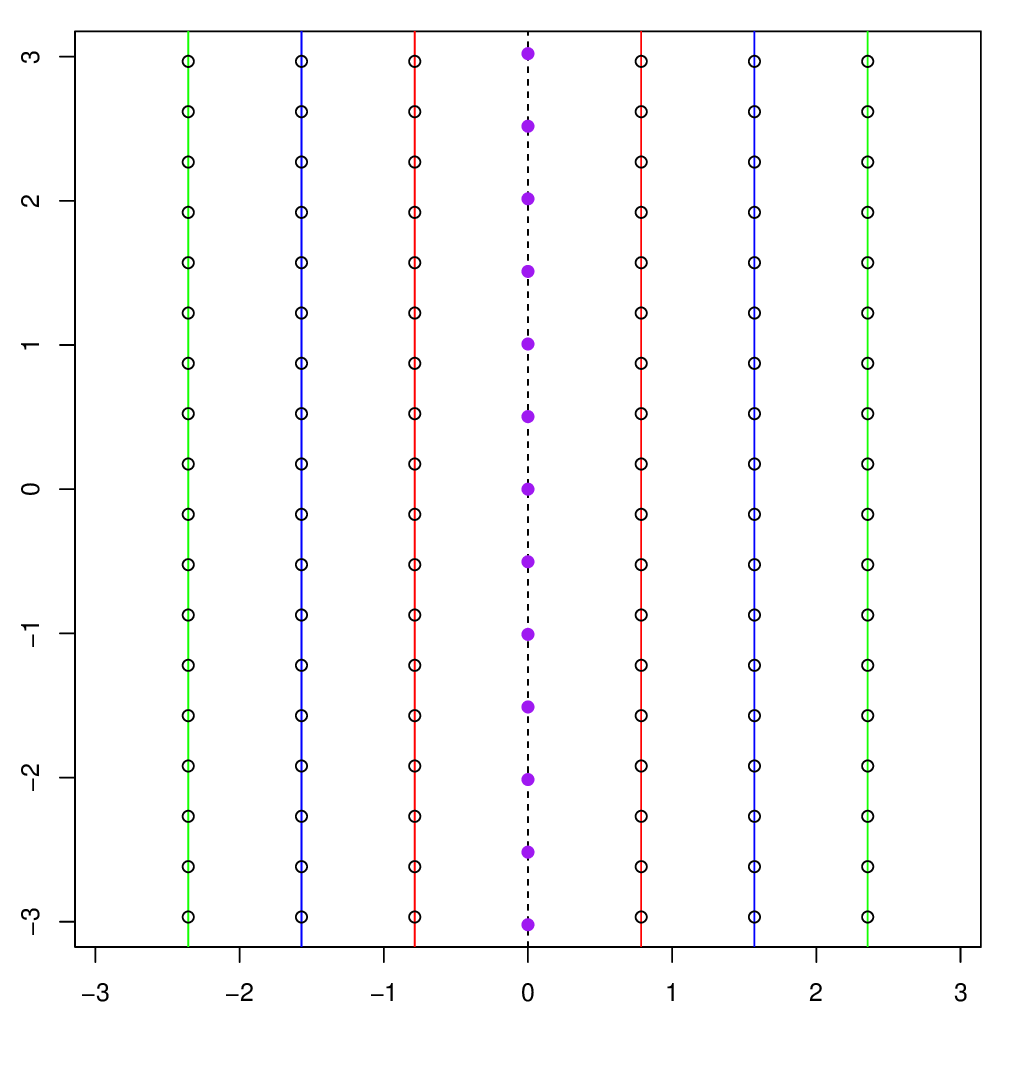}\vspace{-1mm}& 
		\raisebox{7mm}{\hspace{-7mm} \includegraphics[scale=0.15]{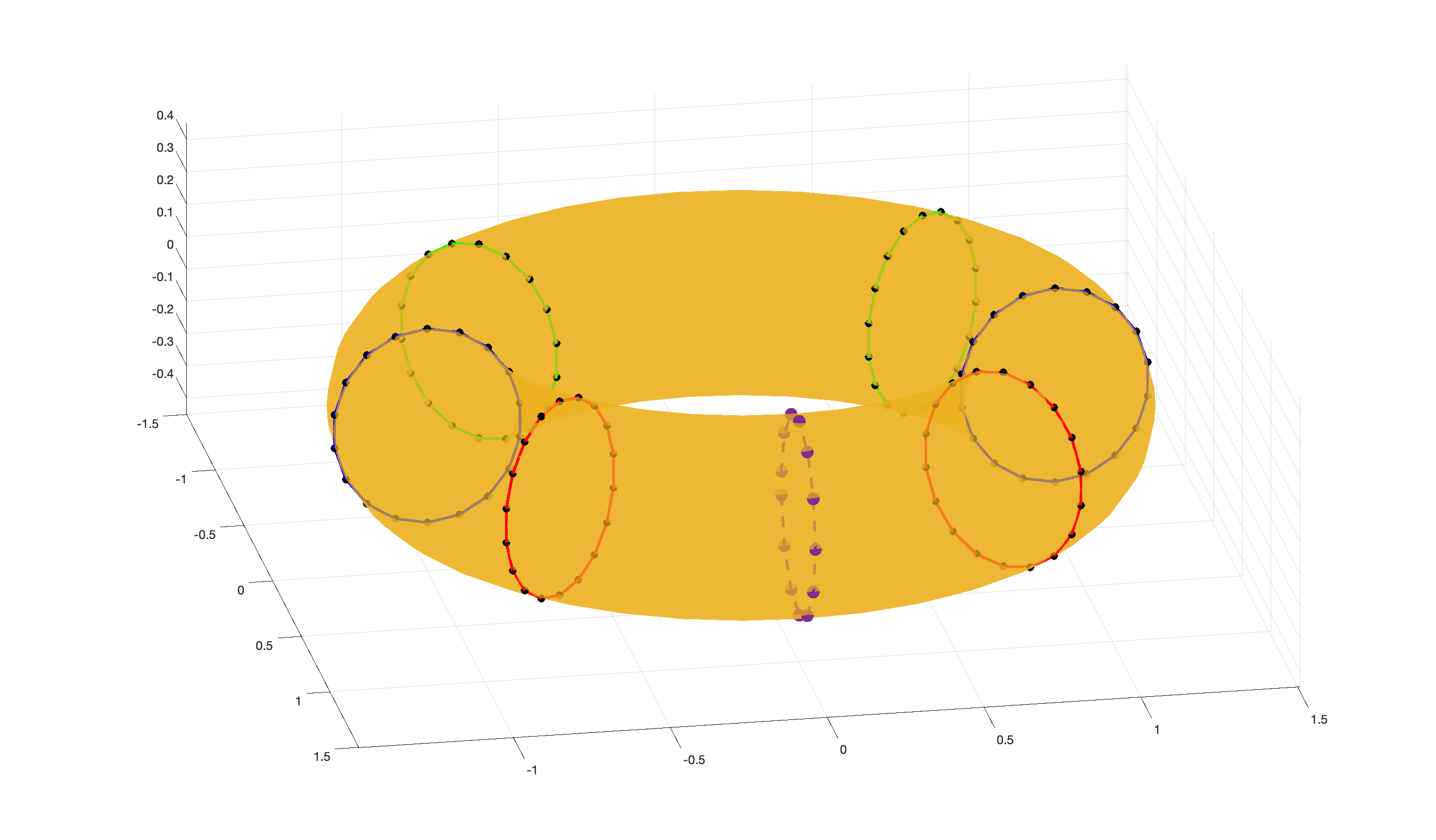}}\vspace{-1mm}
		&\raisebox{-3mm}{\hspace{-18mm}  \includegraphics[trim=50 20 40 50,clip, width=0.45\textwidth]{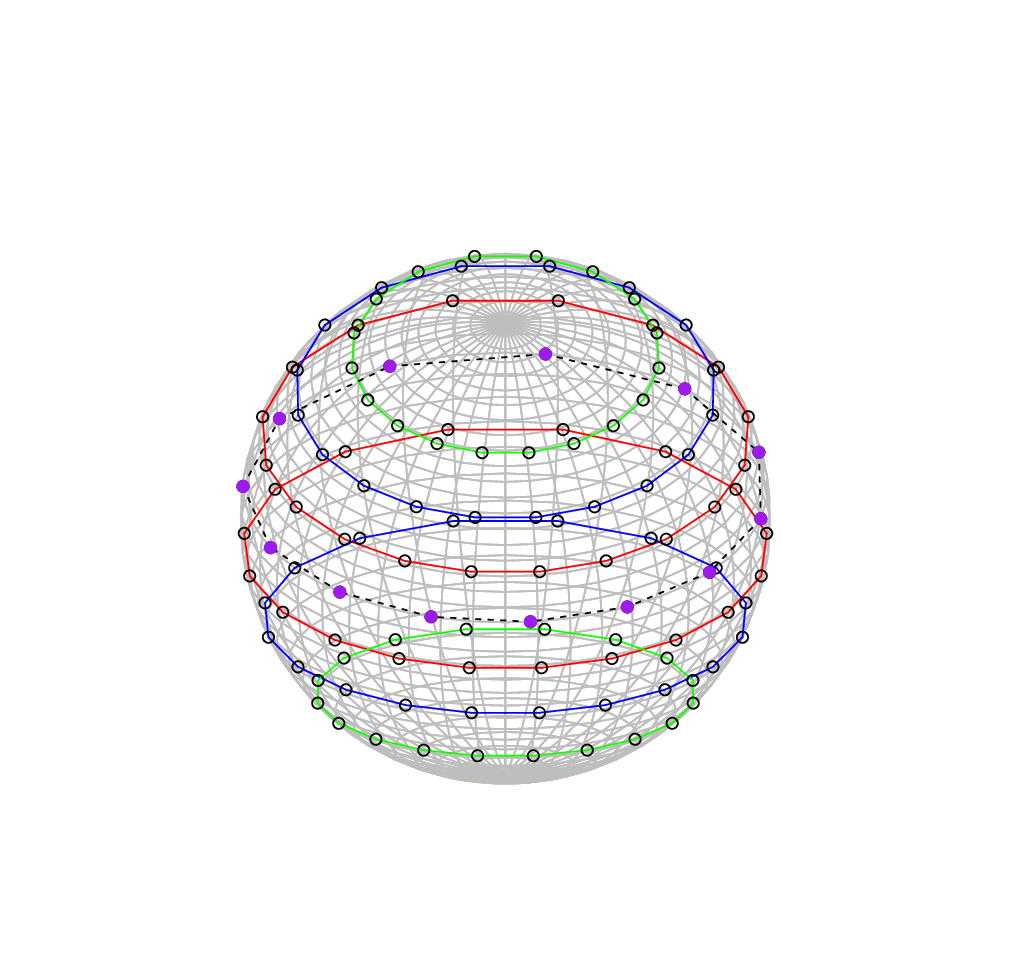}}\vspace{-1mm}
	\end{tabular}\vspace{-4mm}
	\caption{\small Plots of ${\mathfrak{G}}\n_{\widehat{\cal M}_0\n}$, $n=121$ for the 2-torus ${\cal T}^2$  (flat square representation  in the first col~\!\!umn,  ${\mathbb R}^3$-embedded representation in the second column), and the 2-sphere~${\cal S}^2$ (third column). The purple points represent the $n_0$ points on ${\widehat{\cal M}_0\n}$ and the~$n_R=3$ contours~${\mathcal C}^{\rm U}_{\widehat{\cal M}\n_0} (r/(n_R + 1))$ for $r = 1, 2, 3$ are shown in red, blue and, green, respectively. The first row corresponds to the case that ${\cal M}_0$ is a singleton, where we set $n_0=1$, $n_R = 3$, and~$n_S=40$. The second row is an illustration of the case that ${\cal M}_0$ is of dimension $(p-1)$, where we set $n_0=13$, $n_R = 3$, and $n_S=36$, and ${\widehat{\cal M}\n_0}$ is shown as the purple dashed line/loop.}
	\label{Fig:Grid}
\end{figure}

Next, to illustrate the concepts of ranks and empirical quantiles on $\mathcal{T}^2$ and $\mathcal{S}^2$, we simulate samples from the following distributions ((T1)-(T3) for $\mathcal{T}^2$, (S1)-(S3) for $\mathcal{S}^2$). 

\begin{enumerate}
	\item[(T1)] A {\it Bivariate Sine von Mises 
		distribution} BSvM$(\mub= {\bf 0}, \kappab = (3, 3)^\top, \lambda =0)$ with two independent components. The probability density function 
	of the BSvM$(\mub , \kappab , \lambda  )$ distribution, proposed by \cite{singh2002}, takes the form
	\begin{align}
		(\phi_1, \phi_2)^\top\in [-\pi,\pi]^2 \mapsto\,  & C_{\mub, \kappab, \lambda} \exp\big[\kappa_1\cos(\phi_1-\mu_1) \nonumber \\
		& \qquad + \kappa_2\cos(\phi_2-\mu_2)+\lambda\sin(\phi_1-\mu_1)\sin(\phi_2-\mu_2)\big], \label{eq.BSvM}
	\end{align}
	where $C_{\mub, \kappab, \lambda}$ is the normalization constant;
	\item[(T2)] a {\it Bivariate Sine von Mises 
		distribution} BSvM$(\mub= {\bf 0}, \kappab = (3, 3)^\top, \lambda =1.5)$ with two dependent components;
	\item[(T3)] a mixture of three BSvM's, i.e., a random variable $\Yb$ of the form
	\begin{equation}\label{eq.3MixBSvM}
		\Yb = I(U \leq 3/7)\Yb_1 + I(3/7<U\leq 6/7)\Yb_2+ I(U>6/7)\Yb_3,
	\end{equation}
	where $U\sim {\rm U}_{[0, 1]}$ and  $\Yb_i\sim {\rm BSvM}(\mub_i,\kappab_i,\lambda_i)$, $i=1,2,3$  
	are mutually independent; we set $\mu_1 = (-\pi/2, \pi/2)^\top$, $\kappab_1 = (4, 4)^\top$, $\lambda_1 = -2$,  $\mu_2 = (\pi/2, \pi/2)^\top$, $\kappab_2 = (4, 4)^\top$, $\lambda_2 = 2$, $\mu_3 = (0, -\pi/5)^\top$, $\kappab_3 = (6, 6)^\top$, and $\lambda_3 = 0$;
	\item[(S1)] a von Mises-Fisher $\text{vMF}_2(\mub =(0, 0, 1)^\top, \kappa =10)$ distribution, with density of the form
	$\yb \mapsto C_{\kappa}\exp(\kappa\yb^\top\mub)$ 
	where $\mub\in \mathcal{S}^{2}$ and $\kappa\geq 0$ are the location and concentration parameters, respectively, and $C_{\kappa}$ the normalization constant; 
	\item[(S2)] a tangent vMF distribution  \citep{Garcia-Portugues2020} $\text{TvMF}_2(\mub, G, \nub, \kappa)$  with location~$\mub$, angular function~$G$, skewness direction $\nub$, and skewness intensity~$\kappa$, i.e., the distribution of 
	$\Yb = V_G\mub + \sqrt{1 - V_G^2} {\boldsymbol \Gamma}_{\mub} {\bf U}$, 
	where $\Gammab_{\mub}$ is a $3 \times2$ semi-orthogonal matrix  such that $\Gammab_{\mub} \Gammab_{\mub}^\top = \Ib_3 - {\mub} {\mub}^\top$and $\Gammab_{\mub}^\top \Gammab_{\mub} = \Ib_{2}$, and $V_G$ 
	an absolutely continuous scalar r.v. (with values in $[-1,1]$ and density characterized by~$G$) independent of ${\bf U} \sim {\rm vMF}_{1}(\nub, \kappa)$;  we set~$\mub = (0, 0, 1)^\top$,  ${\mbf \nu} = (0.7, \sqrt{0.51})^\top$, $\kappa = 10$, and $V_G= 2\widetilde{V} -1$ with~$\widetilde{V}\sim\text{\rm Beta}(2, 8)$;
	\item[(S3)] a mixture of three vMF's, with random variable 
	\begin{equation}\label{eq.3MixvMF}
		\Yb = I(U \leq 0.3)\Yb_1 + I(0.3<U\leq 0.6)\Yb_2+ I(U>0.6)\Yb_3,
	\end{equation}
	where $U\sim {\rm U}_{[0, 1]}$, and $\Yb_i\sim {\rm vMF}_2(\mub_i, \kappa_i)$, $i=1,2,3$ are mutually independent; we set~$\mub_1 = (0.3, 0.4, \sqrt{0.75})^\top$, $\kappa_1 = 20$, $\mub_2 = (-0.3, -0.4, \sqrt{0.75})^\top$, $\kappa_2 = 20$,\linebreak  $\mub_3 = (-0.3, 0.2, \sqrt{0.87})^\top$, and $\kappa_3 = 20$.
\end{enumerate}

In the simulation, we set $n=4001$, $n_0=1$, $n_R = 40$, and $n_S=100$. Figure~\ref{Fig:Contour} shows the empirical quantile contours ${\mathcal C}\n (r/(n_R + 1))$ ($r = 0, 5, 10, 20, 28$)  for the 2-torus $\mathcal{T}^2$ (flat square representations  in the left panels,  ${\mathbb R}^3$-embedded representations in the central panels) and the 2-sphere ${\cal S}^2$ (right panels). 
\begin{figure}[htbp!]
	\centering
	\begin{tabular}{ccc}
		\hspace{-10mm}    \includegraphics[scale=0.28]{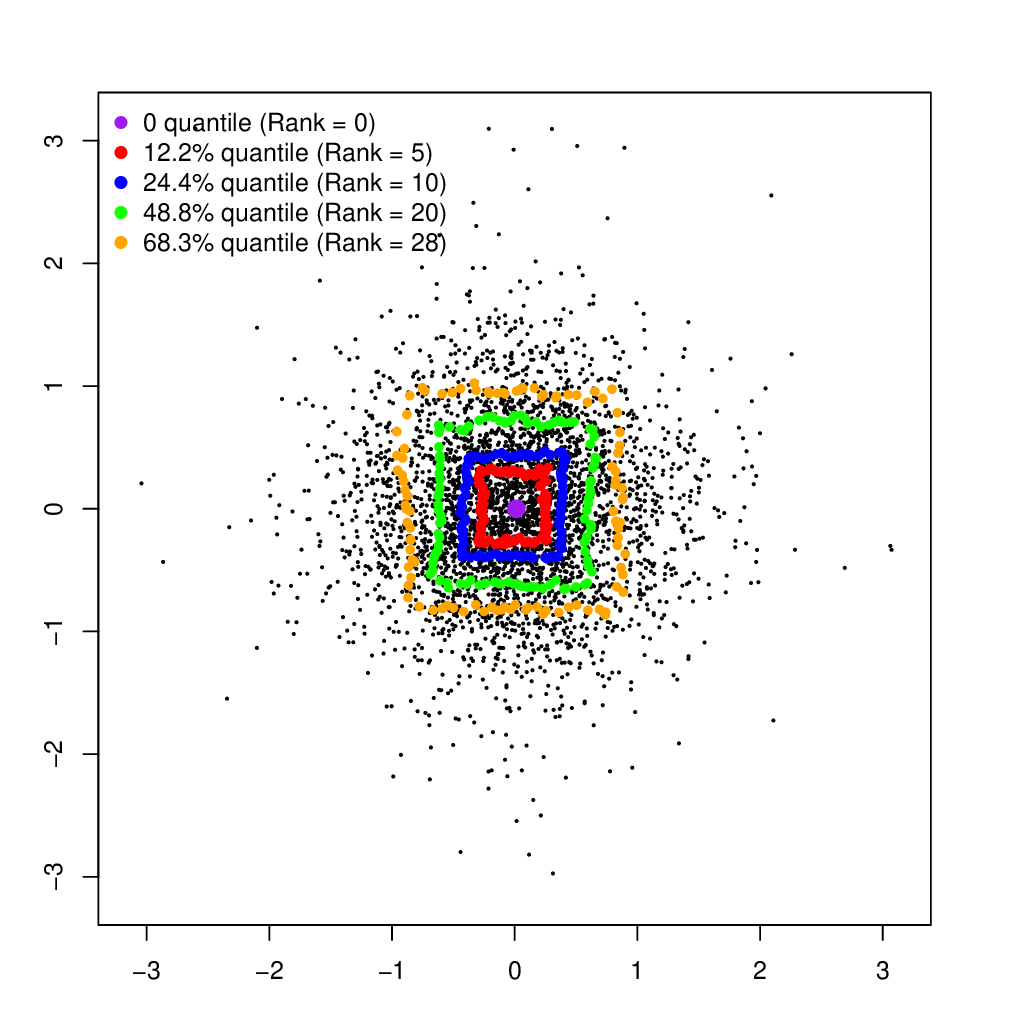} \vspace{0mm}&
		\raisebox{8mm}{ \hspace{-10mm}\includegraphics[scale=0.12]{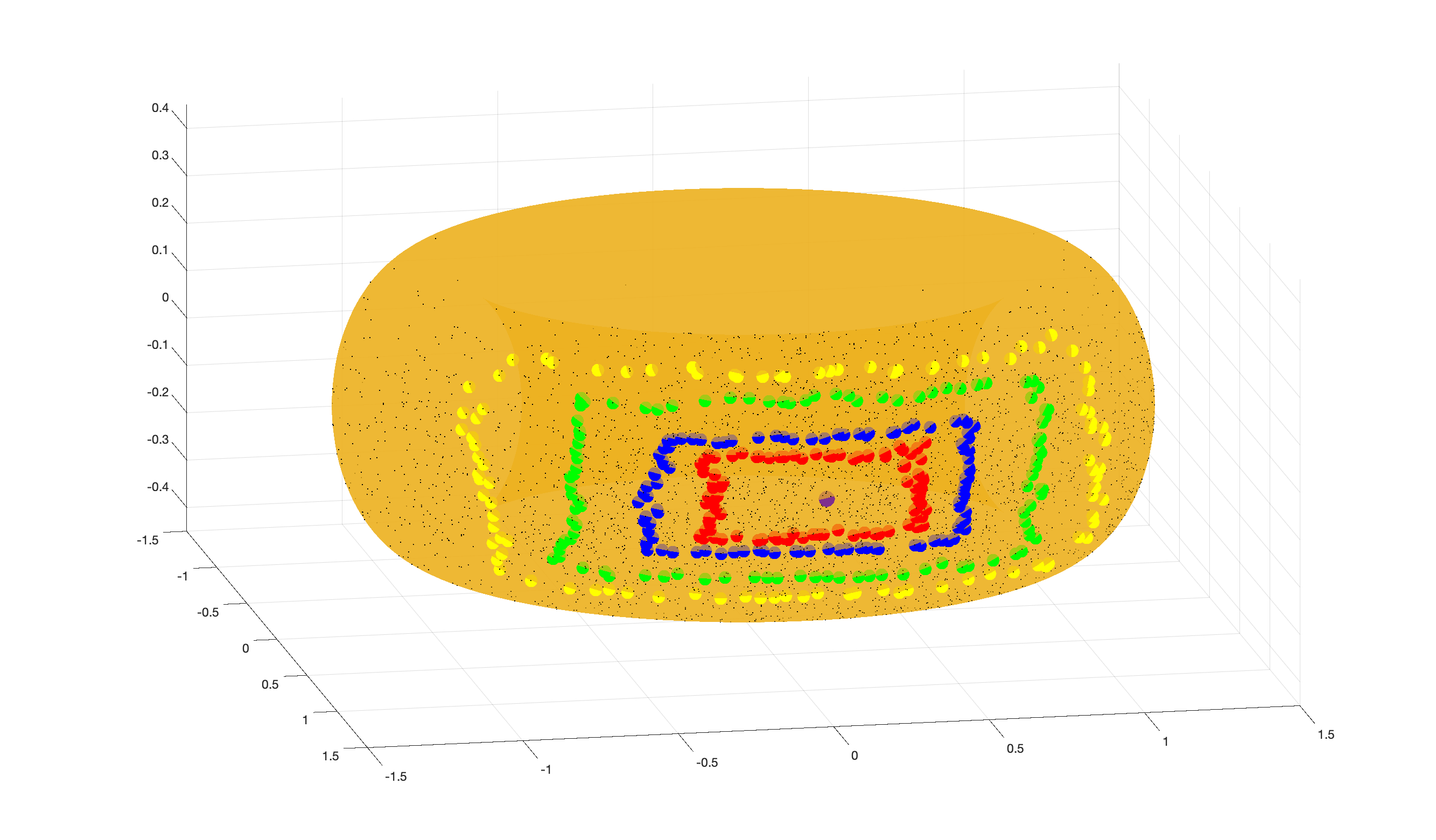}}\vspace{-6mm}&{\hspace{-20mm} \includegraphics[trim=70 20 70 60,clip, width=0.3\textwidth]{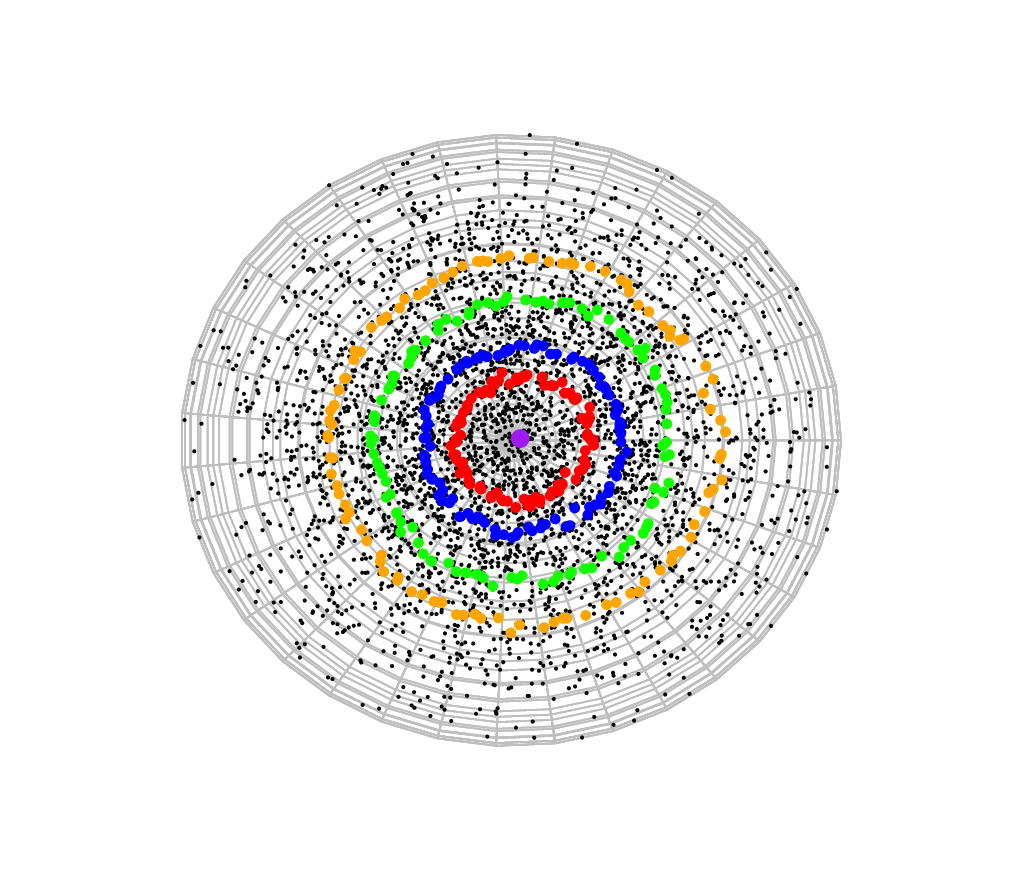}} 
		\\
		
		\hspace{-10mm}        \includegraphics[scale=0.28]{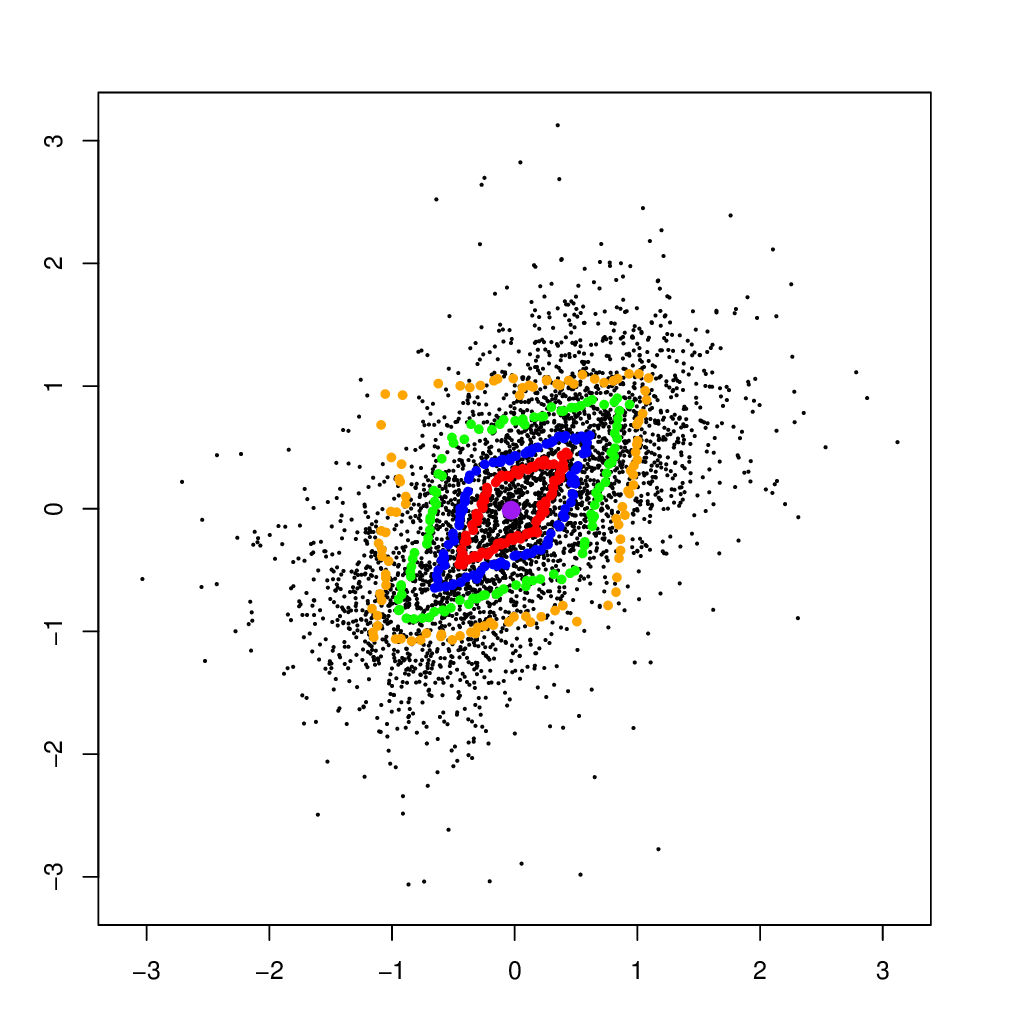}\vspace{0mm}& 
		\raisebox{8mm}{\hspace{-10mm} \includegraphics[scale=0.125]{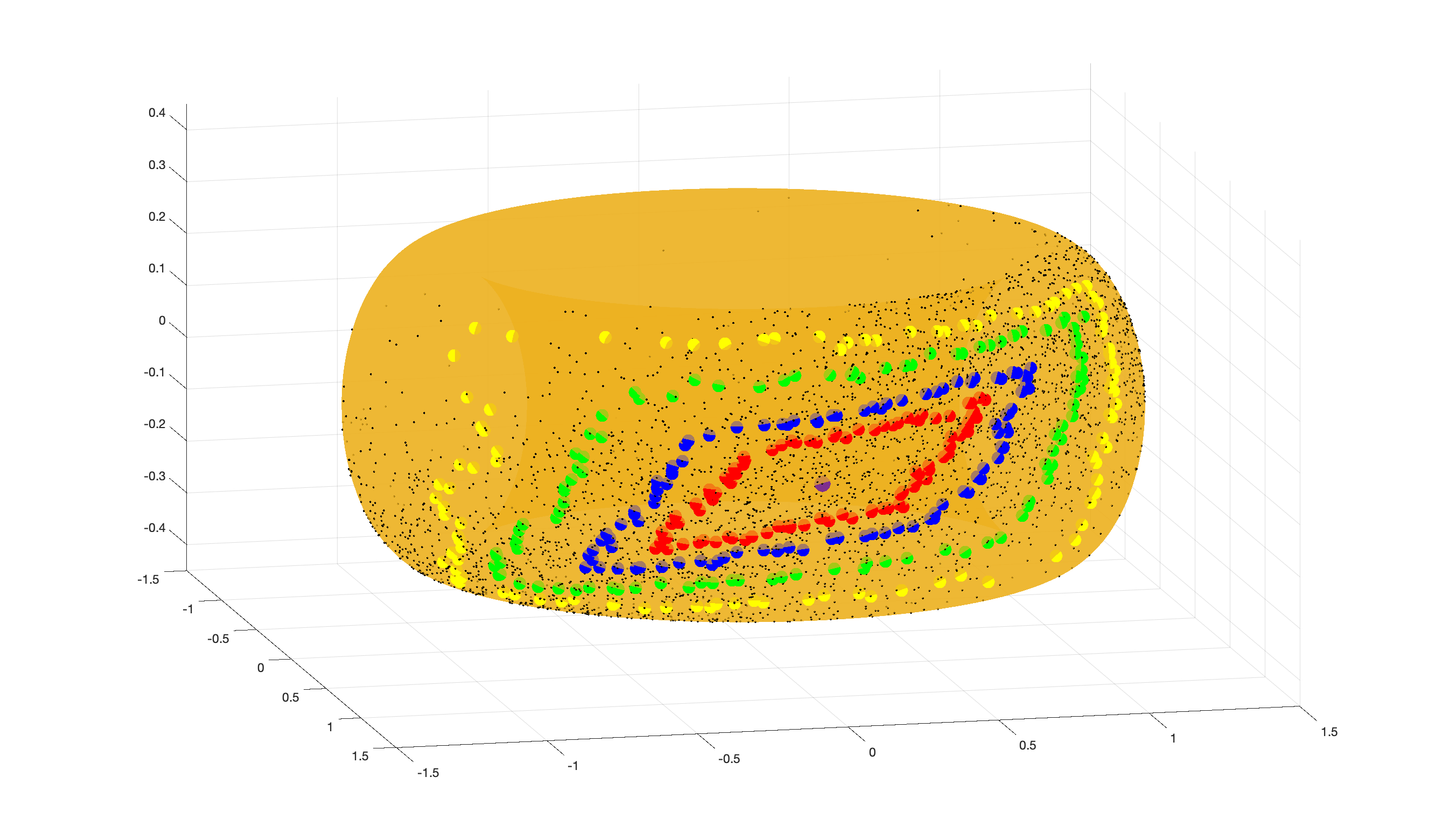}}\vspace{-3mm}
		&{\hspace{-10mm}  \includegraphics[trim=50 0 0 20,clip, width=0.4\textwidth]{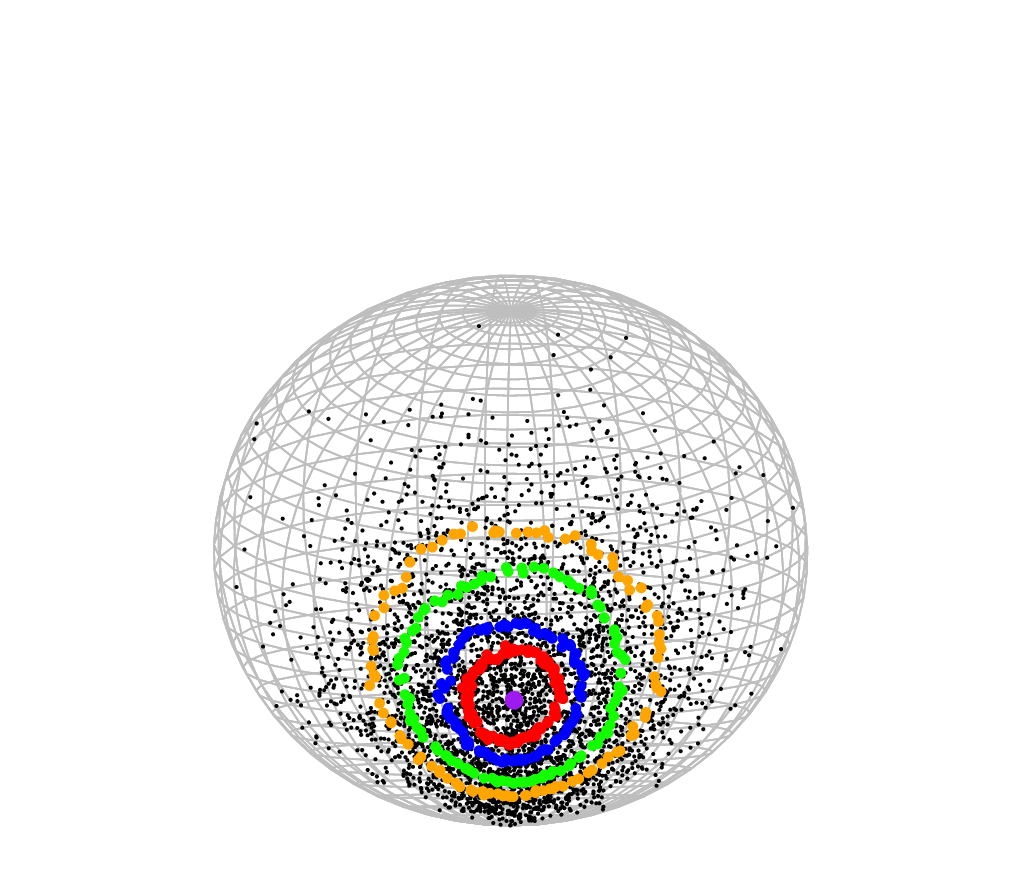}}
		\\
		\hspace{-10mm}    \includegraphics[scale=0.28]{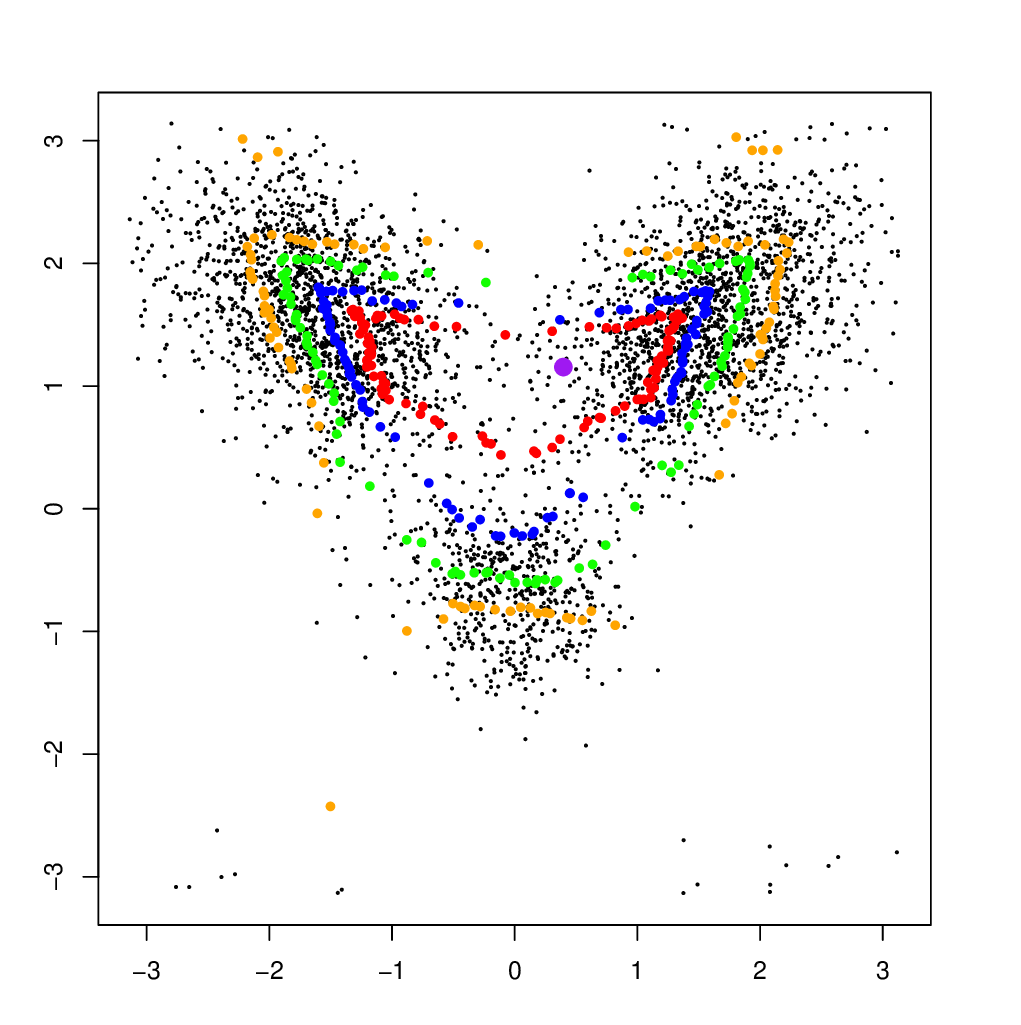}&
		\raisebox{8mm}{\hspace{-10mm}   \includegraphics[scale=0.13, trim={35 0 0 0},clip]{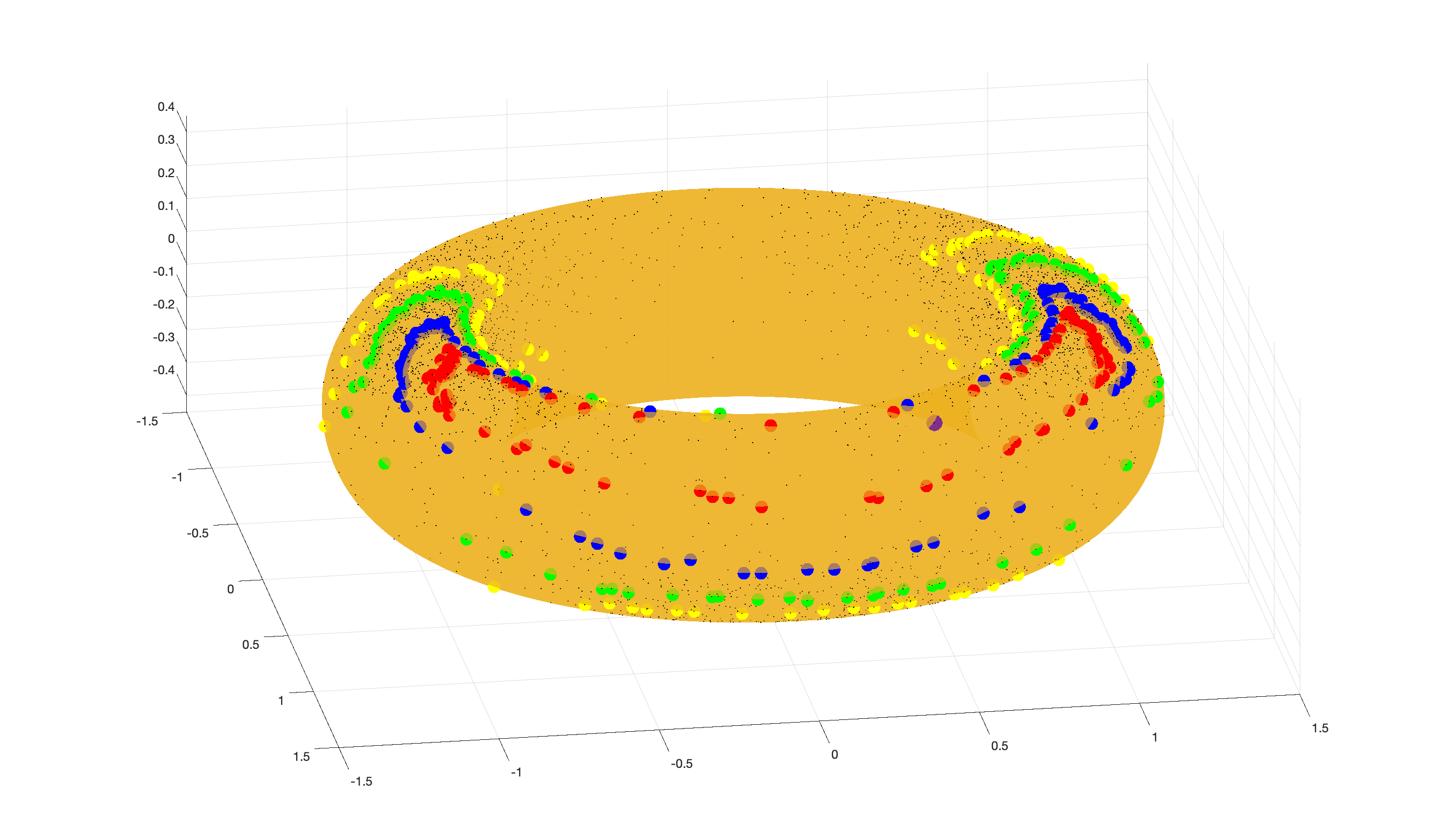}}
		&{\hspace{-17mm} \includegraphics[trim=70 70 70 70,clip,scale=0.7, width=0.41\textwidth]{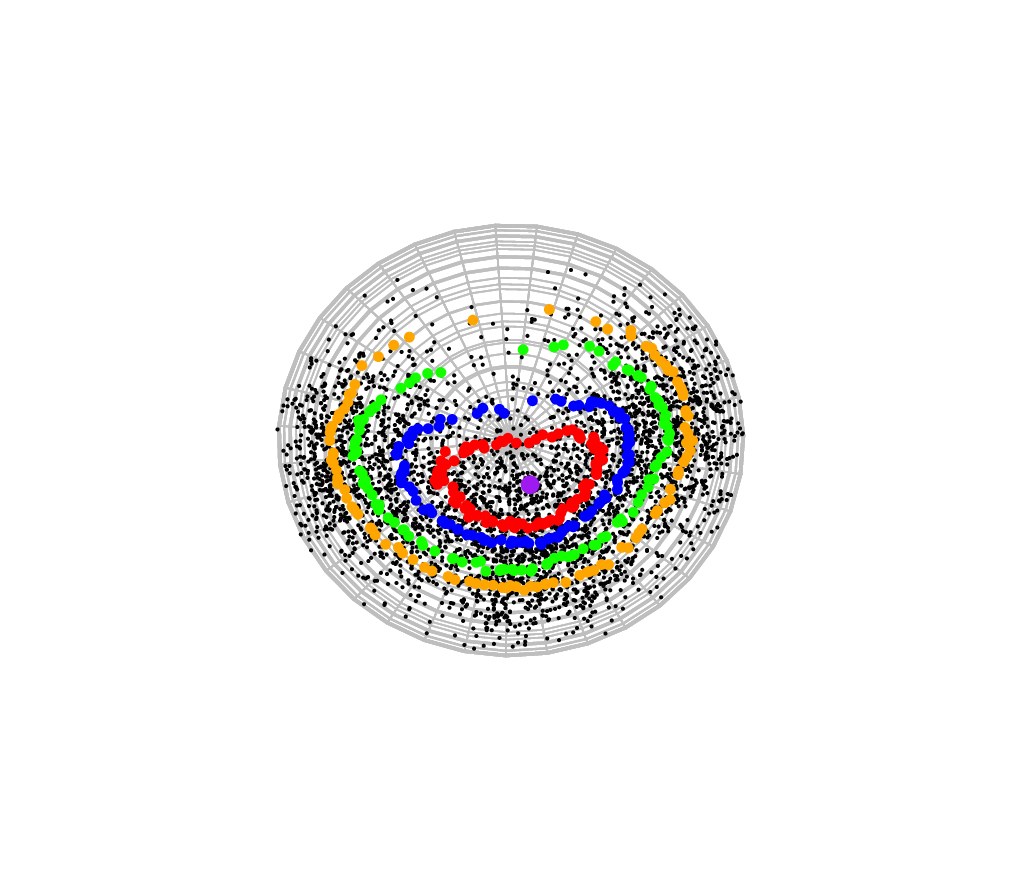}}
	\end{tabular}\vspace{-5mm}
	\caption{\small Empirical quantile contours ${\mathcal C}\n (r/(n_R + 1))$ ($r = 0, 5, 10, 20, 28$), with ${\cal M}_0$ being a singleton, for the 2-torus $\mathcal{T}^2$ (flat square representations  in the left panels,  ${\mathbb R}^3$-embedded representations in the central panels) and the 2-sphere ${\cal S}^2$ (right panels), computed from~$n=~\!4001$ i.i.d.\  observations with distributions (T1) (left and central top panels), (T2) (left and central middle panels), (T3) (left and central bottom panels), (S1) (right top panel), (S2) (middle right panel), and (S3) (bottom right panel); $n_0=1$, $n_R = 40$, and $n_S=100$.}
	\label{Fig:Contour}
\end{figure}

In all cases, the empirical quantile regions are nicely nested around the central point~${\mathcal C}\n (0)$. 
For distribution (T1), a visual inspection of the quantile contours reveals the independence of the two  ${\cal S}^1$ marginals. For (T2), the shape of the contours correctly unveils their mutual dependence. For (T3), the quantile contours nicely pick up the non-convex shape of the mixture distribution. As for the (S1)-(S3) data, the quantile contours also capture the rotational symmetry of (S1), the skewness of (S2),  and the non-convexity of (S3).


\subsubsection{Equatorial strip-type quantile regions}
 Regular grids of the equatorial strip type for $\mathcal{T}^2$ and  $\mathcal{S}^2$ are shown in the bottom row of  Figure~\ref{Fig:Grid}, with~$n=121$, $n_0=13$, $n_R = 3$, and~$n_S=36$. We plot ${\mathfrak{G}}\n_{\widehat{\cal M}_0\n}$, ${\widehat{\cal M}_0\n}$ and ${\mathcal C}^{\rm U}_{\widehat{\cal M}\n_0} (r/(n_R + 1))$, $r = 1, 2, 3$. Unlike the cap-type regions in the first row, which are nested around a central point, the regions enclosed by ${\mathcal C}^{\rm U}_{\widehat{\cal M}\n_0} (r/(n_R + 1))$ for $\mathcal{T}^2$ and  $\mathcal{S}^2\vspace{-2mm}$ are  nested    around, and homotopic-equivalent to, the ``equator'' ${\widehat{\cal M}_0\n}$ (the purple dashed lines/loops).

\begin{figure}[htbp!]
	\centering
	\begin{tabular}{ccc}
		\includegraphics[scale=0.25]{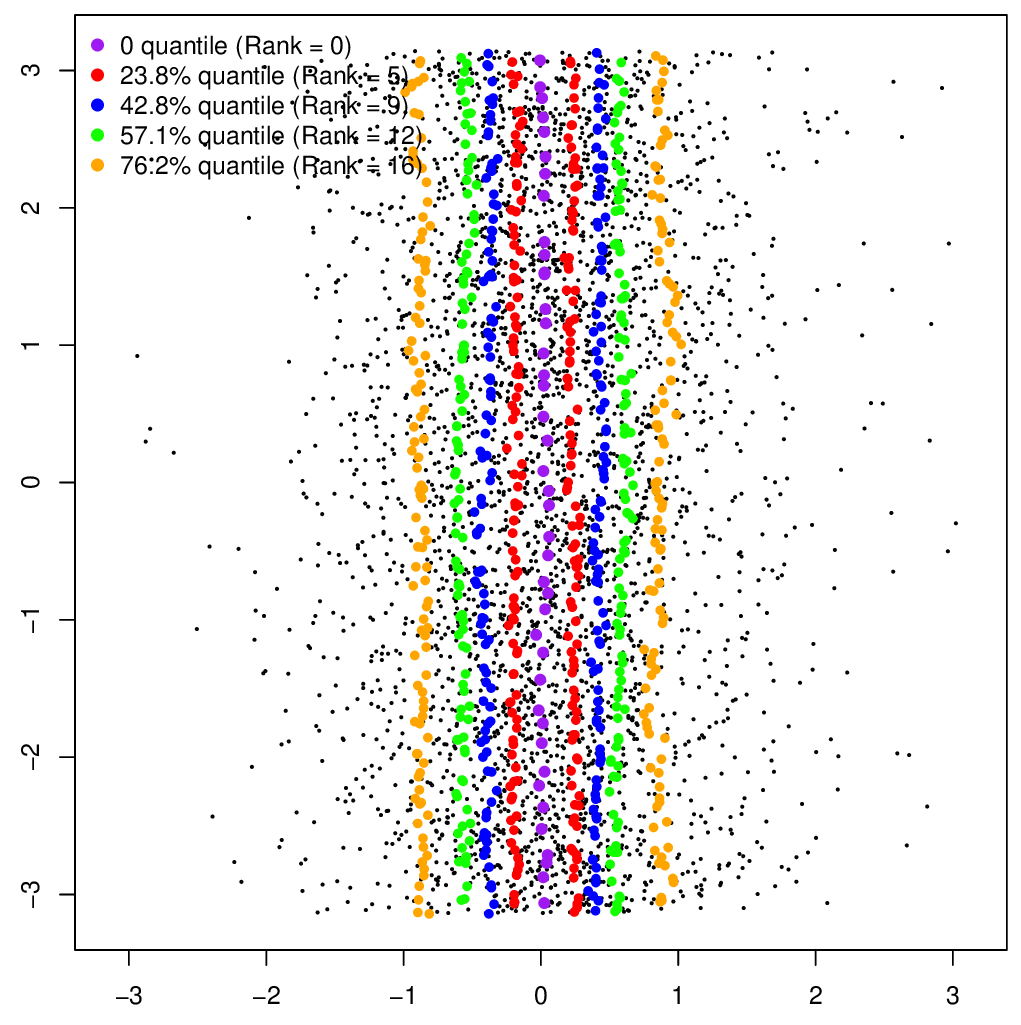} &
		\raisebox{8mm}{ \hspace{-3mm}\includegraphics[scale=0.12]{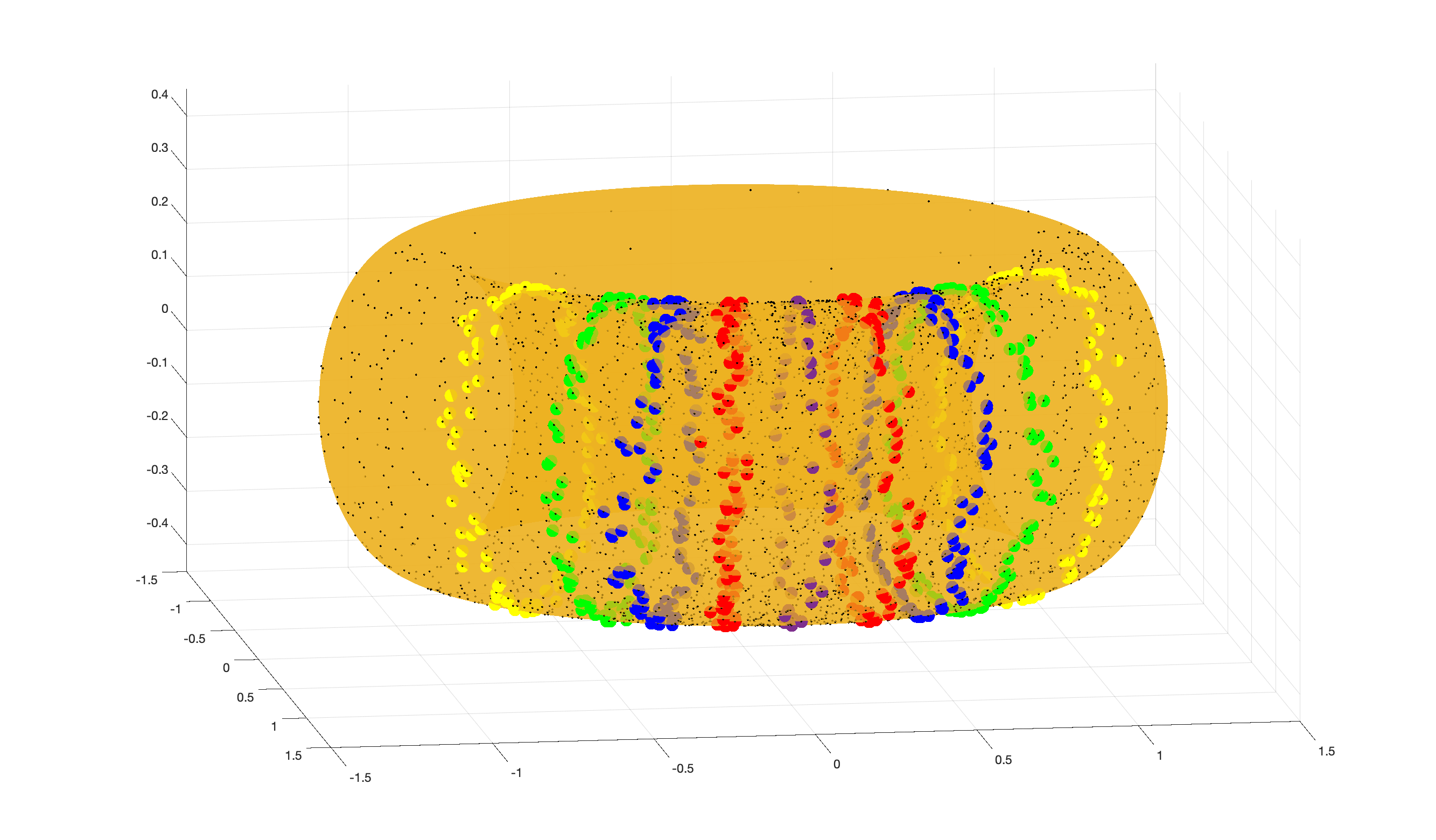}}\vspace{-9mm}&\raisebox{-4mm}{\hspace{-13mm} \includegraphics[trim=20 0 0 20,clip,scale=0.28, width=0.4\textwidth]{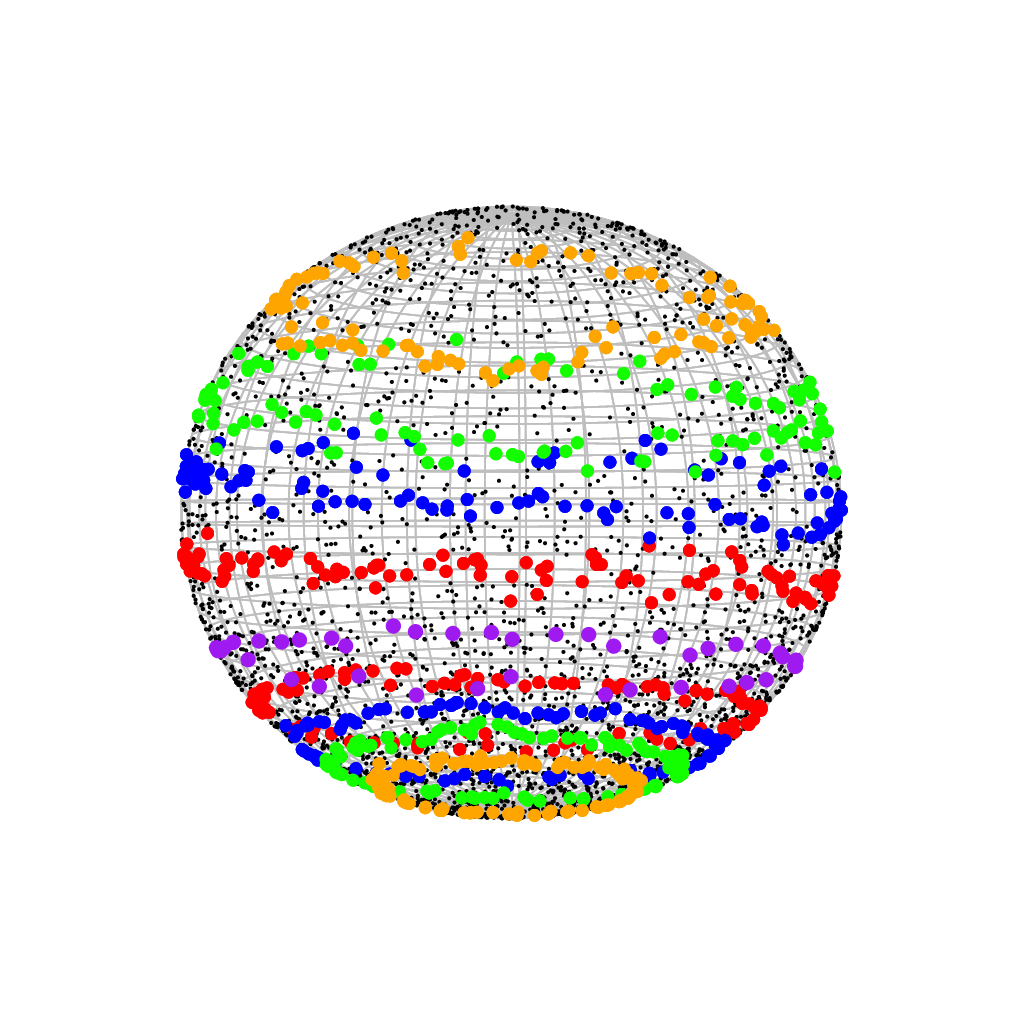}} 
		\\
		
		\includegraphics[scale=0.25]{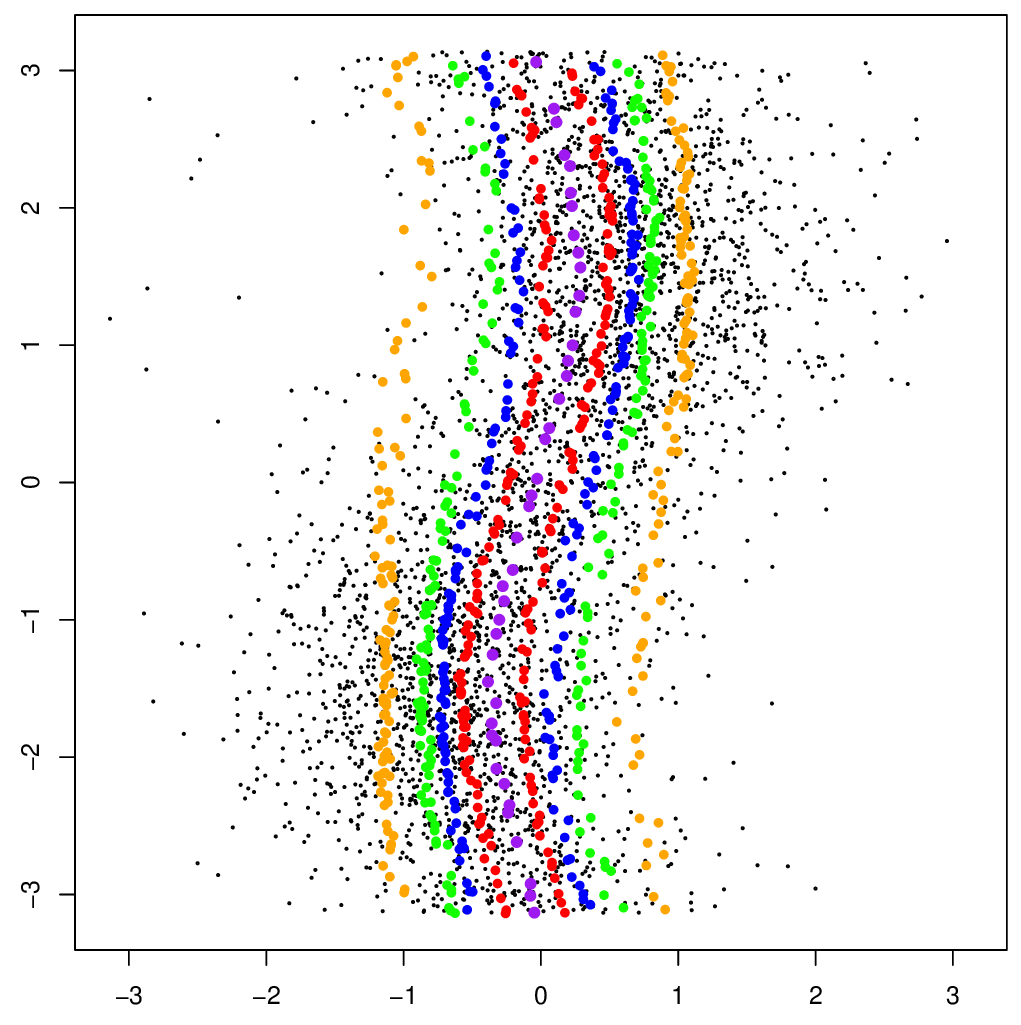}\vspace{-6mm}&
		\raisebox{8mm}{\hspace{-3mm}   \includegraphics[scale=0.12, trim={35 0 0 0},clip]{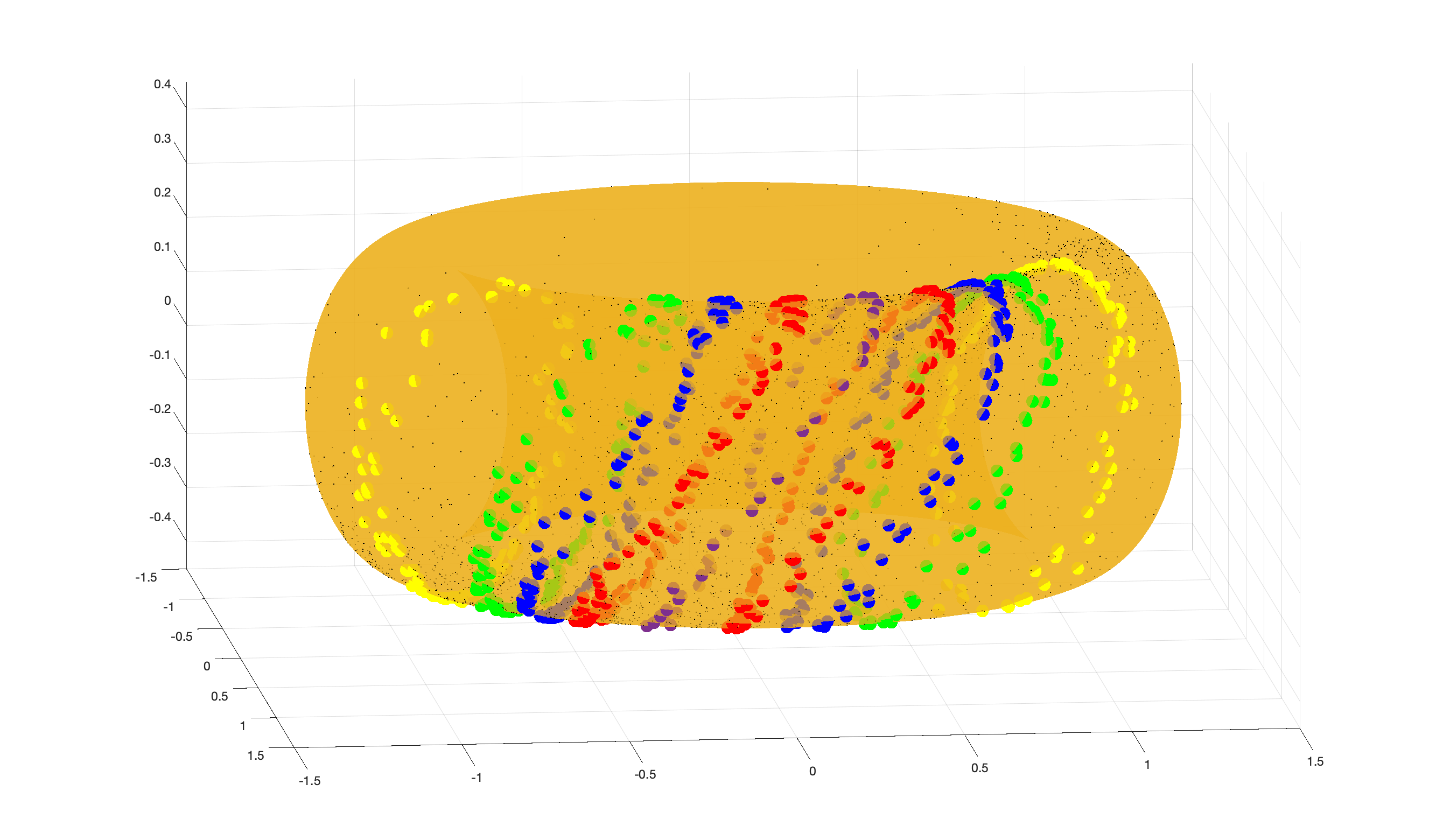}}\vspace{-5mm}
		&\raisebox{-9mm}{\hspace{-13mm} \includegraphics[trim=70 0 70 20,clip,scale=0.2, width=0.3\textwidth]{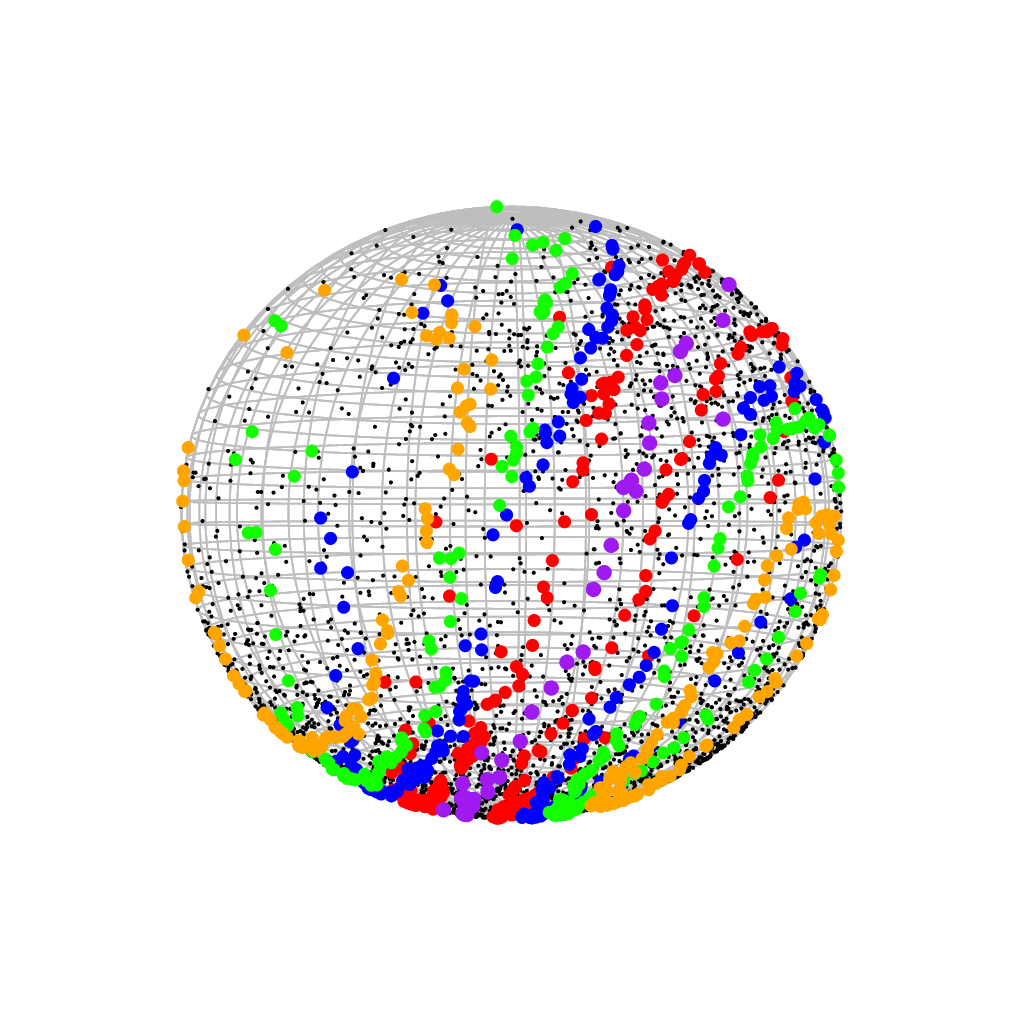}}\vspace{-1mm}\\
		
		\includegraphics[scale=0.25]{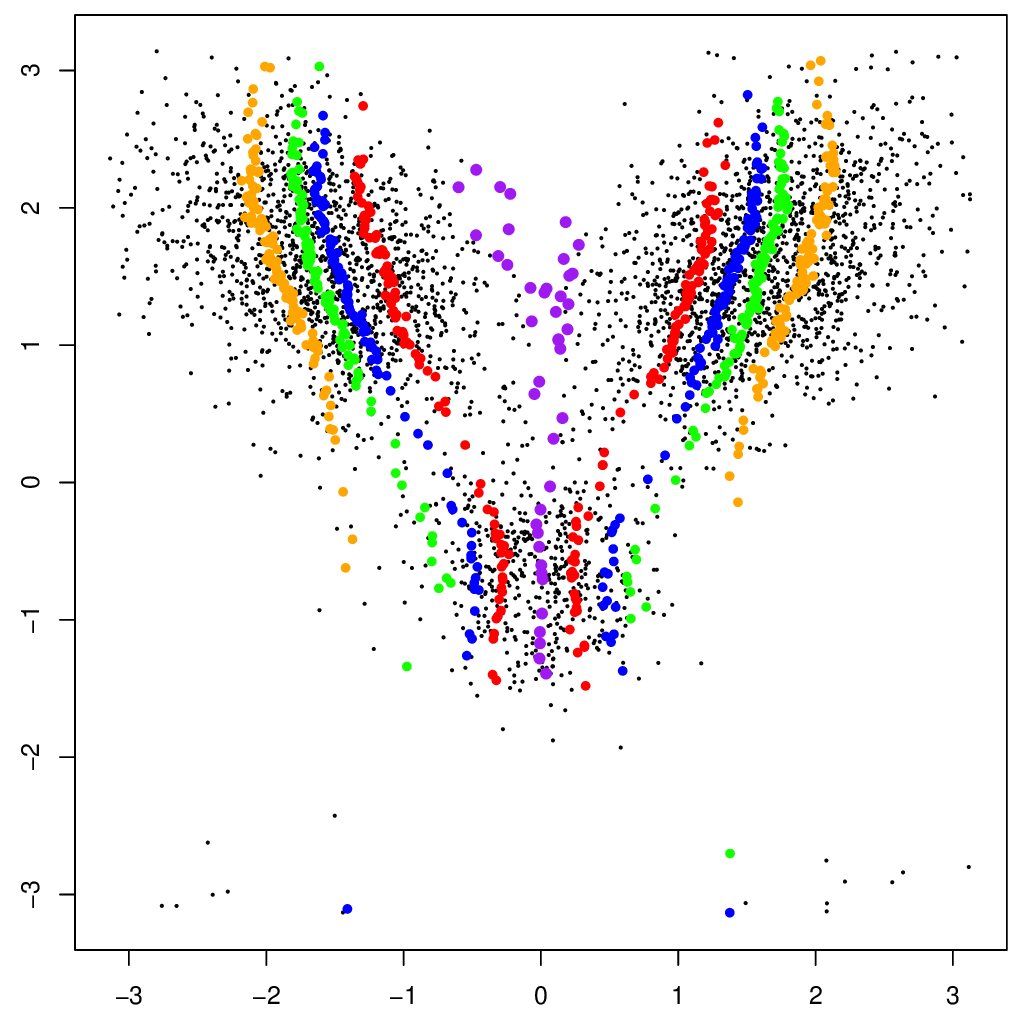}\vspace{-3mm}&
		\raisebox{8mm}{\hspace{-3mm}   \includegraphics[scale=0.13, trim={35 0 0 0},clip]{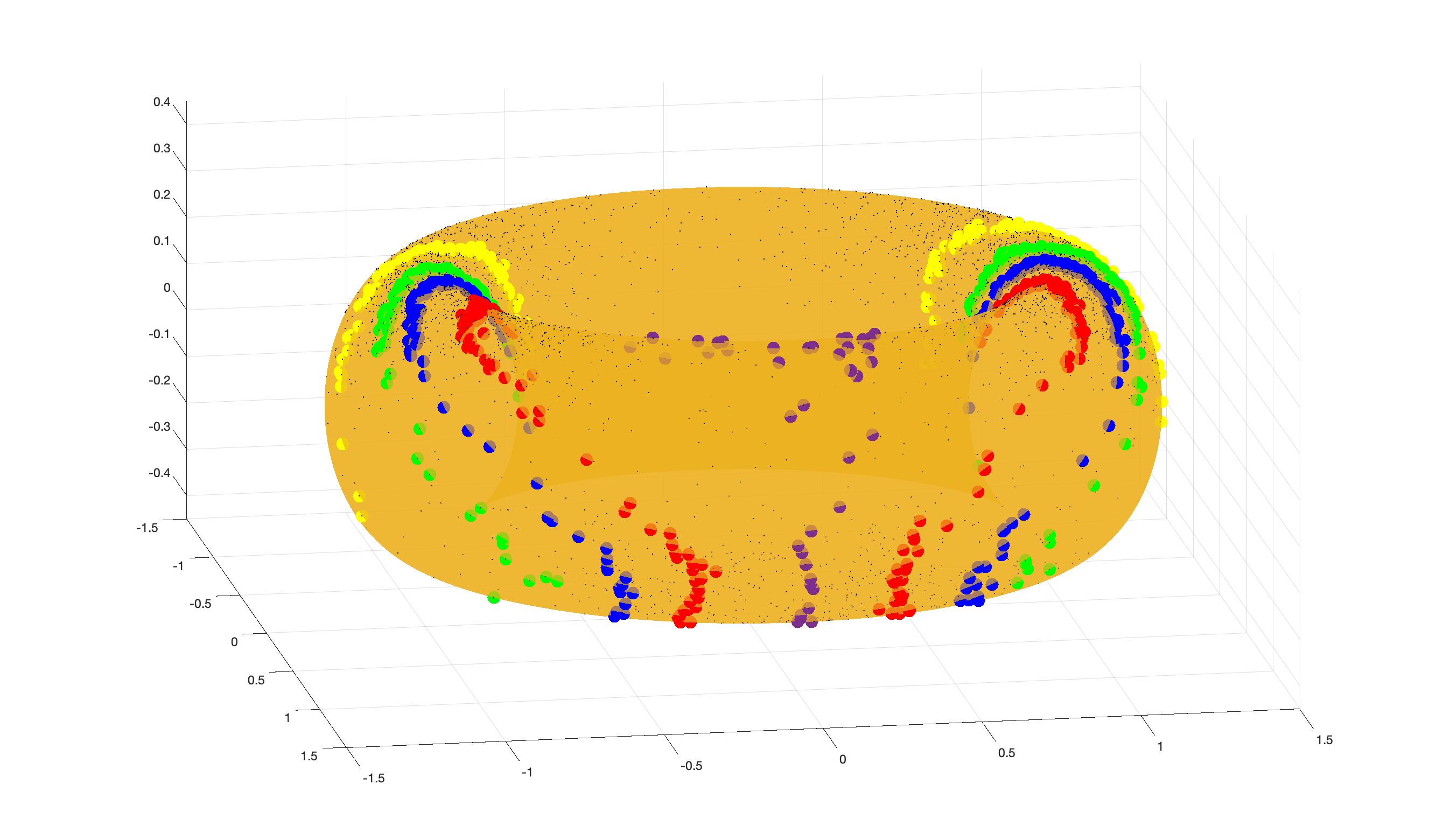}}\vspace{-3mm}
		&\raisebox{-12mm}{\hspace{-13mm} \includegraphics[trim=70 70 70 70,clip,scale=0.7, width=0.41\textwidth]{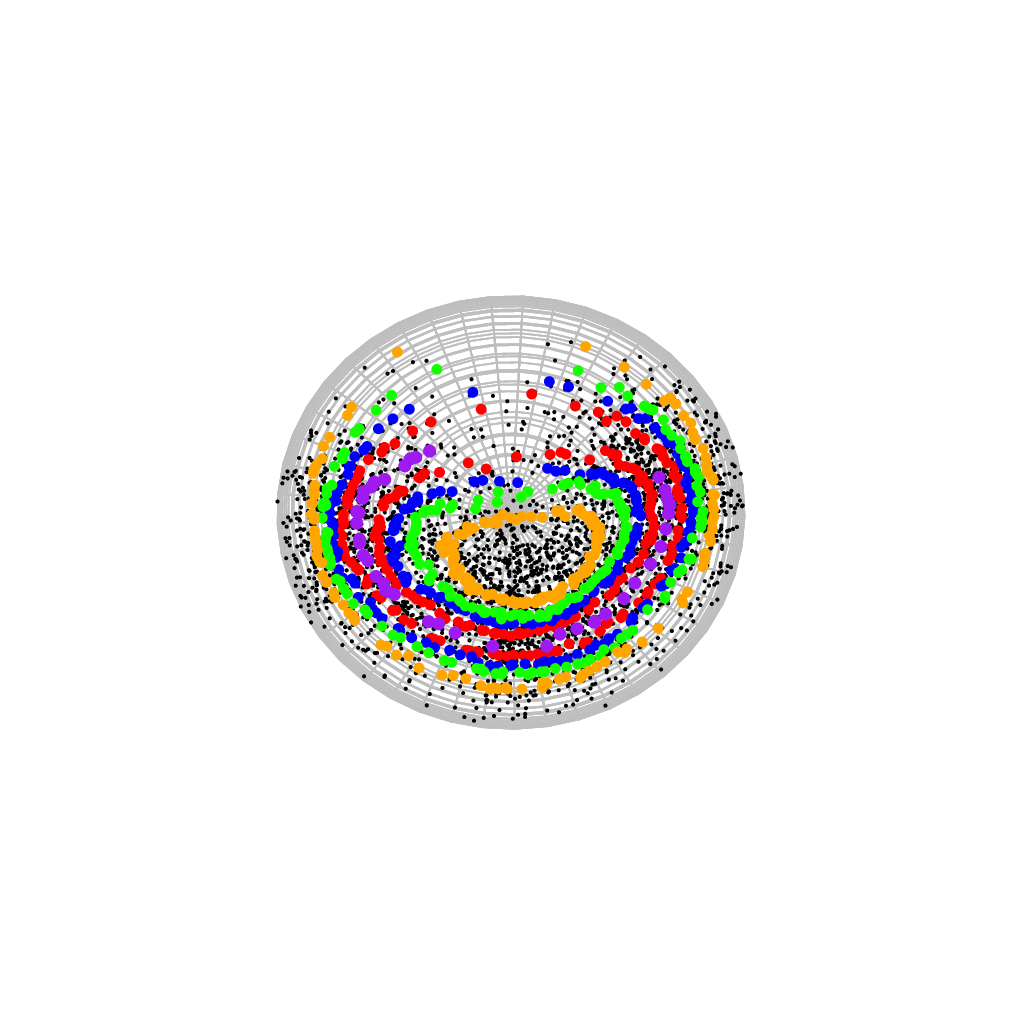}}\vspace{-3mm}
		
	\end{tabular}\vspace{-3mm}
	\caption{\small Empirical quantile contours ${\mathcal C}\n (r/(n_R + 1))$ ($r = 0, 5, 9, 12, 16$), with ${\cal M}_0$ of dimension~$(p-1)$,  for the 2-torus $\mathcal{T}^2$ (flat square representations  in the left panels,  ${\mathbb R}^3$-embedded representations in the central panels) and the 2-sphere ${\cal S}^2$ (right panels), computed from~$n=4001$ i.i.d.\  observations with distributions (Ta) (left and central top panels), (Tb) (left and central middle panels), (Tc) (left and central bottom panels), (Sa) (right top panel), (Sb) (middle right panel), and (Sc) (bottom right panel); $n_0=41$, $n_R = 20$, and $n_S=198$.}
	\label{Fig:StripCont}
\end{figure}

Next, we illustrate the concepts of  ranks and empirical quantile contours for the same equatorial strip-type quantile regions through simulation under 

\begin{enumerate}
	\item[(Ta)] a BSvM$(\mub= {\bf 0}, \kappab = (2.5, 0)^\top, \lambda =0)$ distribution;
	\item[(Tb)] a BSvM$(\mub= {\bf 0}, \kappab = (2.5, 0)^\top, \lambda =2)$ distribution with  two dependent components;
	\item[(Tc)]the mixture of three BSvM's as in (T3);
	\item[(Sa)] a mixture of two vMF's, with random variable 
	\begin{equation}\label{eq.Sph2Mix}
		\Yb = I(U \leq 0.3)\Yb_1 + I(U>0.3)\Yb_2,
	\end{equation}
	where $U\sim {\rm U}_{[0, 1]}$, and $\Yb_i\sim {\rm vMF}_2(\mub_i, \kappa_i)$, $i=1,2$ are mutually independent; we set~$\mub_1 = (0, 0, 1)^\top$, $\kappa_1 = 1$, $\mub_2 = (0, 0, -1)^\top$, and $\kappa_2 = 2$; 	
	\item[(Sb)] a mixture of a tangent vMF and a vMF, with r.v. $\Yb$ of the form \eqref{eq.Sph2Mix}, where $\Yb_1 \sim \text{TvMF}_2(\mub_1, G, \nub, \kappa_1)$ and $\Yb_2\sim {\rm vMF}_2(\mub_2, \kappa_2)$. We set $\mub_1 = (0, 0, 1)^\top$,  $\nub = (0.7, \sqrt{0.51})^\top$, $\kappa_1 = 2$, $V_G = 2\widetilde{V} -1$ with~$\widetilde{V}\sim\text{\rm Beta}(5, 2)$, $\mub_2 = (0, 0, -1)^\top$, and~$\kappa_2 = 3$; 	
	\item[(Sc)] a mixture of three vMF's, with random variable $\Yb$ as in (S3).
\end{enumerate}

We set $n = 4001$, $n_0=41$, $n_R = 20$, $n_S=198$, and plot in Figure~\ref{Fig:StripCont} the empirical quantile contours ${\mathcal C}\n (r/(n_R + 1))$ ($r = 0, 5, 9, 12, 16$). For distribution (Ta), where the variable~$\Yb = (\Yb_1^\top, \Yb_2^\top) ^\top$ has independent components $\Yb_1 \sim \text{vMF}_1((0, 1)^\top, 2.5)$ and~$\Yb_2  \sim~\!\rm{U}_{\mathcal{S}^1}$, the empirical quantile contours distribute symmetrically over the~$\Yb_1$ coordinate and uniformly over the $\Yb_2$ coordinate. For (Tb), the strips nested around~${\mathcal C}\n (0)$ are skewed, reflecting the positive dependence of the two components. For (Tc), the non-convexity of the mixture distribution are well captured by the shape of the strips. Turning to (Sa), the distribution is symmetric and concentrates around the south and north poles, and is more dense in the southern hemisphere. The empi\-rical quantile contours capture all these properties. For (Sb), where the first mixture component of (Sa) is replaced by a tangent vMF one, the skewness is also reflected by the shape of the quantile contours. Finally, as in the torus case,  the non-convexity of the mixture distribution (Sc) is  also is graphically visible.

\section{Riemannian quantile regression}\label{sec.QuanReg}
A natural and substantive application of the concepts of quantile, quantile region, and quantile contour developed in Sections~\ref{sec.DistQuan} and~\ref{sec.ranksign} is the possibility of performing quantile regression for response variables taking values in Riemaniann manifolds. 

   In the sequel, we restrict to response variables $\Yb$ with values in the manifolds $\cal M$ described in Examples~\ref{ex.sphere}-\ref{ex.product}, for which the distribution function $\Fb$ and the quantile function $\Qb$ are well-defined and continuous. As for the covariate $\Xb$, we assume it lives on a Polish metric space $({\cal M}_{\Xb}, d_{\Xb})$ equipped with the corresponding $\sigma$-field.  Assuming that~$({\Xb},\Yb)$ has distribution~${\rm P}_{\Xb,\Yb}$ over~${\cal M}_{\Xb}\times\cal M$ (equipped with the product   sigma-field),  
    denote by ${\rm P}_{\Yb\vert \Xb = \xb}$ the distribution of $\Yb$ conditional on $\Xb = \xb$, $\xb\in{\cal M}_{\Xb}$, by~${\bf y}\mapsto{\bf F} (\yb\vert\,  \xb)$ and ${\bf y}\mapsto{\bf Q}(\yb\vert\,   \xb)$, with~$(\xb,\yb)\in {\cal M}_{\Xb}\times \cal M$, the corresponding conditional distribution and quantile functions.  The conditional distributions ${\rm P}_{\Yb\vert \Xb = \xb}$, moreover, characterize, for each $\xb\in{\mathcal  M}_X$ and $\tau\in [0,1]$, {\it conditional quantile regions} and {\it contours}, which we denote as  $\mathcal{C}^{\rm P}_{\mathcal{M}_0(\xb)}(\tau \vert\,\xb)$   and  ${\mathbb C}^{\rm P}_{\mathcal{M}_0(\xb)}(\tau\vert\, \xb)$, respectively,  for some sub\-manifold~${\mathcal{M}_0(\xb)} \subset~\!\mathcal{M}$ that depends on $\xb$.
   
  The objective of quantile regression, in this context, is the nonparametric estimation, based on a sample~$(\Xb\n, \Yb\n) := \{(\Xb\n_1, \Yb\n_1), \ldots, (\Xb\n_n, \Yb\n_n)\}$ of~$n$ i.i.d.\ copies of~$(\Xb,\Yb)$, of ${\bf F} (\cdot \,\vert\,  \xb)$,  ${\bf Q}(\cdot\,\vert\,   \xb)$, $\mathcal{C}^{\rm P}_{\mathcal{M}_0(\xb)}(\tau \vert\,\xb)$,   and  ${\mathbb C}^{\rm P}_{\mathcal{M}_0(\xb)}(\tau\vert\, \xb)$      for all $\xb \in{\mathcal M}_X$ and~$\tau\in~\![0,1]$.




\subsection{Empirical conditional distributions on $\cal M$}

 Consider a sequence of weight functions $\wb\n: {\cal M}_{\Xb}^{n+1} \rightarrow \mathbb{R}^n$ that are measurable with respect to $\xb \in {\cal M}_{\Xb}$ and the sample~$\Xb\n := (\Xb_1, \ldots, \Xb_n)$, of  the form
\begin{equation*} \wb\n : 
(\xb, \Xb\n) \mapsto \wb\n(\xb, \Xb\n) := (w_1\n(\xb, \Xb\n), \ldots, w_n\n(\xb, \Xb\n))
\end{equation*}
where $w_j\n: {\cal M}_{\Xb}^{n+1} \rightarrow \mathbb{R}$, $j = 1, \ldots, n$ are such that
\begin{equation*}
w_j\n(\xb, \Xb\n) \geq 0  \quad {\rm and} \quad \sum_{j=1}^n w_j\n(\xb, \Xb\n) = 1 \quad \text{a.s.~for all  $n$ and $\xb$}.
\end{equation*}
Define the empirical conditional distribution of $\Yb$ given $\Xb = \xb$ as
\begin{equation}\label{empiricalDF}
{\rm P}\n_{\wb\n(\xb)} := \sum_{j=1}^n w_j\n(\xb, \Xb\n) \delta_{\Yb_j}\qquad \xb\in{\mathcal M}_0.
\end{equation}
Adopting the terminology used by \cite{Stone1977} for $\xb$ in the Euclidean space to our setting where $\xb \in {\cal M}_{\Xb}$, we say that $\wb\n$ is a {\it consistent} weight function if, when\-ever~$(\Xb, Z), (\Xb\n_1, Z\n_1), \ldots, (\Xb\n_n, Z\n_n)$ are i.i.d.\  with $Z$ real-valued  and ${\rm E}\vert Z\vert^r <\infty$ for some~$r >1$,
\begin{equation}\label{eq.cstWeight}
{\rm E} \left\vert  \sum_{j=1}^n w_j\n(\Xb, \Xb\n) Z_j\n - {\rm E}(Z\vert\Xb)  \right\vert^r \rightarrow 0 \quad {\rm as} \quad n\rightarrow \infty;
\end{equation}
see  Section~\ref{sec.Weight} for   examples  and sufficient conditions. By definition, for any consistent weight function $\wb\n$ and  any measurable function $\tilde{f}: {\cal M} \rightarrow \mathbb{R}$ satisfying ${\rm E}\vert \tilde{f}(\Yb)\vert^r <\infty$, 
\begin{equation}\label{eq.cstWeight2}
{\rm E} \left\vert  \sum_{j=1}^n w_j\n(\Xb, \Xb\n) \tilde{f}(\Yb_j\n) - {\rm E}(\tilde{f}(\Yb)\vert\Xb))  \right\vert^r \rightarrow 0 \quad {\rm as} \quad n\rightarrow \infty.
\end{equation}

In order to show convergence of the empirical conditional distribution ${\rm P}\n_{\wb\n(\Xb)}$ to the conditional distribution ${\rm P}_{\Yb\vert\Xb}$, let us  introduce some further notation. Let ${\cal F}_{{\rm BL}_1}^{\cal M}$ denote the set of all  real functions on 
 ${\cal M}$ with  Lipschitz norm less than or equal to one, namely
$${\cal F}_{{\rm BL}_1}^{\cal M} :=  \big\{f: {\cal M} \rightarrow \mathbb{R}: \Vert f \Vert_{\infty} \leq 1 \ {\rm and} \ \vert f(\yb_1) - f(\yb_2) \vert \leq d(\yb_1, \yb_2), \forall \yb_1, \yb_2 \in {\cal M}\big\}.$$
Denoting by $\mu$ and $\nu$  two probability measures on ${\cal M}$,   define the {\it bounded Lipschitz distance} between $\mu$ and $\nu$ as
$$d_{\rm BL}(\mu, \nu) := \sup_{f \in {\cal F}_{{\rm BL}_1}^{\cal M}} \left\vert \int_{{\cal M}} f {\rm d}\mu - \int_{{\cal M}} f {\rm d}\nu\right\vert.$$
Then  (see \citet[Theorem 1.12.4]{vaart1996weak}), weak convergence   on~${\cal M}$ is metrizable by~$d_{\rm BL}$.

Proposition~\ref{prop.weight} states that for a consistent weight function $\wb\n$, ${\rm P}\n_{\wb\n(\Xb)}$ is a consistent estimator of the conditional distribution ${\rm P}_{\Yb\vert\Xb}$. See Appendix~\ref{App.proof} for the proof (recall that $(\Xb\n, \Yb\n) := \{(\Xb\n_1, \Yb\n_1), \ldots, (\Xb\n_n, \Yb\n_n)\}$ stands for an i.i.d.\ sample of~$n$ copies of the ${\mathcal M}_X\times{\mathcal M}$-valued random variable~$(\Xb,\Yb)$).

\begin{proposition}\label{prop.weight}
Let  $\wb\n: {\cal M}_{\Xb}^{n+1} \rightarrow \mathbb{R}$ be a consistent weight function. Then, for any $\epsilon >0$,
{${\rm P}\left[{d}_{\rm BL}\big({\rm P}\n_{\wb\n(\Xb)}, {\rm P}_{\Yb\vert\Xb}\big) > \epsilon \right] \rightarrow 0$
as $n\rightarrow \infty$.} 
\end{proposition}

\subsection{Empirical conditional quantile regions and contours on $\cal M$}\label{subsec.EmpQR}

Once  the  empirical conditional distributions  ${\rm P}\n_{\wb\n(\xb)}$ of $\Yb$,  $\xb \in \mathcal{M}_{\Xb}$ are obtained, the empirical conditional quantiles are  computed by solving the optimal coupling problem between the sample and a structured grid constructed via a two-step approach similar to that of Section~\ref{subsec.empQuan}. That grid depends on the empirical distributions defined in \eqref{empiricalDF} and the preliminary estimation of  $\mathcal{M}_0(\xb)$, a step that is not needed in the multiple-output quantile regression approach of  \cite{delBarrio2022} since $\Yb$ there takes values in ${\mathbb R}^d$.  To this end, as in Section~\ref{Sec31}, we make the following assumption on the submanifold $\mathcal{M}_0(\xb)$.

\begin{assumption}\label{Ass.M0x}
		{\rm 
			For each $\xb \in \mathcal{M}_{\Xb}$,  $\mathcal{M}_0(\xb)$ is of the form ${\cal M}_0(\xb)={\cal M}_{\Fb (\thetab(\xb)\vert \xb)}$ for some $\thetab(\xb) \in{\cal M}$ where (i) the mapping ${\bf f}_{\xb} \mapsto {\cal M}_{{\bf f}_{\xb}}$, ${\bf f}_{\xb} \in~\!{\cal M}$ is continuous 
			in the sense that $d({\bf f}\n_{\xb}, {\bf f}_{\xb})\rightarrow 0$ as $n\to\infty$  implies~$d_{\rm H}({\cal M}_{{\bf f}\n_{\xb}}, {\cal M}_{{\bf f}_{\xb}})\to 0$ as $n\to\infty$;
			(ii) a strongly consistent 
			estimator ${\hat\thetab}\n(\xb)$ of~$\thetab(\xb)$ is available.
		}
\end{assumption}

\begin{Rem}
	{\rm
Suppose that,  in Assumption~\ref{Ass.M0x},  ${\cal M}_0(\xb)={\cal M}_{\Fb (\thetab_{\rm Fr}(\xb) \vert \xb)}$ for some ele\-ment~$\thetab(\xb) = \thetab_{\rm Fr}(\xb)$ 
 of the conditional Fr\'echet mean set 
\begin{equation*}
	\mathcal{A}_{\text{\rm Fr}}(\xb) \coloneqq  {\rm argmin}_{\yb \in  \mathcal{M}} {\rm E}_{{\rm P}_{\Yb\vert\Xb=\xb}}[c(\Yb, \yb)],
	\qquad \xb\in{\mathcal M}_0
\end{equation*}
of ${\rm P}_{\Yb\vert\Xb=\xb}$. Let
\begin{equation*}
	\mathcal{A}\n_{\text{\rm Fr}}(\xb)\coloneqq {\rm argmin}_{\yb \in  \mathcal{M}} \sum_{j=1}^n w_j\n(\xb, \Xb\n) c(\Yb_j\n, \yb)  ,  
	\qquad \xb\in{\mathcal M}_0
\end{equation*}
denote the empirical version of $\mathcal{A}_{\text{\rm Fr}}(\xb)$. Proposition~\ref{prop.FrConsistent} states that, for any $\epsilon >0$,  with probability converging to one,  $\mathcal{A}\n_{\text{\rm Fr}}(\xb)$ is a subset of the $\epsilon$-fattening $\mathcal{A}^{\epsilon}_{\text{\rm Fr}} (\xb)$ of $\mathcal{A}_{\text{\rm Fr}}(\xb)$, which is non-empty due to the compactness of~${\cal M}$ (see Appendix~\ref{App.proof} for the proof). This  implies that  any $\hat{\thetab}_{\rm Fr}\n(\xb)\in 	\mathcal{A}\n_{\text{\rm Fr}}(\xb)$ is a consistent estimator of some  conditional Fr\' echet mean $\thetab_{\rm Fr}(\xb)\in 	\mathcal{A}_{\text{\rm Fr}}(\xb)$. Therefore, Assumption~\ref{Ass.M0x} holds for  $\thetab(\xb) = \thetab_{\rm Fr}(\xb)$. 
	}
\end{Rem}

\begin{proposition}\label{prop.FrConsistent}
	Let $\wb\n$ be a consistent weight function. Then, for any $\epsilon >0$ and $\xb \in \mathcal{M}_{\Xb}$,
	$\underset{n\rightarrow\infty}{\lim} {\rm P}(\mathcal{A}\n_{\text{\rm Fr}}(\xb)\subset \mathcal{A}^{\epsilon}_{\text{\rm Fr}}(\xb)) = 1.$
\end{proposition}

Assumption~\ref{Ass.M0x} allows us to construct empirical conditional quantile regions and contours on $\cal M$ via the following two-steps procedure.\smallskip 

\noindent{\bf Step 1.} Denoting by $N$  an arbitrary large number, construct a consistent estimator~$\widehat{\cal M}^{(N)}_0(\xb)$ of $\mathcal{M}_0(\xb)$. This relies on a consistent estimator of $\Fb (\thetab(\xb)\vert \xb)$, which can be achieved by a transport between the empirical conditional distribution ${\rm P}\n_{\wb\n(\xb)}$ and any empirical distribution that converges weakly to ${\rm U}_{\cal M}$. Denote, for instance,\linebreak  by~$\mathfrak{G}^{(N)}_0(\xb) = \{{\scriptstyle{\mathfrak{G}}}_{0i}^{(N)}(\xb); \ i = 1, \ldots, N\}$   an $N$-point grid  such that the  empirical distribution  ${\rm P}^{{\mathfrak{G}}^{(N)}_0(\xb)}\!\!\coloneqq N^{-1}\sum_{i=1}^N \delta_{{\scriptstyle{\mathfrak{G}}}_{0i}^{(N)}(\xb)}$ converges weakly to ${\rm U}_{\cal M}$ as $N\to\infty$.  Based on this grid, solve the Kantorovich formulation of the optimal transport problem, namely,
\begin{align}\label{eq.OT}
	\begin{split}
		& \min_{{\bf p}:= \{p_{ij}\}} \sum_{i=1}^N \sum_{j=1}^n c(\Yb_j\n, {\scriptstyle{\mathfrak{G}}}_{0i}^{(N)}(\xb)) p_{ij}, \\
		{\rm s.t.} \  &\sum_{j=1}^n p_{ij} = \frac{1}{N}, \  i = 1, \ldots, N,\\
		&  \sum_{i=1}^N p_{ij} = w_j\n(\xb, \Xb\n), \  j = 1, \ldots, n, \\
		& p_{ij} \geq 0,  \ i = 1, \ldots, N,\  j = 1, \ldots, n.
	\end{split}
\end{align}
Denoting by ${\bf p}^*(\xb) = \{p_{ij}^*(\xb);  \ i = 1, \ldots, N,\,  j = 1, \ldots, n\}$ the solution of \eqref{eq.OT}  for~$\xb\in~\!{\mathcal M}_{\Xb}$,  set 
\begin{equation}\label{eq.thetaNx}
	{\scriptstyle\mathfrak{G}}^{(N)}_*(\xb) \coloneqq  {\scriptstyle{\mathfrak{G}}}_{0\ell}^{(N)}(\xb)
\end{equation}
where $\ell \in \{1, \ldots, N\}$ satisfies 
\begin{equation}\label{eq.ConditonalStep1}
	p_{\ell j}^*(\xb) = \max\{p_{ij}^*(\xb); \ i = 1, \ldots, N\}
\end{equation}
with $j \in \{h: w_h\n(\xb, \Xb\n)\neq 0, h = 1, \ldots, n\}$\footnote{This set is nonempty since otherwise $p_{i j}^*(\xb)=0$ for all $i=1, \ldots, N$.}, that is, such that $\Yb_j$ is the sample point closest to $\hat{\thetab}\n(\xb)$ in the~$d(\cdot, \cdot)$  metric (in case \eqref{eq.ConditonalStep1} admits multiple solutions,  choose as~${\scriptstyle{\mathfrak{G}}}_{0\ell}^{(N)}(\xb)$ the one  that is closest to $\hat{\thetab}\n(\xb)$ in terms of the geodesic distance). Use~${\scriptstyle\mathfrak{G}}^{(N)}_*(\xb)$ defined in \eqref{eq.thetaNx}  as an estimator of $\Fb (\thetab(\xb)\vert \xb)$. \smallskip

\noindent{\bf Step 2(i).} Construct a further $N$-points\footnote{{Not necessarily the same $N$ as in Step 1.}}  regular grid $$\mathfrak{G}^{(N)}_{\widehat{\cal M}^{(N)}_0(\xb)} = \{{\scriptstyle{\mathfrak{G}}}^{(N)}_{1i}(\xb) \in {\cal M}, i = 1, \ldots, N\}$$ of the form
$\mathfrak{G}^{(N)}_{\widehat{\cal M}^{(N)}_0(\xb)} \coloneqq \mathring{\mathfrak{G}}^{(N_0)}_{\widehat{\cal M}^{(N)}_0(\xb)} \bigcup {\mathfrak{G}}^{(N_S)}_{\widehat{\cal M}^{(N)}_0(\xb)}(1) \bigcup \cdots \bigcup {\mathfrak{G}}^{(N_S)}_{\widehat{\cal M}^{(N)}_0(\xb)}(N_R)$ 
by factorizing~$N$ into $N_R N_S + N_0$ and proceeding  as in Steps 2(a) and (b) of Section~\ref{subsec.empQuan}. \smallskip 

\noindent{\bf Step 2(ii).} Solve the optimal transport problem \eqref{eq.OT} again, now with the  new grid~$\mathfrak{G}^{(N)}_{\widehat{\cal M}^{(N)}_0(\xb)}$. Denote by ${\bf p}^{**}(\xb) = \{p_{ij}^{**}(\xb);  \ i = 1, \ldots, N, j = 1, \ldots, n\}$ the solution of this problem. Since {for each $i$ there exists at least one $j \in \{1, \ldots, n\}$ such that $p_{ij}^{**}(\xb)>0$,} 
we follow \cite{delBarrio2022} and choose the one which gets the highest mass from the gridpoint ${\scriptstyle{\mathfrak{G}}}^{(N)}_{1i}(\xb) \in \mathfrak{G}^{(N)}_{\widehat{\cal M}^{(N)}_0(\xb)}\vspace{-1.5mm}$; in case of ties,  
pick the one which is closest  to~$\mathring{\mathfrak{G}}^{(N_0)}_{\widehat{\cal M}^{(N)}_0(\xb)}$, that is,
\begin{equation*}
	\Qb^{(N, n)}_{\wb}({\scriptstyle{\mathfrak{G}}}^{(N)}_{1i}(\xb) \vert \xb) := {\arg\min} \left\lbrace d\left(\yb, \mathring{\mathfrak{G}}^{(N_0)}_{\widehat{\cal M}^{(N)}_0(\xb)}\right): \yb \in \{\Yb\n_J: J \in \arg\max_j p_{ij}^{**}(\xb)\}\right\rbrace,
\end{equation*}
where $d\left(\yb, \mathring{\mathfrak{G}}^{(N_0)}_{\widehat{\cal M}^{(N)}_0(\xb)}\right) = \min \left\lbrace d \left(\yb, {\scriptstyle{\mathfrak{G}}}^{(N)}_{1k}(\xb)\right): {\scriptstyle{\mathfrak{G}}}^{(N)}_{1k}(\xb) \in  \mathring{\mathfrak{G}}^{(N_0)}_{\widehat{\cal M}^{(N)}_0(\xb)}\right\rbrace.$
Based on this~empi- rical conditional quantile function $\Qb^{(N, n)}_{\wb}(\cdot \vert \xb)$, we   naturally define the {\it Riemannian  empirical conditional quantile contours} and {\it regions} of order $r/(N_R +1)$  as
\[
\mathcal{C}_{\wb}^{(N, n)}\left(r/(N_R+1) \Big\vert \xb \right) \coloneqq 
\begin{cases}
	\Qb^{(N, n)}_{\wb}\Big(\mathring{\mathfrak{G}}^{(N_0)}_{\widehat{\cal M}^{(N)}_0(\xb)}\Big\vert \xb\Big)  & \text{if}\  \ r=0\\
	\Qb_{\wb}^{(N, n)}\Big({\mathfrak{G}}^{(N_S)}_{\widehat{\cal M}^{(N)}_0(\xb)}(r)\Big\vert \xb\Big) & \text{if}\  \ r\geq 1 \\
\end{cases}
\]
($r = 1, \ldots, N_R$ if $N_0 = 0$,  and~$r = 0, \ldots, N_R$ if $N_0 \neq 0$) and $$\mathbb{C}_{\wb}^{(N, n)}\left(r/(N_R+1) \Big\vert \xb \right) \coloneqq  \Qb_{\wb}^{(N, n)}\Big(\mathring{\mathfrak{G}}^{(N_0)}_{\widehat{\cal M}^{(N)}_0(\xb)} \bigcup {\mathfrak{G}}^{(N_S)}_{\widehat{\cal M}^{(N)}_0(\xb)}(1) \bigcup \cdots \bigcup {\mathfrak{G}}^{(N_S)}_{\widehat{\cal M}^{(N)}_0(\xb)}(r) \Big\vert \xb\Big),$$
respectively.

The following Lemma, as a generalization of Lemma~\ref{lem.Mhat}, states the consistency of~$\widehat{\cal M}^{(N)}_0(\xb)$ constructed in   Step 1; see Appendix~\ref{App.proof} for the proof.

\begin{Lem}\label{lem.Mxhat}
	Let ${\rm P} \in \mathfrak{P}$ have bounded density on ${\cal M}$ and  let Assumptions~\ref{ass.cost}--\ref{Ass.M0x} hold. Then, for any $\xb\in{\mathcal M}_X$, $d_{\rm H}(\widehat{\cal M}^{(N)}_0(\xb), {\cal M}_{\Fb(\thetab(\xb)\vert \xb)}) \rightarrow 0$ a.s. as $n, N \rightarrow \infty$.
\end{Lem}

Proposition~\ref{prop.QwAsymptotic} (See Appendix~\ref{App.proof} for the proof) summarizes the consistency properties of these empirical concepts as $N$ and $n$ tend to infinity.

\begin{proposition}\label{prop.QwAsymptotic}
	Let ${\rm P} \in \mathfrak{P}$ have bounded density on ${\cal M}$. For~$c=d^2/2$, let Assumptions~\ref{ass.cost}--\ref{Ass.M0x} hold. Let $\{(\Xb\n_1, \Yb\n_1), \ldots, (\Xb\n_n, \Yb\n_n)\}$ be an i.i.d. sample from~${\rm P}_{\Xb, \Yb}$ and denote by~$\wb\n$  a consistent weight function. Then, for any $\xb\in{\mathcal M}_{\Xb}$, as $n, N \rightarrow \infty$,
	\begin{enumerate}
		\item[\it (i)] $\max \Big\{d\left(\Qb_{\wb}^{(N, n)}({\scriptstyle{\mathfrak{G}}}^{(N)}_{1i}(\xb) \vert \xb),   \Qb({\scriptstyle{\mathfrak{G}}}^{(N)}_{1i}(\xb) \vert \xb)\right) \Big\vert \,
		{{\scriptstyle{\mathfrak{G}}}^{(N)}_{1i}(\xb) \in \mathfrak{G}^{(N)}_{\widehat{\cal M}^{(N)}_0(\xb)}}\Big\}
		 \rightarrow 0,$
		a.s.;
		\item[\it (ii)] letting \vspace{-2mm}
		\begin{equation*}
			{\cal R}_N \coloneqq 
			\begin{cases}
				\{1, \ldots, N_R\} & \text{if} \ \ N_0 = 0 \\
				\{0, \ldots, N_R\} & \text{if} \ \ N_0 \neq 0
			\end{cases}
			,
		\end{equation*}
		for any $\epsilon>0$ and   any $\xb\in {\mathcal M}_{\Xb}$,  
		$${\rm P}\left\lbrace \max_{r \in {\cal R}_N} d_{\rm H}\left(\mathbb{C}_{\wb}^{(N, n)}(r/(N_R+1) \vert \xb),  \mathbb{C}^{\rm P}_{\mathcal{M}_0(\xb)}(r/(N_R+1) \vert \xb)\right) > \epsilon \right\rbrace \rightarrow 0$$
		and
		$${\rm P}\left\lbrace \max_{r \in {\cal R}_N} d_{\rm H}\left(\mathcal{C}_{\wb}^{(N, n)}(r/(N_R+1) \vert \xb),  \mathcal{C}^{\rm P}_{\mathcal{M}_0(\xb)}(r/(N_R+1) \vert \xb)\right) > \epsilon \right\rbrace \rightarrow 0.$$
	\end{enumerate}
\end{proposition}

\subsection{Examples of consistent weight function}\label{sec.Weight}

Consistency of the empirical conditional quantile functions $\Qb_{\wb}^{(N, n)}(\cdot \vert\, \xb)$ requires consistent weight functions  (see \eqref{eq.cstWeight}). Assuming that the space $({\cal M}_{\Xb}, d_{\Xb})$ of covariates  is a locally compact metric space, \citet[Theorem~3.1]{forzani2012} show that \eqref{eq.cstWeight} holds for $r=2$ when~$\wb\n$ satis\-fies the following three conditions:\smallskip
\begin{compactenum}
\item[\it (i)] for all $a>0$, $\underset{n\rightarrow\infty}{\lim} {\rm E}\left\lbrace\sum_{i=1}^n w_i\n(\Xb, \Xb\n) I[d_{\Xb}(\Xb, \Xb\n_i)>a] \right\rbrace =0$;
\item[\it (ii)] $\underset{n\rightarrow\infty}{\lim} {\rm E}\Big[\underset{1\leq i\leq n}{\max}w_i\n(\Xb, \Xb\n) \Big] =0$;
\item[\it (iii)] there is a constant $c>0$ such that, for every nonnegative measurable function $f$ with ${\rm E}(f(\Xb))<\infty$, 
${\rm E}\Big[ \sum_{i=1}^n w_i\n(\Xb, \Xb\n) f(\Xb\n_i) \Big] \leq c  {\rm E}(f(\Xb)).$
\end{compactenum}

When $({\cal M}_{\Xb}, d_{\Xb})$ is complete and separable, the same authors provide  two classes of mean-square consistent\footnote{That is, satisfies \eqref{eq.cstWeight} for $r=2$.} weight functions:
\begin{enumerate}
\item[\it (a)] $k$-nearest neighbors ($k$-NN) weights, of the form 
\begin{equation}\label{eq.NeighborW}w_j\n(\xb, \Xb\n) \coloneqq k_n^{-1} I\left[\Xb_j\n \in \mathcal{N}_{k_n}(\xb) \right], \ j= 1, \ldots, n,
\end{equation}
where $\mathcal{N}_{k_n}(\xb)$ denotes the set of $k_n$-nearest neighbors of $\xb$ among $\{\Xb_1\n, \ldots, \Xb_n\n\}$ in the~$d_{\Xb}$ metric and the positive integers $k_n$ are such that  $k_n \rightarrow \infty$ and $k_n/n \rightarrow 0$ as $n \rightarrow \infty$;
\item[\it (b)] kernel weights, of the form 
\begin{equation}\label{eq.KernelW}
w_j\n(\xb, \Xb\n) \coloneqq  \Bigg\lbrace \sum_{i=1}^n K\Big(\frac{d_{\Xb}(\xb, \Xb_i\n)}{h_n(\xb)} \Big)\Bigg\rbrace^{-1} K\Big(\frac{d_{\Xb}(\xb, \Xb_j\n)}{h_n(\xb)} \Big), \ j= 1, \ldots, n,
\end{equation}
where the bandwidth $h_n$ satisfies  
$$h_n(\Xb) \rightarrow 0 \ \ \text{a.s. and} \ \ n {\rm P}\big[\bar{B}^{\mathcal{M}_{\Xb}}(\xb, h_n(\Xb))\big]/\log(n) \rightarrow \infty\quad \text{for all $\xb$ as $n\to\infty$}$$ 
with $\bar{B}^{\mathcal{M}_{\Xb}}(\xb, r)$  the closed ball of ${\mathcal{M}_{\Xb}}$ with radius $r$ centered at $\xb$ and   kernel~$K$  such that  $c_1 I\big[u\in [0, 1]\big] \leq K(u) \leq c_2 I\big[u\in [0, 1]\big]$ for some positive constants~$c_1$~and~$c_2$.

\end{enumerate}


\subsection{Numerical illustration: non-Euclidean covariates}

In this section, we   illustrate, via simulations, the Riemannian quantile regression model and the concept of  empirical conditional quantile contours defined in Section~\ref{subsec.EmpQR}   for  ${\cal S}^1$- and~$\mathcal{S}^2$- va\-lued covariates~$\Xb$ and  ${\cal S}^2$- and  $\mathcal{T}^2$- valued responses $\Yb$.  In both cases, we first compare the performance of the $k$-NN and kernel weights (see \eqref{eq.NeighborW} and~\eqref{eq.KernelW}), and then highlight the variation with $\xb$ of the quantile contours.\smallskip

\textbf{Comparison of the $k$-NN and kernel weights.}  We compare the peformance of the $k$-NN and kernel weights for both cap- and equatorial strip-type quantile contours. For the $k$-NN weights, we set $k=N$, and for the kernel weights, we use the trimmed Gaussian kernel $K(u) = e^{-u^2}I\big[u\in [0, 1]\big]$ with bandwidth $h_n(\Xb) = \pi/10$. We set the sample size to~$n=~\!10000$. For the cap-type quantile contours, we choose $N=2001$,  which factorizes into~$N_RN_S+N_0$ with~$N_R=20$, $N_S=100$, and $N_0=1$. Details of the simulation settings are as follows.

\begin{enumerate}
\item[(TS1)] ($\mathcal{S}^2\times \mathcal{T}^2$)-valued $(\Xb,\Yb)$, with $\Xb\sim {\rm vMF}_2({\mub}_\Xb, {\kappa}_\Xb)$, ${\mub}_\Xb = (0.7,0.7,\sqrt{0.02})^\top$ and~${\kappa}_\Xb = 1$; conditional on~$\Xb=(X_1, X_2,X_3)^\top$,     $\Yb\sim{\rm BSvM}(\mub, \kappab (\Xb), \lambda)$  (see~\eqref{eq.BSvM}), with~$\mub = {\bf 0}$, $\kappab = (e^{\vert X_1 \vert + \vert X_2 \vert + \vert X_3 \vert},\,  2)^\top\!$, and $\lambda = 1$.\smallskip
\item[(TS2)] ($\mathcal{S}^2\times \mathcal{T}^2$)-valued $(\Xb,\Yb)$; the covariate $\Xb=(X_1, X_2,X_3)^\top$ is from the same distribution as in~(TS1); conditionally on $\Xb$,  the response $\Yb$ is from the mixture~\eqref{eq.3MixBSvM} of three mutually independent BSvMs, namely,  
 $\Yb_1\sim {\rm BSvM}(\mub_1,\kappab_1,\lambda_1)$, $\Yb_2\sim {\rm BSvM}(\mub_2,\kappab_2,\lambda_2)$, and $\Yb_3\sim {\rm BSvM}(\mub_3,\kappab_3,\lambda_3)$, with $\mub_1 = (-\pi/2, \pi/2)^\top$, \linebreak $\mub_2 = (\pi/2, \pi/2)^\top$,  $\mub_3 = (0, -\pi/5)^\top$, $\kappab_1 = e^{2X_1}(1, 1)^\top$,  $\kappab_2 = e^{2X_2}(1, 1)^\top$, \linebreak $\kappab_3 = e^{\vert X_1 \vert + \vert X_2 \vert + \vert X_3 \vert}(1, 1)^\top$, $\lambda_1 = -2$, $\lambda_2 = 2$, and $\lambda_3 = 0$.\smallskip

\item[(SS1)] ($\mathcal{S}^1\times \mathcal{S}^2$)-valued $(\Xb,\Yb)$; here, $\Xb=(X_1, X_2)^\top \sim {\rm vMF}_1((0, 1)^\top, 3)$ and, conditionally on $\Xb$,  $\Yb \sim {\rm vMF}_2((0, 0, 1)^\top, 5e^{X_1})$.\smallskip

\item[(SS2)] ($\mathcal{S}^1\times \mathcal{S}^2$)-valued $(\Xb,\Yb)$; the covariate  $\Xb$ has the same  ${\rm vMF}_1((0, 1)^\top, 3)$ distribution as in (SS1) and, condi\-tionally~on $\Xb$, $\Yb$ is from the mixture \eqref{eq.3MixvMF} of three mutually independent vMFs, namely,   $\Yb_1\sim {\rm vMF}_2(\mub_1, \kappa_1)$, $\Yb_2\sim {\rm vMF}_2(\mub_2, \kappa_2)$, and $\Yb_2\sim {\rm vMF}_3(\mub_3, \kappa_3)$,  with $\mub_1 = (0.3, 0.4, \sqrt{0.75})^\top$, $\mub_2 = (-0.3, -0.4, \sqrt{0.75})^\top$, $\mub_3 = (-0.3, 0.2, \sqrt{0.87})^\top$, $\kappa_1 = 6e^{X_1}$, $\kappa_2 = 4e^{2X_2}$, and  $\kappa_3 = 2e^{3X_1+2X_2}$.
\end{enumerate}

\begin{figure}[!htbp]
	\centering
	\begin{subfigure}[b]{0.3\textwidth}
		\centering
		\includegraphics[scale=0.25]{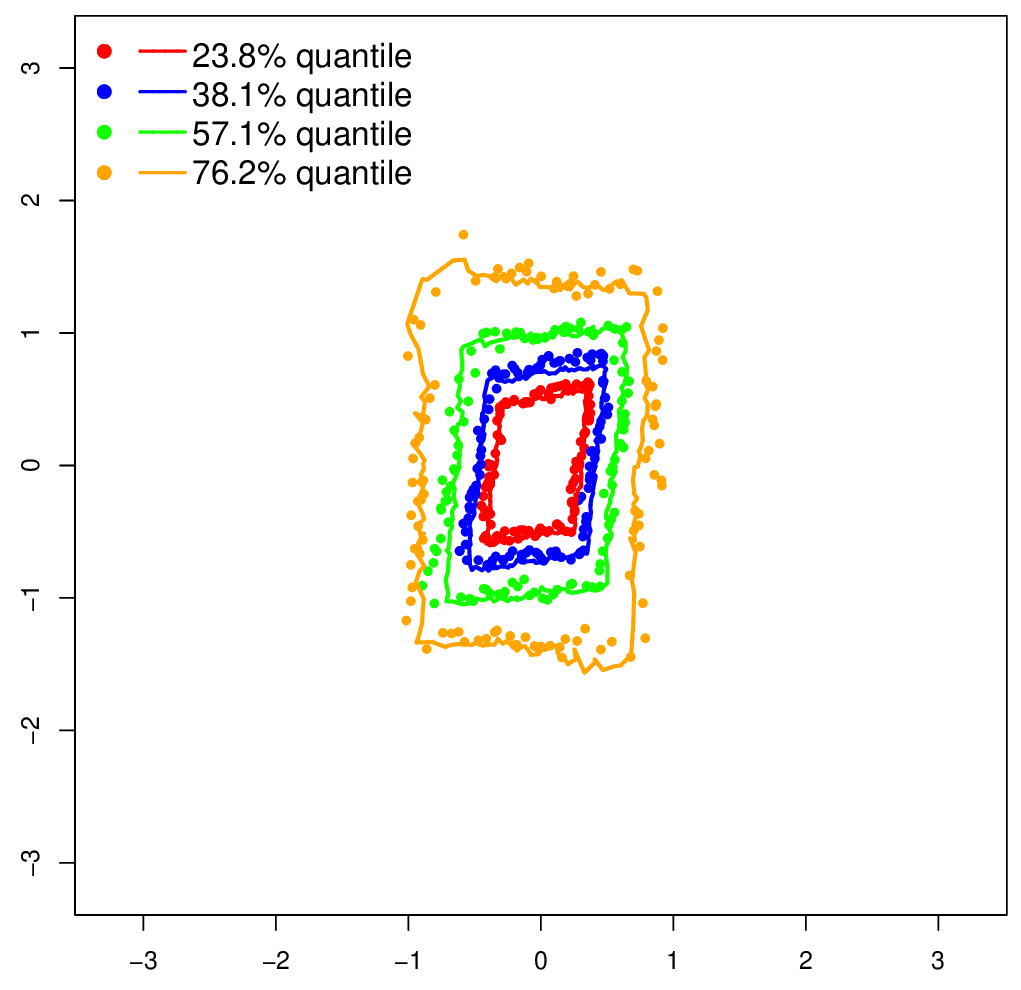}
	\end{subfigure}
	\begin{subfigure}[b]{0.3\textwidth}
		\centering
		\raisebox{8mm}{\includegraphics[scale=0.12]{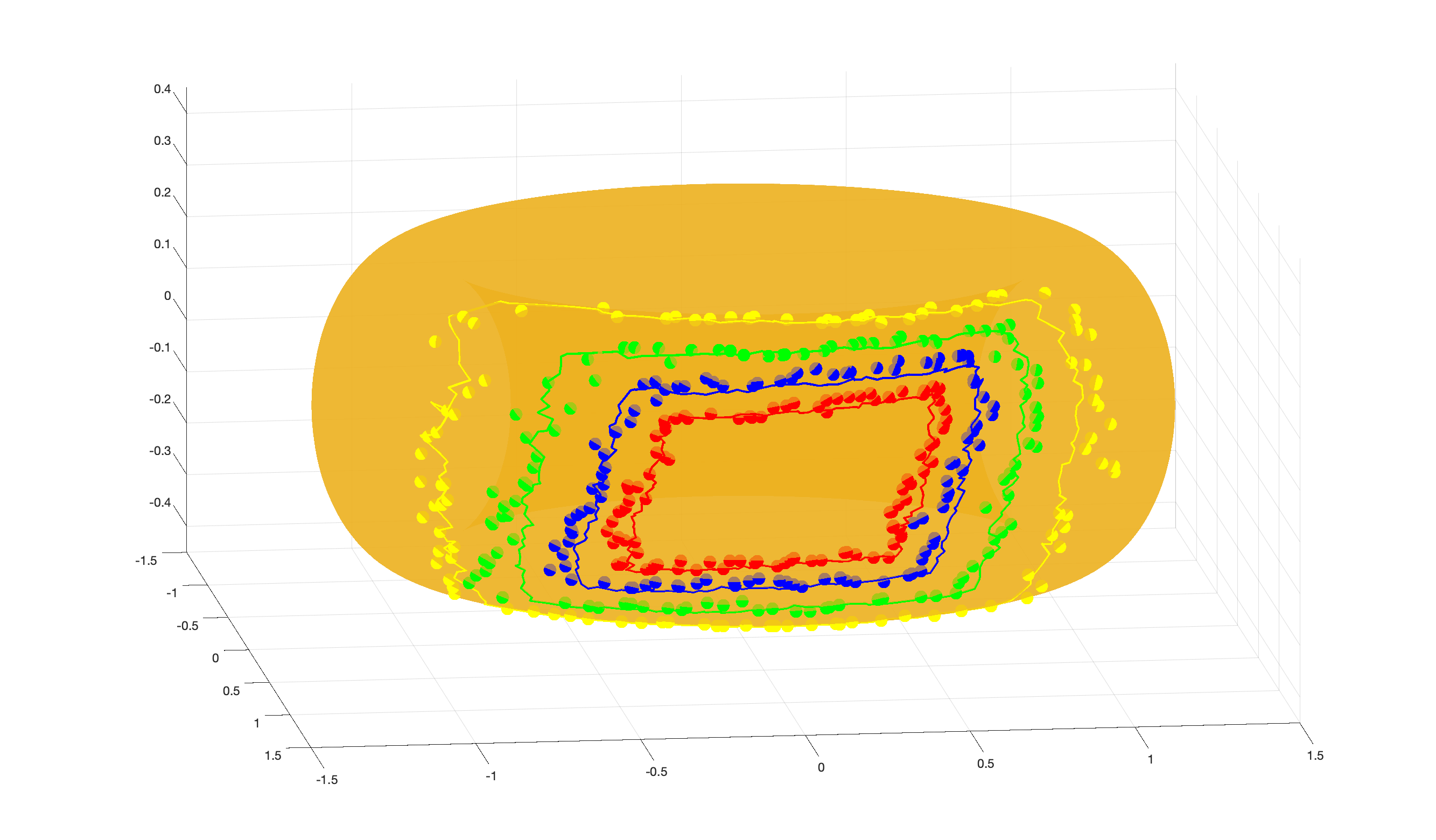}} \vspace{-5mm}
	\end{subfigure}
	\begin{subfigure}[b]{0.3\textwidth}
	\centering
	\hspace*{6mm}\includegraphics[trim=60 70 60 70,clip,scale=0.3]{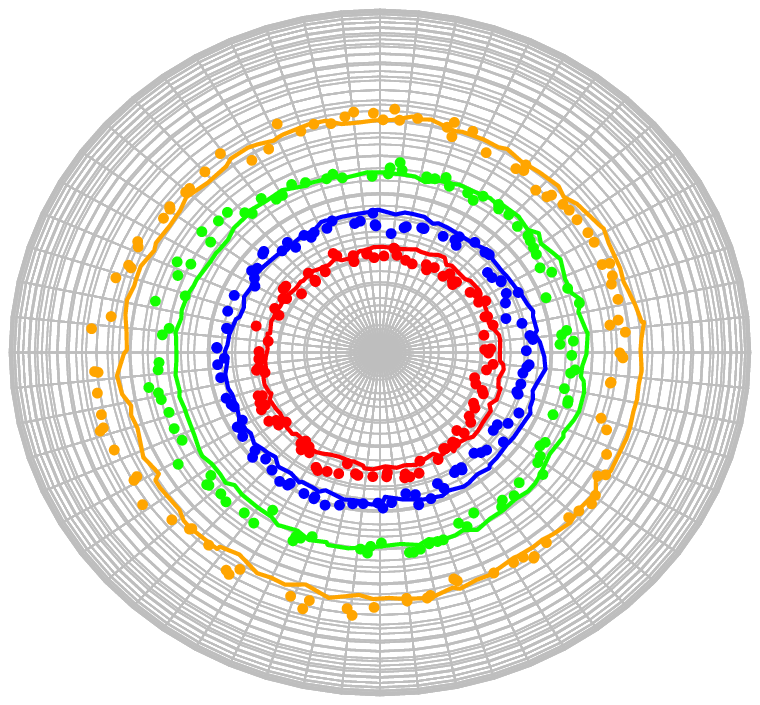}
	\end{subfigure}

	\begin{subfigure}[b]{0.3\textwidth}
		\centering
		\includegraphics[scale=0.25]{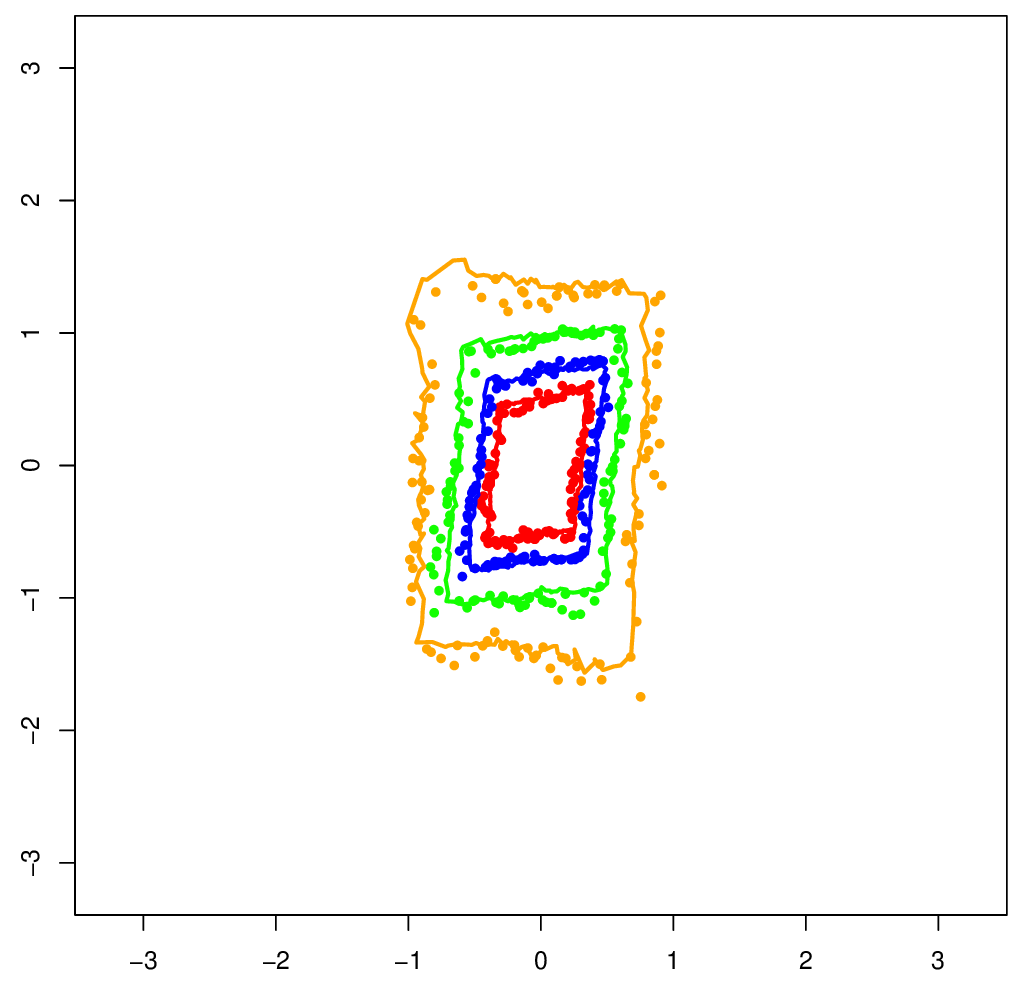}
	\end{subfigure}
	\begin{subfigure}[b]{0.3\textwidth}
		\centering
		\raisebox{8mm}{\includegraphics[scale=0.12]{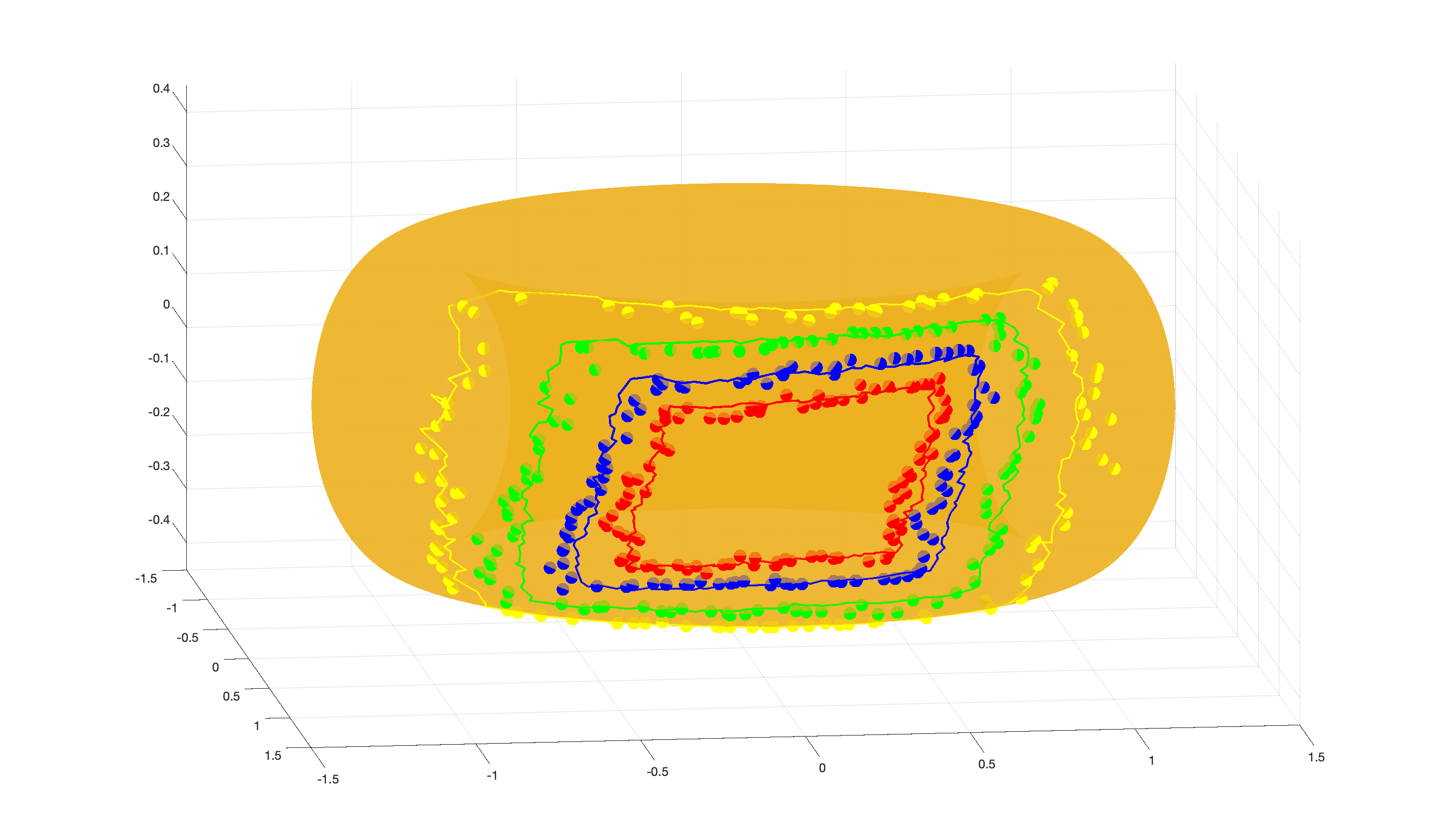}} \vspace{-5mm}
	\end{subfigure}
	\begin{subfigure}[b]{0.3\textwidth}
		\centering
		\hspace*{6mm}\includegraphics[trim=60 70 60 70,clip,scale=0.3]{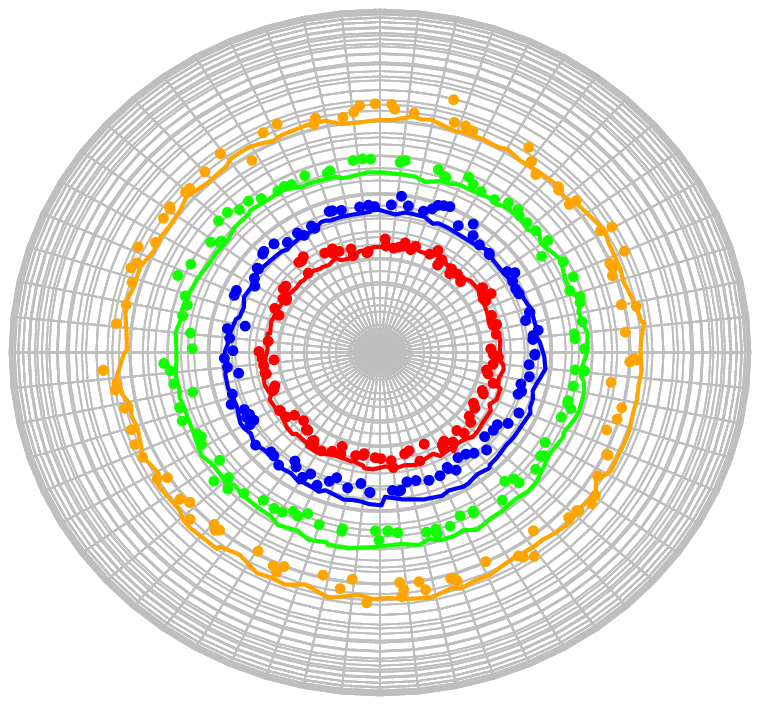}
	\end{subfigure}

	\begin{subfigure}[b]{0.3\textwidth}
		\centering
		\includegraphics[scale=0.25]{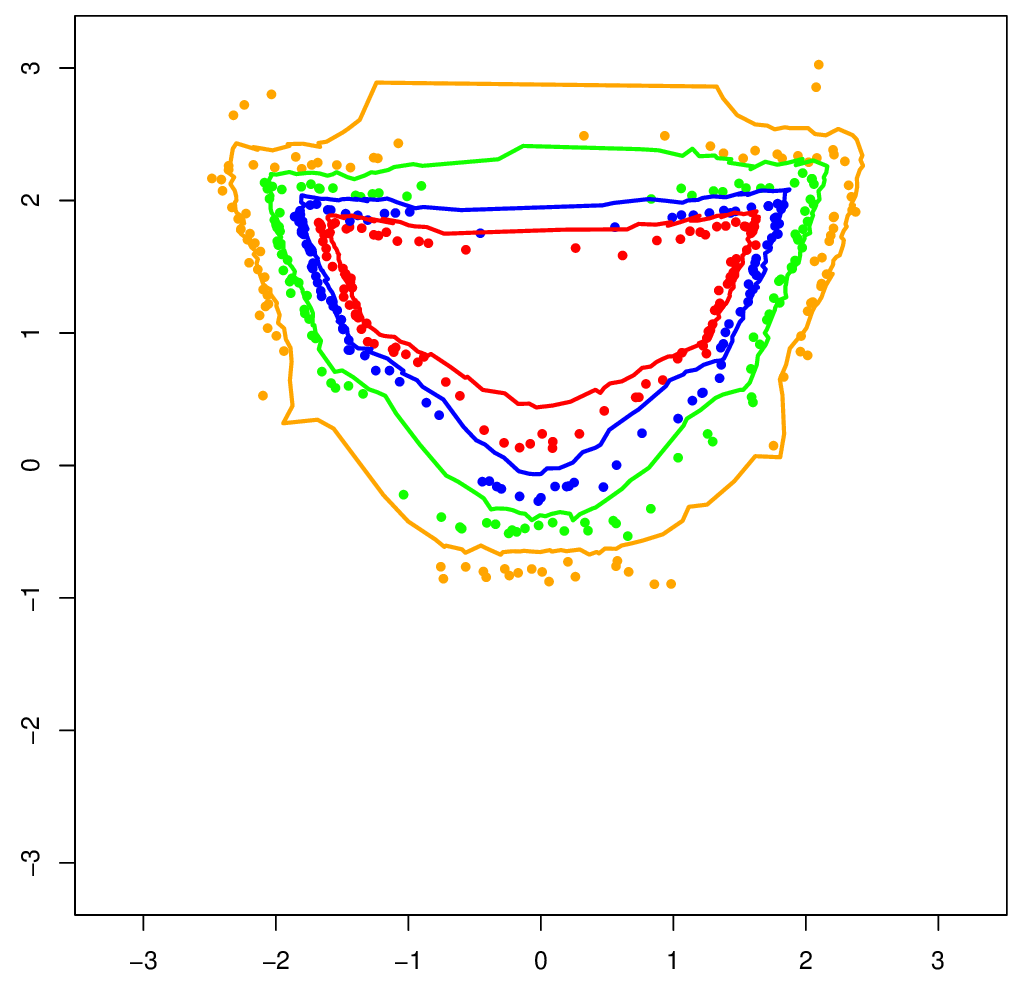} 
	\end{subfigure}
	\begin{subfigure}[b]{0.3\textwidth}
		\centering
		\raisebox{8mm}{\includegraphics[scale=0.12]{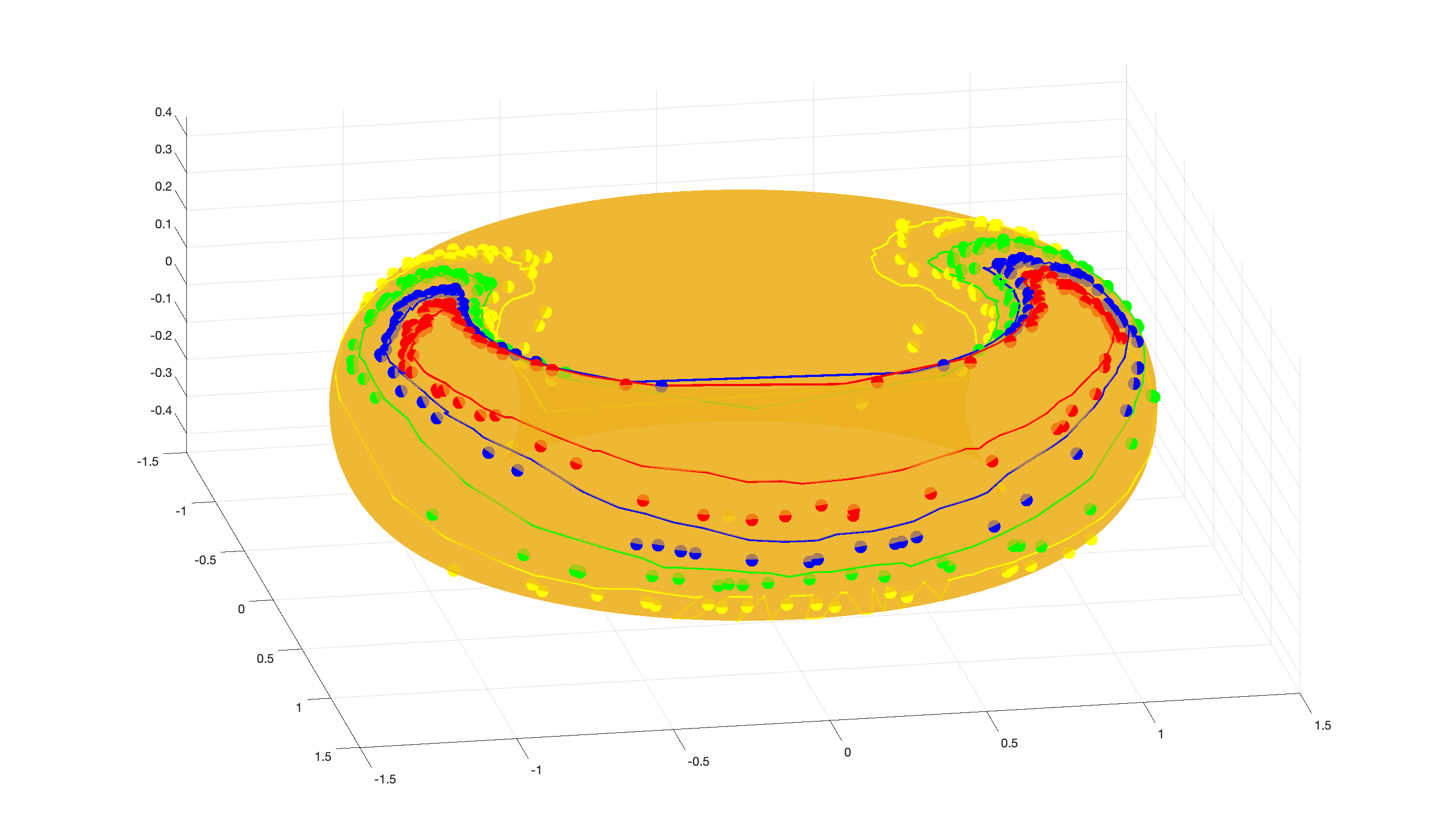}} \vspace{-8mm}
	\end{subfigure}
	\begin{subfigure}[b]{0.3\textwidth}
		\centering
		\includegraphics[trim=70 70 70 110,clip,scale=0.42]{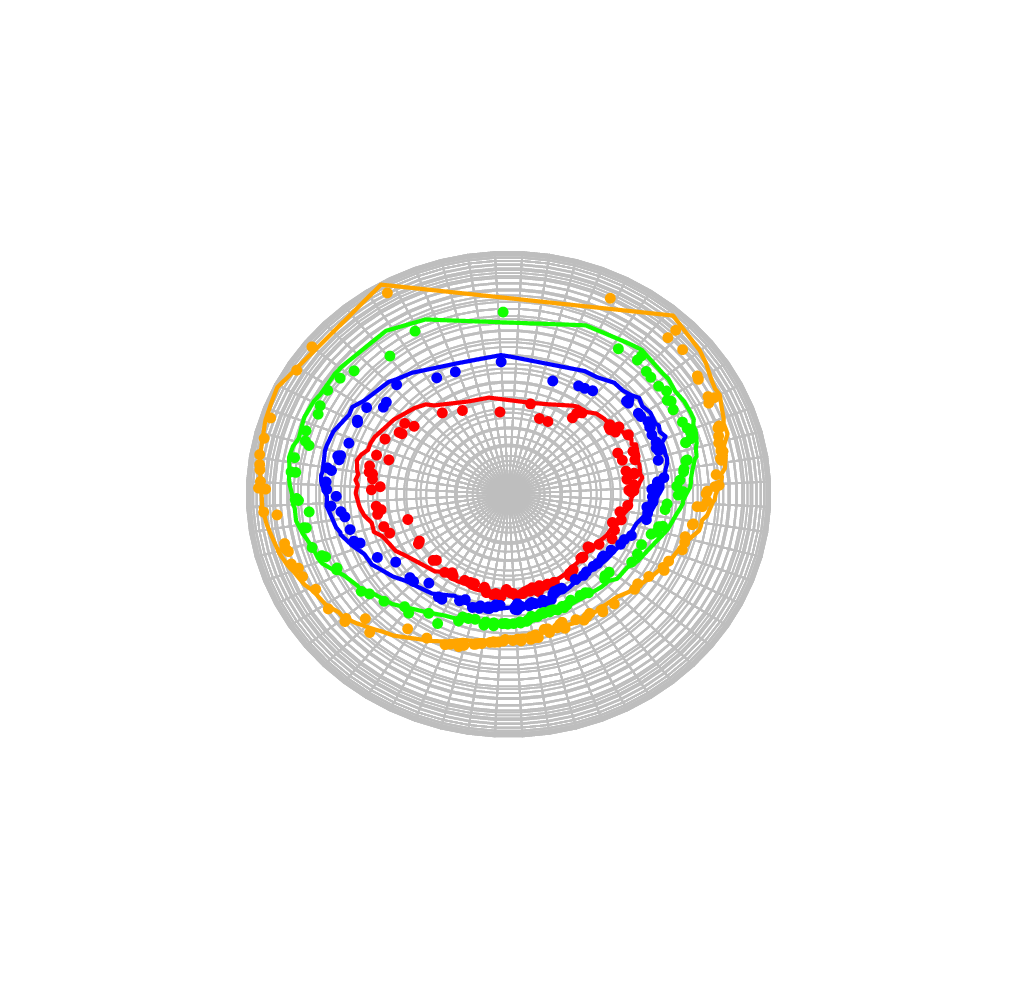}\vspace{-8mm} 
	\end{subfigure}
	
	\begin{subfigure}[b]{0.3\textwidth}
		\centering
		\includegraphics[scale=0.25]{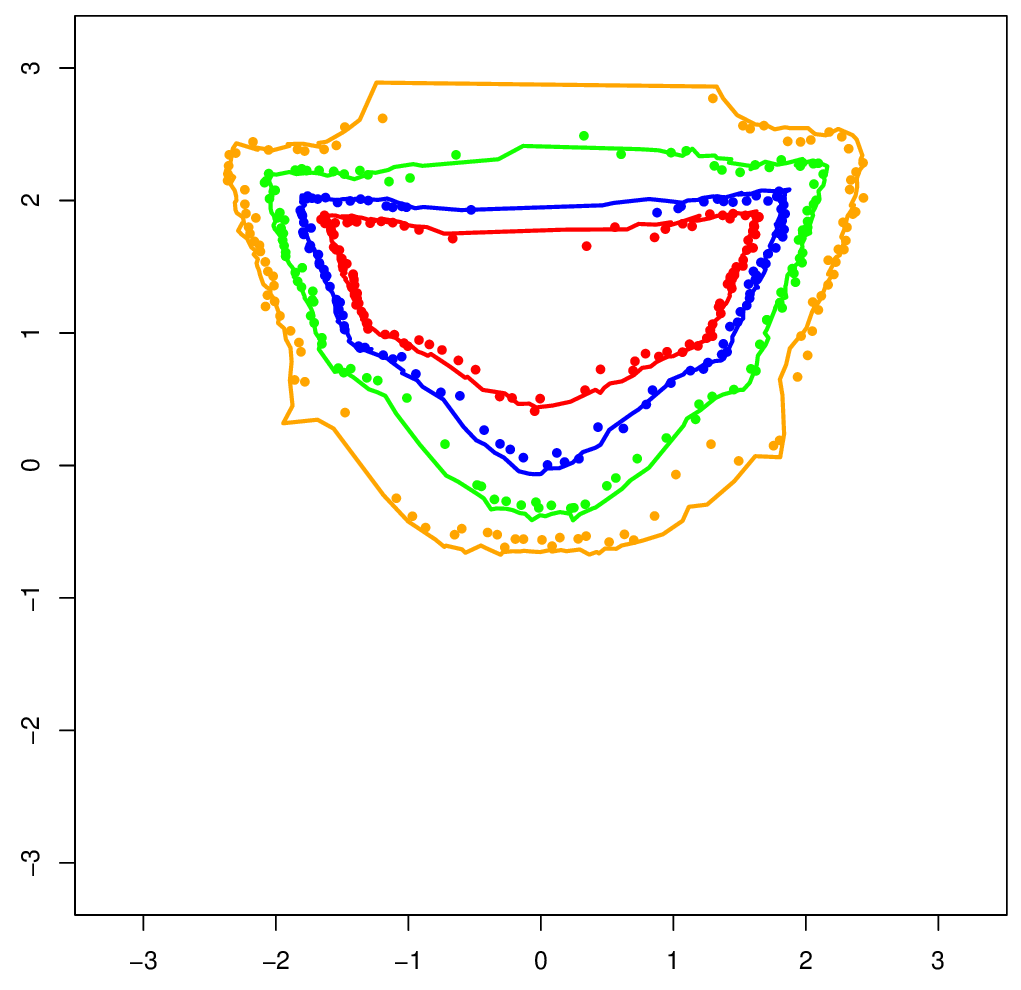} 
	\end{subfigure}
	\begin{subfigure}[b]{0.3\textwidth}
		\centering
		\raisebox{8mm}{\includegraphics[scale=0.12]{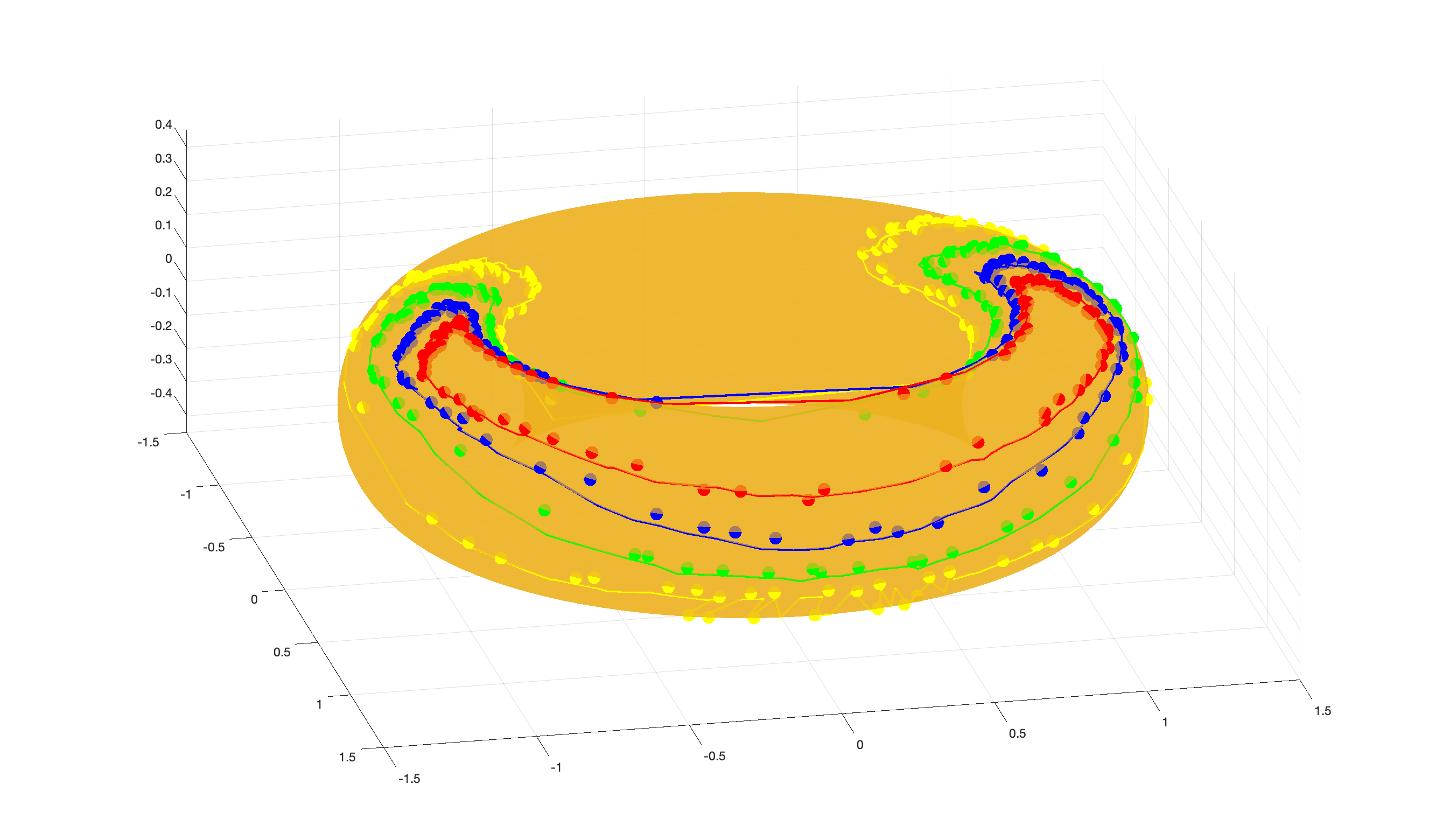}} \vspace{-8mm}
	\end{subfigure}
	\begin{subfigure}[b]{0.3\textwidth}
		\centering
		\includegraphics[trim=70 70 70 110,clip,scale=0.42]{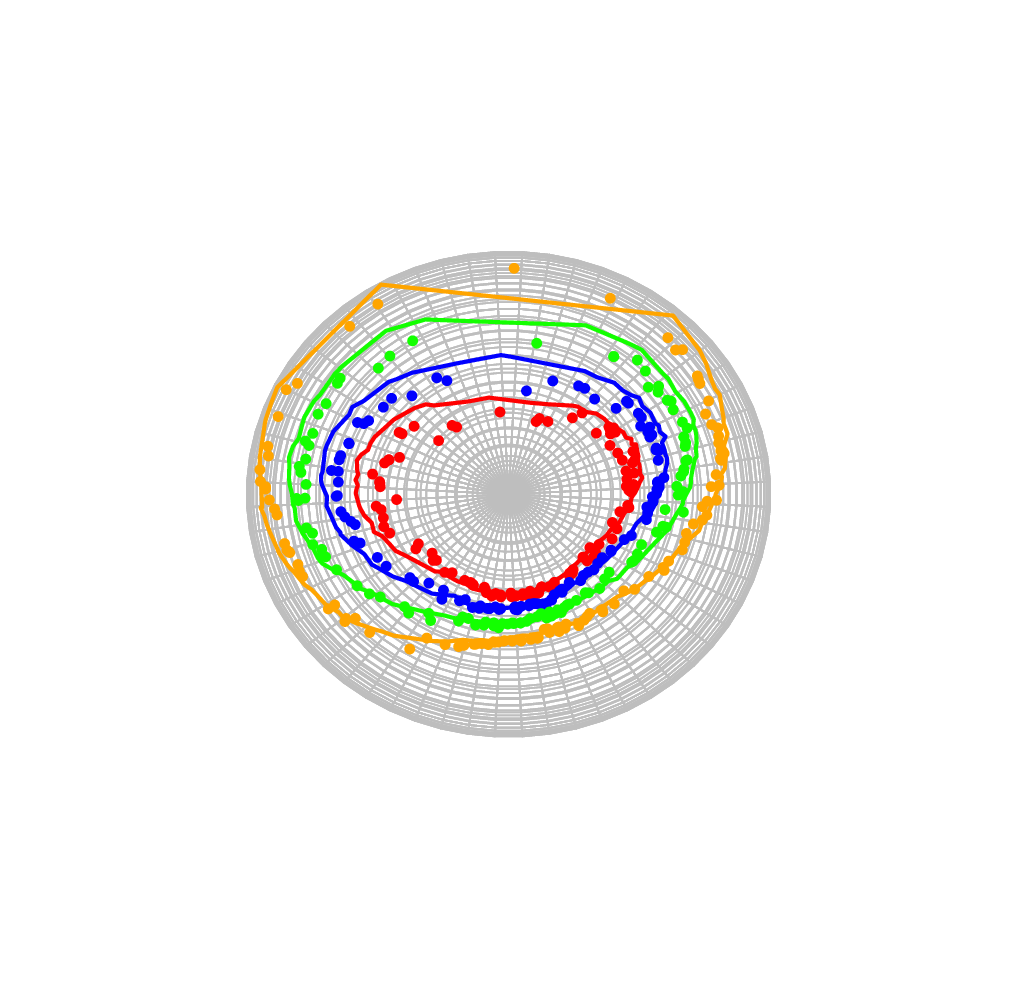}\vspace{-8mm} 
	\end{subfigure}
	
	\caption{\small Plots of empirical (points) and popluation (solid lines) conditional quantile contours when ${\cal M}_0(\xb)$ is a singleton.  1st-2nd rows: (TS1) and (SS1) distributions; 3rd-4th rows: (TS2) and (SS2) distributions. $k$-NN (1st and 3rd rows) and Gaussian kernel (2nd and 4th rows) weights are used;  $n=10000$, $N=2001$, $N_0=1$, $N_R = 20$, $N_S=100$. Quantile\linebreak  levels:~$r/({N_R+1})$, $r=5$ (red),   8  (blue),  12 (green), and  16 (orange);  $\xb=~\!(0.75, 0.65, \sqrt{0.015})^\top$ for $\Yb \in \mathcal{T}^2$ and $\xb=~\!(0.6, 0.8)^\top$ for $\Yb \in \mathcal{S}^2$.}
	\label{Fig:QRContT2S2cap}
\end{figure}  

For the equatorial strip-type quantile contours, we set $N = 2024$, $N_R=20$, $N_S=100$, and $N_0=24$, then  simulate data from 
\begin{enumerate}
\item[(TS3)] ($\mathcal{S}^2\times \mathcal{T}^2$)-valued $(\Xb,\Yb)$ from the same distribution as in (TS1) but with $\kappab_\Xb =~(e^{\vert X_1 \vert + \vert X_2 \vert + \vert X_3 \vert},\,  0)^\top\!$.
\end{enumerate}

\begin{figure}[!htbp]
	\centering
	\begin{subfigure}[b]{0.45\textwidth}
		\centering
		\includegraphics[scale=0.3]{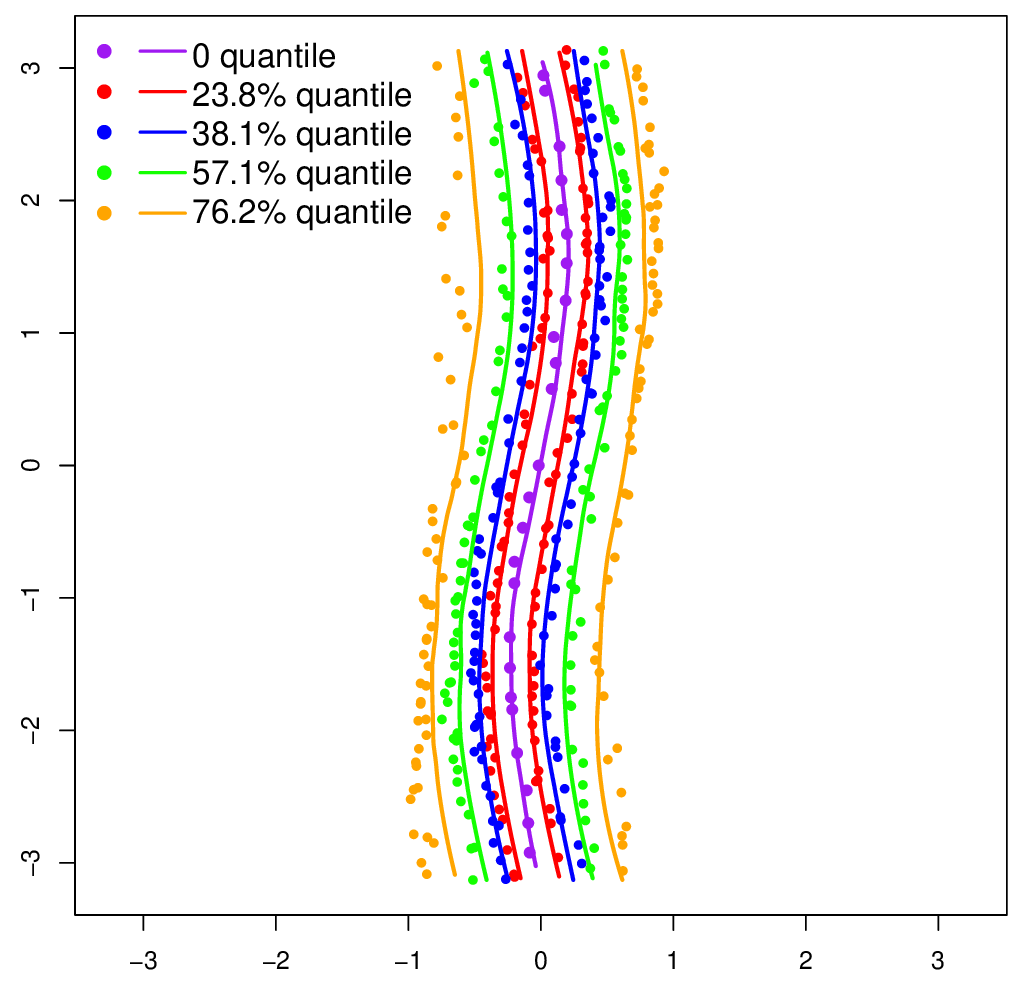}
	\end{subfigure}
	\begin{subfigure}[b]{0.45\textwidth}
		\centering
		\includegraphics[scale=0.3]{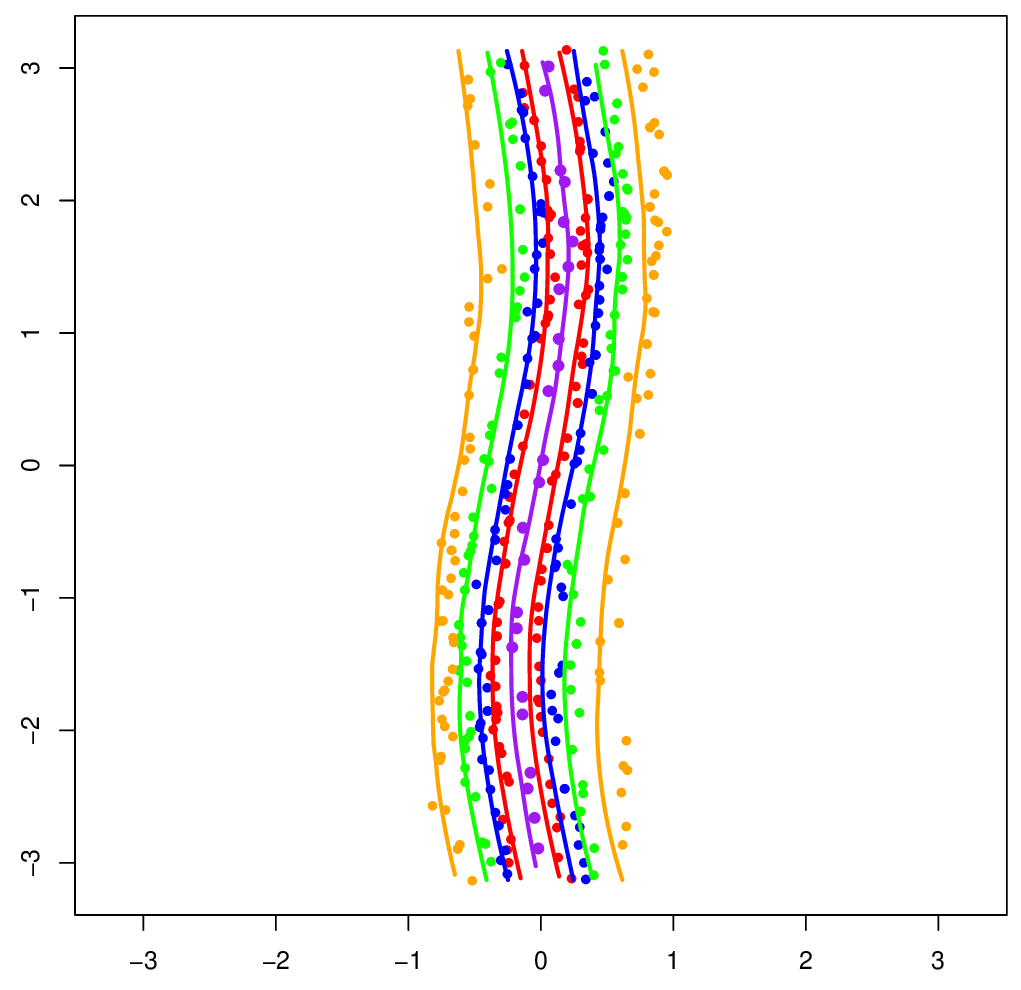}
	\end{subfigure}
	\begin{subfigure}[b]{0.45\textwidth}
		\centering
		\includegraphics[scale=0.15]{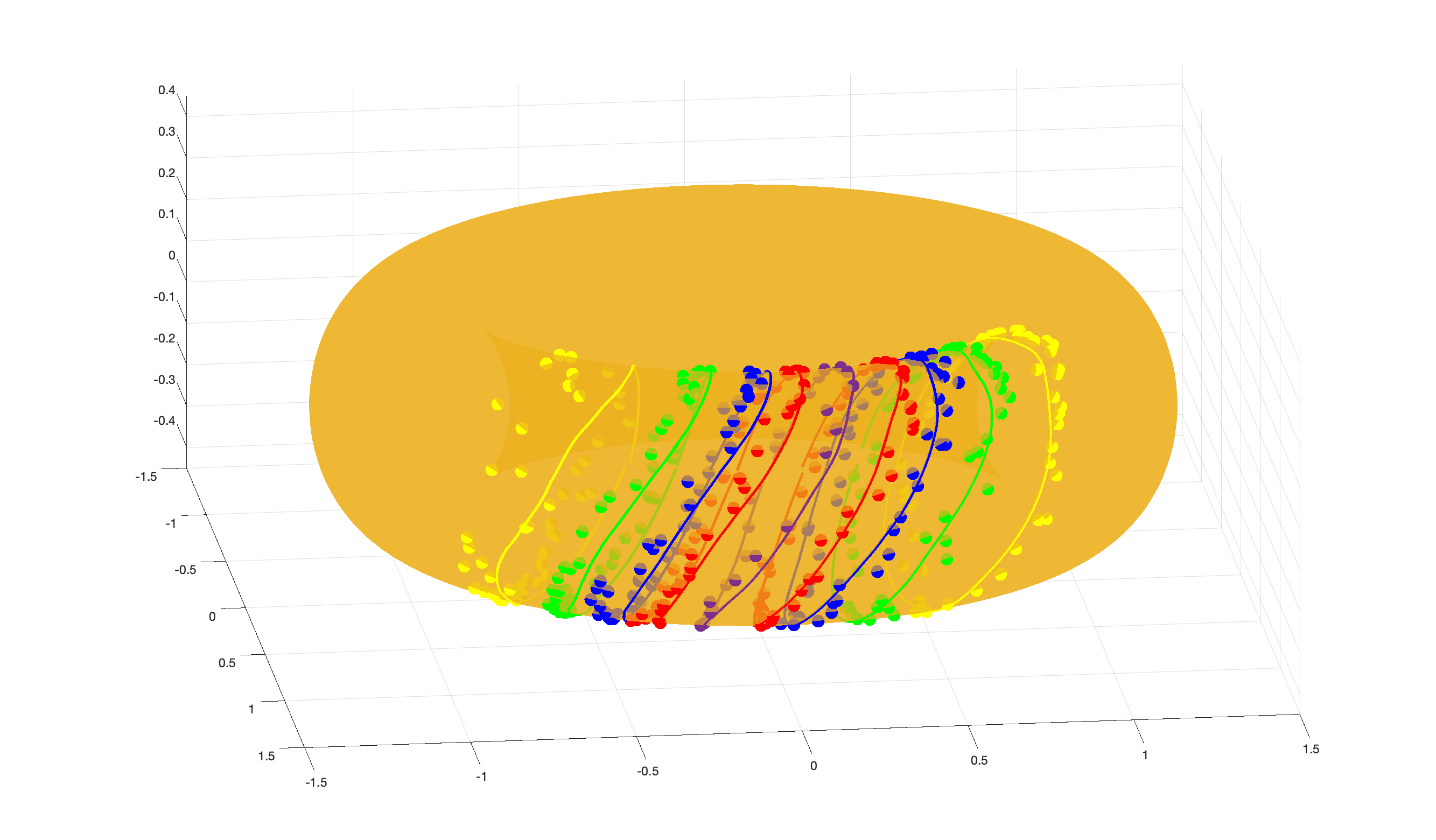} \vspace{-5mm}
	\end{subfigure}
	\begin{subfigure}[b]{0.45\textwidth}
		\centering
		\includegraphics[scale=0.15]{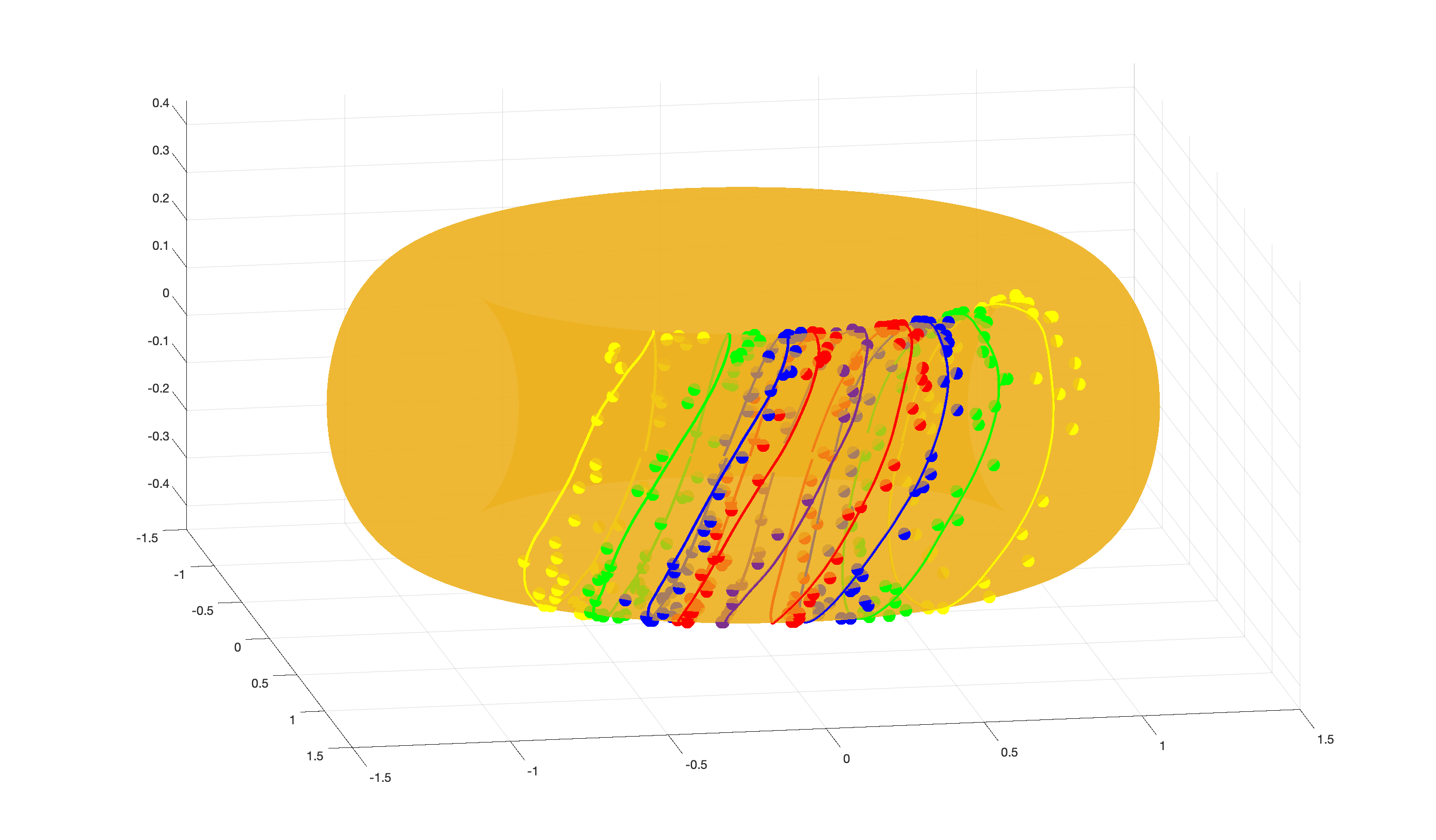} \vspace{-5mm}
	\end{subfigure}
	\caption{\small Plots of strip-type empirical (points) and population (solid lines) conditional quantile contours for (TS3)   when ${\cal M}_0(\xb)$ is of dimension $(p-1)$.  Left column: $k$-NN weights; right column:  Gaussian kernel weights. $n=10000$, $N=2024$, $N_0=24$, $N_R = 20$, $N_S=100$ and~$\xb=~\!(0.75, 0.65, \sqrt{0.015})^\top$. Quantile levels $r/({N_R+1})$, $r=0$ (purple),  5 (red),   8  (blue),~12 (green), and  16 (orange).}
	\label{Fig:QRContT2Strips}
\end{figure}

Figures~\ref{Fig:QRContT2S2cap} (for cap-type contours) and~\ref{Fig:QRContT2Strips} (for equatorial strip-type), provide  a comparison of the performances  of the $k$-NN and kernel weights.  Figure~\ref{Fig:QRContT2S2cap} shows the empirical  conditional quantile contours~$\mathcal{C}^{(N, n)}_{\wb}\big(r/{(N_R+1)}\vert \xb\big)$ for~$r = 5,\, 8,\, 12$, and~$16$ at~$\xb =~\!(0.75, 0.65, \sqrt{0.015})^\top$ for $\Yb \in \mathcal{T}^2$ (1st and 2nd columns) and  for $\Yb \in \mathcal{S}^2$ at~$\xb = (0.6, 0.8)^\top$ (3rd column). The 1st and 3rd rows are for the $k$-NN weights,  the 2nd and~4th rows   for the kernel weights. Along with these empirical contours, we also provide (solid lines) the corresponding ``population  contours'' $\mathcal{C}^{\rm P}_{\mathcal{M}_0(\xb)}\big({r}/({N_R+1})\vert \xb\big)$---actually,  numerical approximations thereof, computed from  ``very large'' samples (size~10000 for each fixed~$\xb$)  generated from the underlying  distributions of $\Yb$ conditional on $\Xb = \xb$. Figure~\ref{Fig:QRContT2Strips} provides plots for  (TS3) based on the $k$-NN and kernel weights, shown in the left and right columns, respectively.

In all plots, the empirical contours  $\mathcal{C}^{(N, n)}_{\wb}(\cdot\vert \xb)$  are very  close to their population counterparts $\mathcal{C}^{\rm P}_{\mathcal{M}_0(\xb)}\left(\cdot \vert \xb\right)$, which is in line with the asymptotic result in Proposition~\ref{prop.QwAsymptotic}(ii).   In particular, they   adapt quite well to the main  features of the conditional distributions---rotational symmetry for (SS1),   multimodality for (SS2) and (TS2), positive dependence of $\Yb$'s two components  for (TS3). In general, the  $k$-NN and kernel weights yield very similar results, but the kernel-weights-based empirical conditional contours fit slightly better than the $k$-NN-based ones for (TS2). In the sequel, we are adopting kernel weights to illustrate the dependence of conditional quantile contours   with $\xb$.\smallskip

\begin{figure}[!htbp]
     \centering
     \begin{subfigure}[b]{0.45\textwidth}
         \centering 
         \includegraphics[scale=0.27]{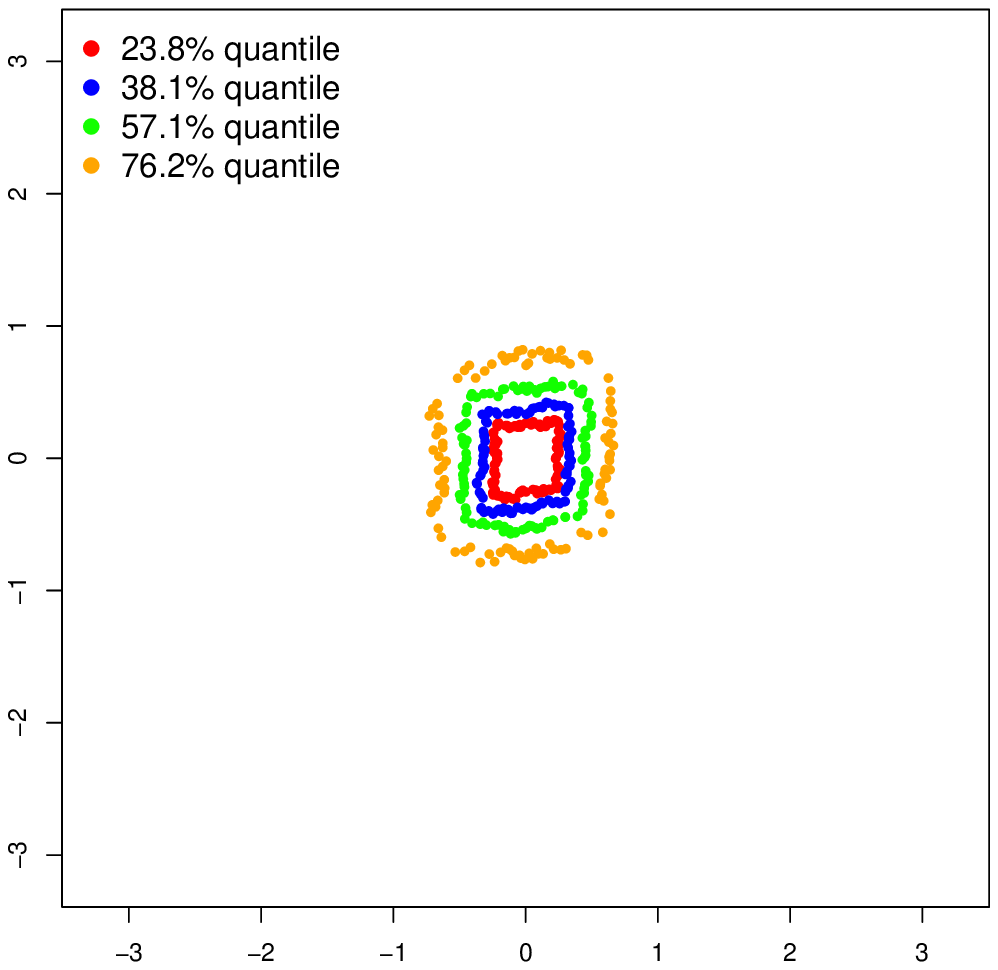}
     \end{subfigure}
     \begin{subfigure}[b]{0.45\textwidth}
         \centering
         \includegraphics[scale=0.27]{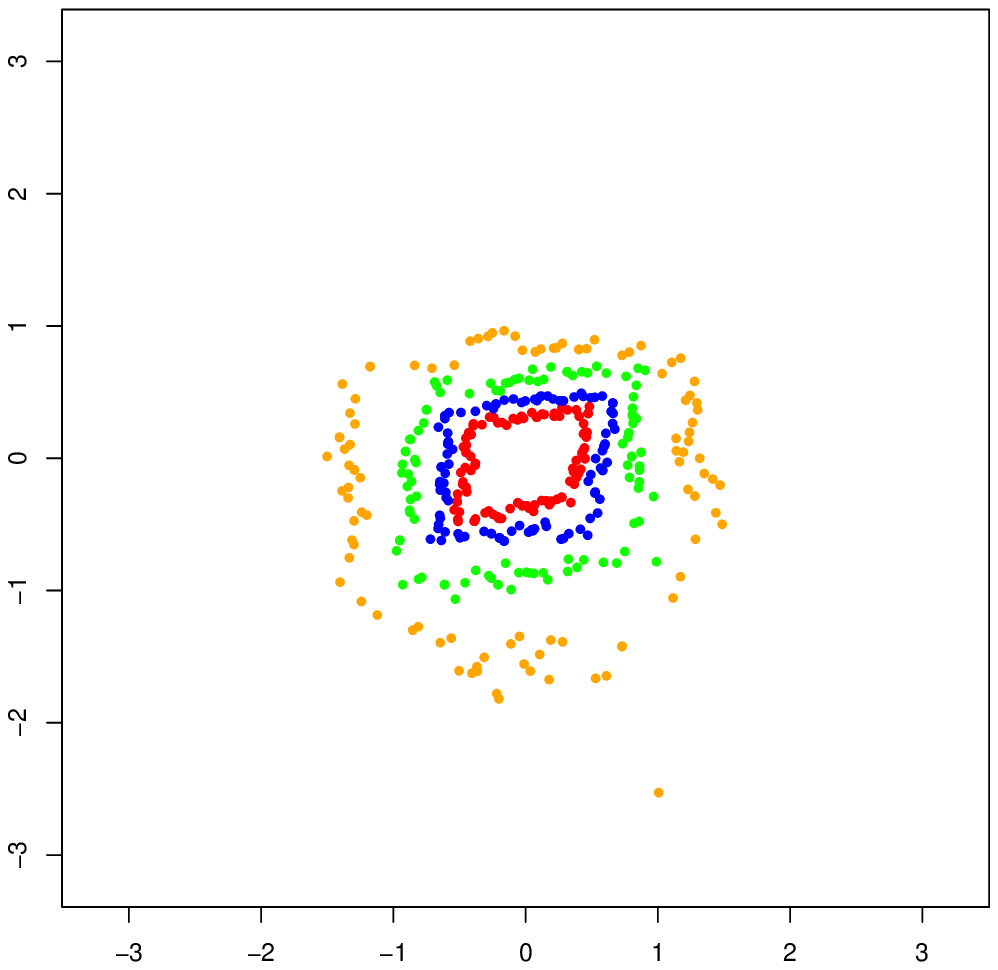}
     \end{subfigure}
     \begin{subfigure}[b]{0.45\textwidth}
         \centering
         \includegraphics[scale=0.27]{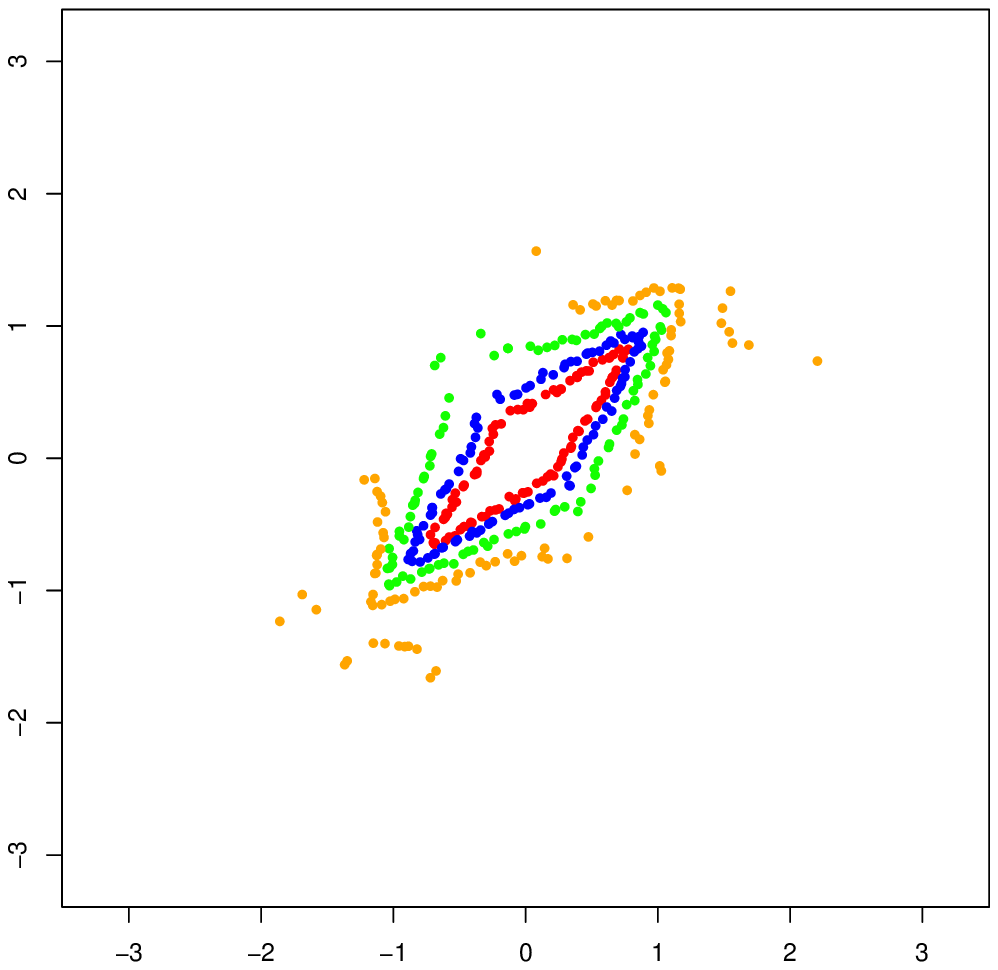} 
     \end{subfigure}
     \begin{subfigure}[b]{0.45\textwidth}
         \centering
         \includegraphics[scale=0.27]{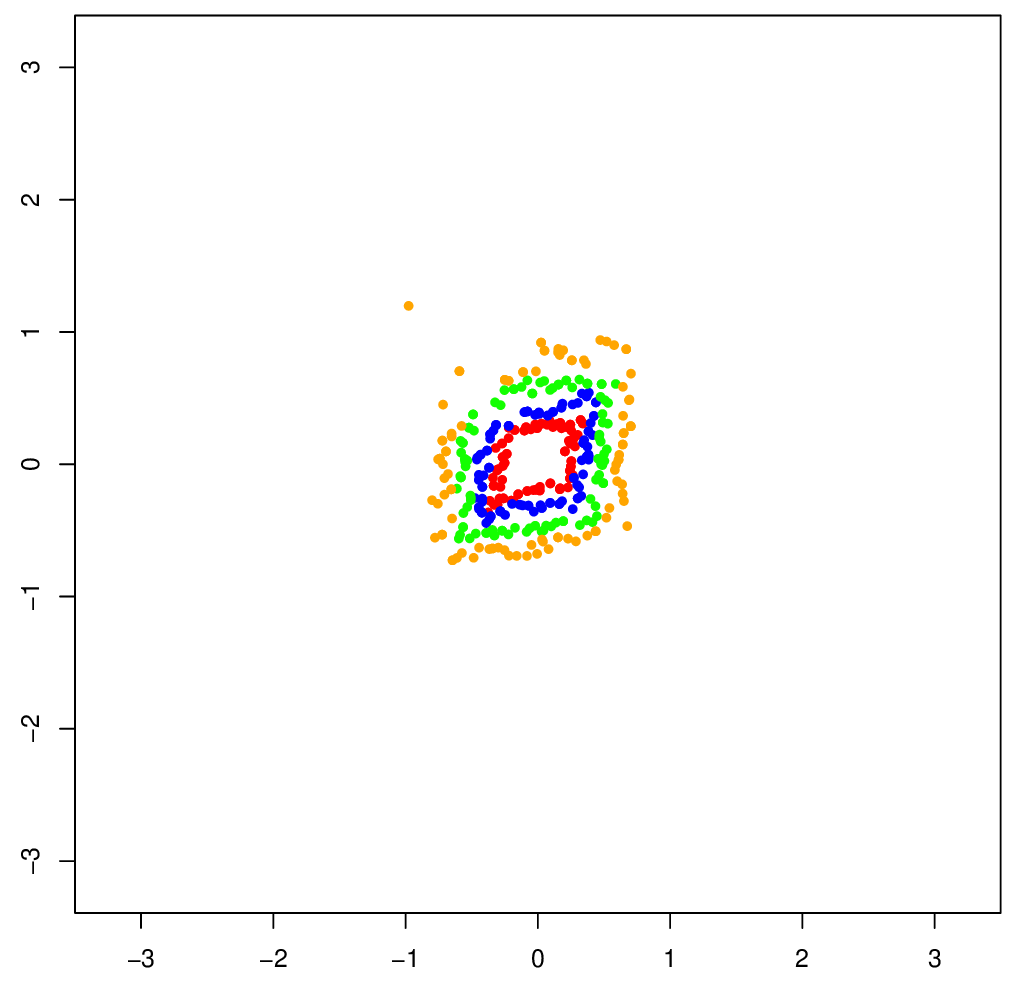} 
     \end{subfigure}
     \begin{subfigure}[b]{0.45\textwidth}
         \centering
         \includegraphics[scale=0.13]{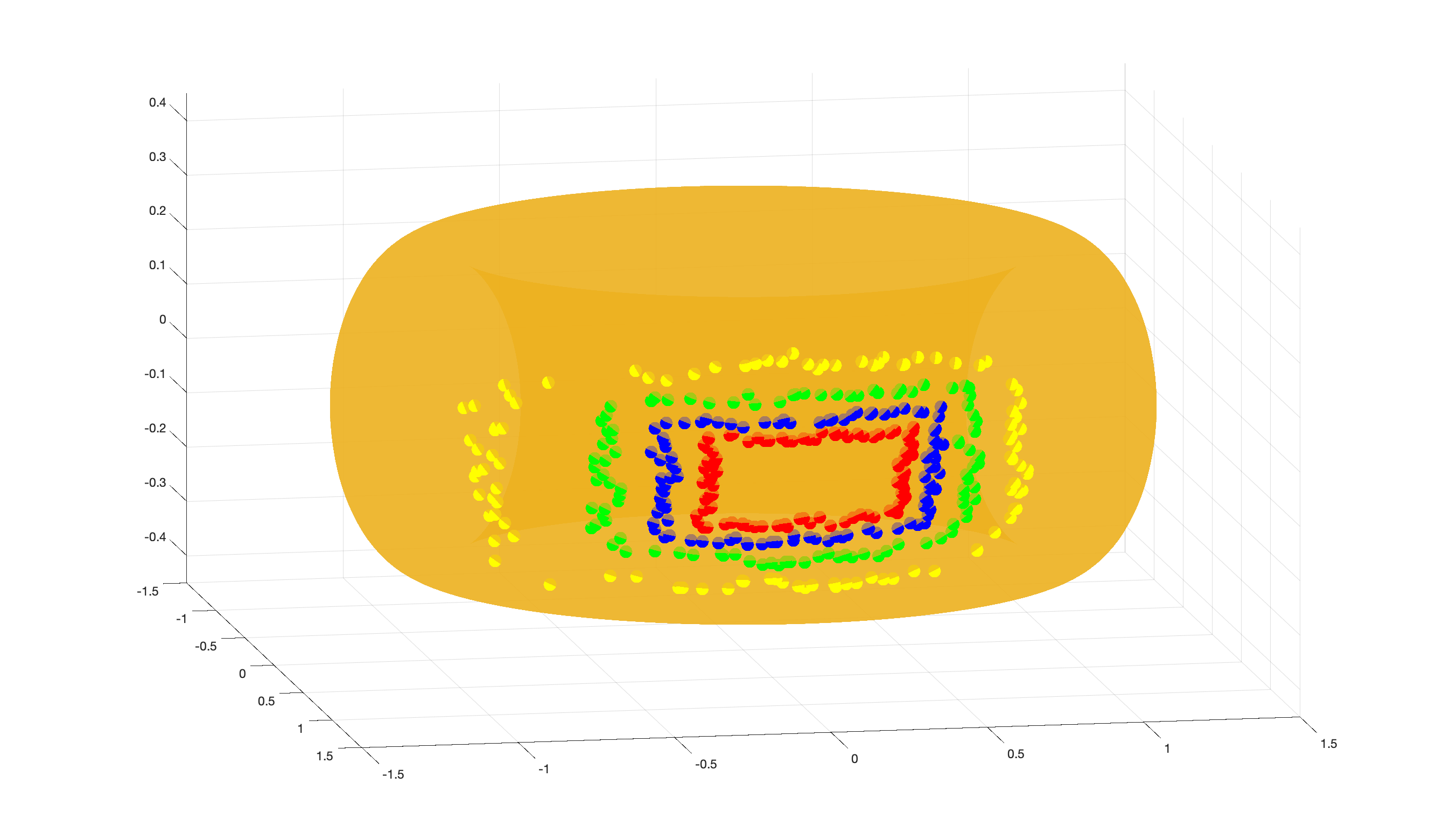} \vspace{-2mm}
     \end{subfigure}
     \begin{subfigure}[b]{0.45\textwidth}
         \centering
         \includegraphics[scale=0.13]{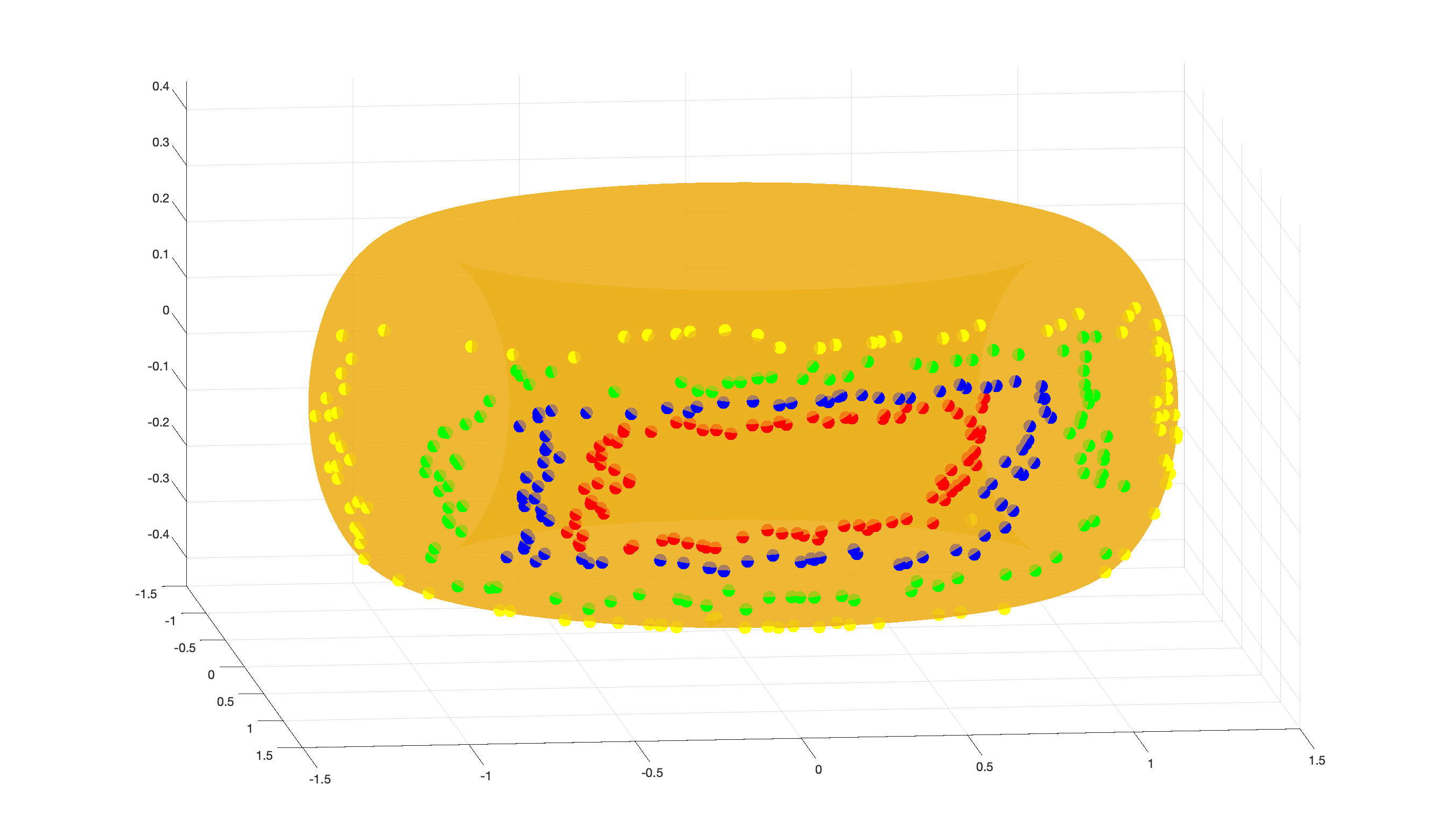} \vspace{-2mm}
     \end{subfigure}
     \begin{subfigure}[b]{0.45\textwidth}
         \centering
         \includegraphics[scale=0.13]{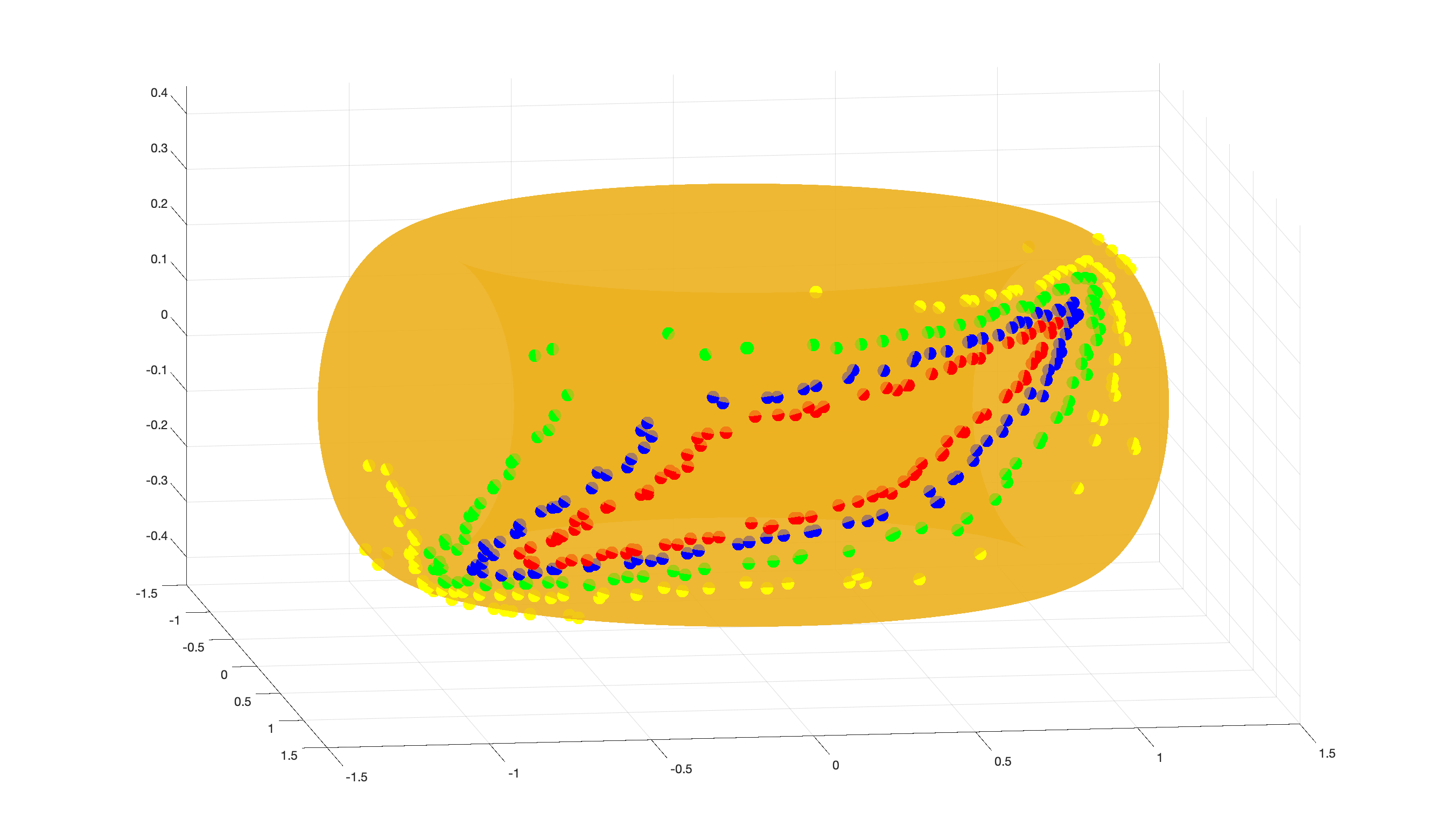} \vspace{-2mm}
     \end{subfigure}
     \begin{subfigure}[b]{0.45\textwidth}
         \centering
         \includegraphics[scale=0.13]{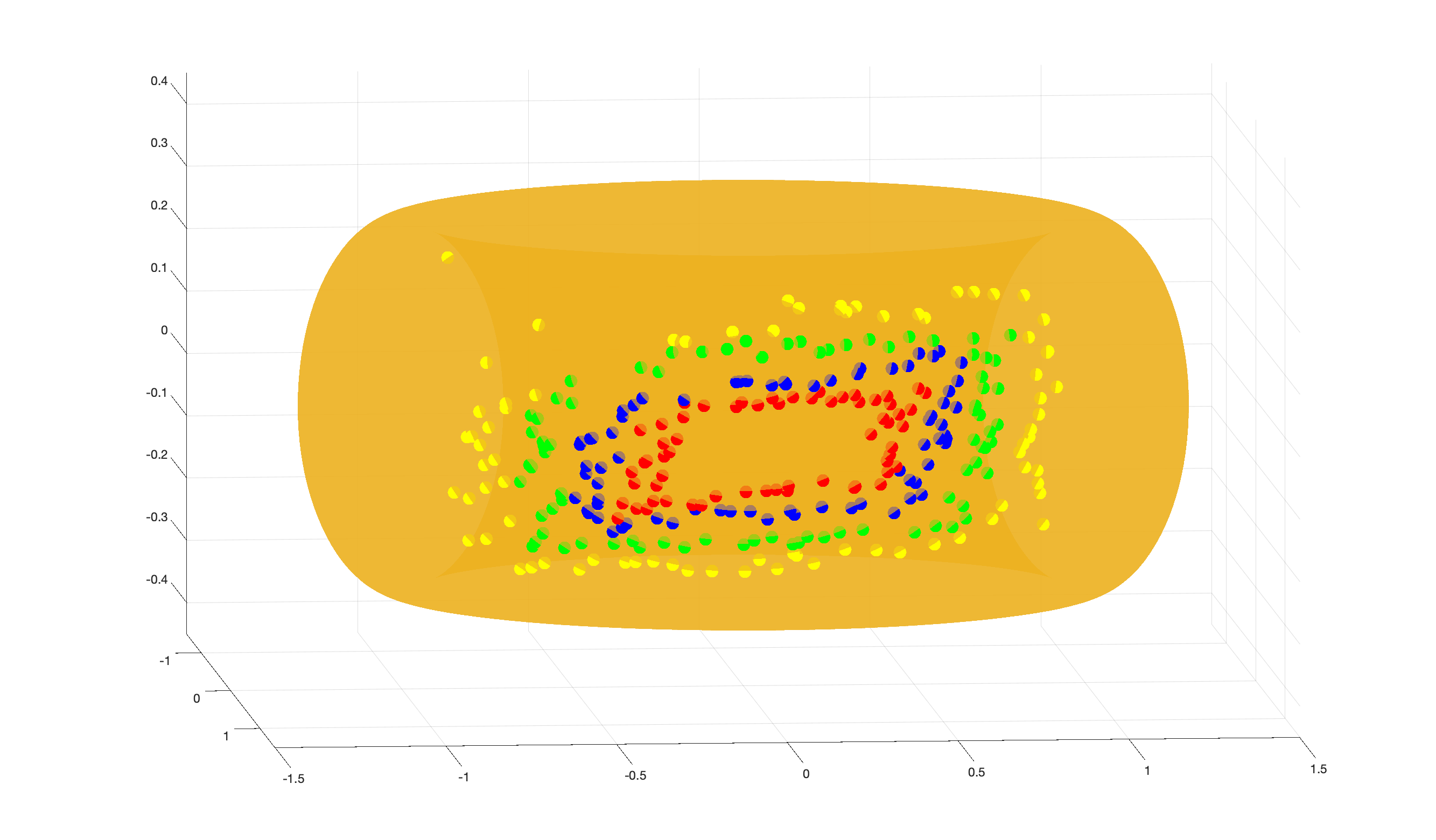} \vspace{-2mm}
     \end{subfigure}
     \begin{subfigure}[b]{0.45\textwidth}
         \centering
         \includegraphics[trim=70 40 70 90,clip, width=0.8\textwidth]{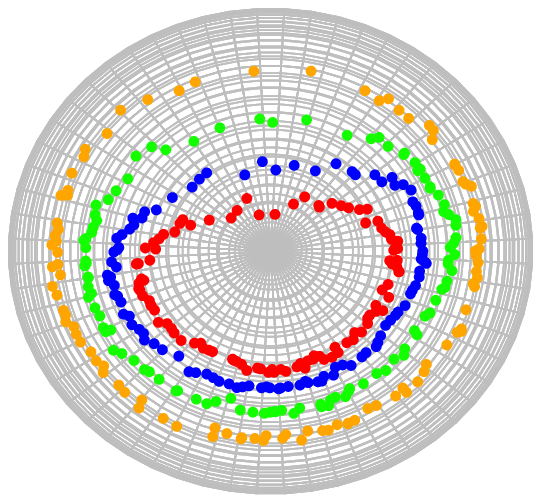} \vspace{-10mm}
     \end{subfigure}
     \begin{subfigure}[b]{0.45\textwidth}
         \centering
         \includegraphics[trim=70 40 70 90,clip, width=0.8\textwidth]{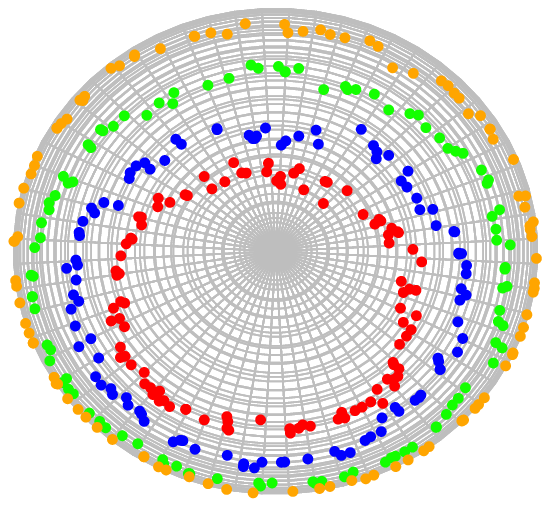} \vspace{-10mm}
     \end{subfigure}
        \caption{\small      Cap-type conditional quantile contours $\mathcal{C}^{(N, n)}_{\wb}\big(\frac{r}{N_R+1}\vert \xb\big)$ ($r=5$ (red),   8  (blue),  12 (green), and~16 (orange)) for 
         (TS1*) (top-left for $\xb = (0.75, 0.65, \sqrt{0.015})^\top$, top-right for~$\xb = (0.2, 0.4, \sqrt{0.8})^\top$), (TS1**) (2nd row: left for $\xb = (0.4, 0.4, \sqrt{0.68})^\top$, right for $\xb = (0, 0, 1)^\top$) and (SS2*) (bottom-left for $\xb = (0, 1)^\top$, bottom-right for $\xb = (0.6, 0.8)^\top$) 
        ($n=10000$,\linebreak $N=2001$, $N_0=1$, $N_R = 20$, $N_S=100$). {Third and fourth rows:  (TS1*) and (TS1**) plots under  $\mathbb{R}^3$-embedded representations.}
        }
 \label{Fig:QRContourT2diffx}
\end{figure}

\textbf{Dependence on $\xb$ of the quantile contours.}
Turning to the main objective of quantile regression, we now analyze how our method is capturing the impact of $\Xb$ on the conditional distribution (the conditional quantile contours)  of $\Yb$.  For cap-type quantile contours, we therefore consider the following variants of (TS1) and (SS2). 

\begin{enumerate}
\item[(TS1*)] $(\Xb,\Yb)$  from the same distribution as in (TS1) but with  $\kappa_1 = e^{3\vert X_1 \vert}$, $\kappa_2 =~\!e^{3\vert X_2 \vert}$, and ${\kappa}_\Xb = 2$;\smallskip

\item[(TS1**)] $(\Xb,\Yb)$  from the same distribution as in (TS1) but with $\lambda = 10(\vert X_1 \vert + \vert X_2 \vert)$, $\kappa_1 = 6 = \kappa_2$,  ${\kappa}_\Xb = 2$ and ${\mub}_\Xb = (0.3,0.3,\sqrt{0.82})^\top$.

\item[(SS2*)] $(\Xb,\Yb)$  from the same distribution as in (SS2) but with  concentration parameters  $\kappa_1, \kappa_2, \kappa_3 = e^{\vert X_2\vert/(\vert X_1\vert +0.25)}$.
\end{enumerate}

\begin{figure}[b!]
	\centering
	\begin{subfigure}[b]{0.45\textwidth}
		\centering
		\includegraphics[scale=0.3]{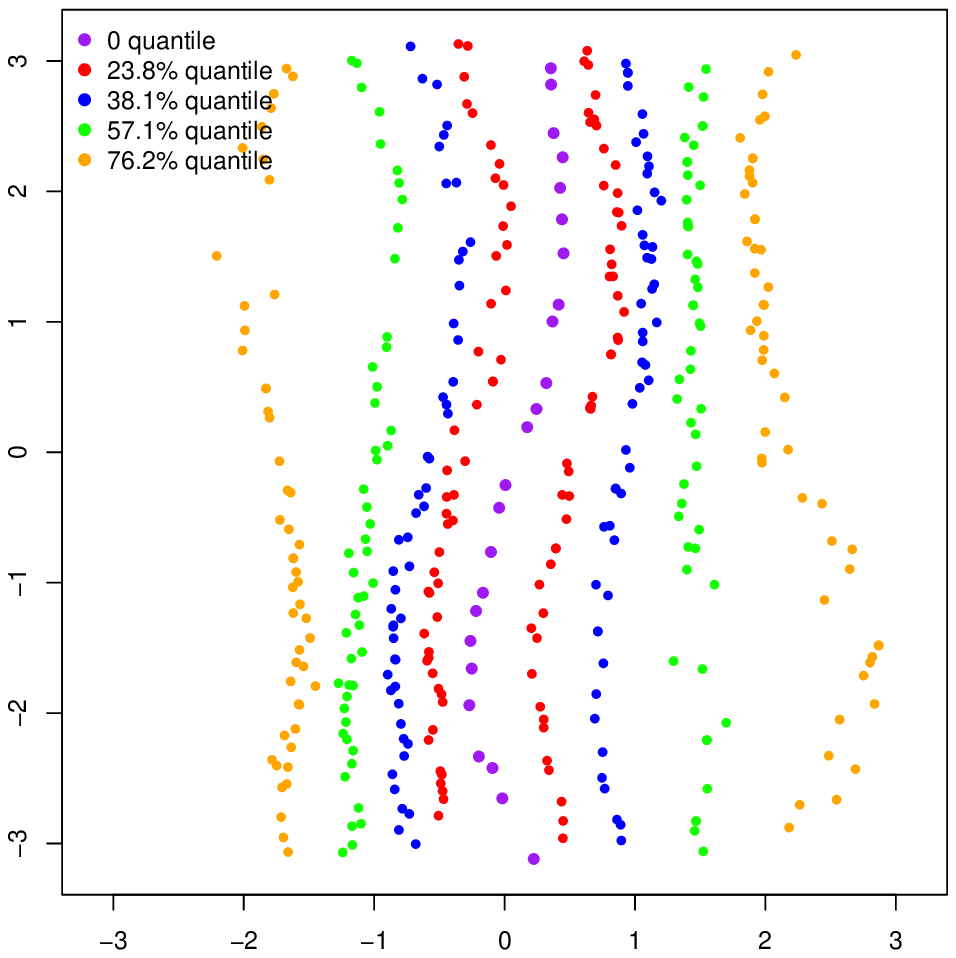}
	\end{subfigure}
	\begin{subfigure}[b]{0.45\textwidth}
		\centering
		\includegraphics[scale=0.3]{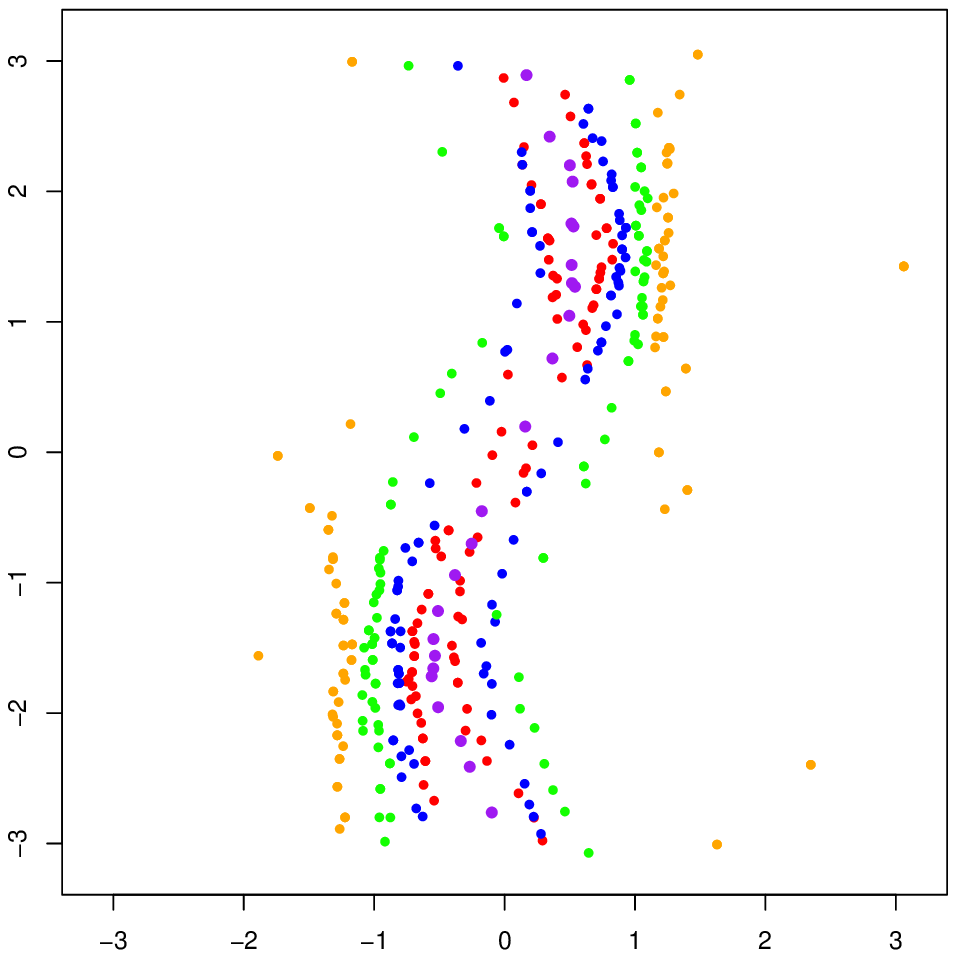}
	\end{subfigure}
	\begin{subfigure}[b]{0.45\textwidth}
		\centering
		\includegraphics[scale=0.15]{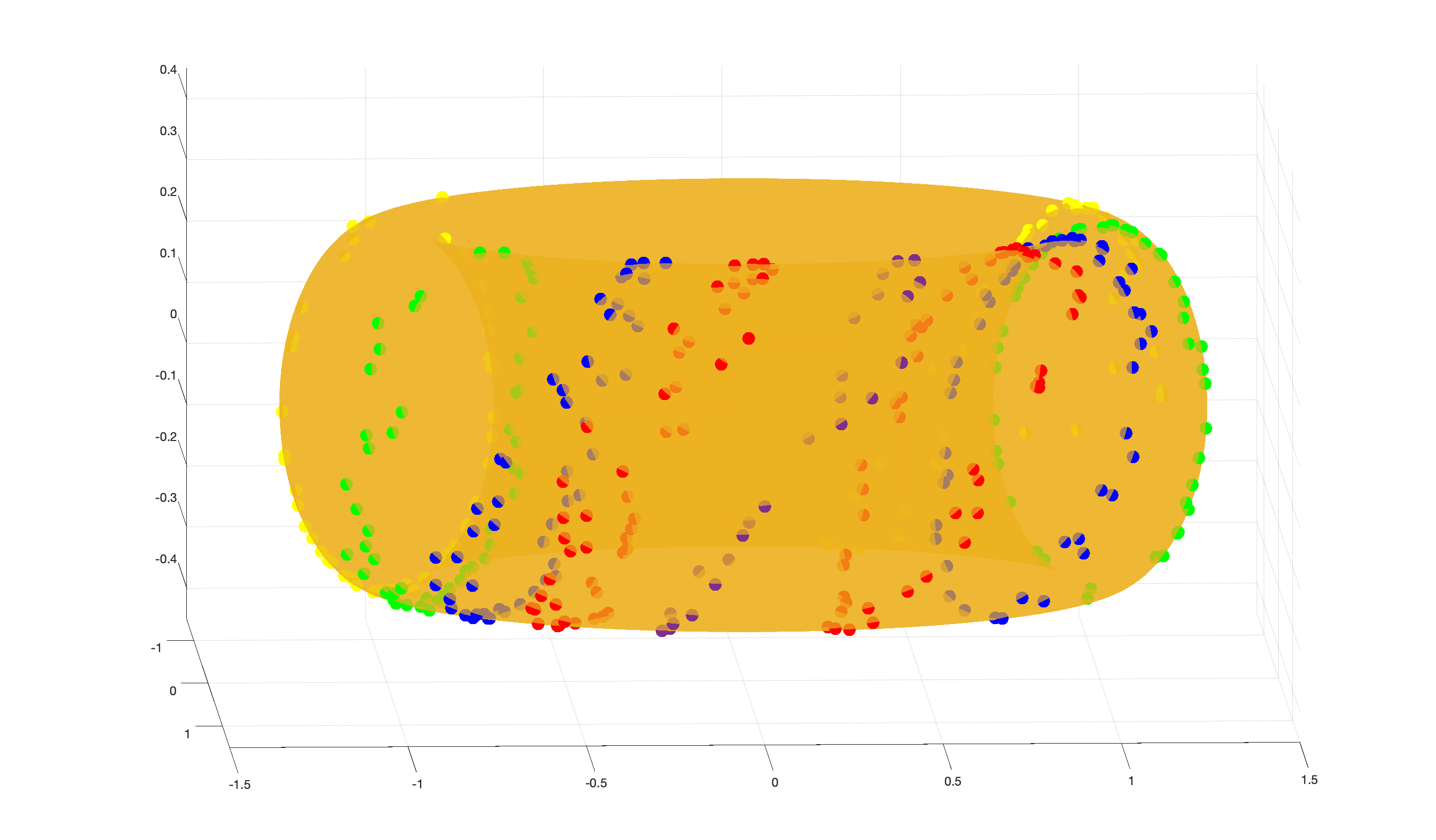} \vspace{-5mm}
	\end{subfigure}
	\begin{subfigure}[b]{0.45\textwidth}
		\centering
		\includegraphics[scale=0.15]{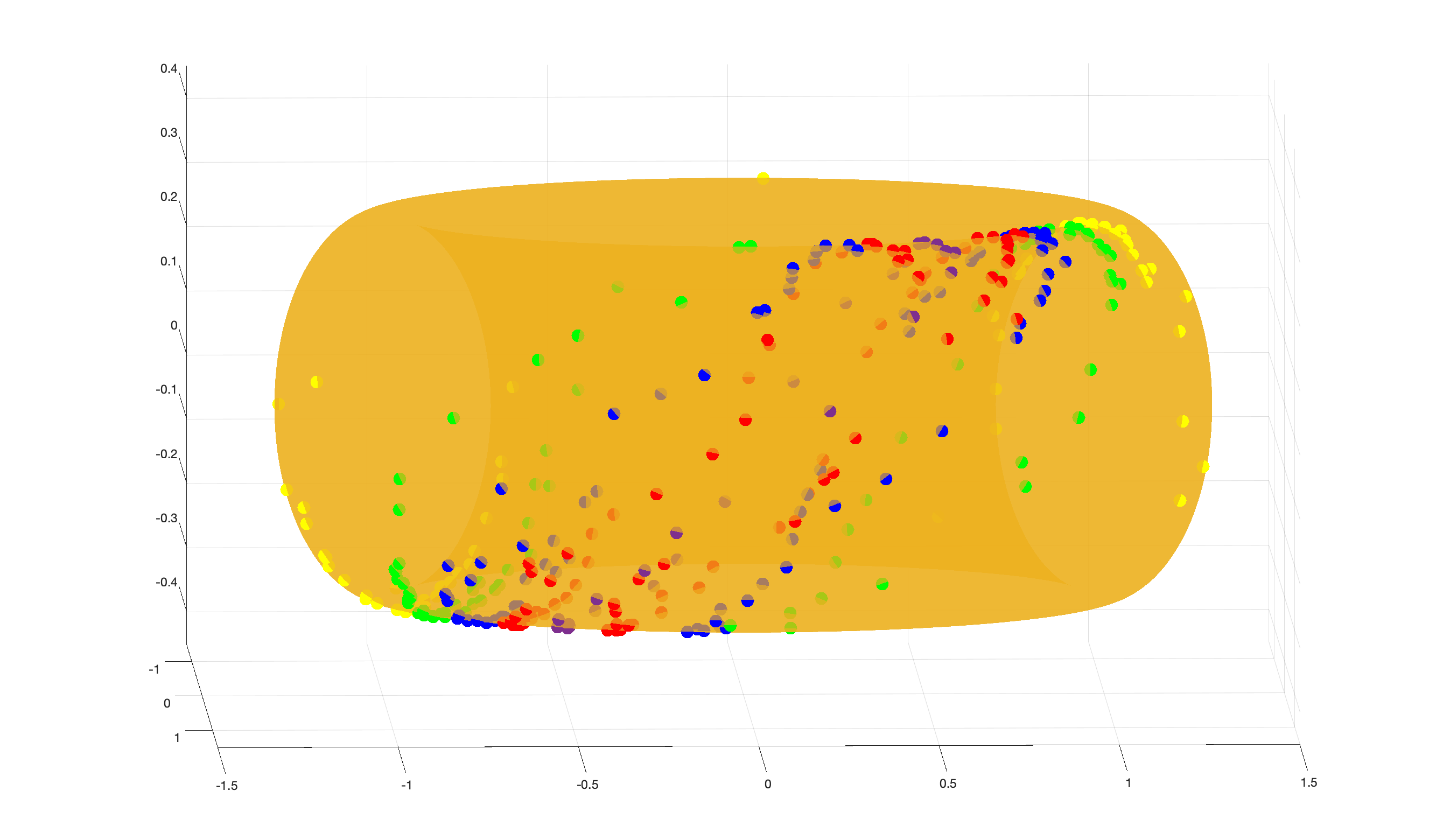} \vspace{-5mm}
	\end{subfigure}
	\caption{\small Equatorial strip-type conditional quantile contours $\mathcal{C}^{(N, n)}_{\wb}\big(\frac{r}{N_R+1}\vert \xb\big)$ ($r=0$ (purple),~5 (red),   8  (blue),  12 (green), and~16 (orange)) for 
		(TS3*) (top-left for $\xb = (0.7, 0.7, \sqrt{0.02})^\top$, top-right for $\xb = (0.7, -0.7, \sqrt{0.02})^\top$, along with their $\mathbb{R}^3$-embedded representations in the second row.}
	\label{Fig:QRContT2diffxStrips}
\end{figure}

Figure~\ref{Fig:QRContourT2diffx} provides plots of $\mathcal{C}^{(N, n)}_{\wb}\big(\frac{r}{N_R+1}\vert \xb\big)$ ($r = 5,\, 8,\, 12$, and 16) for  (TS1*) (top-left for $\xb = (0.75, 0.65, \sqrt{0.015})^\top$ and top-right for $\xb = (0.2, 0.4, \sqrt{0.8})^\top$, corresponding to high and low concentrations of the conditional distributions of $\Yb \in {\cal T}^2$), (TS1**) (second row: left for $\xb = (0.4, 0.4, \sqrt{0.68})^\top$ and right for $\xb = (0, 0, 1)^\top$, corresponding to conditionally strongly positive dependent and independent  marginals of $\Yb\in {\cal T}^2$) and (SS2*) (bottom-left for $\xb = (0, 1)^\top$ and bottom-right for $\xb = (0.6, 0.8)^\top$, corresponding to high and low (hence also more symmetric) conditional concentrations of $\Yb \in {\cal S}^2$). The $\mathbb{R}^3$ representations for (TS1*) and (TS1**) are shown in the third and fourth rows, respectively. The kernel weights are used in the computation. The plots demonstrate how the impact of $\Xb$ on the conditional distributions of $\Yb$---significant differences in dispersion (for (TS1*) and (SS2*)) and dependence (for (TS1**))---are picked up by the conditional quantile contours. 

For equatorial strip-type quantile contours, we consider a variant of (TS3):

\begin{enumerate}
	\item[(TS3*)] ($\mathcal{S}^2\times \mathcal{T}^2$)-valued $(\Xb,\Yb)$ from the same distribution as in (TS3) but\linebreak  with $\kappab_\Xb = (2\vert X_1 - X_2 \vert,\,  0)^\top\!$ and $\lambda = e^{\vert X_1 - X_2 \vert}$.
\end{enumerate}

Figure~\ref{Fig:QRContT2diffxStrips} shows  plots of $\mathcal{C}^{(N, n)}_{\wb}\big(\frac{r}{N_R+1}\vert \xb\big)$ ($r = 0, 5,\, 8,\, 12$, and 16) at $\xb = (0.7, 0.7, \sqrt{0.02})^\top$ (left column) where the two components of $\Yb\in {\cal T}^2$ are both uniform over $\mathcal{S}^1$ and independent, and at $\xb = (0.7, -0.7, \sqrt{0.02})^\top$ (right column) where the first component is more concentrated  and positively  dependent on the other component around the center. Again, the empirical conditional quantile contours correctly reflect these variations. 

\subsection{Numerical illustration: Euclidean covariates}

\begin{figure}[b!]
	\begin{subfigure}[b]{0.48\textwidth}
		\includegraphics[trim=160 40 80 90,clip,scale=0.19]{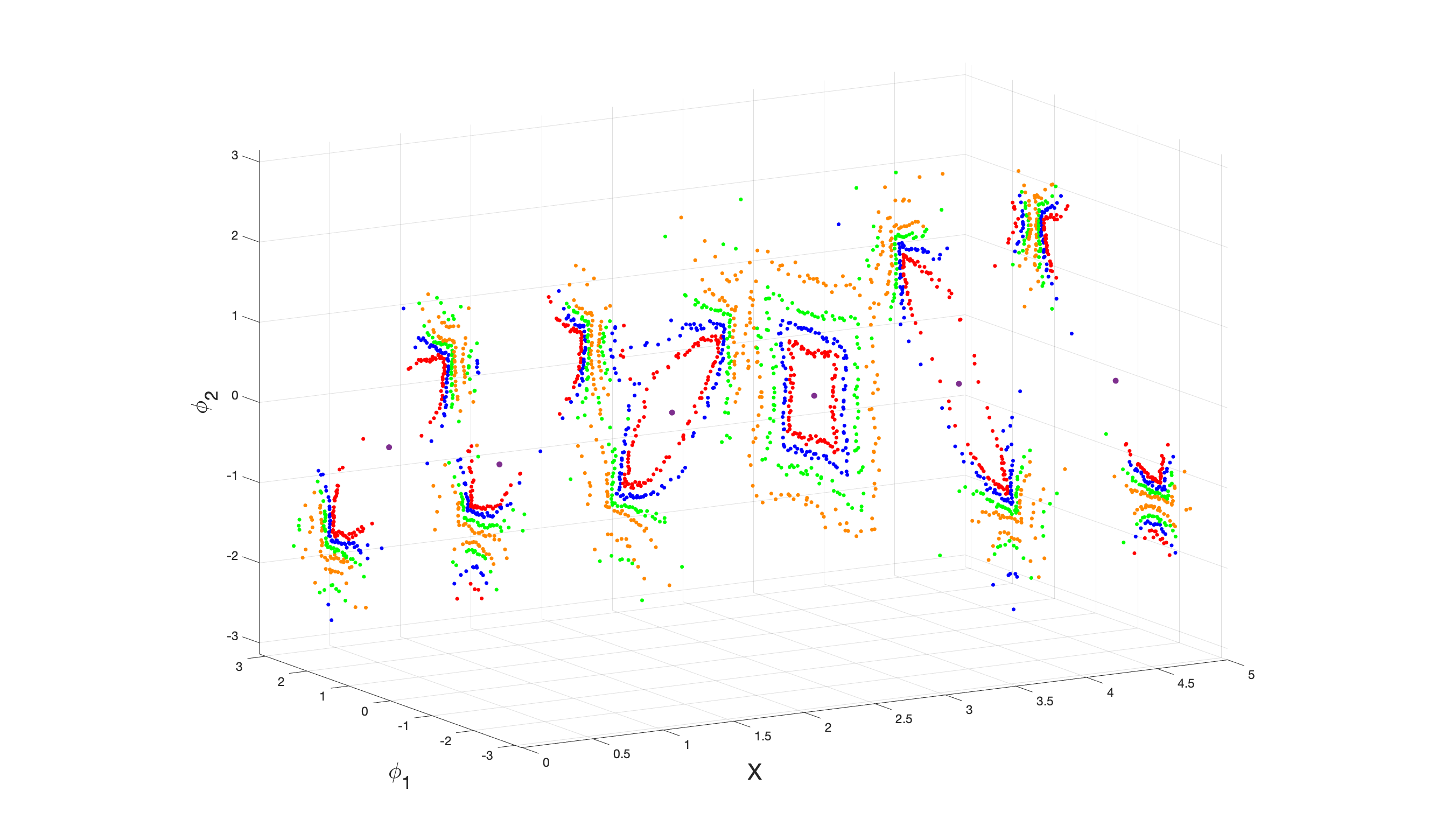} 
	\end{subfigure}
	\begin{subfigure}[b]{0.48\textwidth}
		\includegraphics[trim=160 40 80 90,clip,scale=0.19]{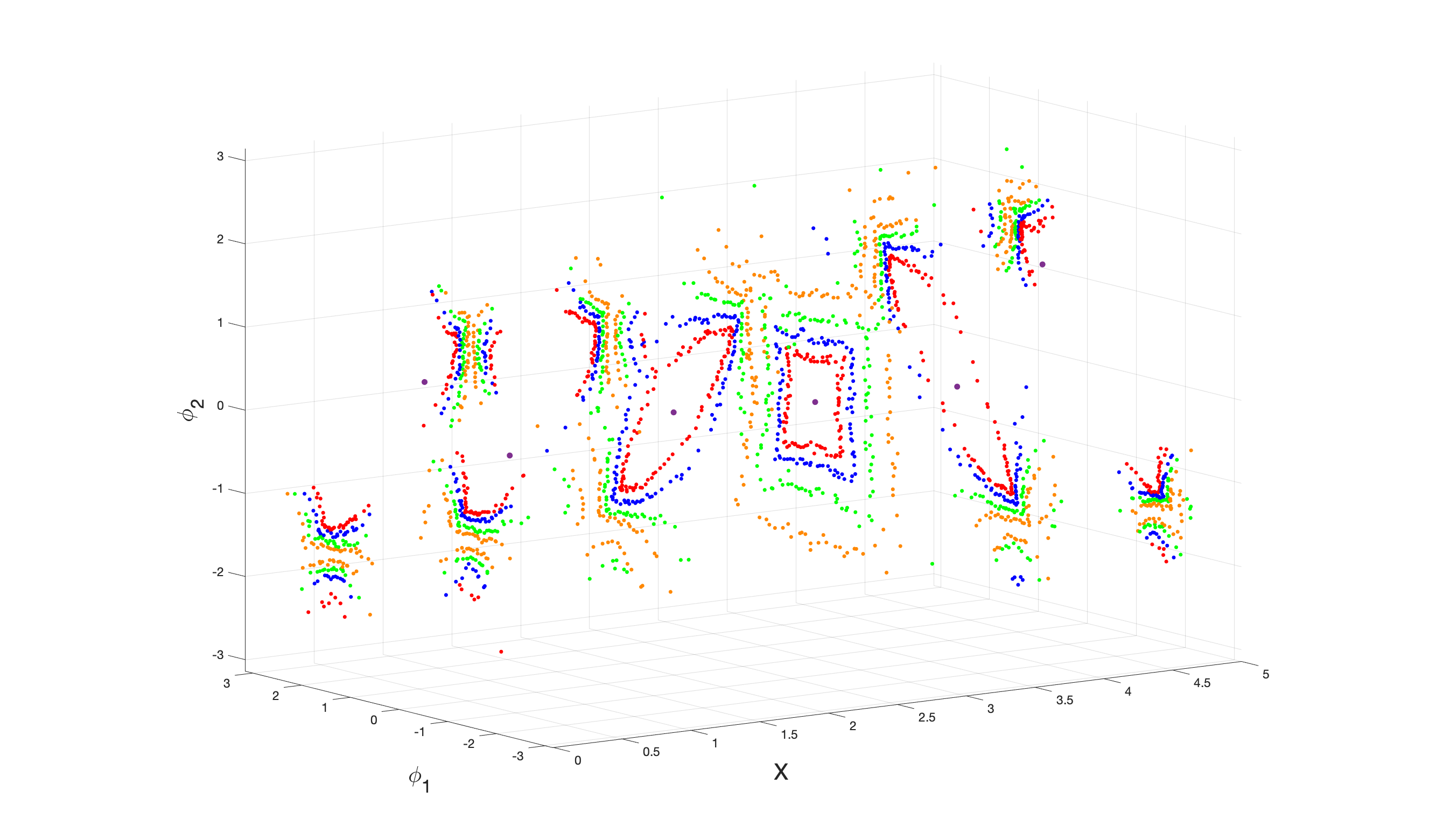} 
	\end{subfigure}
	\begin{subfigure}[b]{0.48\textwidth}
		\includegraphics[trim=140 40 80 90,clip,scale=0.19]{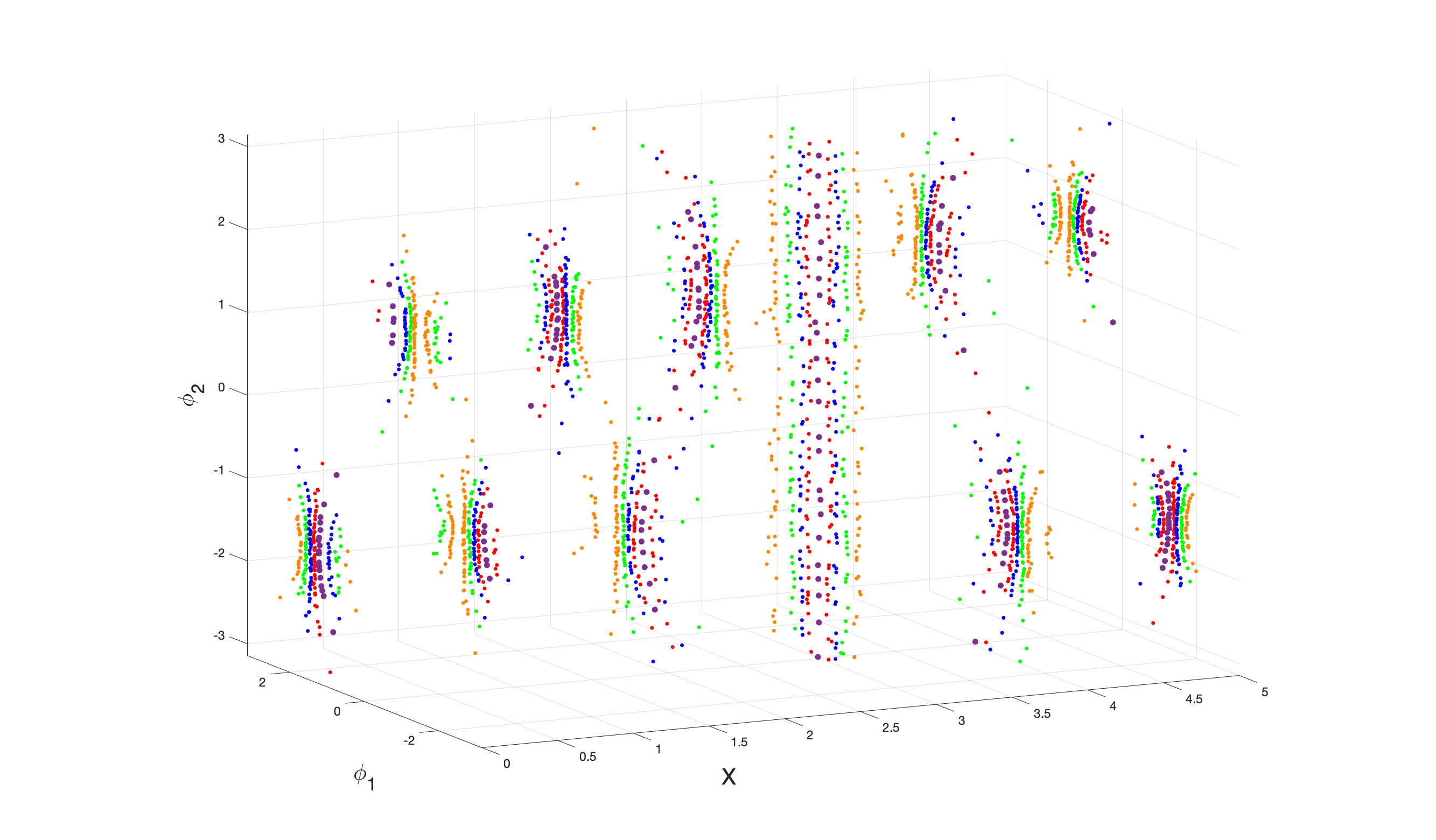} 
	\end{subfigure}
	\begin{subfigure}[b]{0.48\textwidth}
		\includegraphics[trim=140 40 80 90,clip,scale=0.19]{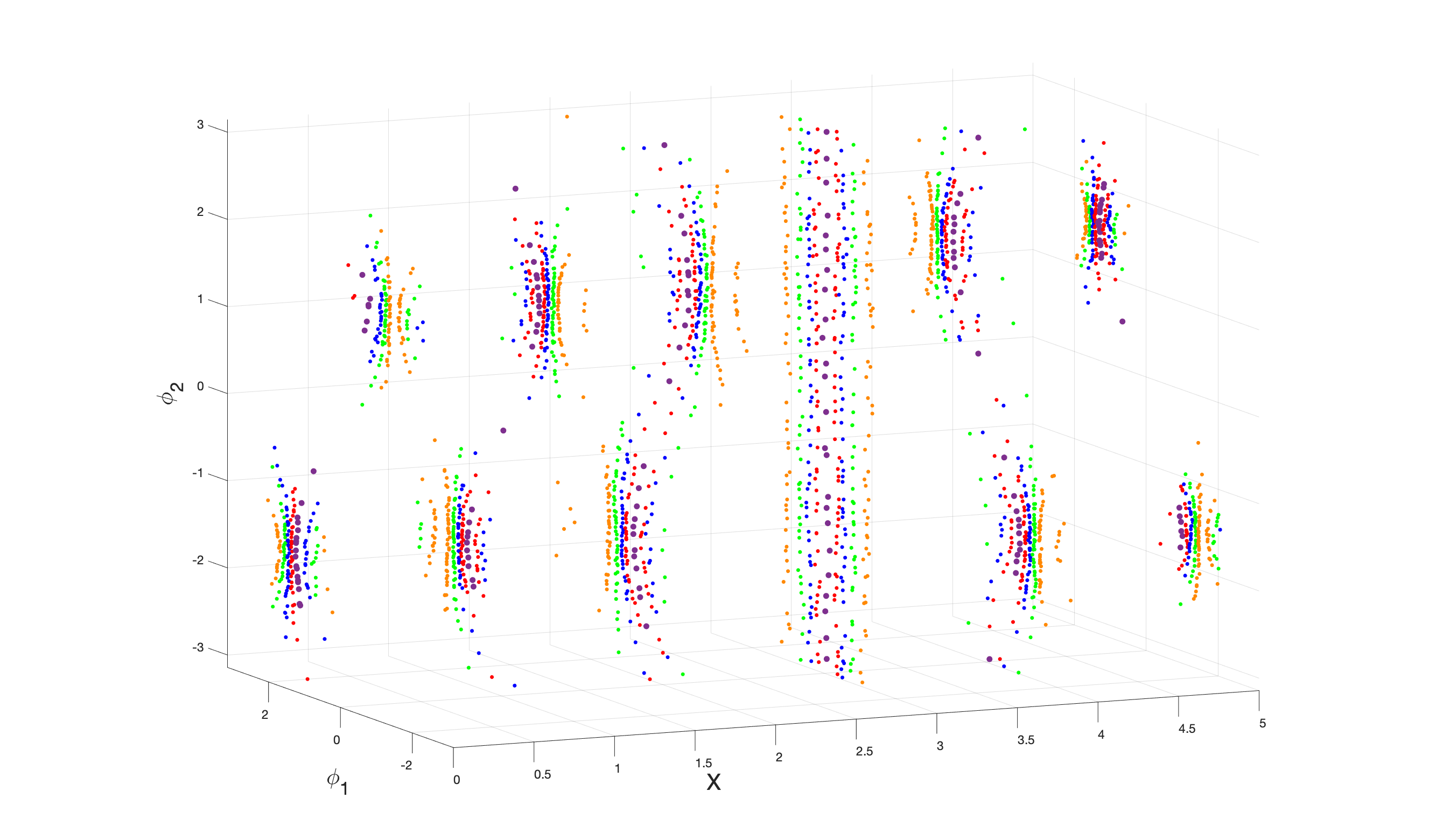} 
	\end{subfigure}
	\begin{subfigure}[b]{0.45\textwidth}
		\includegraphics[trim=190 40 80 90,clip,scale=0.21]{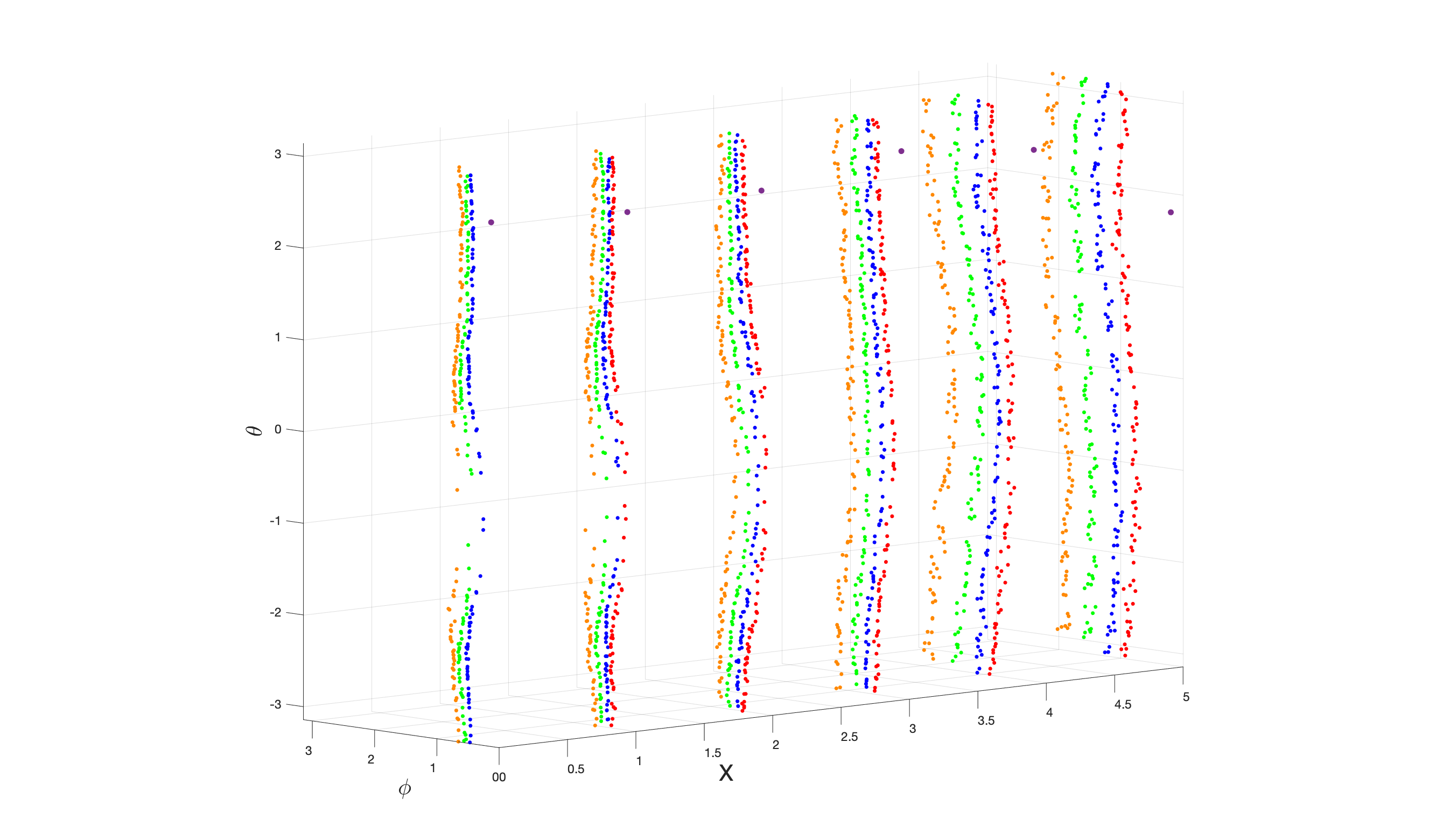} \vspace{-2mm}
	\end{subfigure}
	\begin{subfigure}[b]{0.45\textwidth}
		\includegraphics[trim=190 40 80 90,clip,scale=0.21]{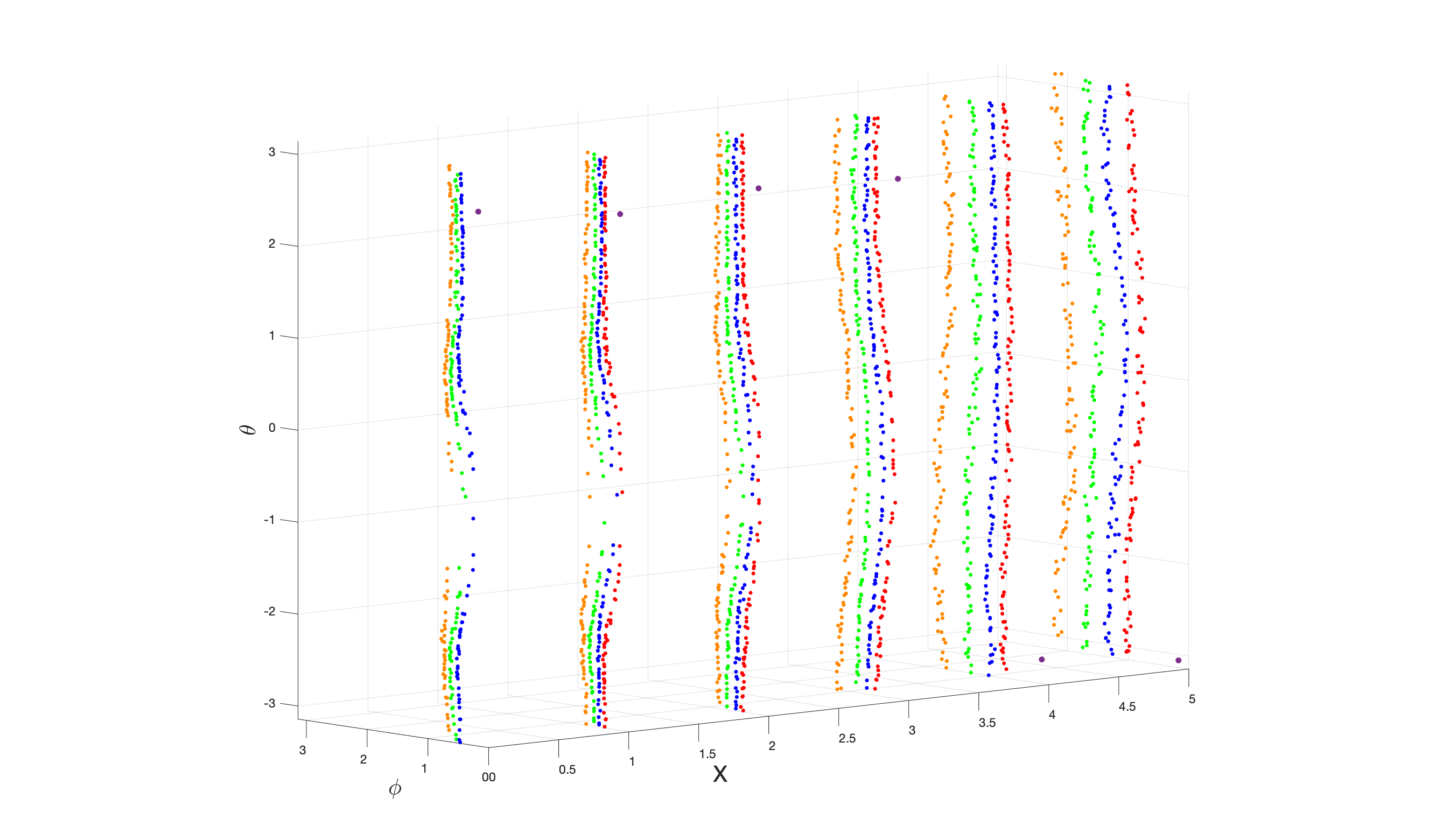} \vspace{-2mm}
	\end{subfigure}
	\caption{\small Conditional quantile contours $\mathcal{C}^{(N, n)}_{\wb}(\frac{r}{N_R+1}\vert X)$ ($j=0$ (purple), 5 (red),   8  (blue),~12 (green), and  16 (orange))  at various   $x$ values of $X  \in \mathbb{R}$. Top row:  (TR1) with cap-type quantile contours; middle row: (TR2) with strip-type quantile contours; bottom row: (SR) with cap-type quantile contours. 
	The left and right columns are for the $k$-NN and kernel weights, respectively.   The coordinate system adopted here is $(\phi_1,\phi_2)\in [-\pi, \pi)^2]$ for the\linebreak  2-torus~${\mathcal T}^2={\mathcal S}^1\times {\mathcal S}^1$, $(\phi, \theta)$ and $(\phi, \theta)\in [0,\pi)\times [-\pi, \pi)$ for the 2-sphere~${\mathcal S}^2$.
	}
	\label{Fig:QRContourXR}
\end{figure} 

Representations are easier, too, when the covariate $X$   is real-valued. In this section, we consider the following  cases of~$({\mathbb R}\times \mathcal{T}^2)$- and $({\mathbb R}\times \mathcal{S}^2)$-valued $(X,\Yb)$. 
\begin{enumerate}
\item[(TR1)] $({\mathbb R}\times \mathcal{T}^2)$-valued $(X,\Yb)$; the real-valued covariate $X$ is ${\rm U}_{[0, 5]}$  and, conditional on $X=x$, the $ \mathcal{T}^2$-valued response   $\Yb$ is  ${\rm BSvM}(\mub, \kappab, \lambda)$ as defined in \eqref{eq.BSvM}, with $\mub = {\bf 0}$, $\kappa_1 = \kappa_2= 2$, and dependence parameter $\lambda = x^2-9$;
{\item[(TR2)] $(X,\Yb)$ from the same distribution as in (TR1) but with $\kappa_2= 0$;}
{\item[(SR)] $({\mathbb R}\times \mathcal{S}^2)$-valued $(X,\Yb)$;  the real-valued covariate $X$ is ${\rm U}_{[0, 5]}$ and, conditional on~$X=x$, the $\mathcal{S}^2$-valued response   $\Yb$ is from the mixture as in (SS2) but with mixture components concentration parameters   $\kappa_1, \kappa_2, \kappa_3 = (x-5)^2$.}
\end{enumerate}

Figure~\ref{Fig:QRContourXR} shows plots of the empirical conditional quantile contours $\mathcal{C}^{(N, n)}_{\wb}(\frac{r}{N_R+1}\vert \xb)$ for~$r = 0, 5, 8, 12$, and 16, at various   $\xb$ values. For (TR1) (top row) and (SR)  (bottom row), we plot the cap-type contours; for (TR2) (middle row), the contours are of the equatorial  strip type. The plots for the $k$-NN and kernel weights are shown in the left and right columns, respectively; again, their respective performances are quite similar. Note that for the (SR) case, in order to visualize the contours, we transform the Euclidean coordinates of $\Yb \in \mathcal{S}^2$ to the two-dimensional spherical coordi\-nates~$(\phi, \theta)^\top$ with values in $[0, \pi]\times [-\pi, \pi)$, with north pole $(0, 0, 1)^\top$. As shown in the bottom row of\linebreak  Figure~\ref{Fig:QRContourXR}, when $x$ decreases, conditional concentrations increase, leading to more   multimodal and nonconvex  contours. For $(X, \Yb)\in ({\mathbb R}\times \mathcal{T}^2)$, as shown in the first and second rows, $x$ has great impact on the dependence between  the two compo\-nents~$(\phi_1, \phi_2) \in [-\pi, \pi)^2$: as $x$ increases, that dependence decreases from positive to nil,  then   to negative.

\section{A real-data analysis}

\begin{figure}[b!]
	\centering
	\begin{subfigure}[b]{0.32\textwidth}
		\centering
		\includegraphics[scale=0.27]{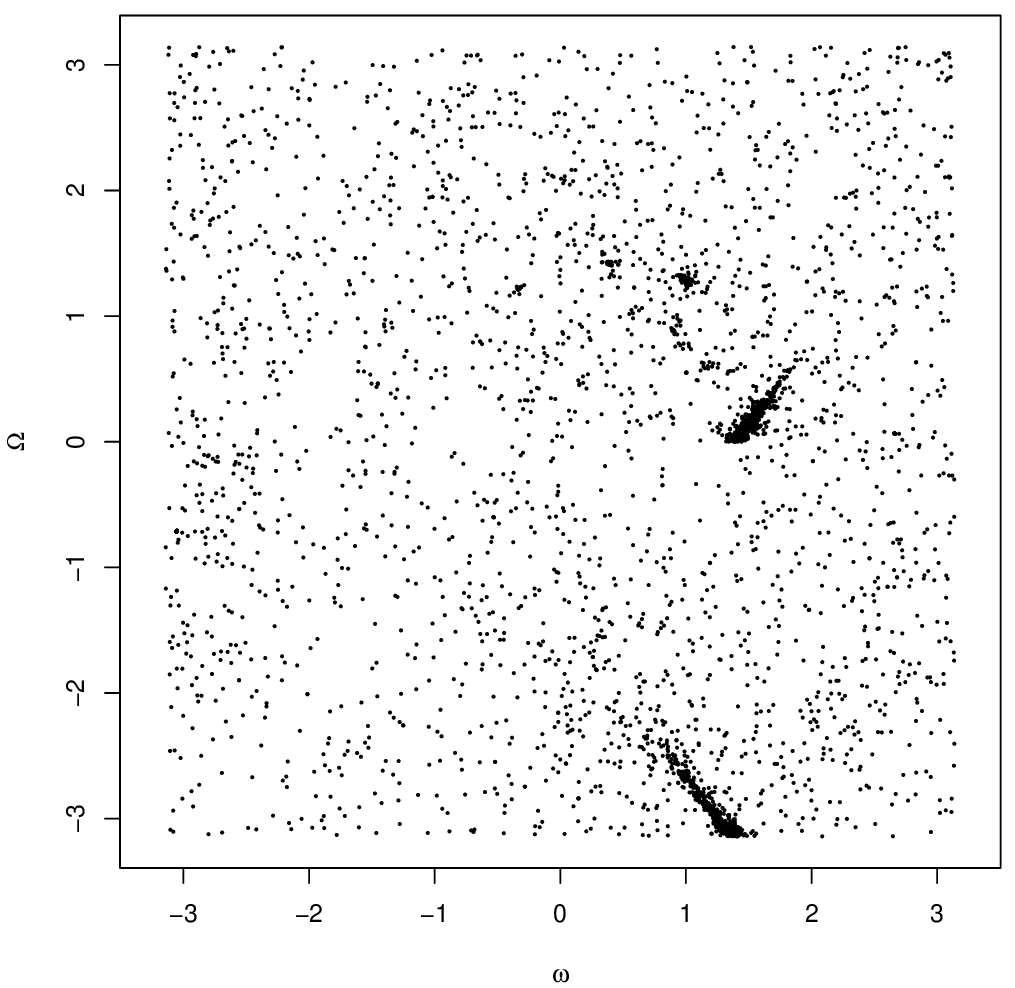}
	\end{subfigure}
	\begin{subfigure}[b]{0.32\textwidth}
		\centering
		\includegraphics[scale=0.27]{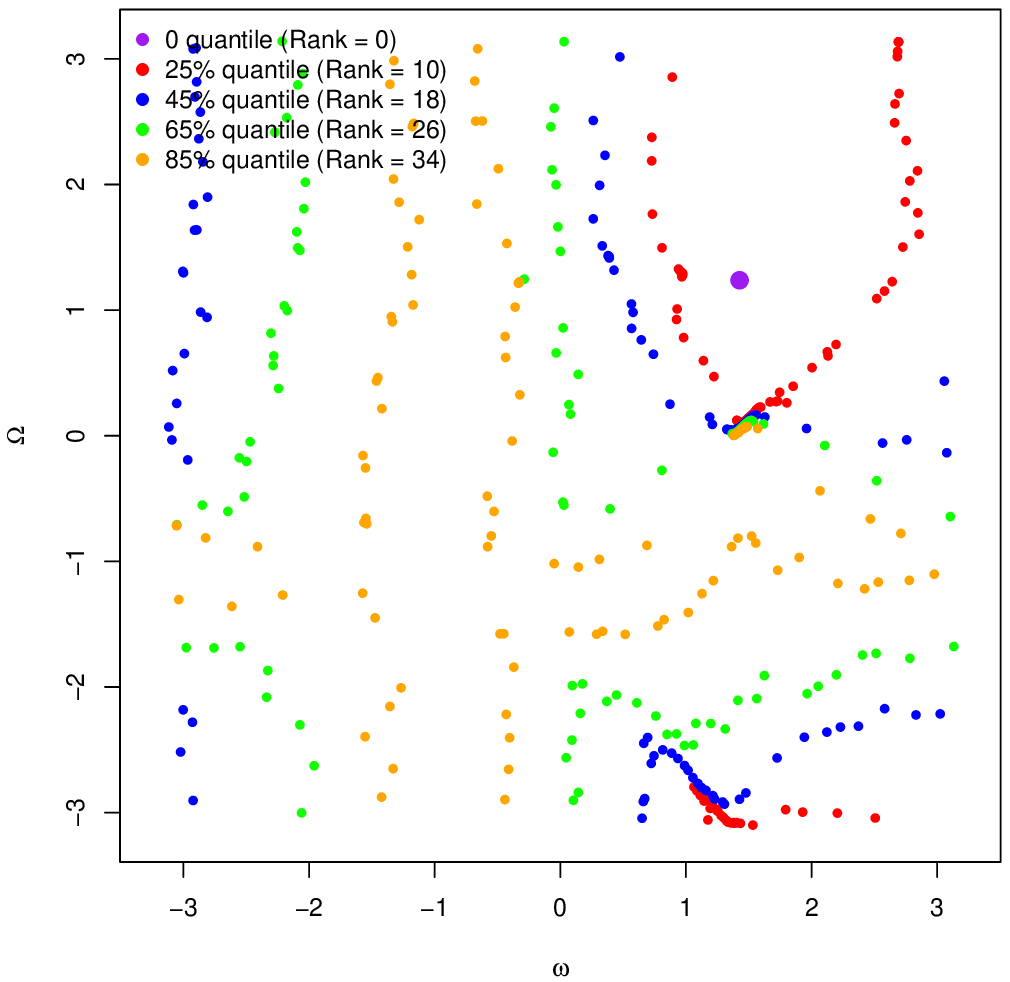}
	\end{subfigure}
	\begin{subfigure}[b]{0.32\textwidth}
		\centering
		\includegraphics[scale=0.27]{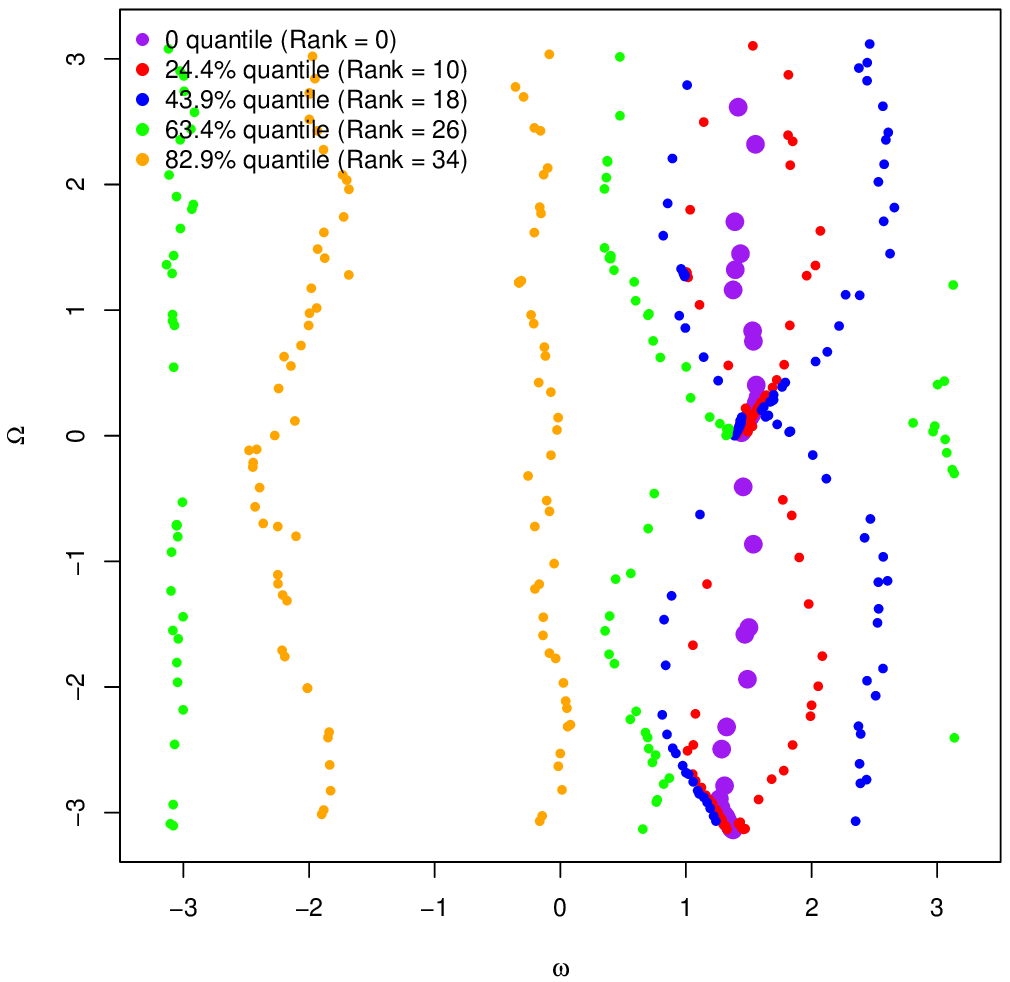}
	\end{subfigure}
	\begin{subfigure}[b]{0.48\textwidth}
		\centering
		\includegraphics[scale=0.18]{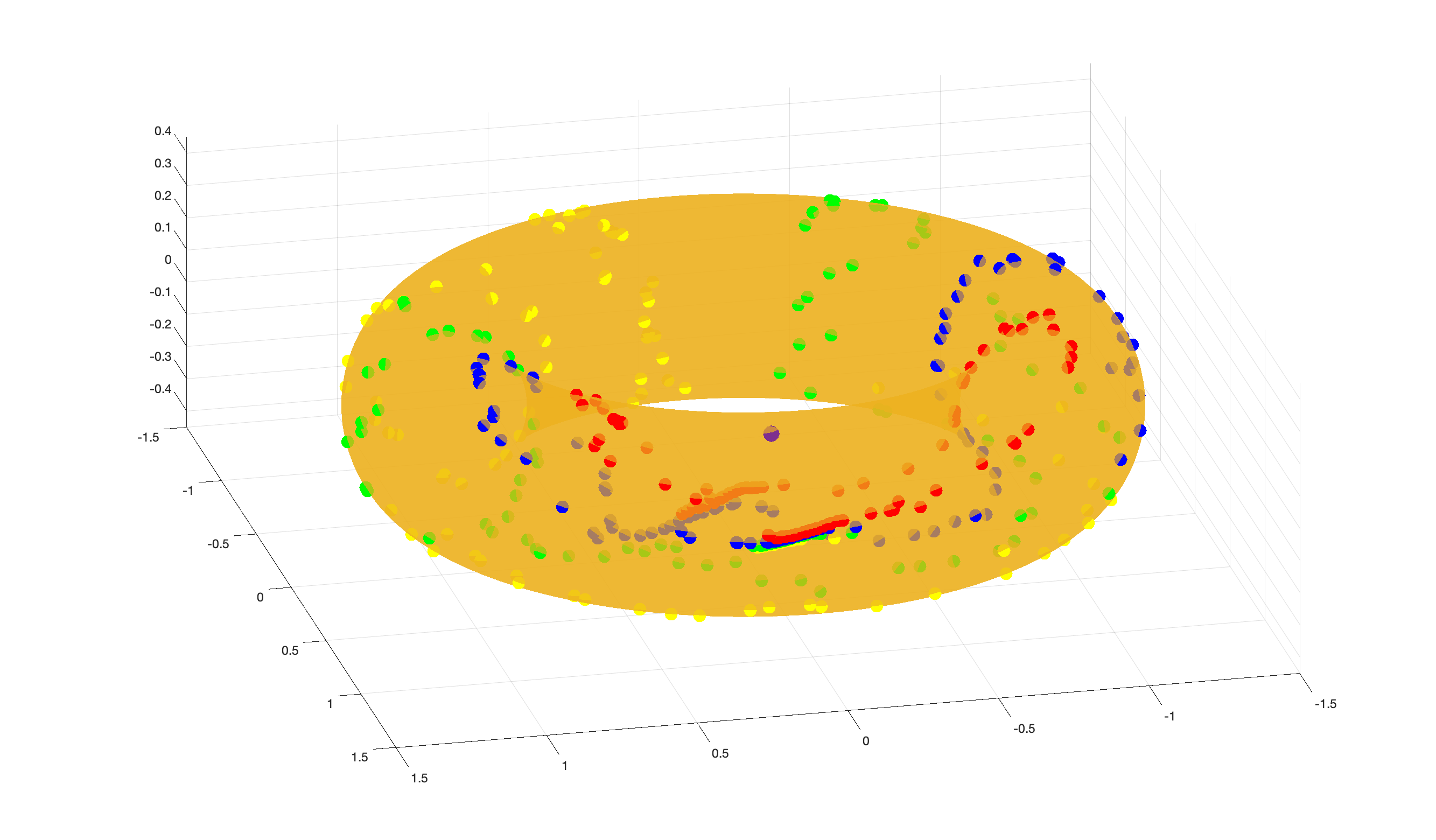}
	\end{subfigure}
	\begin{subfigure}[b]{0.48\textwidth}
		\centering
		\includegraphics[scale=0.18]{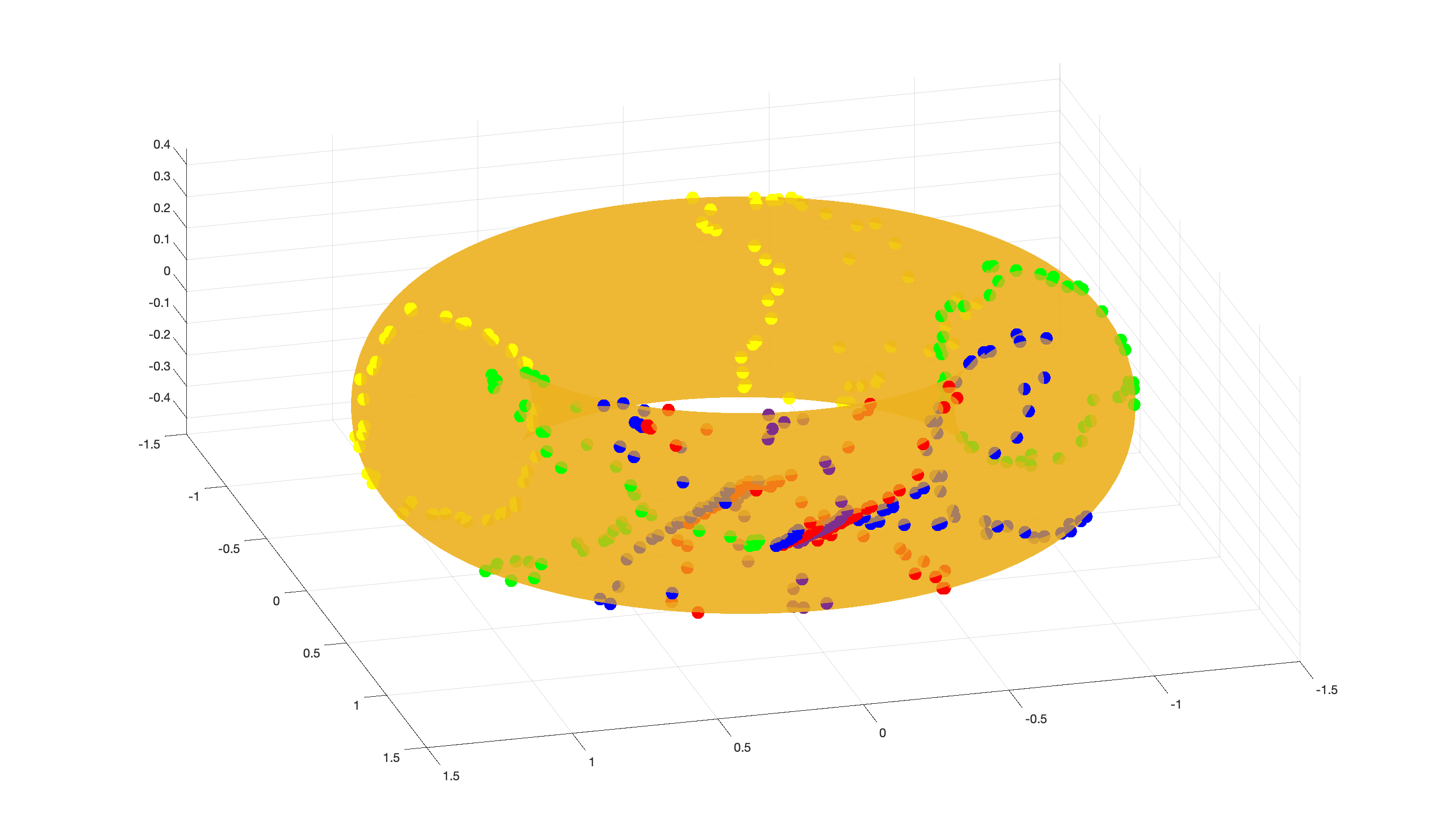}
	\end{subfigure}
	\caption{\small Plots of $(\omega, \Omega)^\top \in \mathcal{T}^2$ (top-left), the cap-type (top-middle) and strip-type  (top-right) empirical quantile contours for the comets dataset, along with their $\mathbb{R}^3$-embedded representations in the 2nd row.}
	\label{Fig:QRComets}
\end{figure}

We apply our novel concept of quantile contour to an analysis of  a dataset on comets, downloaded from the website of Jet Propulsion Laboratory (\url{https://ssd.jpl.nasa.gov/tools/sbdb_query.html}). The dataset contain orbital information of $3901$ comets, including longitude of the ascending node~$\Omega \in [-\pi, \pi)$, argument of perihelion $\omega \in [-\pi, \pi)$, and so on. The orbital planes of a comet and the earth intersect at two points/nodes. The \textit{ascending node} is the point where the comet crosses the ecliptic plane and heads from south to north. The longitude of the ascending node $\Omega$ is the angle, measured eastward (in the direction of the earth's  orbital motion) from the vernal equinox to the ascending node of the comet's orbit. The argument of the perihelion $\omega$ is the angle between 
the line connecting the center of the comet's orbit and the perihelion point and the line of the ascending node (on the earth's  orbital plane). Based on this data, we examine the empirical quantile contours of $(\omega, \Omega)^\top \in \mathcal{T}^2$.

The scatterplot of $(\omega, \Omega)^\top$ in the top-left panel of Figure~\ref{Fig:QRComets} shows that there are two small areas (around $[1, 2]\times\left( [-\pi, -2.5] \bigcup [0, 0.5] \right)$) in which $(\omega, \Omega)^\top$ is quite concentrated, while in other areas the observed points are more dispersed and uniformly distributed. Here we provide two types of quantile contours $\mathcal{C}\n\big(r/{(n_R+1)}\vert \xb\big)$,\linebreak   for~$r = 0, 10, 18, 26, 34$: the cap-type (top-middle panel, with $n_0=1$, $n_R = 39$,\linebreak  and~$n_S=100$) and the equatorial strip-type  (top-right panel, with $n_0=61$, $n_R = 40$ and $n_S=96$). The empirical   quantile contour of order 0 (purple point(s)) is at $\omega \approx 1.5$. For both types of  contours, the low level contours ($r = 10, 18$) are more skewly nested around the   contour of order 0 than the high level ones ($r= 26, 34$). This is consistent with the scatterplot where the points are more concentrated around $\omega \in [1, 2]$.


\bibliographystyle{imsart-nameyear.bst}
\bibliography{QRmanifold.bib}  


\appendix

\section{Some basic concepts of Riemannian geometry}\label{App.Rieman}

This section collects some basic concepts of Riemannian geometry that are used in the main part of the paper. For more details, we refer to the monographs by \cite{petersen2006riemannian}, \cite{jost2008riemannian}, and \cite{lee2012smooth}. \smallskip

To begin with, recall that a {\it Riemannian manifold}   is a pair $({\cal M}, g)$, where ${\cal M}$ is a {\it smooth manifold} and $g$ is a {\it Riemannian metric} on ${\cal M}$. Before  defining what is meant with a {\it smooth}  manifold, the following notions are needed.

\begin{definition}\label{def.manifold}{\rm 
A topological space ${\cal M}$ is a {\it topological $m$-manifold}  if (i) it is a {\it Hausdorff space}, (ii) it is {\it second-countable}, and (iii) it is {\it locally $m$-dimensional Euclidean}, that is, each point in ${\cal M}$ has a neighborhood that is homeomorphic to an open subset of $\mathbb{R}^m$.
}\end{definition}

\begin{definition}\label{def.chart}{\rm 
A {\it chart} on a topological $m$-manifold ${\cal M}$ is a pair $(U, \varphi)$, where~$U$ is an open subset of ${\cal M}$ and $\varphi: U\rightarrow V$ is a homeomorphism from $U$ to an open sub\-set~$V = \varphi(U) \subset \mathbb{R}^m$. An {\it atlas} ${\cal A}$ for ${\cal M}$ is a collection of charts whose domains cover~${\cal M}$. 
}\end{definition}

\begin{definition}\label{def.atlas}{\rm 
Two charts $(U, \varphi)$ and $(W, \psi)$ on a topological $m$-manifold {\cal M} are {\it smoothly compatible} if either $U\cap W=\emptyset$ or the {\it transition map} 
$$\psi\circ \varphi^{-1}: \varphi(U\cap W)\rightarrow \psi(U\cap W)$$ is a diffeomorphism. An {\rm atlas} ${\cal A}$ is called a {\it smooth atlas} if any two charts on ${\cal A}$ are smoothly compatible. A smooth atlas ${\cal A}$ on ${\cal M}$ is {\it maximal} if it is not  contained in any strictly  larger smooth atlas.
}\end{definition} 

Now we are ready to provide the definition of  a {\it smooth manifold}.

\begin{definition}\label{def.SmoothManifold}{\rm 
A {\it smooth manifold} is a pair $({\cal M}, {\cal A})$, where ${\cal M}$ is a topological manifold and ${\cal A}$ is a maximal smooth atlas on ${\cal M}$.
}\end{definition}

As in Euclidean spaces, a Riemannian metric on a smooth manifold ${\cal M}$  is defined via  an  inner product on its tangent spaces (see, e.g., \citet[Chapter~3]{lee2012smooth}).

\begin{definition}\label{def.metric}{\rm 
A {\it Riemannian metric} $g$ on a smooth manifold ${\cal M}$   assigns an inner product
$g_{\yb}(\cdot, \cdot) = \langle\cdot, \cdot\rangle_{\yb}$
on the tangent space $T_{\yb}{\cal M}$ to each point~$\yb \in {\cal M}$.
}\end{definition}

Note that, using a local coordinate chart $(U, \varphi)$ (with $U \ni \yb$), the Riemannian\linebreak  metric~$g$ can be characterized by a positive definite symmetric matrix $(g_{ij}(\yb))_{i, j=1, \ldots, p}$, where~$g_{ij}(\yb)$ is such that  
$$\langle\ub, \vb\rangle_{\yb} = \sum_{i,j=1}^{p} g_{ij}(\yb) u_i v_j$$
for any pair $\ub = (u_1, \ldots, u_{p})^\top$, $\vb = (v_1, \ldots, v_{p})^\top$ of tangent vectors  in $T_{\yb}{\cal M}$. 
   Conse\-quently,~$g_{ij}(\yb) = \langle\xib_i, \xib_j\rangle_{\yb}$ for vectors $\xib_i, \, \xib_j$ of the orthonormal basis of $T_{\yb}{\cal M}$. See, e.g., \citet[Section~1.4]{jost2008riemannian} for   details.
\smallskip

The concept of {\it geodesic} on a Riemannian manifold is   essential for the concepts proposed in this paper.\footnote{Note that $\gamma{''}$ here is the {\it covariant derivative} of $\gamma\pr$ along $\gamma\pr$, where $\gamma\pr(t)$ denotes the \textit{velocity} of~$\gamma$ at~$t$, namely, the vector in $T_{\gamma(t)}{\cal M}$ the coordinates of which are  the derivatives with respect to $t$  of the local coordinates of $\gamma(t)$; see \citet[Chapter~3]{lee2012smooth} and \citet[Chapter 5]{petersen2006riemannian}  for precise definitions of $\gamma\pr(t)$ and~$\frac{{\rm d}^2\gamma}{{\rm d}t^2}$, respectively.}

\begin{definition}\label{def.geodesic}{\rm 
Let $I$ denote a connected subset of $\mathbb{R}$. A {\it geodesic} on $({\cal M},g)$ is a smooth curve $\gamma: I \rightarrow {\cal M}$ with~$t \mapsto \gamma(t)$ such that $\gamma{''} \coloneqq \frac{{\rm d}^2\gamma}{{\rm d}t^2} = \boldsymbol{0}$.  
}\end{definition}

The intuitive interpretation of Definition~\ref{def.geodesic} is that a geodesic is a ``constant speed" curve, that is, a curve that is parameterized proportionally to arc length. 
 A Riemannian manifold ${\cal M}$ is called  {\it geodesic complete} if any geodesic is the subset of a geodesic
defined for all $t\in \mathbb{R}$.   Define  the distance $d(\xb, \yb)$ between $\xb, \yb \in {\cal M}$ as the length of the shortest piecewise smooth curve connecting $\xb$ and $\yb$, namely, 
$$d(\xb, \yb) \coloneqq  \inf\{\ell(\gamma): \gamma: [0, 1] \rightarrow {\cal M}, \gamma \ \text{is piecewise}\ C^\infty ,\ \ \gamma(0) = \xb, \ \text{and}\ \gamma(1) = \yb\},$$
where the length $\ell(\gamma)$ of $\gamma$ is $\ell(\gamma) \coloneqq \int_0^1 g_{\gamma(t)}(\gamma\pr(t),\gamma\pr(t)){\rm d}t.\vspace{1mm}$  By the Hopf-Rinow Theorem (see, e.g., \citet[Chapter~5]{petersen2006riemannian}), if ${\cal M}$, equipped with the distance $d$, is a complete metric space (by Corollary 5.3.2 in \cite{petersen2006riemannian}, a sufficient condition  for this is a compact ${\cal M}$),  then ${\cal M}$ must be geodesically complete, and any two points $\yb_1, \yb_2 \in {\cal M}$ can be connected by a geodesic of length $d(\yb_1, \yb_2)$, i.e., a geodesic of  shortest length. 

Closely related with the concept of geodesic is the exponential map from a tangent space $T_{\yb}{\cal M}$ to ${\cal M}$. For $\vb \in T_{\yb}{\cal M}$, let $\gamma^{\vb}_{\yb}$ be the unique geodesic with $\gamma^{\vb}_{\yb}(0) = \yb$ and~$(\gamma^{\vb}_{\yb})\pr(0) = \vb$;
sufficient conditions  for the existence and uniqueness, for all $\yb$ and~$\vb$,  of such a geodesic 
 are given, e.g., in Theorem 11 and Lemma 7 of \citet[Chapter~5]{petersen2006riemannian}. Let  $[0, \ell_{\vb})$ be the nonnegative part of the maximal interval on which $\gamma^{\vb}_{\yb}$ is defined. Note that $\gamma^{a\vb}_{\yb}(t) = \gamma^{\vb}_{\yb}(at)$ for all $a>0$ and $t<\ell_{a\vb}$. Let $V_{\yb} \subset T_{\yb}{\cal M}$ be the set of vectors $\vb$ such that $\ell_{\vb}>1$, so that $\gamma^{\vb}_{\yb}(t)$ is defined on $[0, 1]$. 
\begin{definition}{\rm 
Call  {\it exponential map} of ${\cal M}$ at~$\yb$ the mapping $\exp_{\yb}: \vb \mapsto \gamma^{\vb}_{\yb}(1)$ from $V_{\yb}$ to ${\cal M}$. 
}\end{definition}


\color{black}

Note that, since the Riemannian manifold ${\cal M}$ of interest is compact, by the Hopf-Rinow Theorem, the exponential map $\exp_{\yb}$ is well defined over $T_{\yb}{\cal M}$ at any $\yb \in {\cal M}$. The associated geodesic $\gamma^{\vb}_{\yb}(t) = \exp_{\yb}(t\vb)$ is a minimizing curve---the shortest curve connecting $\yb$ and any point  on $\gamma^{\vb}_{\yb}(t)$ in the  neighborhood of $\yb$; that minimizing property, however, fails at  the so-called {\it cut locus} of $\yb$.

\begin{definition}{\rm 
Let $({\cal M}, g)$ be a (geodesically) complete Riemannian manifold. The {\em cut locus of} $\yb$ in $T_{\yb}{\cal M}$ is the set 
\begin{align*}
\text{\rm Cut}^T(\yb) \coloneqq &\{\vb \in T_{\yb}{\cal M}: \exp_{\yb}(t\vb) \ \text{is a minimizing geodesic for} \ t \in [0, 1] \\ 
&\ \ \ \ \  \ \ \  \text{but fails to be minimizing for} \  t = 1+\epsilon \ \text{for any} \ \epsilon>0\},
\end{align*}
and the {\it cut locus} of $\yb$ in ${\cal M}$ is
$$\text{\rm Cut}(\yb) \coloneqq \{\exp_{\yb}(\vb): \vb \in \text{\rm Cut}^T(\yb)\}.$$
}\end{definition}


For ${\cal M} = {\cal S}^{p}$ with $p\geq 1$, $\text{\rm Cut}(\yb)$ is the antipodal point of $\yb$; for ${\cal M} = \mathcal{T}^2$, 
$$\text{\rm Cut}(\yb) = \{\text{\rm Cut}(\yb_1)\times {\cal S}^1\} \bigcup \{{\cal S}^1\times \text{\rm Cut}(\yb_2)\},$$
where $\yb = (\yb_1^\top, \yb_2^\top)^\top$ with $\yb_1, \yb_2 \in {\cal S}^1$. See Figure~\ref{Fig:embedT2'} for an illustration.

%

\begin{figure}[!htbp]
         \centering
         \begin{subfigure}[b]{0.48\textwidth}
         \includegraphics[trim=12 10 0 0mm, scale=0.37]{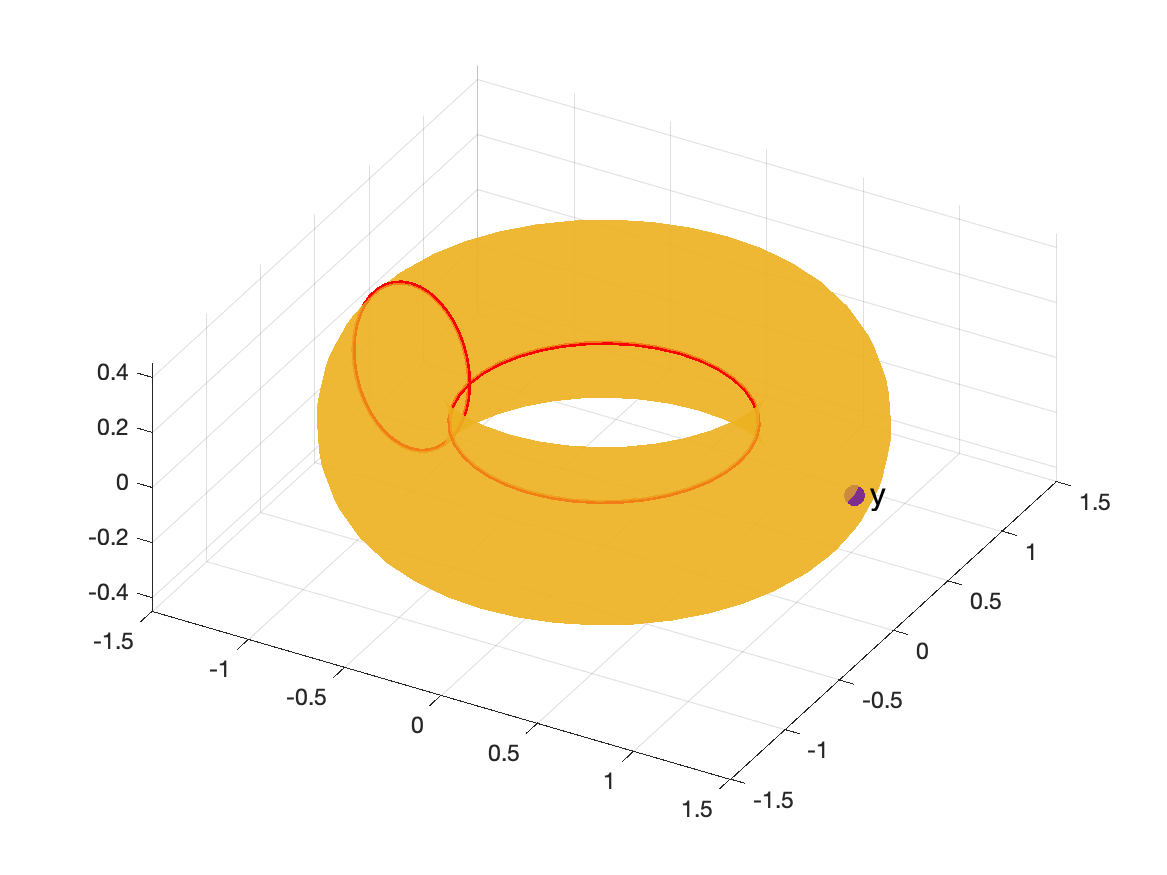} \vspace{-1mm}
         \end{subfigure}
         \begin{subfigure}[b]{0.48\textwidth}
         \includegraphics[trim=12 10 0 0mm, scale=0.37]{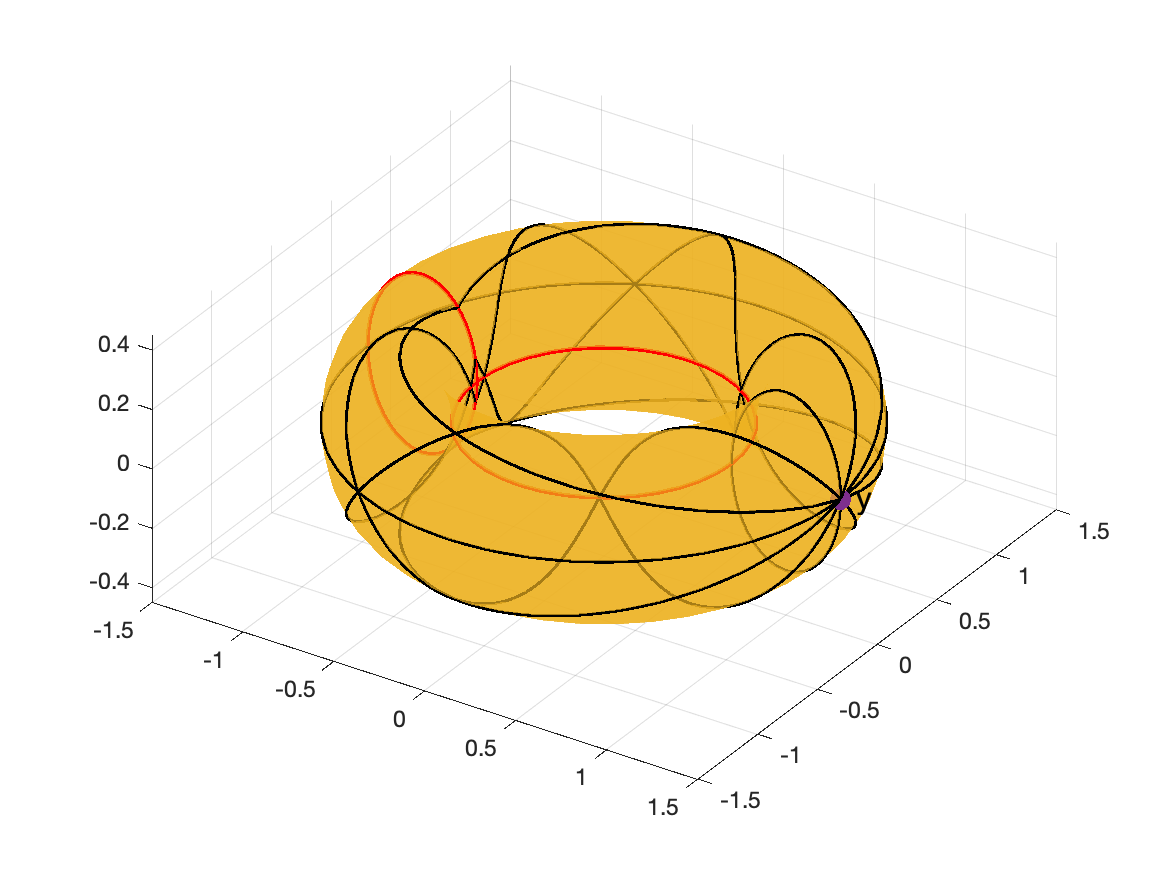} \vspace{-1mm}
         \end{subfigure}
         \caption{\small {Left: The cut locus   $\text{\rm Cut}(\yb)$ (the union of two  red circles)   of $\yb$ (the purple point) in the~2-torus~$\mathcal{T}^2$. Right: A collection of geodesics (the black curves) connecting $\yb$ and $\text{\rm Cut}(\yb)$.}  } 
 \label{Fig:embedT2'}
\end{figure}

Let ${\rm d}f$ denote the differential of a smooth function $f: {\cal M} \rightarrow \mathbb{R}$, which is a smooth 1-form on ${\cal M}$; see, e.g., \citet[Chapters~2 and 14]{lee2012smooth} for more details. The gradient~$\nabla f$ of $f$, which is a $(0, 1)$-tensor field, is defined as follows.

\begin{definition}{\rm 
Let $({\cal M}, g)$ be a Riemannian manifold and $f: {\cal M} \rightarrow ~\!\mathbb{R}$ a mapping from $\cal M$ to $\mathbb R$. Call {\it gradient} of $f$   the vector field $\nabla f \in T{\cal M}$ such that  $g(\nabla f, \vb) = {\rm d}f(\vb)$ for any $\vb \in T{\cal M}$.
}\end{definition}


\section{Proofs}\label{App.proof}

\noindent {\bf Proof of Proposition~\ref{Prop.cont}.} Note that, by construction, $(\yb, {\rm Cut}(\yb))$, for any~$\delta>0$,  is excluded from ${\cal D}_{\delta}$. This allows us to  locally reduce the problem to a  Euclidean one. Then, using the continuity of  optimal transport maps in Euclidean spaces (see \citet[Proposition~5.9]{loeper2009}),  we can prove the desired result as follows. 


Let ${\cal E}_{r\delta} \coloneqq {\cal D}_{\delta} \bigcap \{(\yb,\zb)\in \mathcal{M}\times \mathcal{M}: \yb \in {\rm B}_r^{\cal M}(\yb_0)\}$, where   ${B}^{\cal M}_r(\yb_0)$ denotes the open ball on $\mathcal{M}$ with center $\yb_0$ and radius $r$.  For any $(\yb, \zb) \in {\cal E}_{r\delta}$ and $\xb \in {\rm Cut}(\yb)$, we have
\begin{align*}
d(\zb, {\rm Cut}(\yb_0)) = d(\zb, \xb_0) \geq  d(\zb, \xb) - d(\xb, \xb_0) 
\geq \delta - d(\xb, \xb_0),
\end{align*}
where $\xb_0 \in {\rm Cut}(\yb_0)$ is such that $d(\zb, \xb_0) = d(\zb, {\rm Cut}(\yb_0))$. Taking $\xb \in {\rm Cut}(\yb)$ such that~$d(\xb, \xb_0) = d({\rm Cut}(\yb), \xb_0)$, Assumption~\ref{ass.Cutyy0} then entails
$$d(\zb, {\rm Cut}(\yb_0)) \geq \delta - Cr.$$
Hence, ${\cal E}_{r\delta} \subset {\rm B}_r^{\cal M}(\yb_0) \times \widetilde{\cal D}_{\delta-Cr}(\yb_0)$, where $\widetilde{\cal D}_{\delta}(\yb_0) \coloneqq \{\zb \in \mathcal{M}: d(\zb, {\rm Cut}(\yb_0)) \geq \delta\}$.

Now, for any $(\tilde{\yb}, \zb) \in {\rm B}_r^{\cal M}(\yb_0) \times \widetilde{\cal D}_{\delta-Cr}(\yb_0)$, let $\tilde{\xb} \in {\rm Cut}(\tilde{\yb})$ be such that $$d(\zb, \tilde{\xb}) = d(\zb, {\rm Cut}(\tilde{\yb})).$$ For all $\tilde{\xb}_0 \in {\rm Cut}(\yb_0)$, 
\begin{align*}
d(\zb, {\rm Cut}(\tilde{\yb})) = d(\zb, \tilde{\xb}) \geq d(\zb, \tilde{\xb}_0) - d(\tilde{\xb}_0, \tilde{\xb}) \geq \delta-Cr - d(\tilde{\xb}_0, \tilde{\xb}).
\end{align*}
 Taking $\tilde{\xb}_0$ such that $d(\tilde{\xb}_0, \tilde{\xb}) = d({\rm Cut}(\yb_0), \tilde{\xb})$,   by Assumption~\ref{ass.Cutyy0}, we   have
 $$d(\zb, {\rm Cut}(\tilde{\yb})) \geq \delta-2Cr,$$ which implies ${\rm B}_r^{\cal M}(\yb_0) \times \widetilde{\cal D}_{\delta-Cr}(\yb_0) \subset {\cal D}_{\delta-2Cr}.$ 

Summing up, for any $\yb_0 \in {\cal M}$ and $0<r<\frac{\delta}{2C}$,
$${\cal E}_{r\delta} \subset {\rm B}_r^{\cal M}(\yb_0) \times \widetilde{\cal D}_{\delta-Cr}(\yb_0) \subset {\cal D}_{\delta-2Cr}.$$
In view of Assumption~\ref{Ass.FyCut}, taking $r=\frac{\delta}{3C}$ yields 
$$\{(\yb, \Fb(\yb)), \yb \in {\rm B}_r^{\cal M}(\yb_0)\} \subset {\rm B}_r^{\cal M}(\yb_0) \times \widetilde{\cal D}_{2\delta/3}(\yb_0) \subset {\cal D}_{\delta/3}.$$
Moreover,   $c$ satisfies on ${\cal D}_{\delta/3}$ the regularity conditions in Assumptions~\ref{ass.cost}, \ref{ass.costCurve}, and \ref{ass.convex}. 

Now consider reducing the problem locally at $\yb_0$ to a Euclidean one. Defining 
$${\cal E}^1_r(\yb_0) \coloneqq (\exp_{\yb_0})^{-1}({\rm B}_r^{\cal M}(\yb_0)) \subset \mathbb{R}^{p}\quad\text{and}  \quad {\cal E}^2_{r\delta}(\yb_0) \coloneqq  (\exp_{\yb_0})^{-1}(\widetilde{\cal D}_{\delta-Cr}(\yb_0)) \subset \mathbb{R}^{p},$$
 note that ${\cal E}^1_r(\yb_0)$ and~${\cal E}^2_{r\delta}(\yb_0)$ are diffeomorphic to ${\rm B}_r^{\cal M}(\yb_0)$ and $\widetilde{\cal D}_{\delta-Cr}(\yb_0)$, respectively. Consider  the restriction $\psi^*$ of $\psi$ to ${\cal E}^1_r(\yb_0)$. The following lemma then 
 follows   along the same lines as the proof of Proposition~5.3 in \cite{loeper2011}.



\begin{Lem}\label{lem.psi*}
Let ${\rm P}, {{\rm P}_2 } \in \mathfrak{P}$, where ${{\rm P}_2 }$ has a nonvanishing density on ${\cal M}$. Let Assumptions~\ref{ass.cost} and \ref{ass.costCurve} hold for~$c=~\!d^2/2$. Then,
\begin{compactenum}
  \item[(i)] $\psi^*$ is $c$-convex on ${\cal E}^1_r(\yb_0)$ with respect to the restriction  to ${\cal E}^1_r(\yb_0) \times {\cal E}^2_{r\delta}(\yb_0)$ of $c$;
 \item[(ii)] $\nabla \psi^*$ coincides with $\nabla \psi$ on ${\cal E}^1_r(\yb_0)$;\vspace{.5mm} 
 \item[(iii)]  denoting by  {\rm Vol}  the Riemannian volume measure on ${\cal M}$, {$\dfrac{{\rm d}(\nabla \psi^* \# {\rm Vol})}{{\rm dVol}}\vspace{-1mm}  $ is uniformly bounded over ${\cal E}^1_r(\yb_0)$.}
   \end{compactenum}
\end{Lem}
Proposition~5.9 of \citet{loeper2009}  is a direct consequence of  Lemma~\ref{lem.psi*} and implies the continuity on ${\cal E}^1_r(\yb_0)$ of  $\nabla \psi$. The continuity of $\Fb$  follows. \cqfd
\vspace{4mm}

\color{black}

\noindent {\bf Proof of Proposition~\ref{ProHom}.} 
From Corollary 10 of \cite{McCann2001}, $\Fb =~\!\exp_{\yb}(\nabla \psi (\yb))$ is a bijection from $\mathcal{M}^a$ to $\mathcal{M}^b$, where $\mathcal{M} \setminus \mathcal{M}^a$ is a set of ${\rm P}$-measure zero and $\mathcal{M} \setminus \mathcal{M}^b$ is a set of ${\rm U}_{\mathcal M}$-measure zero. 
Moreover, in view of Proposition~\ref{Prop.cont}, $\Fb$ is continuous over~$\mathcal{M}$.
Since $\mathcal{M}^a$ and $\mathcal{M}^b$, equipped with $d$, are compact metric spaces, it follows (see e.g. Theorem 4.17 of \cite{Rudin}) that $\Fb^{-1}$ is continuous. The result follows.\cqfd
\vspace{4mm}

\noindent {\bf Proof of Lemma~\ref{lem.homotopy}.} 
Part (i). By Assumption~\ref{ass.Wincrease}, for each $\tau \in [0, 1]$ and $\yb \in~\!\mathbb{C}_{{{\cal M}_0}}^{\rm U}(\tau)$, there exists a function $h_{\tau}: \mathbb{C}_{{{\cal M}_0}}^{\rm U}(\tau) \rightarrow \partial{{\cal M}_0}$ with $h_{\tau}(\yb) = \xb$ and $\xb \in \partial{{\cal M}_0}$ such that~$\gamma^{\vb}_{\xb}(cs(\tau)) = \yb$ for some $\vb \in {V}^{\perp}_{\xb}$ and $c\in [0, 1]$. Let $g_{\tau}: \partial{{\cal M}_0} \rightarrow \mathbb{C}_{{{\cal M}_0}}^{\rm U}(\tau)$\linebreak with~$g_{\tau}(\xb) = \gamma^{\vb}_{\xb}(cs(\tau)) = \yb$. Then $(h_{\tau}\circ g_{\tau})(\xb) = \xb$, that is, $h_{\tau}\circ g_{\tau}$ is the identity map on $\partial{{\cal M}_0}$. Also, we have $(g_{\tau}\circ h_{\tau})(\yb) = \yb$, so that $g_{\tau}\circ h_{\tau}$ is the identity map on~$\mathbb{C}_{{{\cal M}_0}}^{\rm U}(\tau)$. Hence, $\partial{{\cal M}_0}$ and $\mathbb{C}_{{{\cal M}_0}}^{\rm U}(\tau)$ are homotopic equivalent.

Part (ii). Consider a function $f: {\cal M} \rightarrow [0, 1]$ such that $f(\yb) = {\rm U}_{\mathcal M}(\mathbb{C}_{{{\cal M}_0}}^{\rm U}(\tau)) = \tau$ if~$\yb\in~\!\mathcal{C}_{{{\cal M}_0}}^{\rm U}(\tau)$. The corresponding sublevel set is $  f^{-1}([0, \tau]) = \mathbb{C}_{{{\cal M}_0}}^{\rm U}(\tau)$ where,  for~$A\subset~\![0, 1]$,  
$f^{-1}(A) \coloneqq \{\xb~\!\in~\!{\cal M}\big\vert f(\xb)~\!\in~\!A\}$.
By construction, $f$ has no critical value (see, e.g., \citet[Chapter 5]{lee2012smooth} for the definition) on $(0, 1)$. It follows from the Morse theory (\citet{milnor1963morse}) that
for any $\tau_1, \tau_2 \in (0, 1)$, $\mathbb{C}_{{\cal M}_0}^{\rm U}(\tau_1)$ is  diffeomorphic to $\mathbb{C}_{{\cal M}_0}^{\rm U}(\tau_2)$.\cqfd
\vspace{4mm}

%


The proof  of  Proposition~\ref{Gliv} relies on two lemmas, the proofs of which follow   along similar lines as the proof of Lemma 11 in \cite{HLV2022} where we refer to for details. 

\begin{Lem}\label{wconv}  Denote by  
 $\gamma \coloneqq({\bf I}_{\rm d} \times {\bf F}) \# {\rm P}$ and $\gamma\n\! \coloneqq({\bf I}_{\rm d} \times {\bf F}\n) \# {\rm P}\n$ 
  the optimal transport plans over $\Gamma({\rm P}, {\rm U}_{\cal M})$ and~$\Gamma( {\rm P}\n, {\rm P}^{{\mathfrak{G}}\n})$, respectively.  Then, 
$\gamma \n$ converges weakly to $\gamma$ as $\ny$. 
\end{Lem}

Before turning to the second lemma, recall that the {\it $c$-subdifferential} of a $c$-convex function $g$ is defined as
 $\partial_c g \coloneqq \{(\yb, \zb) \in \mathcal{M}\times \mathcal{M}: g^c(\zb)-g(\yb)=c(\yb, \zb)\},$ 
where $g^c$ is the $c$-transform of $g$. The $c$-subdifferential of $g$ at  $\yb\in \mathcal{M}$ then is \citep[Definition~5.2]{Villani2009} 
 $\partial_c g(\yb) \coloneqq  \{\zb \in \mathcal{M}: (\yb, \zb) \in \partial_c g\}$. The a.s.~uniqueness of $\Fb$ implies that 
$${\bf F}(\yb) =\exp_{\yb}(\nabla \psi (\yb))=\partial_c \psi (\yb)$$ is $\rm P$-a.s. a singleton. Similarly, associated with $\Fb\n$, there exists a $c$-convex function~$\psi\n$ such that $\Fb\n = \exp(\nabla \psi\n)$ (see, e.g., Theorem 5.10 (ii) in \cite{Villani2009}).  

Combined with the same arguments as in the proof of Lemma~12 in \cite{HLV2022}, Lemma~\ref{wconv} entails the 
convergence of $\partial_c \psi \n$ to $\partial_c \psi$ in the Painlev\' e-Kuratowski sense. Recall that a sequence of set $\mathcal{A}\n$ converges to a set $\mathcal{A}$ in the Painlev\' e-Kuratowski sense (notation: $\mathcal{A}\n \to_{\rm PK} \mathcal{A}$) if (i) the {\it outer limit set} $\lim\sup_{n\to\infty}\mathcal{A}\n$---namely, the set of points $\xb$ for which there exists a sequence~$\xb_n \in \mathcal{A}\n$  admitting a subsequence converging to $\xb$---and (ii) the {\it inner limit set}~$\lim \inf_{n\to\infty} \mathcal{A}\n$---the  set of points $\xb$  which are limits of sequences $\xb_n \in \mathcal{A}\n$---coincide with  $\mathcal{A}$ (see, e.g. Chapter~4 in \cite{RockafellarWets97}).

\begin{Lem}\label{PKconvergence}
Denote by $\partial_c \psi \n$ and $\partial_c \psi$ the $c$-subdifferentials of the empirical and population potentials $\psi\n$ and $\psi$, respectively, Then, 
 $\partial_c \psi \n \to_{\rm PK} \partial_c \psi$ as $\ny$.
\end{Lem}

We now can  conclude with the proof of Proposition~\ref{Gliv}.\medskip 

\noindent {\bf Proof of Proposition~\ref{Gliv}.} 
We first prove the uniform convergency of $\Fb^{(n)}$ to $\Fb$. Assume that $ {\cal M} \times \partial_c\psi\n (\yb_n)$ contains a   sequence $(\yb_n, \zb_n)$ such that, for all $n$ and some~$\delta >0$,
\begin{equation} \label{contradictioneq}
c(\zb_n, \partial_c\psi (\yb_n))>\delta > 0.
\end{equation}
 Since ${\cal M}$ is compact, there exists a subsequence of $\yb_n$ that converges to some $\yb \in {\cal M}$. The continuity of $\partial_c\psi$ together with Lemma \ref{PKconvergence}  (which entails the convergence of $\zb_n$   to~$\partial_c\psi (\yb)= {\bf F}(\yb)$) contradict \eqref{contradictioneq}. The result follows. 
 
 The uniform convergence of $\Qb^{(n)}$ to $\Qb$ follows as a direct consequence of Propo\-sition~\ref{prop.QwAsymptotic}(i).\cqfd
\vspace{4mm}

\noindent {\bf Proof of Proposition~\ref{PropDistFree}.} 
For part (i), noting that the order statistic generates the sigma-field of permutationally invariant measurable functions of $\Yb\n$, it follows from the factorization theorem that it is sufficient for $\mathfrak{P}\n$. 

For part (ii), since ${\rm P} \in \mathfrak{P}$, the probability that $\Yb_i = \Yb_j$ for some  $1\leq i\neq j\leq n$ is zero,  $\Fb\n(\Yb\n)$  with probability one is a permutation of $\mathfrak{G}\n$. Since $\Yb_1\n, ..., \Yb_n\n$ are i.i.d., 
 all $n!$ permutations of $\mathfrak{G}\n$ have the same probability. The result follows.

As for part (iii), first consider  the case that $\tilde{n}_0 =0$ or $1$. Since $\hat{\thetab}\n$is a {sym\-metric} function of $\Yb\n$, $\hat{\thetab}\n$and hence also ${\scriptstyle\mathfrak{G}}\n_*$ are measurable with respect to $\Yb_{(\cdot)}\n$. Thus, conditional on $\Yb_{(\cdot)}\n$, the $n$-tuple of vectors $\Fb\n_*(\Yb\n)$ is uniform over the $n!$ permutations of the grid $\mathfrak{G}\n_{\widehat{\cal M}_0\n}$---irrespective of the value of~$\Yb_{(\cdot)}\n$. Therefore, it is also unconditionally uniform over these permutations.  If $\tilde{n}_0>~\!1$, the situation is exactly the same, except that the~$\tilde{n}_0!$ permutations are  involving    $\tilde{n}_0$    copies in $\mathring{\mathfrak{G}}^{(n_0)}_{\widehat{\cal M}_0\n}$ that are indistinguishable. Therefore, the $n!$ permutations   reduce to $n!/\tilde{n}_0!$ permutations with repetitions. This, however, can be avoided by breaking the $\tilde{n}_0$ ties 
 as in \cite{Hallin2021}.

Parts (iv)-(vi) are direct consequences of part (iii). 

Finally, for part (vii),   assume that the grid does not exhibit any multiplicity at~$\hat{\thetab}\,\!\n$ (either $\tilde{n}_0=0$ or 1, or the tie-breaking device has  been performed). In view of \cite{Basu1959} (see Proposition~E.4 in Appendix~E of \cite{Hallin2021}), since $\Yb_{(\cdot)}\n$ is sufficient and complete, we only have to show that the sigma-field ${\sigma}(\Yb_{(\cdot)}\n, {\bf R}\n, {\Sb}\n)$ is  essentially equivalent to the Borel $\sigma$-field~${\frak F}_{\cal M}$. Since~$\hat{\thetab}\,\!\n$ is $\Yb_{(\cdot)}\n$-measurable,  by defi\-nition of~${{\mathfrak{G}}}\n_{\widehat{\cal M}_0\n}$,  the mapping~$\Yb\n \in {\cal M} \mapsto \left(\Yb_{(\cdot)}\n, {\bf R}\n(\Yb\n), \Sb\n(\Yb\n)\right)$ clearly is injective. The result  follows.\cqfd
\vspace{4mm}

\noindent {\bf Proof of Proposition~\ref{prop.weight}.} 
Since $({\cal M}, d)$ is a compact metric space, for any $f \in~\!{\cal F}_{{\rm BL}_1}^{\cal M}$, by the Arzel\'a–Ascoli theorem, there exists a sequence $f_1, \ldots, f_{m_\epsilon}$ such that, for any~$\epsilon>~\!0$,
$$\sup_{f \in {\cal F}_{{\rm BL}_1}^{\cal M}} \inf_{k=1, \ldots, m_\epsilon} \Vert f - f_k \Vert_\infty \leq \epsilon/4.$$
This entails, for any $f \in {\cal F}_{{\rm BL}_1}^{\cal M}$, 
\begin{align*}
& \left\vert  \sum_{j=1}^n w_j\n(\Xb, \Xb\n) f(\Yb_j\n) - {\rm E}(f(\Yb)\vert\Xb))  \right\vert \\
&\quad\quad\leq \sup_{k=1, \ldots, m_\epsilon} \left\vert  \sum_{j=1}^n w_j\n(\Xb, \Xb\n) f_k(\Yb_j\n) - {\rm E}(f_k(\Yb)\vert\Xb))  \right\vert + \epsilon/2.
\end{align*}
Hence, there exists $n_0$ such that, for $n\geq n_0$,
\begin{align}
{\rm P}\left[{d}_{\rm BL}\big({\rm P}\n_{\wb\n(\Xb)}, {\rm P}_{\Yb\vert\Xb}\big) > \epsilon \right] 
=& \, {\rm P}\left[\sup_{f \in {\cal F}_{{\rm BL}_1}^{\cal M}} \left\vert  \sum_{j=1}^n w_j\n(\Xb, \Xb\n) f(\Yb_j\n) - {\rm E}(f(\Yb)\vert\Xb))  \right\vert > \epsilon  \right] \nonumber \\
\leq & \, {\rm P}\left[\sup_{k=1, \ldots, m_\epsilon} \left\vert  \sum_{j=1}^n w_j\n(\Xb, \Xb\n) f_k(\Yb_j\n) - {\rm E}(f_k(\Yb)\vert\Xb))  \right\vert > \epsilon/2  \right] \nonumber \\
\leq & \, \sum_{k=1}^{m_\epsilon} {\rm P}\left[\left\vert  \sum_{j=1}^n w_j\n(\Xb, \Xb\n) f_k(\Yb_j\n) - {\rm E}(f_k(\Yb)\vert\Xb))  \right\vert > \epsilon/2  \right]. \nonumber
\end{align}
Since $\wb\n$ is a consistent weight function, the RHS of the above inequality is bounded by a positive number $\delta$ that can be arbitrarily small. The result follows. \cqfd
\vspace{4mm}

\noindent {\bf Proof of Proposition~\ref{prop.FrConsistent}.} 
Define  $$G_{\xb}(\yb) \coloneqq  {\rm E}_{{\rm P}_{\Yb\vert\Xb=\xb}}[c(\Yb, \yb)],\, G_{\xb}\n(\yb) = \sum_{j=1}^n w_j\n(\xb, \Xb\n) c(\Yb_j, \yb)$$ 
and let $a \coloneqq  \min\{G_{\xb}(\yb): \yb \in {\cal M}\}$ and $\delta(\epsilon)>0$  be such that
$$a + \delta(\epsilon) = \min\{G_{\xb}(\yb): \yb \in {\cal M}_1\}$$
where ${\cal M}_1 \coloneqq  {\cal M}\setminus \mathcal{A}_{\rm Fr}^{\epsilon}(\xb)$ for some $\epsilon>0$. It is sufficient to prove that
\begin{equation}\label{eq.GGn}
\max_{\yb\in {\cal M}} \vert G_{\xb}(\yb) - G_{\xb}\n(\yb) \vert = o_{\rm P}(1)
\end{equation}
as $n\rightarrow \infty$. To see this, note that, if \eqref{eq.GGn} holds, then
$$\min\{G_{\xb}\n(\yb): \yb\in \mathcal{A}_{\rm Fr}^{\epsilon}(\xb)\} \leq a + \frac{\delta(\epsilon)}{3} + o_{\rm P}(1)$$ 
and
$$\min\{G_{\xb}\n(\yb): \yb\in {\cal M}_1\} \geq a + \frac{\delta(\epsilon)}{2} + o_{\rm P}(1),$$
from which the desired result in Proposition~\ref{prop.FrConsistent} follows.

Now, for any $\yb, \tilde{\yb} \in {\cal M}$, by definition of $G_{\xb}(\cdot)$, we have
\begin{align*}
\vert G_{\xb}(\yb) - G_{\xb}(\tilde{\yb}) \vert & \leq \max\{\vert c(\zb, \yb) -  c(\zb, \tilde{\yb})\vert: \zb \in {\cal M}\} \\
&= \max\{\vert d^2(\zb, \yb) -  d^2(\zb, \tilde{\yb})\vert/2: \zb \in {\cal M}\} \\
&\leq \max\{\vert u^2 -  \tilde{u}^2\vert/2: \vert u -\tilde{u}\vert \leq d(\yb, \tilde{\yb}),\quad  u, \tilde{u} \in [0, {\rm diam}({\cal M})]\}
\end{align*}
where ${\rm diam}({\cal M})$ denotes the diameter of ${\cal M}$ (the last inequality above is due to the fact that $\vert d(\zb, \yb) -  d(\zb, \tilde{\yb})\vert \leq d(\yb, \tilde{\yb})$). In the same manner, we have
$$\vert G\n_{\xb}(\yb) - G\n_{\xb}(\tilde{\yb}) \vert \leq \max\{\vert u^2 -  \tilde{u}^2\vert/2: \vert u -\tilde{u}\vert \leq d(\yb, \tilde{\yb}), \quad  u, \tilde{u} \in [0, {\rm diam}({\cal M})]\}.$$
Therefore, $G_{\xb}(\cdot)$ and $G\n_{\xb}(\cdot)$ are uniformly continuous on ${\cal M}$, and hence there  \linebreak exists~$\delta(\tilde{\epsilon})>0$ such that, for $d(\yb, \tilde{\yb})<\delta(\tilde{\epsilon})$,
$$\vert G_{\xb}(\yb) - G_{\xb}(\tilde{\yb}) \vert \leq \frac{\tilde{\epsilon}}{4} \ \ \text{and} \ \ \vert G\n_{\xb}(\yb) - G\n_{\xb}(\tilde{\yb}) \vert \leq {\tilde{\epsilon}}/{4}.$$
Let $\{\yb_1, \ldots, \yb_k\}$ be a $\delta(\tilde{\epsilon})$-net of ${\cal M}$---namely, be such that,  for all $\yb\in {\cal M}$, there exists a~$\tilde{\yb}\in \{\yb_1, \ldots, \yb_k\}$ such that $d(\yb, \tilde{\yb})<\delta(\tilde{\epsilon})$. It follows from \eqref{eq.cstWeight} that, for any $\tilde{\epsilon}>0$,
$$\underset{n\rightarrow\infty}{\lim} {\rm P}[\vert G_{\xb}(\yb_{\ell}) -  G\n_{\xb}(\yb_{\ell}) \vert > \tilde{\epsilon}/4] = 0, \ \forall \ell = 1, \ldots, k.$$
Then, for all $\yb \in {\cal M}$, 
\begin{align*}
\vert G_{\xb}(\yb) - G_{\xb}\n(\yb) \vert &\leq \vert G_{\xb}(\yb) - G_{\xb}(\tilde{\yb}) \vert + \vert G_{\xb}(\tilde{\yb}) - G_{\xb}\n(\tilde{\yb}) \vert + \vert G\n_{\xb}(\tilde{\yb}) - G_{\xb}\n(\yb) \vert 
\leq \frac{3\tilde{\epsilon}}{4}
\end{align*}
with probability tending to one as $n\rightarrow \infty$. Therefore, \eqref{eq.GGn} holds and the desired result follows.\cqfd
\vspace{4mm}

\noindent {\bf Proof of Lemma~\ref{lem.Mxhat}.} Denote by~$\Yb_j$, $j \in \{1, \ldots, n\}$ the sample point closest to~$\hat{\thetab}\n(\xb)$ in the~$d(\cdot, \cdot)$  metric. By Assumption~\ref{Ass.M0x}, $\Yb_j$ is a consistent estimator of $\thetab(\xb)$. For any fixed $\xb\in {\cal M}_{\Xb}$, by the definition of ${\scriptstyle\mathfrak{G}}^{(N)}_*(\xb)$, there exists a $c$-convex function $\psi^{(N, n)}_{\xb}$ such that
$${\scriptstyle\mathfrak{G}}^{(N)}_*(\xb) = \exp_{\Yb_j}(\nabla \psi_{\xb}^{(N, n)} (\Yb_j)).$$   
Repeating the same arguments as in the proof of Proposition~\ref{Gliv}, we obtain that~$\exp_{\Yb_j}(\nabla \psi_{\xb}^{(N, n)} (\Yb_j))$ converges almost surely to $\Fb(\Yb_j\vert \xb)$. Therefore, due to the continuity of $\Fb$ and the convergence in probability of $\Yb_j$ to $\thetab(\xb)$, ${\scriptstyle\mathfrak{G}}^{(N)}_*(\xb)$ converges in probability to $\Fb(\thetab(\xb) \vert \xb)$ as $n, N\rightarrow \infty$. The result then follows from Assumption~\ref{Ass.M0x}. \cqfd
\vspace{4mm}

\noindent {\bf Proof of Proposition~\ref{prop.QwAsymptotic}.} Part {(i)}. The a.s.~uniqueness of $\Fb$ and Proposition~\ref{ProHom} imply that, for each $\xb \in {\cal M}_0$, there exists a $c$-convex function $\phi_{\xb}$ such that
$$\Qb({\scriptstyle{\mathfrak{G}}}^{(N)}_{1i}(\xb) \vert \xb) =\exp_{{\scriptstyle{\mathfrak{G}}}^{(N)}_{1i}}(\nabla \phi_{\xb} ({\scriptstyle{\mathfrak{G}}}^{(N)}_{1i}(\xb)))=\partial_c \phi_{\xb} ({\scriptstyle{\mathfrak{G}}}^{(N)}_{1i}(\xb))$$ is $\rm P$-a.s. a singleton. Also,  there exists a $c$-convex function $\phi_{\xb}^{(N, n)}$ such that
$$\Qb_{\wb}^{(N, n)}({\scriptstyle{\mathfrak{G}}}^{(N)}_{1i}(\xb) \vert \xb) =\exp_{{\scriptstyle{\mathfrak{G}}}^{(N)}_{1i}}(\nabla \phi_{\xb}^{(N, n)} ({\scriptstyle{\mathfrak{G}}}^{(N)}_{1i}(\xb))).$$
Due to Proposition~\ref{prop.weight}, it follows from the same arguments as in  the proof of Lemma 12 in \cite{HLV2022} that $\partial_c \phi_{\xb}^{(N, n)} \to_{\rm PK} \partial_c \phi_{\xb}$ as $N, n\rightarrow \infty$. To conclude the proof of part {(i)}, assume that there exists a sequence  $(\yb_n, \zb_n) \in {\cal M} \times \partial_c\phi_{\xb}^{(N, n)} (\yb_n)$ such that, for all $n$ and some $\delta >0$,
\begin{equation} \label{contradictioneq2}
c(\zb_n, \partial_c\phi_{\xb} (\yb_n))>\delta > 0.
\end{equation}
 Since ${\cal M}$ is compact, there exists a subsequence of $\yb_n$ that converges to some $\yb \in {\cal M}$. The continuity of $\partial_c\phi_{\xb}$ together with the fact that $\partial_c \phi_{\xb}^{(N, n)} \to_{\rm PK} \partial_c \phi_{\xb}$ imply that $\zb_n$ converges to $\partial_c\phi_{\xb} (\yb)= \Qb(\yb\vert \xb)$, which contradicts \eqref{contradictioneq2}. Part (i) follows.

Part {(ii)}. Let 
\[
\mathfrak{H}_r^{(N)}(\xb) \coloneqq 
\begin{cases}
	\mathring{\mathfrak{G}}^{(N_0)}_{\widehat{\cal M}^{(N)}_0(\xb)}  & \text{if}\  \ r=0\\
	{\mathfrak{G}}^{(N_S)}_{\widehat{\cal M}^{(N)}_0(\xb)}(1) \bigcup \cdots \bigcup {\mathfrak{G}}^{(N_S)}_{\widehat{\cal M}^{(N)}_0(\xb)}(r)  & \text{if}\  \ r\geq 1 \\
\end{cases}
\]
denote the subset of $\mathfrak{G}^{(N)}_{\widehat{\cal M}^{(N)}_0(\xb)}$ that is used in the construction of  the conditional (on~$\Xb=~\!\xb$) quantile regions. 
We first show that, for all $r\in {\cal R}_N$, the Hausdorff distance between $\mathfrak{H}_r^{(N)}(\xb)$ and $\mathbb{C}_{{\cal M}_0(\xb)}^{\rm U}(\tau_r\vert \xb)$ ($\tau_r = r/(N_R+1)$) converges to $0$.


 The proof of Lemma~\ref{lem.Mxhat} implies that, with probability one, ${\scriptstyle\mathfrak{G}}^{(N)}_*(\xb)$ belongs to the $\epsilon_2$-fattening of $\{\Fb(\thetab(\xb) \vert \xb)\}$ for any $\epsilon_2>0$. Indeed, it implies that, as $n, N \rightarrow\infty$,
$$\max_{r\in {\cal R}_N} d_{\rm H}\left(  \mathfrak{H}_r^{(N)}(\xb), \tilde{\mathfrak{H}}_r^{(N)}(\xb)\right) = o_{\rm P}(1),$$
where $\tilde{\mathfrak{H}}_r^{(N)}(\xb)$ denotes the set of regular grids with   center  ${\scriptstyle\mathfrak{G}}^{(N)}_*(\xb)$    replaced by~$\Fb(\thetab(\xb) \vert \xb)$ in the definition of $\mathfrak{H}_r^{(N)}(\xb)$. By this construction, we have
$$\max_{r\in {\cal R}_N} d_{\rm H}\left(  \tilde{\mathfrak{H}}_r^{(N)}(\xb), \mathbb{C}_{{\cal M}_0(\xb)}^{\rm U}(\tau_r\vert \xb)\right)  = o_{\rm P}(1) \quad \text{as} \ \ N\rightarrow \infty.$$
Therefore, as $n$ and $N\rightarrow \infty$, 
$$\max_{r\in {\cal R}_N} d_{\rm H}\left({\mathfrak{H}}_r^{(N)}(\xb), \mathbb{C}_{{\cal M}_0(\xb)}^{\rm U}(\tau_r\vert \xb)\right) = o_{\rm P}(1).$$

Now, the continuity of $\Qb$ entails that
$$\max_{r\in {\cal R}_N} d_{\rm H}\left(\{\Qb(\yb_1 \vert \xb): \yb_1 \in \mathfrak{H}_r^{(N)}(\xb)\}, \{\Qb(\yb_2 \vert \xb): \yb_2 \in \mathbb{C}_{{\cal M}_0(\xb)}^{\rm U}(\tau_r\vert \xb) \right) = o_{\rm P}(1)$$
as $n, N\rightarrow \infty$. The desired uniform convergence of $\mathbb{C}_{\wb}^{(N, n)}(\cdot \vert \xb)$ to $\mathbb{C}_{{\cal M}_0(\xb)}^{\rm P}(\cdot \vert \xb)$ then follows from the uniform convergence of $\Qb_{\wb}^{(N, n)}(\cdot \vert \xb)$ to $\Qb(\cdot \vert \xb)$ in part {(i)}. The same arguments yield the desired uniform convergence of $\mathcal{C}_{\wb}^{(N, n)}(\cdot \vert \xb)$ to $\mathcal{C}_{{\cal M}_0(\xb)}^{\rm P}(\cdot \vert \xb)$. \cqfd
\vspace{4mm}

\end{document}